\documentclass[reqno,10pt]{amsart}

\usepackage[hmargin=1.5cm,top=1.5cm,bottom=1.5cm,includehead]{geometry}
\usepackage{amsmath,mathrsfs}
\usepackage{amssymb, amsmath}
\usepackage{graphicx}
\usepackage{enumerate,color,xcolor}
\usepackage{upgreek}
\usepackage{hyperref}

\newtheorem{thm}{Theorem}[section]
\newtheorem{lem}{Lemma}[section]
\newtheorem{rmk}{Remark}[section]

\newtheorem{prop}{Proposition}[section]

\newcommand{\vv}[1]{\boldsymbol{#1}}

\numberwithin{equation}{section}
\let\f=\frac

\newcommand{\bit}{\begin{itemize}}
	\newcommand{\eit}{\end{itemize}}

\newcommand{\N}{{\mathbb N}}
\newcommand{\Z}{{\mathbb Z}}
\newcommand{\R}{{\mathbb R}}
\newcommand{\T}{{\mathbb T}}
\newcommand{\pa}{\partial}
\newcommand{\vf}{\varphi}
\newcommand{\eps}{\varepsilon}
\newcommand{\ov}{\overline}

\newcommand{\ba}{\begin{aligned}}
	\newcommand{\ea}{\end{aligned}}
\def\na{\nabla}

\let\wt=\widetilde
\let\wh=\widehat
\let\udl=\underline
\def\fA{{\mathfrak A}}
\def\fB{{\mathfrak B}}
\def\fC{{\mathfrak C}}
\def\fD{{\mathfrak D}}
\def\fR{{\mathfrak R}}
\def\fE{{\mathfrak E}}
\def\fM{{\mathfrak M}}
\def\cD{{\mathcal D}}
\def\cE{{\mathcal E}}
\def\sF{{\mathscr F}}
\def\cL{{\mathcal L}}
\def\cN{{\mathcal N}}

\def\sR{{\mathscr R}}

\def\sss{\mathtt{s}}
\def\bfe{\vec{\mathbf{e}}}
\def\fz{\mathrm{z}}
\def\fv{\mathrm{v}}
\def\mrv{\mathring{v}}
\def\mrxi{\mathring{\xi}}
\def\virgp{\raise 2pt\hbox{,}}
\def\cdotpv{\raise 2pt\hbox{;}}

\def\<{\langle}
\def\>{\rangle}

\def\gs{\gtrsim}
\def\ls{\lesssim}
\newcommand{\rr}[1]{\left( #1 \right)}
\newcommand{\nr}[1]{\left| #1 \right|}
\newcommand{\nnr}[1]{\left\| #1 \right\|}

\def\lta{\ell^2(|\alpha|\leq2)}
\def\ltaa{\ell^2_\alpha}
\def\etmu{\mathsf{e}_\theta(\upmu)}
\def\Etmu{\mathsf{E}_\theta(\upmu)}
\def\Etnu{\mathsf{E}_\theta(\upnu)}
\begin{document}

	\title[Plasma Echo Cascades at the Gevrey-3 Threshold]{Infinitely Many Plasma Echo Cascades and\\ Sobolev Norm Divergence at the Gevrey-3 Threshold}

	\author{Yue Luo}
	\address[Yue Luo]{Tsinghua University}
	\email{luo-y21@tsinghua.org.cn}

	\begin{abstract}
We construct a global solution to the confined Vlasov--Poisson system near a fixed homogeneous Penrose-stable equilibrium that exhibits infinite-time growth driven by nonlinear plasma echoes. For every sufficiently small \(\epsilon>0\), the initial perturbation has size \(O(\epsilon)\) in a Gevrey-3 norm with radius \(\lambda_0\simeq\epsilon^{1/3}\). The corresponding solution is real and neutral, and the full distribution remains nonnegative. Its kinetic profile scatters in a lower Sobolev topology, while the spatial density and electric field decay. Nevertheless, the kinetic profile satisfies

\[
\lim_{t\to\infty}\|f(t)\|_{H^{1000}_{z,v}}=\infty.
\]

The construction combines infinitely many well-separated echo cascades within a single global nonlinear solution. Each cascade generates a high-frequency output that persists under the subsequent evolution, and the accumulated outputs produce divergence in a fixed Sobolev norm. Thus, unlike finite-time echo amplification, the instability persists to arbitrarily large times despite scattering and macroscopic Landau damping. The construction identifies an instability mechanism at the critical amplitude-dependent Gevrey-3 radius \(\lambda_0\simeq\epsilon^{1/3}\), where cumulative echo amplification balances the initial Fourier decay.
\end{abstract}
	
	\maketitle

	


	\setcounter{tocdepth}{1}
	\tableofcontents

	\section{Introduction}
	Consider the Vlasov-Poisson system for the fluctuation $h(t,x,v)$ defined on $I\times\T^3_x\times\R^3_v$ around an equilibrium $f^0=f^0(v)$, where $I$ is a time interval and $\T={\R}/{(2\pi \Z)}$:
	\begin{equation}\label{VP}
		\left\{\begin{aligned}
			&\pa_th+v\cdot\na_xh+E(t,x)\cdot\na_v(f^0(v)+h)=0, \\
			&E(t,x)=\na_x\Delta_x^{-1}\rho(t,x), \ \rho(t,x)=\int_{\R^3} h(t,x,v)dv
		\end{aligned}\right.
	\end{equation}
	with the assumption
	\[\int_{\T^3\times\R^3}h(t,x,v)dxdv=0.\]
	Taking out the linear transport effect, we define $f(t,z,v)=h(t,z+vt,v)$, which then satisfies the reduced system
	\begin{equation}\label{1}
		\left\{\begin{aligned}
			&\pa_tf+E(t,z+vt)\cdot(\na_v-t\na_z)(f^0(v)+f)=0, \\
			&E(t,x)=\na_x\Delta_x^{-1}\rho(t,x), \ \rho(t,x)=\int_{\R^3} f(t,x-vt,v)dv.
		\end{aligned}\right.
	\end{equation}

	\subsection{Landau damping, global stability, and plasma echoes}
	
	Landau damping describes the decay of macroscopic electric fields in collisionless plasmas through phase mixing. Following Landau's linear analysis \cite{Landau} and the Penrose stability criterion \cite{Penrose}, the nonlinear stability theory was developed by Mouhot and Villani \cite{MV}, who established Landau damping and scattering for small analytic perturbations of homogeneous equilibria. Subsequent works extended these results to Gevrey regularity: Bedrossian, Masmoudi, and Mouhot \cite{BMM}, and Grenier, Nguyen, and Rodnianski \cite{GNR}, established nonlinear Landau damping for Gevrey weights \(e^{\lambda\langle\nabla\rangle^s}\) with \(s>1/3\). The endpoint \(s=1/3\), corresponding to Gevrey-3 regularity, was obtained by Ionescu, Pausader, Wang, and Widmayer \cite{waveop}.
	
	A central difficulty in nonlinear Landau damping is the plasma echo mechanism; see \cite{echoes,review}. Although phase mixing leads to decay of the electric field, nonlinear interactions can regenerate resonances at later times. These resonances may occur successively, producing echo chains with substantial amplification. Bedrossian \cite{echoes} constructed arbitrarily long nonlinear echo chains, demonstrating the failure of a general extension of nonlinear Landau damping to Sobolev regularity. Zillinger \cite{Zillinger} studied echo chains in the linearized Vlasov--Poisson equation around small, spatially inhomogeneous, time-dependent solutions, identifying Gevrey-3 as a critical regularity threshold. These works demonstrate that the echo mechanism can produce significant growth in strong norms, even when macroscopic damping occurs.
	
	In this paper, we construct a global nonlinear solution near a fixed homogeneous Penrose-stable equilibrium that exhibits persistent echo-driven growth over arbitrarily large times. Unlike constructions based on a single finite echo chain, our solution contains infinitely many well-separated cascades. Each cascade generates a high-frequency output that persists under the subsequent nonlinear evolution, and the outputs accumulated over increasingly large frequency scales force the kinetic profile to diverge in a fixed high Sobolev norm as \(t\to\infty\). At the same time, the kinetic profile scatters in a lower Sobolev topology, while the spatial density and electric field decay. Thus, a single global nonlinear solution exhibits both macroscopic Landau damping and unbounded growth in a strong Sobolev norm.

	The construction takes place at the amplitude-dependent Gevrey-3 radius \(\lambda_0\simeq\epsilon^{1/3}\). This scale reflects the balance between the initial Fourier decay and the cumulative amplification along echo chains. The precise balance is described in Section 1.3.

\subsection{Main result and its interpretation}

	We work around a fixed homogeneous equilibrium $f^0(v)$. We assume that there
	exists $0<\theta<\f12$ such that
		\begin{equation}\label{Penrose}
			\begin{aligned}
				\nr{\wh{f^0}(\eta)}+\nr{\na_\eta\wh{f^0}(\eta)}&\leq\theta^{-2}e^{-2|\eta|},&\qquad f^0(v)&\geq\theta e^{-\f12|v|^2}, \\
				\inf_{\substack{\mu\in\Z^3\backslash\{0\}\\\Im z\ge-\f12|\mu|}}\nr{1+\int_0^\infty t\wh{f^0}(\mu t)e^{i tz}dt}&\geq\theta,&\qquad \nnr{e^{2\<\na\>}\rr{v^\alpha f^0(v)}}_{\lta L^2_v}&\leq\theta^{-2}.
			\end{aligned}
		\end{equation}  
		Here and below, velocity moments are summed over $|\alpha|\leq2$. The first
		line gives quantitative analyticity and a Gaussian lower bound, while the
		second line contains the Penrose stability condition and the weighted
		regularity used in the construction.
		\begin{rmk}
		The Gaussian lower bound is used only to preserve nonnegativity of the full
		distribution $f^0+f$. Its particular variance, as well as the numerical
		analytic radius in \eqref{Penrose}, is inessential.
		\end{rmk}

			For later use, define
			\[
			\|g\|_{\mathcal X^m}:=
			\nnr{\<\na\>^m\rr{v^\alpha g}}_{\ell^2(|\alpha|\leq2)L^2_{z,v}}.
			\]
			and
			\[
			\|g\|_{\mathcal G^{3,m}_{\lambda}}
			:=\nnr{\<\na\>^m e^{\lambda\<\na\>^{1/3}}
				\rr{v^\alpha g}}_{\ell^2(|\alpha|\leq2)L^2_{z,v}}.
			\]
			For each perturbation size $\epsilon>0$, denote by
			$X_*=X_*(\epsilon)$ the critical initial-data space with norm
			\[
			\|g\|_{X_*}:=\nnr{\<\na\>^3
				e^{\<\epsilon\na\>^{1/3}}
				\rr{v^\alpha g}}_{\ell^2(|\alpha|\leq2)L^2_{z,v}},
			\]
			The main theorem constructs a global solution in which echo-driven growth in
				a strong topology coexists with scattering and macroscopic damping.
			 \begin{thm}[Global echo-driven instability at the Gevrey-3 threshold]\label{T1}
				There exist $\epsilon_{\rm th}=\epsilon_{\rm th}(\theta)>0$ and $C=C(\theta)>0$ such that, for every $0<\epsilon<\epsilon_{\rm th}$, there is a real-valued, neutral solution $f(t,z,v)$ to \eqref{1} on $[0,\infty)\times\T^3\times\R^3$ such that, for every $R>0$,
				\[
				f\in C([0,R];\mathcal X^{20})\cap C^1([0,R];\mathcal X^{17}),
				\qquad \int_{\T^3\times\R^3}f(t,z,v)\,dzdv=0,
				\]
				and $f^0(v)+f(t,z,v)\geq0$. Its initial profile satisfies
				\[
				\|f(0)\|_{X_*}\leq C\epsilon,
				\]
				where $X_*=X_*(\epsilon)$ is the initial-data space defined above.
			Moreover, there exists a profile $f_\infty(z,v)$ such that
			\[\lim_{t\to\infty}\nnr{v^\alpha\rr{f(t)-f_\infty}}_{\ell^2(|\alpha|\leq2)H^{18}_{z,v}}=0,
			\qquad \lim_{t\to\infty}\rr{\nnr{\rho_f(t)}_{H^{10}_x}+\nnr{E_f(t)}_{H^{11}_x}}=0,\]
			while
			\[\lim_{t\rightarrow\infty}\nnr{\<\na\>^{1000}f(t)}_{L^2_{z,v}}=\infty.\]
		\end{thm}
			\begin{rmk}
				Theorem \ref{T1} exhibits two distinct aspects of the long-time behavior of a single global nonlinear solution. On the one hand, the kinetic profile scatters in \(H^{18}\), and both the spatial density and electric field decay. On the other hand, its \(H^{1000}\) norm tends to infinity as \(t\to\infty\). Thus, convergence of the kinetic profile in a lower Sobolev topology is compatible with divergence in a stronger topology. In particular, the conclusion describes genuine infinite-time growth of a fixed Sobolev norm, rather than merely a loss of regularity in the estimates or a large amplification over a finite time interval.				
			\end{rmk}
			\begin{rmk}[A fixed-space consequence]
				The space $X_*(\epsilon)$ in Theorem \ref{T1} depends on the size of the
				perturbation. To formulate a consequence in a fixed Banach space, define
				\[
				\|g\|_{X_{\log}}:=\nnr{\<\na\>^3
				e^{\f{\<\na\>^{1/3}}{\log^+\log^+|\na|}}
				\rr{v^\alpha g}}_{\ell^2(|\alpha|\leq2)L^2_{z,v}},
				\qquad \log^+x:=\log(10+x).
				\]
				From scale separation in the construction we can see that
				\[
				\|f(0)\|_{X_{\log}}\leq C\epsilon.
				\]
			\end{rmk}
		\begin{rmk}
			We use the standard Gevrey-3 radius convention \(e^{\lambda_0\langle\nabla\rangle^{1/3}}\). The critical amplitude-dependent radius is \(\lambda_0\sim\epsilon^{1/3}\), which reflects the balance between the initial Fourier decay and the cumulative amplification along echo chains. The Gevrey-3 stability theory of \cite{waveop} is consistent with stability when \(\lambda_0=M\epsilon^{1/3}\) for sufficiently large \(M\), while Theorem~\ref{T1} establishes instability when \(\lambda_0/\epsilon^{1/3}\) is of order one. Together, these results reveal the critical amplitude--radius balance governing the nonlinear echo mechanism at the Gevrey-3 threshold.
			
		\end{rmk}
		\begin{rmk}
			The construction is carried out in the phase space \(\mathbb T^3\times\mathbb R^3\), within the invariant class of solutions depending only on the first spatial coordinate. The echo cascades, including the generation and persistence of their high-frequency outputs, take place entirely within this invariant class. In particular, the instability mechanism does not require genuinely multidimensional resonance interactions.		
		\end{rmk}

	\subsection{The critical echo balance and new ingredients}
	The exponent $1/3$ arises directly from the cumulative echo gain. Insert a
	packet centered near $(k,\mathsf A_k)$, where
	\[
	\mathsf A_k=\f{(1+\eps_0^6)k^3}{\eps_0},\qquad
	t_\upmu=\f{\mathsf A_k}{\upmu+\f12}.
	\]
	Here \(\epsilon_0\) denotes the amplitude of the persistent low-frequency background wave. The inserted perturbation is initially concentrated near the frequency \((k,\mathsf A_k)\). During the echo cascade, its spatial Fourier frequency decreases successively from \(k\) to \(0\), while its velocity frequency remains close to \(\mathsf A_k\). The critical times \(t_\upmu\) describe the successive resonant interactions within this cascade. In the global construction, we choose an increasing sequence \(k_j\) so that the corresponding frequency scales \(\mathsf A_{k_j}\) are sufficiently separated.	Across the critical interval $[t_\upmu,t_{\upmu-1}]$, the leading transfer
	from spatial mode $\upmu$ to $\upmu-1$ has size
	\[
	\f{\eps_0\mathsf A_k}{\upmu^3}.
	\]
		Ignoring for the moment the controlled errors in the recurrence, the formal
		leading gain is
		\[
		\prod_{\upmu=1}^k\f{\eps_0\mathsf A_k}{\upmu^3}
		=\f{\rr{(1+\eps_0^6)k^3}^k}{(k!)^3}
		=(2\pi k)^{-3/2}
		\exp\!\rr{\rr{3+\log(1+\eps_0^6)}k}\rr{1+O(k^{-1})}.
		\]
		Sections 3--5 prove that the omitted transfer errors and nonlinear corrections
		are smaller than this principal term. The gain balances a packet of size about
		$e^{-3k}$. Since
	$\mathsf A_k\simeq k^3/\eps_0$, the exponential scale $e^{-ck}$ becomes the
	critical Fourier decay $e^{-c\eps_0^{1/3}\mathsf A_k^{1/3}}$. This balance between echo amplification and the initial Fourier decay explains the critical Gevrey-3 regularity threshold.

	The proof introduces six ingredients needed to turn this balance into a
	global nonlinear solution.
	\begin{enumerate}[(i)]
	\item We establish weighted energy estimates for the linearized equation around a time-dependent, spatially inhomogeneous background, including equations with source terms. These estimates give coarse Fourier localization and control the nonresonant frequency regions, providing the basis for the refined analysis in the subsequent sections.	
	\item We study the homogeneous linearized equation for a carefully chosen high-frequency perturbation around a prescribed time-dependent background containing a persistent low-frequency wave. Using the Duhamel principle, we analyze the successive resonant interactions and establish matching upper and lower bounds for the amplification along the complete echo chain.
	\item We establish refined estimates for the inhomogeneous linearized equation with source terms, which can generate additional echo chains. Using the Fourier localization obtained from the energy estimates, we apply the Duhamel principle to control the amplification in the resonant frequency regions. These estimates are essential for controlling higher-order nonlinear interactions.
	\item We introduce a Taylor hierarchy \(f^{(n,q)}\) that records the nonlinear interaction order and the Fourier frequencies generated at each order. Combined with the estimates above, it allows us to control the additional echo chains arising from higher-order nonlinear interactions and prove that the first-order echo output dominates the nonlinear corrections.
	\item A scale-separation induction inserts infinitely many cascades into the already constructed nonlinear solution. Summability ensures convergence, while frequency separation prevents cancellation between the outputs of different cascades.
	\item We establish uniform estimates for the sequence of approximate solutions that allow us to pass to the global limit while preserving the lower bounds produced by individual echo cascades. The same estimates also control the long-time behavior of the limiting solution, yielding scattering and decay of the macroscopic observables.
	
	\end{enumerate}

	\subsection{Strategy of the proof}
	The notation for the spatial and velocity components and their Fourier frequencies is introduced in Section \ref{s1.5}.
	    Proposition \ref{propreductionto1d} allows us to work with solutions depending only on the first spatial coordinate. The Fourier form \eqref{xv} of equation \eqref{1} then reduces to
	\[\pa_t\wh{f}(t,\upmu\bfe,\xi)+\sum_{\upnu\in\Z\backslash\{0\}}\wh{\rho}(t,\upnu\bfe)\f{\upxi-\upmu t}{\upnu}\sF\rr{f^0+f}(t,(\upmu-\upnu)\bfe,\xi-\upnu t\bfe)=0.\]

		The proof has five main steps.

		\smallskip
		\noindent\emph{Step 1: a background wave that activates the echo chain.}
		We first use the nonlinear wave operator from \cite{waveop} to construct a small global background of the form
		\[
		\udl f(t,z,v)=\eps_0\vf(v)\cos\fz+f_{*,\infty}(z,v)+f_{*,err}(t,z,v),
		\]
		where $f_{*,err}$ decays rapidly. The persistent spatial Fourier modes \(\pm1\) interact with the high-frequency perturbation and generate components at neighboring spatial frequencies. These components become resonant at different times, allowing the interaction to repeat and produce an echo chain. The estimates are formulated so that the construction remains valid when the background already contains the outputs of all earlier cascades.

		\smallskip
		\noindent\emph{Step 2: one sharp echo cascade.}
		At time $\mathsf{T}_k$ we insert an exponentially small packet localized near the frequency $(k,\mathsf{A}_k)$. In the forward regime, its density is successively resonant near $t_\upmu=\mathsf{A}_k/(\upmu+1/2)$. The source-free linearized analysis gives an approximate recurrence from mode $\upmu$ to mode $\upmu-1$, with leading factor $\eps_0\mathsf{A}_k/\upmu^3$, and proves matching upper and lower bounds through the complete chain. After the last resonance, the cascade leaves a nontrivial output near spatial mode $0$ and velocity frequency $\mathsf{A}_k$. The backward analysis shows that evolving the packet to the fixed initial time creates no comparable reverse cascade, so its Gevrey-3 size remains exponentially small.

			\smallskip
			\noindent\emph{Step 3: control of higher-order nonlinear echo interactions.}
		We introduce a Taylor expansion of the perturbation around the evolving background. The main difficulty is that higher-order nonlinear interactions can generate new echo chains and undergo further resonant amplification. We show that the combined contributions of these higher-order echoes and the nonlinear remainder remain smaller than the leading first-order echo output, thereby establishing the nonlinear amplification required in Proposition~\ref{propudlf}.

		\smallskip
		\noindent\emph{Step 4: iteration over infinitely many scales.}
		Starting from the background wave, we apply Proposition \ref{propudlf} at a rapidly increasing sequence $(k_j,\mathsf{A}_{k_j})$. Each new cascade is constructed around the nonlinear solution containing all previous ones. The scale separation makes the increments summable in a fixed Sobolev space and prevents their outputs from canceling at the distinguished frequency bands. Consequently, the partial solutions converge to a genuine global solution, the equation passes to the limit, and nonnegativity is preserved. For every $j$ and all $t\geq2\mathsf{A}_{k_j}$, the output of the $j$-th cascade gives
		\[
		\|f(t)\|_{H^{1000}_{z,v}}\gtrsim \mathsf{A}_{k_j}^{0.9}.
		\]
		Since $\mathsf{A}_{k_j}\to\infty$, this proves the asserted infinite-time divergence rather than growth only along a finite cascade.

		\smallskip
		\noindent\emph{Step 5: scattering and damping.}
		Every increment is decomposed into a limiting profile and a time-decaying error. Their summability yields convergence of $f(t)$ to $f_\infty$ in $H^{18}$. The identity $\wh\rho_f(t,\mu)=\wh f(t,\mu,\mu t)$ and a Fourier trace estimate then imply decay of $\rho_f$ and $E_f$. This completes the circle: the same solution exists globally and undergoes macroscopic Landau damping, yet the infinite sequence of echo outputs rules out global stability in a strong topology.

		\medskip

	\subsection{Notation and conventions}\label{s1.5}
	For a function $f(z,v)$ on $\T^3\times\R^3$, we use the Fourier transform
	\[
	\wh f(\mu,\xi)=\f{1}{(2\pi)^3}\int_{\T^3\times\R^3}
	f(z,v)e^{-i\mu\cdot z-i\xi\cdot v}\,dzdv,\qquad
	(\mu,\xi)\in\Z^3\times\R^3,
	\]
	and use the same hat notation for Fourier transforms in fewer variables. We
	also write $\sF$ for the Fourier transform.

	The construction lies in the invariant class depending only on the first
	spatial coordinate. Accordingly, we write
	\[
	z=(\fz,\mathring z)\in\T\times\T^2,\qquad
	v=(\fv,\mrv)\in\R\times\R^2,
	\]
	and decompose the corresponding Fourier variables as
	\[
	\mu=(\upmu,\mathring\mu)\in\Z\times\Z^2,\qquad
	\xi=(\upxi,\mrxi)\in\R\times\R^2.
	\]
	We set $\bfe=(1,0,0)$, use $\alpha$ for a three-dimensional multi-index,
	and sum all norms containing $v^\alpha$ over $|\alpha|\leq2$. We also fix
	\[
	\sss:=1-10^{-10},\qquad \mathtt{N}_i:=10^i.
	\]

	For a spatial frequency $k\in\N$, define the cascade scales and critical times
	by
	\begin{equation}\label{DefTk}
		\mathfrak T_k:=(\log k)^{1000};\quad
	\mathsf A_k:=\f{(1+\eps_0^6)k^3}{\eps_0};\quad
	\mathsf T_k:=\f{\mathsf A_k}{k+\f12};\quad
	t_\upmu:=\f{\mathsf A_k}{\upmu+\f12},\quad \upmu\geq0.
	\end{equation}
	Thus $t_k=\mathsf T_k$ and $t_0=2\mathsf A_k$.

	We use the Gaussian
	\[
	\vf(v):=\f{1}{2\pi^2}e^{-|v|^2/2},
	\]
	and choose nonnegative functions $\psi(\fv)$ and $\phi(\mrv)$ such that, for
	$a=0,1,2$ and $|\beta|\leq2$,
	\begin{equation}\label{Defofpsi}
		\mathrm{supp}\,\psi\subset[-1,1],\qquad \nnr{\psi}_{L^2}=1,
	\qquad \nr{\pa_\upxi^a\wh\psi(\upxi)}\ls e^{-10|\upxi|^\sss},
	\end{equation}
	and
	\begin{equation}\label{Defofmrv}
		\mathrm{supp}\,\phi\subset\{|\mrv|\leq1\},\qquad \nnr{\phi}_{L^2}=1,
	\qquad \nr{\pa_{\mrxi}^\beta\wh{\phi}(\mrxi)}\ls e^{-10|\mrxi|^\sss}.
	\end{equation}

	For profiles in the free-transport coordinates, define
	\begin{equation}\label{DefNErho}
		\cN(f,g):=-E_f(t,z+vt)\cdot(\na_v-t\na_z)g,\qquad
	E_f:=\na_x\Delta_x^{-1}\rho_f,\qquad
	\rho_f(t,x):=\int_{\R^3}f(t,x-vt,v)\,dv.
	\end{equation}
	The underlined quantities $\udl E$ and $\udl\rho$ are defined from $\udl f$
	in the same way.

	If $A=A(t,\mu,\xi)$ is a Fourier multiplier, we set
	\begin{equation}\label{multiplier}
		A(t,\na)f:=\sF^{-1}\rr{A(t,\mu,\xi)\wh f(t,\mu,\xi)},\qquad
	A(t,\na)\rho:=\sF^{-1}\rr{A(t,\mu,\mu t)\wh\rho(t,\mu)}.
	\end{equation}
	The notations $A(t,\na_{\fz,\fv})$ and $A(t,\na_{\mrv})$ denote the
	corresponding multipliers in the longitudinal variables and transverse
	velocity variables. Finally,
	\[
	\<x\>:=(1+|x|^2)^{1/2},\qquad
	|\mu,\eta|:=(|\mu|^2+|\eta|^2)^{1/2},\qquad
	\<\mu,\eta\>:=\<|\mu,\eta|\>.
	\]
	We write $f\ls g$ if $f\leq Cg$ for an absolute constant, and $C_\theta$
	denotes a constant depending only on $\theta$ whose value may change from line
	to line.

	\subsection{The single-cascade bootstrap}
	To construct a global solution containing infinitely many echo cascades, we need to successively add high-frequency perturbations to a time-dependent, spatially inhomogeneous background solution. At each step, the perturbation must undergo a complete echo cascade and produce a sufficiently large high-frequency output. At the same time, its evolution must satisfy suitable regularity and decay estimates so that the resulting nonlinear solution can serve as the background for the next cascade. Proposition~\ref{propudlf} provides the key step in this iterative construction.
	
	More precisely, the background contains persistent spatial Fourier modes at frequencies \(\pm1\) and is controlled in a rescaled Gevrey norm with frequency scale \(A_*\). At the reference time \(\mathsf T_k\), we introduce an exponentially small perturbation concentrated near the frequency \(\pm(k,\mathsf A_k)\), where \(\mathsf A_k\gg A_*\). The persistent low-frequency modes of the background couple neighboring spatial frequencies of the perturbation, producing successive resonant transfers from mode \(\upmu\) to \(\upmu-1\) near the critical times \(t_\upmu\). The resulting echo cascade generates an output near spatial frequency zero and velocity frequency \(\mathsf A_k\), which provides a lower bound in a high Sobolev norm. Backward evolution produces no comparable amplification, so the perturbation remains exponentially small at the fixed initial time. Moreover, in a rescaled Gevrey norm with frequency scale \(\mathsf A_k\), the perturbation scatters to a limiting profile, with an error that decays in its rescaled Gevrey norm. Its spatial density also decays as \(t\to\infty\). These properties allow the cascade to be incorporated into the global nonlinear construction without losing control of its long-time evolution.

	\begin{prop}\label{propudlf}
		Let  $\udl{f}(t,z,v)$, with $(t,z,v)\in[10,\infty)\times\T^3\times\R^3, z=(\fz,\mathring{z}), v=(\fv,\mathring{v})$,  be a real-valued solution to 
		\[\pa_t\udl{f}=\cN(\udl{f},f^0+\udl{f}),\] satisfying that 
		\begin{itemize}
			\item $\udl{f}(t,z,v)=\eps_0\vf(v)\cos \fz+f_{*}(t,z,v)$, $f_{*}(t,z,v)=f_{*,\infty}(z,v)+f_{*,err}(t,z,v)$, where $f_{*,\infty}(z,v)=f_{*,\infty}(\fz\bfe,v)$, $f_{*,err}(t,z,v)=f_{*,err}(t,\fz\bfe,v)$. 
			\item For $\sss:=1-10^{-10}$, there exist constants \[ \f{1}{\eps_*}\geq\exp\exp\rr{\f{1}{\eps_0}}, A_*\geq\exp\exp\rr{\f{1}{\eps_0}},\] such that
		\begin{equation}\label{udlf}
			\begin{aligned}
				\nnr{\<\na_{\fz,\fv}\>^{20}e^{\<\f{\na_{\fz,\fv}}{A_*}\>^\sss+\<\na_{\mathring{v}}\>^\sss}\rr{v^\alpha f_{*,\infty}} }_{\lta L^2_{z,v}}&\leq\eps_0\eps_{*}, \\
				\nnr{\<\na_{\fz,\fv}\>^{20}e^{\<\f{\na_{\fz,\fv}}{A_*}\>^\sss+\<\na_{\mathring{v}}\>^\sss}\rr{v^\alpha f_{*,err}(t)} }_{\lta L^2_{z,v}}&\leq C_\theta\eps_0e^{-\<\f{t}{A_*}\>^\sss}.
			\end{aligned}
		\end{equation}
	\end{itemize}
		Choose  $k\in\Z$ satisfying $k\geq\max\{\exp\exp(\f{1}{\eps_*}),\exp\exp(A_*)\}$ and define $\mathsf{A}_k=\f{(1+\eps_0^6)k^3}{\eps_0}$ and  $\mathsf{T}_k=\f{\mathsf{A}_k}{k+\f12}$. Then there exists $\eps_k\ll1$  with  $e^{-3.1k}<\eps_k<e^{-3k}$ such that the initial-value problem
		\begin{equation}\label{perturbationeq}
			\left\{\begin{aligned}
				&\pa_tf=\cN(f,f^0+\udl{f})+\cN(\udl{f},f)+\cN(f,f), \\
				&f|_{t=\mathsf{T}_k}=\eps_k\psi(k^\f13\fv)\phi(\mrv)\cos(k\fz+\mathsf{A}_k\fv),
			\end{aligned}\right.
		\end{equation}
		admits a  solution $f=f(t,z,v)$ for $t\ge10$ satisfying the following properties:
		\begin{equation}\label{finfty8}
			f(t,z,v)=f_{\infty}(z,v)+f_{err}(t,z,v),\quad \nnr{\<\na_{\fz,\fv}\>^{20}e^{\<\f{\na_{\fz,\fv}}{\mathsf{A}_k}\>^\sss+\<\na_{\mrv}\>^\sss}\rr{v^\alpha f_{\infty}} }_{\lta L^2_{z,v}}\leq \mathsf{A}_k^{-\mathtt{N}_3+28},
		\end{equation}
		\begin{equation}\label{fperturb}
			\nnr{\<\na_{\fz,\fv}\>^{20}e^{\<\f{\na_{\fz,\fv}}{\mathsf{A}_k}\>^\sss+\<\na_{\mrv}\>^\sss}\rr{v^\alpha f_{err}(t)} }_{\lta L^2_{z,v}}\leq \mathsf{A}_k^{-\mathtt{N}_3+28}e^{-\<\f{t}{\mathsf{A}_k}\>^\sss}.
		\end{equation}
		\begin{equation}\label{initialdata}
			\nnr{\<\na\>^3e^{1.1\<\eps_0\na\>^\f13}\rr{v^\alpha f(t)}}_{\lta L^2_{z,v}}\leq e^{-\f14k},\quad\text{for any}\ t\in[10,\mathsf{T}_k],
		\end{equation}
		\begin{equation}\label{lowerbound}
			\rr{\int_{0.9\mathsf{A}_k\leq|\upxi|\leq1.1\mathsf{A}_k}\nr{\wh{f}(t,0,\xi)}^2d\xi}^\f12\geq \mathsf{A}_k^{-\mathtt{N}_3+0.9},\quad\text{for any}\ t\geq2\mathsf{A}_k.
		\end{equation}
		Here $\xi=(\upxi,\mathring{\xi})\in\R^1\times\R^2$.
	\end{prop}
	\begin{rmk}
		From \eqref{udlf} we get
		\begin{equation}\label{udlf1}
			\nnr{\<\na_{\fz,\fv}\>^{20}e^{\<\f{\na_{\fz,\fv}}{A_*}\>^\sss+\<\na_{\mathring{v}}\>^\sss}\rr{v^\alpha \udl{f} }}_{\lta L^2_{z,v}}\ls\eps_0,
		\end{equation}
		\begin{equation}\label{udlrho}
			\begin{aligned}
				&\nnr{\<\na_{\fz,\fv}\>^{10}e^{\f12\<\f{\na_{\fz,\fv}}{A_*}\>^\sss}\udl{\rho}(t)}_{L^2_x}^2\ls\sum_{\upmu\neq0}\nr{\<\upmu,\upmu t\>^{10}e^{\f12\<\f{|\upmu,\upmu t|}{A_*}\>^\sss}\wh{\udl{f}}(t,\upmu\bfe,\upmu t\bfe)}^2 \\
				& \ls\<t\>^{-20}e^{-\<\f{t}{A_*}\>^\sss} \sum_{\upmu\neq0} \sup_{\upxi}\nr{\<\upmu,\upxi\>^{20}e^{\<\f{|\upmu,\upxi|}{A_*}\>^\sss}\wh{\udl{f}}(t,\upmu\bfe,\upxi)}^2\ls\eps^2_0 \<t\>^{-20}e^{-\<\f{t}{A_*}\>^\sss}.
			\end{aligned}
		\end{equation}
		We recall from \eqref{multiplier} the definition of the multiplier used for $\rho$.
	\end{rmk}
	
	\medskip

	\subsection{Taylor hierarchy for the nonlinear closure}
The purpose of this subsection is to explain the nonlinear decomposition
used to prove Proposition~\ref{propudlf}. We consider the perturbation
equation \eqref{perturbationeq} around the evolving background $\udl f$.
The leading linearized dynamics produces an echo cascade, but nonlinear
self-interactions of the inserted high-frequency perturbation generate
additional components. To obtain a lower bound for the full nonlinear
solution, we must show that these contributions do not cancel the leading
echo output. We use a Taylor expansion in its amplitude, following an idea employed by Deng and Masmoudi \cite{DengMasmoudi} in their study of the long-time instability of the Couette flow. We organize the perturbation by nonlinear interaction order
and then refine the decomposition according to its Fourier frequencies.

Introduce the linearized operator around the evolving background,
\begin{equation}\label{linearizedOP}
	\cL f:=\cN(f,f^0+\udl f)+\cN(\udl f,f).
\end{equation}
The perturbation equation can therefore be written as
\[
\pa_t f=\cL f+\cN(f,f).
\]
We formally expand
\[
f=\sum_{n=1}^\infty f^{(n)},
\]
where $f^{(n)}$ collects terms of order $n$ in the amplitude $\eps_k$
of the inserted packet. Since $\cL$ already includes interactions with the
time-dependent background, the first-order term captures the echo cascade
generated by those interactions:
\begin{equation}\label{fone}
	\pa_tf^{(1)}=\cL f^{(1)},\qquad
	f^{(1)}(\mathsf T_k)=\eps_k\psi(k^\f13\fv)\phi(\mrv)
	\cos(k\fz+\mathsf A_k\fv).
\end{equation}
For $n\geq2$, the higher-order terms are generated by nonlinear
interactions of the lower-order ones and satisfy
\begin{equation}\label{fn}
	\pa_tf^{(n)}=\cL f^{(n)}+
	\sum_{n_1+n_2=n}\cN(f^{(n_1)},f^{(n_2)}),\qquad
		f^{(n)}(\mathsf T_k)=0.
\end{equation}
This recursion is triangular: once the terms of orders below $n$ are
known, the source in the equation for $f^{(n)}$ is determined.

The interaction order \(n\) alone is not sufficient to control the higher-order terms. Indeed, nonlinear interactions between lower-order components can generate new echo chains, which may undergo further resonant amplification. To analyze these interactions, we also need to track the frequencies generated at each order. The initial cosine contains two Fourier components centered at \(\pm(k,\mathsf A_k)\). Their nonlinear interactions generate components with different combinations of these two phases. We therefore refine the Taylor expansion by introducing a second index \(q\):
\begin{equation}\label{Deffn}
	f^{(n)}=\sum_{\substack{|q|\leq n\\2\mid n-q}}f^{(n,q)}.
\end{equation}
The index $q$ records this signed difference; the component $f^{(n,q)}$
is localized near $(qk,q\mathsf A_k)$, with the localization propagated
under the linearized evolution. The conditions $|q|\leq n$ and
$n-q\in2\Z$ follow from counting the signs of $n$ factors of the
initial packet. At first order,
\begin{equation}\label{fonemp}
	\partial_t f^{(1,\pm1)}=\cL f^{(1,\pm1)},\qquad
	f^{(1,\pm1)}(\mathsf T_k)=\f{\eps_k}{2}\psi(k^\f13\fv)\phi(\mrv)
	e^{\pm i(k\fz+\mathsf A_k\fv)}.
\end{equation}
For $n\geq2$, we define the components recursively by
\begin{equation}\label{fnq}
	\partial_t f^{(n,q)}=\cL f^{(n,q)}+
	\sum_{\substack{n_1+n_2=n,\ q_1+q_2=q\\
	2\mid(q_1-n_1),\ 2\mid(q_2-n_2)}}
	\cN(f^{(n_1,q_1)},f^{(n_2,q_2)}),\qquad
		f^{(n,q)}(\mathsf T_k)=0.
\end{equation}
Since the solution is real-valued,
\begin{equation}\label{qmq}
	f^{(n,-q)}=\ov{f^{(n,q)}}.
\end{equation}
The weighted energy estimates provide Fourier localization for \(f^{(n,q)}\) and control the nonresonant contributions. In the resonant regions, the refined Duhamel estimates control the additional amplification generated by the nonlinear source terms. Combining these bounds, we show in Section 5 that the higher-order contributions remain smaller than the first-order echo output.

To justify this expansion for the exact nonlinear solution, truncate it
at $M_{\rm Tay}=\mathsf A_k^{\mathtt{N}_1}$ (denoted by $M$ in Section~5).
Writing
\[
f^\#=\sum_{n=1}^{M_{\rm Tay}}f^{(n)},\qquad \wt f=f-f^\#,
\]
the equation for $\wt f$ has forcing terms involving at least
$M_{\rm Tay}+1$ factors of the inserted packet amplitude. The smallness
of $\eps_k$, together with energy estimates in the different time regimes,
controls this remainder. In particular, neither higher-order echo chains
nor the nonlinear remainder can destroy the first-order lower bound.

\subsection{Organization of the paper}
Sections 2--5 are devoted to the proof of Proposition~\ref{propudlf}, which constructs a single nonlinear echo cascade around an evolving background solution. Section 2 establishes weighted energy estimates that provide Fourier localization and control the nonresonant contributions. Section 3 uses the Duhamel principle to analyze the resonant interactions in the homogeneous linearized equation, proving matching upper and lower bounds for the echo amplification. Section 4 establishes refined estimates for the inhomogeneous linearized equation, controlling the additional echo growth generated by source terms.

Section 5 combines these estimates to control the higher-order echo chains generated by nonlinear interactions and the nonlinear remainder. This shows that the leading echo amplification survives in the full nonlinear solution and completes the proof of Proposition~\ref{propudlf}.

Finally, Section 6 applies Proposition~\ref{propudlf} recursively, inserting high-frequency perturbations at increasingly separated frequency scales. The resulting sequence converges to a global solution containing infinitely many echo cascades, whose persistent high-frequency outputs yield the instability asserted in Theorem~\ref{T1}.

\section{Weighted estimates for the forced linearized equation}

	This section develops the weighted estimates for the forced linearized equation
    \begin{equation}\label{lineaeq}
    	\left\{\begin{aligned}
    		&\pa_tf=\cL f+Q(t,z,v),\quad z=(\fz,\mathring{z})\in \T^3, \\
    		&\int_{\T^3\times\R^3} f(t,z,v)dzdv=0, \\
    		&f|_{t=\mathsf{T}_k}=f(\mathsf{T}_k,\fz\bfe,v), \,Q(t,z,v)=Q(t,\fz\bfe,v).
    	\end{aligned}\right.
    \end{equation}
	We establish weighted energy estimates for the linearized equation with source terms. These estimates control the size and Fourier localization of the solutions and will be used in Sections 3--5 to bound the nonresonant contributions. We first derive the Fourier--Volterra formulation of \eqref{lineaeq} and prove the forward and backward kernel bounds needed for the Schur estimates. We then establish weighted estimates on the long-time interval \([\mathfrak{T}_k,\infty)\) and the short-time interval \([10,\mathfrak{T}_k]\), which provide the starting point for the refined analysis of resonant interactions in the subsequent sections.

	Throughout Sections 2--4 we use the following conventions.
	\begin{itemize}
		\item We write $z=(\fz,\mathring{z})\in\T\times\T^2$ and $v=(\fv,\mrv)\in\R\times\R^2$. The operators $\na_{\fz,\fv}$, $\na_{\mrv}$, and $\na_{\fz,v}$ act on the longitudinal variables $(\fz,\fv)$, the transverse velocity variable $\mrv$, and all reduced variables $(\fz,v)$, respectively.
		\item $\wh f$ denotes the full Fourier transform, while $\sF_z$ and $\sF_\fv$ denote the partial Fourier transforms in $z$ and $\fv$, respectively. When a hat is applied to a function of $\fv$ alone, it denotes its one-dimensional Fourier transform.
		\item All norms containing velocity moments are summed over $|\alpha|\leq2$; we write this sum uniformly as $\ell^2(|\alpha|\leq2)$, abbreviated by $\ltaa$. The symbol $\alpha_1$ denotes the first component of the multi-index $\alpha$.
		\item Density norms are taken in the physical variable $x$. The critical times $t_\upmu$, including $t_0$, are those defined in \eqref{DefTk}.
	\end{itemize}
	\begin{rmk}
		The zero-mass condition $\int_{\T^3\times\R^3}f(t,z,v)\,dzdv=0$ is automatically propagated by \eqref{fone}, \eqref{fn}, \eqref{fonemp}, and \eqref{fnq}, owing to the divergence structure of $\mathcal{N}$ in \eqref{DefNErho}.
	\end{rmk}
	\begin{rmk}
		The nearby decimal exponents and constants used below record the finite losses between successive estimates. Their precise numerical values are immaterial; only the strict hierarchy between them is used. We keep the numerical choices to make every loss explicit.
	\end{rmk}

	\subsection{Fourier--Volterra formulation}
	Whenever a spatial Fourier index appears in a denominator, the corresponding sum is restricted to nonzero modes. The zero spatial mode of the electric field is set to zero, and the mass condition in \eqref{lineaeq} gives $\wh{\rho}(t,0)=0$.
	For the nonzero spatial modes, applying the Fourier transform to \eqref{lineaeq} gives
		\begin{equation}\label{whfeq}
			\begin{aligned}
				&\pa_t\wh{f}(t,\upmu\bfe,\xi)+\wh{\rho}(t,\upmu\bfe)\f{\upxi-\upmu t}{\upmu}\wh{f^0}(\xi-\upmu t\bfe) \\
				&\quad +\sum_{\upnu\in\Z\backslash\{0\}}\wh{\rho}(t,\upnu\bfe)\f{\upxi-\upmu t}{\upnu}\wh{\udl{f}}(t,(\upmu-\upnu)\bfe,\xi-\upnu t\bfe) \\
				&\quad +\sum_{\upnu\in\Z\backslash\{0\}}\wh{\udl{\rho}}(t,\upnu\bfe)\f{\upxi-\upmu t}{\upnu}\wh{f}(t,(\upmu-\upnu)\bfe,\xi-\upnu t\bfe)=\wh{Q}(t,\upmu\bfe,\xi).
			\end{aligned}
		\end{equation}	
		Duhamel's principle yields the following integral equation for the density:
		\begin{equation}\label{rhohateq}
			\begin{aligned}
				\wh{\rho}(t,\upmu\bfe)&=\wh{f}(\mathsf{T}_k,\upmu\bfe,\upmu t\bfe)+\int_{\mathsf{T}_k}^t\wh{Q}(\tau,\upmu\bfe,\upmu t\bfe)d\tau \\
				&\quad -\int_{\mathsf{T}_k}^t\wh{\rho}(\tau,\upmu\bfe)(t-\tau)\wh{f^0}(\upmu(t-\tau)\bfe)d\tau \\
				&\quad -\sum_{\upnu\neq0}\int_{\mathsf{T}_k}^t\wh{\rho}(\tau,\upnu\bfe)\f{\upmu(t-\tau)}{\upnu}\wh{\udl{f}}(\tau,(\upmu-\upnu)\bfe,(\upmu t-\upnu\tau)\bfe)d\tau \\
				&\quad -\sum_{\upnu\neq0}\int_{\mathsf{T}_k}^t\wh{\udl{\rho}}(\tau,\upnu\bfe)\f{\upmu(t-\tau)}{\upnu}\wh{f}(\tau,(\upmu-\upnu)\bfe,(\upmu t-\upnu\tau)\bfe)d\tau.
			\end{aligned}
		\end{equation}
		For $\upmu\neq0$, let $R^0(t,\upmu\bfe)$ solve
		\begin{equation}\label{Rz}
			R^0(t,\upmu\bfe)=t\wh{f^0}(\upmu t\bfe)\vv{1}_{t>0}-\int_0^tR^0(\tau,\upmu\bfe)(t-\tau)\wh{f^0}(\upmu(t-\tau)\bfe)d\tau.
		\end{equation}
		The standard resolvent estimate in \cite{GNR,waveop} gives
		\begin{equation}\label{eRz}
			\nr{R^0(t,\upmu\bfe)}\ls_\theta\f{1}{|\upmu|+1}e^{-\f12|\upmu t|}.
		\end{equation}
		Consequently, the density equation can be written as
		\begin{equation}\label{eqrho}
			\wh{\rho}(t,\upmu\bfe)=H(t,\upmu\bfe)-\int_{\mathsf{T}_k}^tH(\tau,\upmu\bfe)R^0(t-\tau,\upmu\bfe)d\tau,
		\end{equation}
		where
		\begin{equation}\label{DefH}
			\begin{aligned}
				H(t,\upmu\bfe)&=\wh{f}(\mathsf{T}_k,\upmu\bfe,\upmu t\bfe)+\int_{\mathsf{T}_k}^t\wh{Q}(\tau,\upmu\bfe,\upmu t\bfe)d\tau \\
				&\quad -\sum_{\upnu\neq0}\int_{\mathsf{T}_k}^t\wh{\rho}(\tau,\upnu\bfe)\f{\upmu(t-\tau)}{\upnu}\wh{\udl{f}}(\tau,(\upmu-\upnu)\bfe,(\upmu t-\upnu\tau)\bfe)d\tau \\
				&\quad -\sum_{\upnu\neq0}\int_{\mathsf{T}_k}^t\wh{\udl{\rho}}(\tau,\upnu\bfe)\f{\upmu(t-\tau)}{\upnu}\wh{f}(\tau,(\upmu-\upnu)\bfe,(\upmu t-\upnu\tau)\bfe)d\tau \\
				&=\sum_{i=1}^4\mathscr{J}_i.
			\end{aligned}
		\end{equation}
	This representation separates the stable Volterra response generated by $f^0$ from the two perturbative interactions with the background solution $\underline f$. The next lemma controls the time--frequency kernels of these interactions in both time directions.

	\subsection{Kernel bounds for forward and backward evolution}
	We first record the kernel bounds used in the Schur estimates below.
	\begin{lem}\label{auxi}  Let $k,  \mathsf{T}_k, \mathfrak{T}_k$  be the parameters defined in \eqref{DefTk} and Proposition \ref{propudlf}.  Then  
		\[\ba&\sup_{\substack{\upmu\in\Z\backslash\{0\}\\t\geq \mathsf{T}_k} }\sum_{\upnu\in\Z\backslash\{0\}}\int_{\mathsf{T}_k}^t\f{\tau\vv{1}_{\upmu\neq\upnu}+1}{|\upnu|}\<\upmu-\upnu,\upmu t-\upnu\tau\>^{-6}e^{\rr{\rr{\f{\mathsf{T}_k}{t}}^{0.01}-\rr{\f{\mathsf{T}_k}{\tau}}^{0.01}}\<\upmu,\upmu t\>^{0.34}}\f{\<t\>^3}{\<\tau\>^3}d\tau\ls1,\\
		&\sup_{\substack{\upnu\in\Z\backslash\{0\}\\\tau\geq \mathsf{T}_k} }\sum_{\upmu\in\Z\backslash\{0\}}\int_\tau^\infty\f{\tau\vv{1}_{\upmu\neq\upnu}+1}{|\upnu|}\<\upmu-\upnu,\upmu t-\upnu\tau\>^{-6}e^{\rr{\rr{\f{\mathsf{T}_k}{t}}^{0.01}-\rr{\f{\mathsf{T}_k}{\tau}}^{0.01}}\<\upmu,\upmu t\>^{0.34}}\f{\<t\>^3}{\<\tau\>^3}dt\ls1,\\
		&\sup_{\substack{\upmu\in\Z\backslash\{0\}\\\mathfrak{T}_k\leq t\leq \mathsf{T}_k} }\sum_{\upnu\in\Z\backslash\{0\}}\int_t^{\mathsf{T}_k}\f{\tau\vv{1}_{\upmu\neq\upnu}+1}{|\upnu|}\<\upmu-\upnu,\upmu t-\upnu\tau\>^{-6}e^{\f{\log t-\log \tau}{\log \mathsf{T}_k}\<\upmu,\upmu t\>^{0.34}}d\tau\ls1,\\
		&\sup_{\substack{\upnu\in\Z\backslash\{0\}\\\mathfrak{T}_k\leq\tau\leq \mathsf{T}_k} }\sum_{\upmu\in\Z\backslash\{0\}}\int_{\mathfrak{T}_k}^\tau\f{\tau\vv{1}_{\upmu\neq\upnu}+1}{|\upnu|}\<\upmu-\upnu,\upmu t-\upnu\tau\>^{-6}e^{\f{\log t-\log \tau}{\log \mathsf{T}_k}\<\upmu,\upmu t\>^{0.34}}dt\ls1.	
		\ea\]
	\end{lem}
	\begin{proof}
		We first prove the first inequality. Without loss of generality we only consider the $\upmu>0$ case; otherwise simply replace $(\upmu,\upnu)$ by $(-\upmu,-\upnu)$. We begin with the summation over $\{\upnu\neq\upmu,0\}$. The diagonal contribution $\upnu=\upmu$ is bounded by
		\[\sup_{\upmu\neq0}\sup_{t\geq\mathsf{T}_k}\int_{\mathsf{T}_k}^t\f1{|\upmu|}\<\upmu(t-\tau)\>^{-6}
		e^{\rr{\rr{\f{\mathsf{T}_k}{t}}^{0.01}-\rr{\f{\mathsf{T}_k}{\tau}}^{0.01}}\<\upmu,\upmu t\>^{0.34}}\f{\<t\>^3}{\<\tau\>^3}d\tau\lesssim1,\]
		using the monotonicity of the exponential weight and a direct split into $\tau\leq t/2$ and $\tau\geq t/2$.
		
		For the off-diagonal contribution, split the integration domain into $[\mathsf{T}_k,t]\cap\{\tau\leq t/2\}$ and $[\mathsf{T}_k,t]\cap\{\tau\geq t/2\}$.

		\noindent\textit{Region I: $[\mathsf{T}_k,t]\cap\{\tau\leq t/2\}$.}
		In this region,
		\[\rr{-\rr{\f{\mathsf{T}_k}{t}}^{0.01}+\rr{\f{\mathsf{T}_k}{\tau}}^{0.01}}\<\upmu,\upmu t\>^{0.34}\geq\f1{t^{0.01}}\<\upmu,\upmu t\>^{0.34}\geq\<\upmu,\upmu t\>^{0.33},\] which implies that
		\[\ba
		&\sum_{\upnu\neq\upmu,0}\int_{[\mathsf{T}_k,t]\cap \{\tau\leq\f t2\}}\f{\tau}{|\upnu|}\<\upmu-\upnu,\upmu t-\upnu\tau\>^{-6}e^{\rr{\rr{\f{\mathsf{T}_k}{t}}^{0.01}-\rr{\f{\mathsf{T}_k}{\tau}}^{0.01}}\<\upmu,\upmu t\>^{0.34}}\f{\<t\>^3}{\<\tau\>^3}d\tau\\
		\leq&\sum_{\upnu\neq\upmu,0}\int_{[\mathsf{T}_k,t]\cap \{\tau\leq\f t2\}}\<\upmu-\upnu,\upmu t-\upnu\tau\>^{-6}e^{-\<\upmu,\upmu t\>^{0.33}}\<t\>^4d\tau\ls\sum_{\upnu\neq\upmu,0}\int_{[\mathsf{T}_k,t]\cap \{\tau\leq\f t2\}}\<\upmu-\upnu,\upmu t-\upnu\tau\>^{-6}d\tau\ls1.\ea\]

        \noindent\textit{Region II: $[\mathsf{T}_k,t]\cap\{\tau\ge t/2\}$.}
		Here $\langle t\rangle/\langle\tau\rangle\lesssim1$. We further decompose
		\[\ba\sum_{\upnu\neq\upmu,0}\int_{[\mathsf{T}_k,t]\cap \{\tau\geq\f t2\}}\f{\tau}{|\upnu|}\<\upmu-\upnu,\upmu t-\upnu\tau\>^{-6}e^{\rr{\rr{\f{\mathsf{T}_k}{t}}^{0.01}-\rr{\f{\mathsf{T}_k}{\tau}}^{0.01}}\<\upmu,\upmu t\>^{0.34}}\f{\<t\>^3}{\<\tau\>^3}d\tau
		=\mathcal{I}_{\{|\upmu t-\upnu\tau|\geq\f t2\}}+\mathcal{I}_{\{|\upmu t-\upnu\tau|\leq\f t2\}}, 
		\ea\]
		where  $\mathcal{I}_{\{|\upmu t-\upnu\tau|\geq\f t2\}}$ and $\mathcal{I}_{\{|\upmu t-\upnu\tau|\leq\f t2\}}$ denote the term in the left-hand side with the additional restriction on the integration domain $\{|\upmu t-\upnu\tau|\geq\f t2\}$ and $\{|\upmu t-\upnu\tau|\leq\f t2\}$ respectively.
		On the off-resonant region, the polynomial kernel directly gives
		\[\mathcal{I}_{\{|\upmu t-\upnu\tau|\geq\f t2\}}\ls\sum_\upnu\int_{\mathsf{T}_k}^t\<\upmu-\upnu,\upmu t-\upnu\tau\>^{-5}d\tau\ls1,\]
		so it remains to estimate $\mathcal{I}_{\{|\upmu t-\upnu\tau|\leq t/2\}}$. The constraints $|\upmu t-\upnu\tau|\leq t/2$, $\tau\geq t/2$, $\upmu>0$, and $\upnu\neq0,\upmu$ imply that $0<\upmu<\upnu<4\upmu$ and $\frac{\upmu-1/2}{\upnu}t\leq\tau\leq\frac{\upmu+1/2}{\upnu}t$. Hence
		\[\ba &\rr{-\rr{\f{\mathsf{T}_k}{t}}^{0.01}+\rr{\f{\mathsf{T}_k}{\tau}}^{0.01}}\<\upmu,\upmu t\>^{0.34}\ge \bigg(\bigg(\f{2\upnu}{2\upmu+1}\bigg)^{0.01}-1\bigg) \rr{\f{\mathsf{T}_k}{t}}^{0.01}\<\upmu,\upmu t\>^{0.34}\\&\gtrsim \f{\upnu-\upmu}{\upmu}\<\upmu,\upmu t\>^{0.34}\f{1}{t^{0.01}}\gtrsim\rr{\upnu-\upmu}\rr{\f{t}{\upmu^2}}^{0.33}.\ea\]
	This implies that  
		\[\ba\mathcal{I}_{\{|\upmu t-\upnu\tau|\leq\f t2\}}\leq&\sum_{\upmu<\upnu<4\upmu}\int_{\mathsf{T}_k}^t\f{\tau}{\upnu}\<\upmu-\upnu,\upmu t-\upnu\tau\>^{-6} \ e^{-\rr{\upnu-\upmu}\rr{\f{t}{\upmu^2}}^{0.33}}d\tau\\& \ls\sum_{\upmu<\upnu<4\upmu}\f{t}{\upmu^2}\<\upnu-\upmu\>^{-4.5}e^{-\rr{\upnu-\upmu}\rr{\f{t}{\upmu^2}}^{0.33}}\ls1.
		\ea\]
        And from the proof we also get 
		\begin{equation}\label{canc}
			e^{\rr{\rr{\f{\mathsf{T}_k}{t}}^{0.01}-\rr{\f{\mathsf{T}_k}{\tau}}^{0.01}}\<\upmu,\upmu t\>^{0.34}}\f{\<t\>^3}{\<\tau\>^3}\mathrm{1}_{t\geq\tau}\ls1.
		\end{equation}

		We next prove the third inequality. Again, it suffices to consider $\upmu>0$. We first restrict to $\{\upnu\neq\upmu,0\}$; the diagonal contribution follows directly from $\langle\upmu(t-\tau)\rangle^{-6}$. Split the integration domain into $[t,\mathsf{T}_k]\cap\{|\upmu t-\upnu\tau|\geq\tau/2\}$ and $[t,\mathsf{T}_k]\cap\{|\upmu t-\upnu\tau|\leq\tau/2\}$. The polynomial kernel controls the first region. On the second region, $0<\upnu<\upmu$ and $\frac{\upmu t}{\upnu+1/2}\leq\tau\leq\frac{\upmu t}{\upnu-1/2}$. Since $\log x\geq(x-1)/x$ for $x\geq1$, the assumptions $t\geq\mathfrak{T}_k=(\log k)^{1000}$ and $\log\mathsf{T}_k\sim\log k$ imply, for sufficiently large $k$, that
		\[\f{\log\tau-\log t}{\log \mathsf{T}_k}\<\upmu,\upmu t\>^{0.34}\ge\f{\tau-t}{\tau}\f{\<\upmu,\upmu t\>^{0.34}}{\log \mathsf{T}_k}\ge\f{\upmu-\upnu-\f12}{\upmu}\f{\<\upmu,\upmu t\>^{0.34}}{\log \mathsf{T}_k}\ge \f{\upmu-\upnu}{\upmu} \<\upmu,\upmu t\>^\f13.\]
		This yields 
		\[\sum_{\upnu\neq\upmu,0}\int_{[t,\mathsf{T}_k]\cap\{|\upmu t-\upnu\tau|\leq\f \tau2\}}\f{\tau}{|\upnu|}\<\upmu-\upnu,\upmu t-\upnu\tau\>^{-6}e^{\f{\log t-\log \tau}{\log \mathsf{T}_k}\<\upmu,\upmu t\>^{0.34}}d\tau\ls\sum_{\upnu\neq\upmu,0}\f{\upmu t}{|\upnu|^3}\<\upmu-\upnu\>^{-4}e^{-\f{\upmu-\upnu}{\upmu}\<\upmu,\upmu t\>^\f13}\ls1.\]
		For the second inequality, fix $(\upnu,\tau)$ and split $t\geq\tau$ into $t\geq2\tau$ and $\tau\leq t\leq2\tau$; in the latter region, split further according to $|\upmu t-\upnu\tau|\gtrless t/2$. The same off-resonant and resonant estimates used above are uniform after summing in $\upmu$. The fourth inequality follows analogously by fixing $(\upnu,\tau)$, integrating over $t\leq\tau$, and using the backward exponential weight. The diagonal terms $\upmu=\upnu$ in both estimates are bounded directly by integrating $\<\upmu(t-\tau)\>^{-6}$. This completes the proof of all four bounds.
	\end{proof}

	\subsection{Weighted long-time estimates}
	We now combine the Volterra representation with Lemma~\ref{auxi} to establish weighted energy estimates that control both the Gevrey regularity of the solution and its Fourier localization near a prescribed frequency center \((k',A')\). The estimates apply to both forward and backward evolution from \(\mathsf T_k\).
	
	\begin{prop}\label{longenergyl} Let $k,A_*,\mathsf{A}_k, \mathsf{T}_k, \mathfrak{T}_k,\sss$ be the parameters defined in \eqref{DefTk} and Proposition \ref{propudlf}. Let $k'\in\Z$, $A'\in\R$, $\mathsf{K}\geq0$, and assume that
		  $0\leq\mathsf{r}\leq1$, $0\le \mathtt{r}\leq A_*^{-1}$ and $|k',A'|\leq \mathsf{A}_k^{\mathtt{N}_1+2}$. All the norms below contain the $\ell^2$ sum over $|\alpha|\leq2$, and the implicit constants depend only on $\theta$.
		  \begin{itemize}
		  	\item[(i).]\textnormal{(Forward evolution, i.e., $t\in[\mathsf{T}_k,\infty)$)} If $f(\mathsf{T}_k)$ and $Q(t)$ satisfy
			\[\ba&\nnr{\<\na_{\fz,\fv}\>^6e^{\<\na_{\fz,\fv}\>^{0.34}}e^{\f{\mathsf{r}}{4}(\f34+\f{1}{\mathsf{T}_k})\<\mathtt{r}\nr{\na_\fz-k',\na_\fv-A'}\>^\sss+\<\na_{\mrv}\>^\sss} \rr{v^\alpha f(\mathsf{T}_k)}}_{\lta L^2_{z,v}}\leq \mathsf{K},\\	   	
			&\int_{\mathsf{T}_k}^t\nnr{\<\na_{\fz,\fv}\>^6e^{(\f{\mathsf{T}_k}{\tau})^{0.01}\<\na_{\fz,\fv}\>^{0.34}}e^{\f{\mathsf{r}}{4}(\f34+\f{1}{\tau})\<\mathtt{r}\nr{\na_\fz-k',\na_\fv-A'}\>^\sss+\<\na_{\mrv}\>^\sss} \rr{v^\alpha Q(\tau)}}_{\lta L^2_{z,v}}d\tau\leq Z(t),\ea\]
			where $Z(t)\geq0$ is an increasing function, then the solution $f$ to \eqref{lineaeq} satisfies
			\[\nnr{\<\na_{\fz,\fv}\>^4e^{(\f{\mathsf{T}_k}{t})^{0.01}\<\na_{\fz,\fv}\>^{0.34}}e^{\f{\mathsf{r}}{4}(\f34+\f{1}{t})\<\mathtt{r}\nr{\na_\fz-k',\na_\fv-A'}\>^\sss+\<\na_{\mrv}\>^\sss} \rr{v^\alpha f(t)}}_{\lta L^2_{z,v}}\ls  \mathsf{K}+Z(t).\]
		  	\item[(ii).]\textnormal{(Backward evolution, i.e., $t\in[\mathfrak{T}_k,  \mathsf{T}_k]$)} If $f(\mathsf{T}_k)$ and $Q(t)$ satisfy
			\[\ba&\nnr{\<\na_{\fz,\fv}\>^6e^{\<\na_{\fz,\fv}\>^{0.34}}e^{\f{\mathsf{r}}{4}(1-\f{1}{\mathsf{T}_k})\<\mathtt{r}\nr{\na_\fz-k',\na_\fv-A'}\>^\sss+\<\na_{\mrv}\>^\sss} \rr{v^\alpha f(\mathsf{T}_k)}}_{\lta L^2_{z,v}}\leq \mathsf{K},\\
			&\int_t^{\mathsf{T}_k}\nnr{\<\na_{\fz,\fv}\>^6e^{\f{\log \tau}{\log \mathsf{T}_k}\<\na_{\fz,\fv}\>^{0.34}}e^{\f{\mathsf{r}}{4}(1-\f{1}{\tau})\<\mathtt{r}\nr{\na_\fz-k',\na_\fv-A'}\>^\sss+\<\na_{\mrv}\>^\sss} \rr{v^\alpha Q(\tau)}}_{\lta L^2_{z,v}}d\tau\leq Z(t),\ea\]
			where $Z(t)\geq0$ is a decreasing function, then the solution $f$ to \eqref{lineaeq} satisfies
			\[\nnr{\<\na_{\fz,\fv}\>^4e^{\f{\log t}{\log \mathsf{T}_k}\<\na_{\fz,\fv}\>^{0.34}}e^{\f{\mathsf{r}}{4}(1-\f{1}{t})\<\mathtt{r}\nr{\na_\fz-k',\na_\fv-A'}\>^\sss+\<\na_{\mrv}\>^\sss} \rr{v^\alpha f(t)}}_{\lta L^2_{z,v}}\ls  \mathsf{K}+Z(t).\]
		  \end{itemize}
		 \end{prop}
		 \begin{rmk}
		 	The restriction on $\mathtt r$ is dictated by the background estimate \eqref{udlf}. The bounds on $k'$ and $A'$ reflect the maximal Taylor order $n\leq\mathsf A_k^{\mathtt{N}_1}$ in \eqref{fn} and ensure that every localization center arising in the expansion is admissible.
		 \end{rmk}
	\begin{proof}
		To prove the propagation of the regularity, we introduce the multiplier as follows:
		\[
		\fA_a(t,\mu,\xi)=\left\{\ba&\<\upmu,\upxi\>^{3+a}e^{(\f{\mathsf{T}_k}{t})^{0.01}\<\upmu,\upxi\>^{0.34}}e^{\f{\mathsf{r}}{4}(\f34+\f{1}{t})\<\mathtt{r}\nr{\upmu-k',\upxi-A'}\>^\sss+\<\mrxi\>^\sss}\quad &t\geq \mathsf{T}_k,\\
		&\<\upmu,\upxi\>^{3+a}e^{\f{\log t}{\log \mathsf{T}_k}\<\upmu,\upxi\>^{0.34}}e^{\f{\mathsf{r}}{4}(1-\f{1}{t})\<\mathtt{r}\nr{\upmu-k',\upxi-A'}\>^\sss+\<\mrxi\>^\sss}\quad &\mathfrak{T}_k\leq t<\mathsf{T}_k.
		\ea\right.
		\]
		We have $\fA_a(t,\mu,\xi)=\fA_a(t,\upmu\bfe,\xi)$ and   set $\sF(\fA_a(t,\na_{\fz,\fv})\rho)(t,\upmu\bfe):=\fA_a(t,\upmu,\upmu t\bfe) \wh{\rho}(t,\upmu\bfe)$ for simplicity. 

		\noindent\textit{Proof of (i).} We divide the argument into three steps.
		
		\noindent\textit{Step 1. Density estimate.}
		Fix $T\ge\mathsf{T}_k$. The Volterra equation \eqref{eqrho}, the resolvent bound \eqref{eRz}, Schur's test, and the inequality
		  \begin{equation}\label{xyaineq}
		  	\<x\>^a\le\<y\>^a+ \<x-y\>^a, \forall a\in[0,1],
		  \end{equation} it holds that
		\[\ba
		&\nnr{\<t\>^3\fA_0(t,\upmu\bfe,\upmu t\bfe)\wh{\rho}(t,\upmu\bfe)}_{L^2_t([\mathsf{T}_k,T];L^2_\upmu)}\leq \nnr{\<\tau\>^3\fA_0(\tau,\upmu\bfe,\upmu \tau\bfe)H(\tau,\upmu\bfe)}_{L^2_\tau([\mathsf{T}_k,T];L^2_\upmu)}\\&\times\left(1+\sup_\upmu\int_0^\infty \nr{\<\upmu t\>^3e^{\f14\<\upmu t\>} R^0(t,\upmu\bfe)}\<t\>^3dt\right)\ls_\theta \nnr{\<\tau\>^3\fA_0(\tau,\upmu\bfe,\upmu \tau\bfe)H(\tau,\upmu\bfe)}_{L^2_\tau([\mathsf{T}_k,T];L^2_\upmu)}. \ea\]
		It remains to estimate $\nnr{\<t\>^3\fA(t,\upmu\bfe,\upmu t\bfe)H(t,\upmu\bfe)}_{L^2_t([\mathsf{T}_k,T];L^2_\upmu)}$. The initial-data contribution in \eqref{DefH} satisfies
		\begin{equation}\label{estimateJ1}
			\begin{aligned}
				&\nnr{\<t\>^3\fA_0(t,\upmu\bfe,\upmu t\bfe)\mathscr{J}_1(t,\upmu\bfe)}^2_{L^2_t([\mathsf{T}_k,T];L^2_\upmu)}=\sum_{\upmu\neq0}\int_{\mathsf{T}_k}^{T}\nr{\fA_3(t,\upmu\bfe,\upmu t\bfe)\wh{f}(\mathsf{T}_k,\upmu\bfe,\upmu t\bfe)}^2dt \\
				& \le\sum_{\upmu\neq0}\int_{-\infty}^{+\infty} \nr{\fA_3(\mathsf{T}_k,\upmu\bfe,\upxi\bfe)\wh{f}(\mathsf{T}_k,\upmu\bfe,\upxi\bfe)}^2d\upxi\le\sum_{\upmu\neq0}\sum_{|\alpha|\le2}\int_{\R^3} \nr{\pa^{\alpha}_\xi \fA_3(\mathsf{T}_k,\upmu\bfe,\xi)\wh{f}(\mathsf{T}_k,\upmu\bfe,\xi)}^2d\xi \\
				&\ls \sum_{|\alpha|\le2}\|\fA_3(\mathsf{T}_k,\na_\fz,\na_v)(v^\alpha f(\mathsf{T}_k))\|^2_{L^2_{z,v}} \leq\mathsf{K}^2.
			\end{aligned}
		\end{equation}
		Here we use the change of variables $\upxi=\upmu t$ in the first inequality and the trace theorem in the second.
		
		The source contribution is treated by the same trace argument and Minkowski's inequality:
		\[\ba
		&\nnr{\<t\>^3\fA(t,\upmu\bfe,\upmu t\bfe)\mathscr{J}_2(t,\upmu\bfe)}^2_{L^2_t([\mathsf{T}_k,T];L^2_\upmu)}\leq\sum_{\upmu\neq0}\int_{\mathsf{T}_k}^{T}\rr{\int_{\mathsf{T}_k}^t\nr{  \fA_3(\tau,\upmu\bfe,\upmu t\bfe)\wh{Q}(\tau,\upmu\bfe,\upmu t\bfe)}d\tau}^2dt\\
		\leq&\rr{\int_{\mathsf{T}_k}^T\rr{\sum_{\upmu\neq0}\int_{\mathsf{T}_k}^{T} \nr{ \fA_3(\tau,\upmu\bfe,\upmu t\bfe)\wh{Q}(\tau,\upmu\bfe,\upmu t\bfe)}^2dt }^\f12d\tau}^2 \leq Z(T)^2.\ea\]
		
		For the term $\mathscr{J}_3$  in \eqref{DefH}, for $A_*^5\ll \mathsf{T}_k\le \tau\le t$,  we first claim 
		 \begin{equation}\label{claimflen1}
		 	\begin{aligned}
		 		\rr{\f{\mathsf{T}_k}{\tau}}^{0.01}\<\upmu,\upmu t\>^{0.34}&\leq\rr{\f{\mathsf{T}_k}{\tau}}^{0.01}\<\upnu,\upnu\tau\>^{0.34}+\f14\<\f{\upmu-\upnu,\upmu t-\upnu\tau}{A_*}\>^\sss, \\
		 		\f{\mathsf{r}}{4}(\f34+\f{1}{t})\<\mathtt{r}\nr{\upmu-k',\upmu t-A'}\>^\sss&\leq\f{\mathsf{r}}{4}(\f34+\f{1}{t})\<\mathtt{r}\nr{\upnu-k',\upnu\tau-A'}\>^\sss+\f14\<\f{\upmu-\upnu,\upmu t-\upnu\tau}{A_*}\>^\sss.
		 	\end{aligned}
		 \end{equation} 
	       To prove the claim, observe that for $a\in[0,1]$ and $x\ge y\ge A_*^5$, $\<x\>^a-\<y\>^a\le a\<y\>^{a-1}|x-y|$.
        By interpolation with \eqref{xyaineq}, we deduce that for $a\in]0,\sss[$,
       \[\<x\>^a-\<y\>^a\le \f{a^{\f{1-\sss}{1-a}}}{\<y\>^{\sss-a}}\<x-y\>^\sss\le \f14\<\f{x-y}{A_*}\>^\sss,\]
		This proves the first inequality because $|\upnu,\upnu\tau|\ge A_*^5$. The second follows from \eqref{xyaineq} and $\mathtt{r}\leq A_*^{-1}$. Thus \eqref{claimflen1} holds. Using this comparison,
  we may derive that
 \[ \fA_0(t,\upmu\bfe,\upmu t\bfe)\le 2\fA_0(\tau,\upnu\bfe,\upnu \tau\bfe)e^{\rr{\rr{\f{\mathsf{T}_k}{t}}^{0.01}-\rr{\f{\mathsf{T}_k}{\tau}}^{0.01}}\<\upmu,\upmu t\>^{0.34}}\<\upmu-\upnu,\upmu t-\upnu \tau\>^3e^{\f12\<\f{\upmu-\upnu,\upmu t-\upnu\tau}{A_*}\>^\sss}. \] From this together with the fact that $|\upmu(t-\tau)|\leq|\upmu t-\upnu\tau|+|\upmu-\upnu|\tau\leq\<\upmu-\upnu,\upmu t-\upnu\tau\>\tau$, we have 
		\[\ba
		&|\<t\>^3\fA_0(t,\upmu\bfe,\upmu t\bfe)\mathscr{J}_3(t,\upmu\bfe)|
		\leq2\sum_{\upnu\neq0}\int_{\mathsf{T}_k}^t\<\tau\>^3\nr{\fA_0(\tau,\upnu\bfe,\upnu \tau\bfe)\wh{\rho}(\tau,\upnu\bfe)}\f{\tau\vv{1}_{\upnu\neq\upmu}+1}{|\upnu|}\\
		&\times\nr{\<\upmu-\upnu,\upmu t-\upnu \tau\>^4e^{\f12\<\f{\upmu-\upnu,\upmu t-\upnu\tau}{A_*}\>^\sss}\wh{\udl{f}}(\tau,(\upmu-\upnu)\bfe,(\upmu t-\upnu\tau)\bfe)}e^{\rr{\rr{\f{\mathsf{T}_k}{t}}^{0.01}-\rr{\f{\mathsf{T}_k}{\tau}}^{0.01}}\<\upmu,\upmu t\>^{0.34}}\f{\<t\>^3}{\<\tau\>^3}d\tau\\
		\leq& \eps_{0}^{0.95}\sum_{\upnu\neq0}\int_{\mathsf{T}_k}^t\<\tau\>^3\nr{\fA_0(\tau,\upnu\bfe,\upnu \tau\bfe)\wh{\rho}(\tau,\upnu\bfe)}\f{\tau\vv{1}_{\upnu\neq\upmu}+1}{|\upnu|}\<\upmu-\upnu,\upmu t-\upnu\tau\>^{-6}e^{\rr{\rr{\f{\mathsf{T}_k}{t}}^{0.01}-\rr{\f{\mathsf{T}_k}{\tau}}^{0.01}}\<\upmu,\upmu t\>^{0.34}}\f{\<t\>^3}{\<\tau\>^3}d\tau.
		\ea\]
		Here we use \eqref{udlf1} for $\udl{f}$ in the last inequality.	Now applying the first two inequalities in  Lemma \ref{auxi}, we get from  Schur  test   that
		\[\nnr{\<t\>^3\fA_0(t,\upmu\bfe,\upmu t\bfe)\mathscr{J}_3(t,\upmu\bfe)}_{L^2_t([\mathsf{T}_k,T];L^2_\upmu)}\leq\eps_{0}^{0.94}\nnr{\<\tau\>^3\fA_0(\tau,\upnu\bfe,\upnu\tau\bfe)\wh{\rho}(\tau,\upnu\bfe)}_{L^2_\tau([\mathsf{T}_k,T];L^2_\upnu)}.\]
		
		Finally we estimate $\mathscr{J}_4$. Noting $\<\upnu,\upnu\tau\>^{0.34}\leq\f14\<\f{\upnu,\upnu\tau}{A_*}\>^{\sss}$ since $|\upnu,\upnu\tau|\geq A_*^5$, by \eqref{xyaineq}, we obtain that
		\begin{equation}\label{claimflen2}
			\begin{aligned}
				\rr{\f{\mathsf{T}_k}{\tau}}^{0.01}\<\upmu,\upmu t\>^{0.34}&\leq\f14\rr{\f{\mathsf{T}_k}{\tau}}^{0.01}\<\f{\upnu,\upnu\tau}{A_*}\>^{\sss}+\rr{\f{\mathsf{T}_k}{\tau}}^{0.01}\<\upmu-\upnu,\upmu t-\upnu\tau\>^{0.34}, \\
				\f{\mathsf{r}}{4}(\f34+\f{1}{t})\<\mathtt{r}\nr{\upmu-k',\upmu t-A'}\>^\sss&\leq\f{1}{4}\<\f{\upnu,\upnu\tau}{A_*}\>^\sss+\f{\mathsf{r}}{4}(\f34+\f{1}{t})\<\mathtt{r}\nr{\upmu-\upnu-k',\upmu t-\upnu\tau-A'}\>^\sss.
			\end{aligned}
		\end{equation} These yield that
		  \[ \fA_0(t,\upmu\bfe,\upmu t\bfe)\le 2\fA_0(\tau,(\upmu-\upnu)\bfe,(\upmu t-\upnu \tau)\bfe)e^{\rr{\rr{\f{\mathsf{T}_k}{t}}^{0.01}-\rr{\f{\mathsf{T}_k}{\tau}}^{0.01}}\<\upmu,\upmu t\>^{0.34}}\< \upnu, \upnu \tau\>^3e^{\f12\<\f{ \upnu, \upnu\tau}{A_*}\>^\sss}. \]
		From this together with \eqref{DefH} and the fact $|\upmu(t-\tau)|\leq \<\upmu-\upnu,\upmu t-\upnu\tau\>\tau$, we have
		\[\ba
		&|\<t\>^3\fA(t,\upmu\bfe,\upmu t\bfe)\mathscr{J}_4(t,\upmu\bfe)|\leq2\sum_{\upnu\neq0}\int_{\mathsf{T}_k}^t\<\tau\>^4
		\nr{\< \upnu, \upnu \tau\>^3e^{\f12\<\f{ \upnu, \upnu\tau}{A_*}\>^\sss}\wh{\udl{\rho}} (\tau,\upnu\bfe)}\nr{ (\fA_1\wh{f})(\tau,(\upmu-\upnu)\bfe,(\upmu t-\upnu\tau)\bfe)}\\
		&\times e^{\rr{\rr{\f{\mathsf{T}_k}{t}}^{0.01}-\rr{\f{\mathsf{T}_k}{\tau}}^{0.01}}\<\upmu,\upmu t\>^{0.34}}\f{\<t\>^3}{\<\tau\>^3}d\tau
		\leq\eps_0^{0.95}\sum_{\upnu\neq0}\int_{\mathsf{T}_k}^t\<\upnu,\upnu\tau\>^{-3}\nr{(\fA_1\wh{f})(\tau,(\upmu-\upnu)\bfe,(\upmu t-\upnu\tau)\bfe)}d\tau.
		\ea\]
		In the last inequality we use \eqref{canc} and \eqref{udlrho}.
		Minkowski's inequality then yields
		\[\ba
		\nnr{\<t\>^3\fA(t,\upmu\bfe,\upmu t\bfe)\mathscr{J}_4(t,\upmu\bfe)}_{L^2_t([\mathsf{T}_k,T];L^2_\upmu)}
		&\leq\eps_0^{0.95}\sum_{\upnu\neq0}\int_{\mathsf{T}_k}^{T}\<\upnu,\upnu\tau\>^{-3}\nnr{ (\fA_1\wh{f})(\tau,(\upmu-\upnu)\bfe,(\upmu t-\upnu\tau)\bfe)}_{L^2_\upmu L^2_t}d\tau\\
		&\leq\eps_0^{0.94}\sup_{\mathsf{T}_k\leq\tau\leq T}\nnr{\fA_1(\tau,\na_{\fz,v})(v^\alpha f(\tau))}_{\lta L^2_{z,v}}.
		\ea\] The last inequality follows from the same change of variables and trace estimate used in \eqref{estimateJ1}.
		
		Patching together all the estimates, we arrive at 
		\begin{equation}\label{lerho}
			\begin{aligned}
				\nnr{\<t\>^3\fA(t,\na_{\fz,\fv})\rho(t)}_{L^2_t([\mathsf{T}_k,T];L^2_x)}&\leq C_\theta\rr{\mathsf{K}+Z(T)} \\
				&\quad +\eps_0^{0.9}\rr{\nnr{\<t\>^3\fA(t,\na_{\fz,\fv})\rho(t)}_{L^2_t([\mathsf{T}_k,T];L^2_x)}+\sup_{\mathsf{T}_k\leq\tau\leq T}\nnr{\fA_1(\tau,\na_{\fz,v})(v^\alpha f(\tau))}_{\lta L^2_{z,v}}}.
			\end{aligned}
		\end{equation}
		
		\noindent\textit{Step 2. Profile energy estimate.}
		The basic weighted energy identity gives
	
		\begin{equation}\label{u1vfestimate1}
			\begin{aligned}
				&\f12\f{d}{dt}\nnr{\fA_1(t,\na_{\fz,v})\rr{v^\alpha f(t)}}_{\lta L^2_{z,v}}^2+0.01\rr{\f{\mathsf{T}_k}{t}}^{0.01}\f1t\nnr{\fA_{1.17}(t,\na_{\fz,v})\rr{v^\alpha f(t)}}_{\lta L^2_{z,v}}^2 \\
				&\qquad\qquad+\f{\mathsf{r}}{4t^2}\nnr{\<\mathtt{r}\nr{\na_\fz-k',\na_\fv-A'}\>^{\f\sss2}\fA_{1}(t,\na_{\fz,v})\rr{v^\alpha f(t)}}_{\lta L^2_{z,v}}^2 =\sum_{i=1}^{6}\mathscr{A}_i.
			\end{aligned}
		\end{equation}
where
		\[\ba&\mathscr{A}_1:=-\sum_{|\alpha|\leq2}\sum_{\upmu\neq0}\int\wh{\rho}(t,\upmu\bfe)\f1\upmu\pa_{\xi}^\alpha\rr{ (\upxi-\upmu t)\wh{f^0}(\xi-\upmu t\bfe)} \fA_1(t,\upmu\bfe,\xi)
		\ov{\sF\rr{\fA_1(t,\na_{\fz,v})(v^\alpha f)}(t,\upmu\bfe,\xi)}d\xi,\\
		&\mathscr{A}_2:=-\sum_{|\alpha|\leq2}\sum_{\substack{\upmu\in\Z\\\upnu\neq0}}\int\wh{\rho}(t,\upnu\bfe)\f{\upxi-\upmu t}\upnu\pa_{\xi}^\alpha\wh{\udl{f}}(t,(\upmu-\upnu)\bfe,\xi-\upnu t\bfe)\fA_1(t,\upmu\bfe,\xi)\ov{\sF\rr{\fA_1(t,\na_{\fz,v})(v^\alpha f)}(t,\upmu\bfe,\xi)}d\xi,\\
		 &\mathscr{A}_3:=-\sum_{|\alpha|\leq2}\sum_{\substack{\upmu\in\Z\\\upnu\neq0}}\int\wh{\udl{\rho}}(t,\upnu\bfe)\f{\upxi-\upmu t}\upnu\pa_{\xi}^\alpha\wh{f}(t,(\upmu-\upnu)\bfe,\xi-\upnu t\bfe)\\&\qquad\qquad\times\rr{\fA_1(t,\upmu\bfe,\xi)-\fA_1(t,(\upmu-\upnu)\bfe,\xi-\upnu t\bfe)}\ov{\sF\rr{\fA_1(t,\na_{\fz,v})(v^\alpha f)}(t,\upmu\bfe,\xi)}d\xi,\\
		&\mathscr{A}_4:=-\sum_{\substack{|\alpha|\leq2\\\alpha_1>0}}C_\alpha\sum_{\substack{\upmu\in\Z\\\upnu\neq0}}\int\wh{\rho}(t,\upnu\bfe)\f{1}\upnu\pa_{\xi}^{\alpha-\bfe}\wh{\udl{f}}(t,(\upmu-\upnu)\bfe,\xi-\upnu t\bfe)\fA_1(t,\upmu\bfe,\xi)\ov{\sF\rr{\fA_1(t,\na_{\fz,v})(v^{\alpha}f)}(t,\upmu\bfe,\xi)}d\xi,\\
	&\mathscr{A}_5:=-\sum_{\substack{|\alpha|\leq2\\\alpha_1>0}}C_\alpha\sum_{\substack{\upmu\in\Z\\\upnu\neq0}}\int\wh{\udl{\rho}}(t,\upnu\bfe)\f1\upnu\pa_{\xi}^{\alpha-\bfe}\wh{f}(t,(\upmu-\upnu)\bfe,\xi-\upnu t\bfe)\fA_1(t,\upmu\bfe,\xi)\ov{\sF\rr{\fA_1(t,\na_{\fz,v})(v^\alpha f)}(t,\upmu\bfe,\xi)}d\xi,\\
	&\mathscr{A}_6:=\sum_{|\alpha|\leq2}\sum_\upmu\int\sF\rr{\fA_1(t,\na_{\fz,v})v^\alpha Q}(t,\upmu\bfe,\xi)\ov{\sF\rr{\fA_1(t,\na_{\fz,v})(v^\alpha f)}(t,\upmu\bfe,\xi)}d\xi.\ea\]

		The same subadditivity and interpolation arguments used for \eqref{claimflen1} and \eqref{claimflen2} give
		\begin{equation}\label{claimflen3}
			\begin{aligned}
				\rr{\f{\mathsf{T}_k}{t}}^{0.01}\<\upmu,\upxi\>^{0.34}&\leq\rr{\f{\mathsf{T}_k}{t}}^{0.01}\<\upnu,\upnu t\>^{0.34}+\f14\<\f{\upmu-\upnu,\upxi-\upnu t}{A_*}\>^\sss, \\
				\rr{\f{\mathsf{T}_k}{t}}^{0.01}\<\upmu,\upxi\>^{0.34}&\leq\f14\<\f{\upnu,\upnu t}{A_*}\>^{\sss}+\rr{\f{\mathsf{T}_k}{t}}^{0.01}\<\upmu-\upnu,\upxi-\upnu t\>^{0.34}, \\
				\f{\mathsf{r}}{4}(\f34+\f{1}{t})\<\mathtt{r}\nr{\upmu-k',\upxi-A'}\>^\sss&\leq\f{\mathsf{r}}{4}(\f34+\f{1}{t})\<\mathtt{r}\nr{\upnu-k',\upnu t-A'}\>^\sss+\f14\<\f{\upmu-\upnu,\upxi-\upnu t}{A_*}\>^\sss, \\
				\f{\mathsf{r}}{4}(\f34+\f{1}{t})\<\mathtt{r}\nr{\upmu-k',\upxi-A'}\>^\sss&\leq\f{1}{4}\<\f{\upnu,\upnu t}{A_*}\>^\sss+\f{\mathsf{r}}{4}(\f34+\f{1}{t})\<\mathtt{r}\nr{\upmu-\upnu-k',\upxi-\upnu t-A'}\>^\sss .
			\end{aligned}
		\end{equation}
These imply that 
		\[\fA_1(t,\upmu\bfe,\xi)\ls  \fA_1(t,\upnu\bfe,\upnu t\bfe)e^{\f12\<\f{\upmu-\upnu,\upxi-\upnu t}{A_*}\>^\sss+\<\mrxi\>^\sss}\<\upmu-\upnu,\upxi-\upnu t\>^4.\]
From  these together with the estimates $|\xi-\upmu t|\leq\<\upmu-\upnu,\xi-\upnu t\>t$ and $\f{\<\upnu,\upnu t\>}{|\upnu|}\ls\<t\>$, we get that
\[\ba|\mathscr{A}_1|\ls\sum_{|\alpha|\leq2}\nnr{\fA_0(t,\na_{\fz,\fv})\rho(t)}_{L^2_x}\<t\>\nnr{\<\na_v\>^5e^{\<\na_v\>^\sss}\rr{v^\alpha f^0(v)}}_{L^2_v}\nnr{\fA_1(t,\na_{\fz,v})\rr{v^\alpha f(t)}}_{L^2_{z,v}};\\
		|\mathscr{A}_2|\ls\sum_{|\alpha|\leq2}\nnr{\fA_0(t,\na_{\fz,\fv})\rho(t)}_{L^2_x}\<t\>^2\nnr{\<\na_{\fz,\fv}\>^{6.5}e^{\f12\<\f{\na_{\fz,\fv}}{A_*}\>^\sss+\<\na_{\mrv}\>^\sss}\rr{v^\alpha\udl{f}(t)}}_{L^2_{z,v}}\nnr{\fA_1(t,\na_{\fz,v})\rr{v^\alpha f(t)}}_{L^2_{z,v}}.\ea\]
		 
		 To estimate $\mathscr{A}_3$, by \eqref{claimflen3} and the mean value theorem, we observe that  
		\begin{equation}\label{commusymbol1}
			\begin{aligned}
				&\qquad\nr{\fA_1(t,\upmu\bfe,\xi)-\fA_1(t,(\upmu-\upnu)\bfe,\xi-\upnu t\bfe)} \\
				&\leq e^{\<\mrxi\>^\sss}\nr{\<\upmu,\xi\>^4e^{\rr{\f{\mathsf{T}_k}{t}}^{0.01}\<\upmu,\upxi\>^{0.34}}-\<\upmu-\upnu,\xi-\upnu t\bfe\>^4e^{\rr{\f{\mathsf{T}_k}{t}}^{0.01}\<\upmu-\upnu,\upxi-\upnu t\>^{0.34}}} e^{\f{\mathsf{r}}{4}(\f34+\f{1}{t})\<\mathtt{r}\nr{\upmu-k',\upxi-A'}\>^\sss} \\
				&\quad+e^{\<\mrxi\>^\sss}\nr{e^{\f{\mathsf{r}}{4}(\f34+\f{1}{t})\<\mathtt{r}\nr{\upmu-k',\upxi-A'}\>^\sss}-e^{\f{\mathsf{r}}{4}(\f34+\f{1}{t})\<\mathtt{r}\nr{\upmu-\upnu-k',\upxi-\upnu t-A'}\>^\sss}}\<\upmu-\upnu,\xi-\upnu t\bfe\>^4 e^{\rr{\f{\mathsf{T}_k}{t}}^{0.01}\<\upmu-\upnu,\upxi-\upnu t\>^{0.34}} \\
				&\ls e^{\<\mrxi\>^\sss}\max\big\{\<\upmu,\xi\>^{-0.66}\<\upmu-\upnu,\xi-\upnu t\bfe\>^4,\<\upmu-\upnu,\xi-\upnu t\bfe\>^{-0.66}\<\upmu,\xi\>^4\big\}e^{\f12\<\f{\upnu,\upnu t}{A_*}\>^\sss} e^{\f{\mathsf{r}}{4}(\f34+\f{1}{t})\<\mathtt{r}\nr{\upmu-\upnu-k',\upxi-\upnu t-A'}\>^\sss} \\
				&\quad\times e^{\rr{\f{\mathsf{T}_k}{t}}^{0.01}\<\upmu-\upnu,\upxi-\upnu t\>^{0.34}}\<\upnu,\upnu t\> +\mathtt{r}\mathsf{r}e^{\<\mrxi\>^\sss}\<\upmu-\upnu,\xi-\upnu t\bfe\>^4e^{\rr{\f{\mathsf{T}_k}{t}}^{0.01}\<\upmu-\upnu,\upxi-\upnu t\>^{0.34}}e^{\f{\mathsf{r}}{4}(\f34+\f{1}{t})\<\mathtt{r}\nr{\upmu-\upnu-k',\upxi-\upnu t-A'}\>^\sss} \\
				&\quad\times\max\{\<\mathtt{r}\nr{\upmu-\upnu-k',\upxi-\upnu t-A'}\>^{\sss-1},\<\mathtt{r}\nr{\upmu-k',\upxi-A'}\>^{\sss-1}\} e^{\f14\<\f{\upnu,\upnu t}{A_*}\>^\sss}\<\upnu,\upnu t\>
			\end{aligned}
		\end{equation}
		and $|\xi-\upmu t|\ls\min\{|\upmu,\xi|,|\upmu-\upnu,\xi-\upnu t|\}t\leq\min\{|\upmu-k',\xi-A'|,|\upmu-\upnu-k',\xi-\upnu t-A'|\}t+|k',A'|t$. From these together with the fact that $\max\{x^{-a},y^{-a}\}\times\min\{x,y\}\le \min\{x^{1-a},y^{1-a}\}$ with $a\in]0,1[$ and $x,y>0$, we get  that 
      \[\ba
		&\qquad\nr{\fA_1(t,\upmu\bfe,\xi)-\fA_1(t,(\upmu-\upnu)\bfe,\xi-\upnu t\bfe)}|\xi-\upmu t|\\
		&\ls  e^{\f12\<\f{\upnu,\upnu t}{A_*}\>^\sss}\<\upnu,\upnu t\>^{5}\fA_{1.17}(t,(\upmu-\upnu)\bfe,\xi-\upnu t\bfe)\<\upmu,\xi\>^{0.17}t+\mathsf{r}e^{\f14\<\f{\upnu,\upnu t}{A_*}\>^\sss}\<\upnu,\upnu t\>^{1+\f{\sss}2} \<\mathtt{r}\nr{\upmu-\upnu-k',\upxi-\upnu t-A'}\>^{\f{\sss}2}\\
		&\times\fA_{1}(t,(\upmu-\upnu)\bfe,\xi-\upnu t\bfe) \<\mathtt{r}\nr{\upmu-k',\upxi-A'}\>^{\f{\sss}2}t+\mathsf{r}\mathtt{r}|k',A'|e^{\f14\<\f{\upnu,\upnu t}{A_*}\>^\sss}\<\upnu,\upnu t\>^{}\fA_{1}(t,(\upmu-\upnu)\bfe,\xi-\upnu t\bfe)t.
      		\ea \]
      This gives
		\[\ba
		&|\mathscr{A}_3|\ls\sum_{|\alpha|\leq2}\sum_{\upmu,\upnu}\int\nr{\sF\rr{\<\na_{\fz,\fv}\>^5e^{\f12\<\f{\na_{\fz,\fv}}{A_*}\>^\sss}\udl{\rho}}(t,\upnu\bfe)}t\nr{\sF\rr{\fA_{1.17}(t,\na_{\fz,v})v^\alpha f}(t,(\upmu-\upnu)\bfe,\xi-\upnu t\bfe)} \\
		&\times|\sF\rr{\fA_{1.17}(t,\na_{\fz,v})v^\alpha f}(t,\upmu\bfe,\xi)|d\xi
		+\mathsf{r}\sum_{|\alpha|\leq2}\sum_{\upmu,\upnu}\int\nr{\sF\rr{\<\na_{\fz,\fv}\>^2e^{\f14\<\f{\na_{\fz,\fv}}{A_*}\>^\sss}\udl{\rho}}(t,\upnu\bfe)}t\\
		&\times\left|\sF\rr{\<\mathtt{r}\nr{\na_\fz-k',\na_\fv-A'}\>^{\f\sss2}\fA_{1}(t,\na_{\fz,v})v^\alpha f}\right. (t,(\upmu-\upnu)\bfe,\xi-\upnu t\bfe)|\\
		&\times
		\nr{\sF\rr{\<\mathtt{r}\nr{\na_\fz-k',\na_\fv-A'}\>^{\f\sss2}\fA_{1}(t,\na_{\fz,v})v^\alpha f}(t,\upmu\bfe,\xi)}d\xi
		+\mathsf{r}\mathtt{r}|k',A'|\sum_{|\alpha|\leq2}\sum_{\upmu,\upnu}\int\nr{\sF\rr{\<\na_{\fz,\fv}\>e^{\f14\<\f{\na_{\fz,\fv}}{A_*}\>^\sss}\udl{\rho}}(t,\upnu\bfe)}\\
		&\times t\nr{\sF\rr{\fA_{1}(t,\na_{\fz,v})v^\alpha f}(t,(\upmu-\upnu)\bfe,\xi-\upnu t\bfe)}\nr{\sF\rr{\fA_{1}(t,\na_{\fz,v})v^\alpha f}(t,\upmu\bfe,\xi)}d\xi \\
		&\ls\eps_0\<t\>^{-2.5}
		\nnr{\fA_{1.17}(t,\na_{\fz,v})\rr{v^\alpha f(t)}}^2_{\lta L^2_{z,v}}\\
		&\quad+\eps_0\mathsf{r}\<t\>^{-6.5}e^{-\f12\<\f{t}{A_*}\>^\sss}
		\nnr{\<\mathtt{r}\nr{\na_\fz-k',\na_\fv-A'}\>^{\f\sss2}\fA_{1}(t,\na_{\fz,v})\rr{v^\alpha f(t)}}^2_{\lta L^2_{z,v}}\\
		&\quad+\eps_0\mathsf{r}\mathtt{r}|k',A'|\<t\>^{-6.5}e^{-\f12\<\f{t}{A_*}\>^\sss}
		\nnr{\fA_{1}(t,\na_{\fz,v})\rr{v^\alpha f(t)}}^2_{\lta L^2_{z,v}}.
		\ea\]
		In the last inequality we use \eqref{udlrho}. All the terms on the right-hand side can be absorbed by the damping terms on the left-hand side of \eqref{u1vfestimate1}, since $t\geq \mathfrak{T}_k=(\log k)^{1000}\gg A_*$ and
		\[|k',A'|\leq \mathsf{A}_k^{\mathtt{N}_1+2}\ll e^{\f12\<\f{\mathfrak{T}_k}{A_*}\>^\sss}.\]

		Applying the estimates used for $\mathscr{A}_1$ and $\mathscr{A}_2$ to $\mathscr{A}_4$ and $\mathscr{A}_5$ yields
		\[\ba
		&|\mathscr{A}_4|+|\mathscr{A}_5|\ls\nnr{\fA(t,\na_{\fz,\fv})\rho(t)}_{L^2_x}\<t\>\nnr{\<\na_{\fz,\fv}\>^{5.5}e^{\f12\<\f{\na_{\fz,\fv}}{A_*}\>^\sss+\<\na_{\mrv}\>^\sss}\rr{v^\alpha\udl{f}(t)}}_{\lta L^2_{z,v}}\nnr{\fA_1(t,\na_{\fz,v})\rr{v^\alpha f(t)}}_{\lta L^2_{z,v}}\\
		&\qquad\qquad\qquad+\nnr{\<\na_{\fz,\fv}\>^5e^{\f12\<\f{\na_{\fz,\fv}}{A_*}\>^\sss}\udl{\rho}(t)}_{L^2_x}\nnr{\fA_{1}(t,\na_{\fz,v})\rr{v^\alpha f(t)}}^2_{\lta L^2_{z,v}}.
		\ea\]
		For $\mathscr{A}_6$, it is straightforward to see that
		\[\mathscr{A}_6\leq\nnr{\fA_{1}(t,\na_{\fz,v})\rr{v^\alpha f(t)}}_{\lta L^2_{z,v}}\nnr{\fA_{1}(t,\na_{\fz,v})\rr{v^\alpha Q(t)}}_{\lta L^2_{z,v}}.\]

		Now summing up all the estimates and integrating over the interval $[\mathsf{T}_k,t]$, we get
		\begin{equation}\label{lf}
			\begin{aligned}
				\nnr{\fA_{1}(t,\na_{\fz,v})\rr{v^\alpha f(t)}}_{\lta L^2_{z,v}}&\leq\nnr{\fA_{1}(\mathsf{T}_k,\na_{\fz,v})\rr{v^\alpha f(\mathsf{T}_k)}}_{\lta L^2_{z,v}} \\
				&\quad +\int^t_{\mathsf{T}_k} C_\theta\nnr{\fA(\tau,\na_{\fz,\fv})\rho(\tau)}_{L^2_x}\<\tau\>d\tau \\
				&\quad +\int_{\mathsf{T}_k}^t\eps_0^{0.9}\nnr{\fA(\tau,\na_{\fz,\fv})\rho(\tau)}_{L^2_x}\<\tau\>^2d\tau \\
				&\quad +\int_{\mathsf{T}_k}^t\nnr{\fA_{1}(\tau,\na_{\fz,v})\rr{v^\alpha Q(\tau)}}_{\lta L^2_{z,v}}d\tau \\
				&\leq \mathsf{K}+Z(t)+C_\theta\nnr{\<\tau\>^3\fA(\tau,\na_{\fz,\fv})\rho(\tau)}_{L^2_\tau([\mathsf{T}_k,t];L^2_x)} \\
				&\quad +\eps_0^{0.9}\nnr{\<\tau\>^3\fA(\tau,\na_{\fz,\fv})\rho(\tau)}_{L^2_\tau([\mathsf{T}_k,t];L^2_x)}.
			\end{aligned}
		\end{equation}
		
		\noindent\textit{Step 3. Closing the estimates.}
		Combining \eqref{lerho} and \eqref{lf}, we obtain
		\[\left\{\ba&\nnr{\<\tau\>^3\fA(\tau,\na_{\fz,\fv})\rho(\tau)}_{L^2_\tau([\mathsf{T}_k,t];L^2_x)}\leq C_\theta\rr{\mathsf{K}+Z(t)}
		+\eps_0^{0.9}\sup_{\mathsf{T}_k\leq\tau\leq t}\nnr{\fA_1(\tau,\na_{\fz,v})v^\alpha f}_{\lta L^2_{z,v}}\\
		&\sup_{\mathsf{T}_k\leq \tau\leq t}\nnr{\fA_{1}(\tau,\na_{\fz,v})\rr{v^\alpha f(\tau)}}_{\lta L^2_{z,v}}\leq \mathsf{K}+Z(t)+\rr{C_\theta+\eps_0^{0.9}}\nnr{\<\tau\>^3\fA(\tau,\na_{\fz,\fv})\rho(\tau)}_{L^2_\tau([\mathsf{T}_k,t];L^2_x)}
		\ea\right.\]
		By a standard continuity argument we get that
		\[\sup_{\mathsf{T}_k\leq \tau\leq t}\nnr{\fA_{1}(\tau,\na_{\fz,v})\rr{v^\alpha f(\tau)}}_{\lta L^2_{z,v}}\leq 2C_\theta^2 \rr{\mathsf{K}+Z(t)}.\]
   \smallskip
		
		\noindent\textit{Proof of (ii).}
		We retain the three-step structure of part (i), but use the backward Volterra equation and the last two kernel bounds in Lemma \ref{auxi}. In this time direction, $\wh{\rho}(t,\upmu\bfe)$ satisfies
		\[\ba
		&\wh{\rho}(t,\upmu\bfe)=H(t,\upmu\bfe)-\int_t^{\mathsf{T}_k}H(\tau,\upmu\bfe)R^0(-(t-\tau),-\upmu\bfe)d\tau,\quad
		H(t,\upmu\bfe)=\wh{f}(\mathsf{T}_k,\upmu\bfe,\upmu t\bfe)-\int_t^{\mathsf{T}_k}\wh{Q}(\tau,\upmu\bfe,\upmu t\bfe)d\tau\\&+\sum_{\upnu\neq0}\int_t^{\mathsf{T}_k}\wh{\rho}(\tau,\upnu\bfe)\f{\upmu(t-\tau)}{\upnu}\wh{\udl{f}}(\tau,(\upmu-\upnu)\bfe,(\upmu t-\upnu\tau)\bfe)d\tau
		+\sum_{\upnu\neq0}\int_t^{\mathsf{T}_k}\wh{\udl{\rho}}(\tau,\upnu\bfe)\f{\upmu(t-\tau)}{\upnu}\wh{f}(\tau,(\upmu-\upnu)\bfe,(\upmu t-\upnu\tau)\bfe)d\tau.
		\ea\]
		\noindent\textit{Step 1. Backward density estimate.}
		Fix $T\in[\mathfrak{T}_k,\mathsf{T}_k]$. The representative background-interaction term is
 \[\nnr{\<t\>^3\fA_0(t,\upmu\bfe,\upmu t\bfe)\sum_{\upnu\neq0}\int^{\mathsf{T}_k}_t\wh{\rho}(\tau,\upnu\bfe)\f{\upmu(t-\tau)}{\upnu}\wh{\udl{f}}(\tau,(\upmu-\upnu)\bfe,(\upmu t-\upnu\tau)\bfe)d\tau}_{L^2_t([T,\mathsf{T}_k];L^2_\upmu)}.\]
		The proof of \eqref{claimflen1} applies with the forward analytic factor $(\mathsf{T}_k/t)^{0.01}$ replaced by $\log t/\log\mathsf{T}_k$, and with the localization factor $(3/4+1/t)$ replaced by $(1-1/t)$. Therefore the analogue of the $\mathscr{J}_3$ estimate is
		\[\ba
		&\nr{\sum_{\upnu\neq0}\int_t^{\mathsf{T}_k}\<t\>^3\fA_0(t,\upmu\bfe,\upmu t\bfe)\wh{\rho}(\tau,\upnu\bfe)\f{\upmu(t-\tau)}{\upnu}\wh{\udl{f}}(\tau,(\upmu-\upnu)\bfe,(\upmu t-\upnu\tau)\bfe)d\tau}
		\leq\sum_{\upnu\neq0}\int_t^{\mathsf{T}_k}\<\tau\>^3\nr{\sF\rr{\fA_0(\tau,\na_{\fz,\fv})\rho}(\tau,\upnu\bfe)}\\
		&\f{\tau\vv{1}_{\upnu\neq\upmu}+1}{|\upnu|}\nr{\sF\rr{\<\na_{\fz,\fv}\>^4e^{\f12\<\f{\na_{\fz,\fv}}{A_*}\>^\sss}\udl{f}}(\tau,(\upmu-\upnu)\bfe,(\upmu t-\upnu\tau)\bfe)}e^{\f{\log t-\log \tau}{\log \mathsf{T}_k}\<\upmu,\upmu t\>^{0.34}}d\tau\\
		\leq& \eps_{0}^{0.95}\sum_{\upnu\neq0}\int_t^{\mathsf{T}_k}\<\tau\>^3\nr{\sF\rr{\fA(\tau,\na_{\fz,\fv})\rho}(\tau,\upnu\bfe)}\f{\tau\vv{1}_{\upnu\neq\upmu}+1}{|\upnu|}\<\upmu-\upnu,\upmu t-\upnu\tau\>^{-6}e^{\f{\log t-\log \tau}{\log \mathsf{T}_k}\<\upmu,\upmu t\>^{0.34}}d\tau.
		\ea\]
		Now applying the last two inequalities in Lemma \ref{auxi}, we get from Schur test that
		\[\ba&\nnr{\<t\>^3\fA_0(t,\upmu\bfe,\upmu t\bfe)\sum_{\upnu\neq0}\int^{\mathsf{T}_k}_t\wh{\rho}(\tau,\upnu\bfe)\f{\upmu(t-\tau)}{\upnu}\wh{\udl{f}}(\tau,(\upmu-\upnu)\bfe,(\upmu t-\upnu\tau)\bfe)d\tau}_{L^2_t([T,\mathsf{T}_k];L^2_\upmu)}\\&
		\le \eps_{0}^{0.94}\nnr{\<\tau\>^3\fA_0(\tau,\upnu\bfe,\upnu\tau\bfe)\wh{\rho}(\tau,\upnu\bfe)}_{L^2_\tau([T,\mathsf{T}_k];L^2_\upnu)}.\ea\]
		 
		\noindent\textit{Step 2. Backward profile energy estimate.}
		The weighted energy identity now has the form
		\[\ba
		&-\f12\f{d}{dt}\nnr{\fA_1(t,\na_{\fz,v})\rr{v^\alpha f(t)}}_{\lta L^2_{z,v}}^2+\f{1}{2t\log \mathsf{T}_k}\nnr{\fA_{1.17}(t,\na_{\fz,v})\rr{v^\alpha f(t)}}_{\lta L^2_{z,v}}^2\\
		&+\f{\mathsf{r}}{8t^2}\nnr{\<\mathtt{r}\nr{\na_\fz-k',\na_\fv-A'}\>^{\f\sss2}\fA_{1}(t,\na_{\fz,v})\rr{v^\alpha f(t)}}_{\lta L^2_{z,v}}^2\leq \left(C_\theta\nnr{\fA(t,\na_{\fz,\fv})\rho(t)}_{L^2_x}\<t\>\right.\\&\left.
		+\eps_0^{0.9}\nnr{\fA_0(t,\na_{\fz,\fv})\rho(t)}_{L^2_x}\<t\>^2+\nnr{\fA_{1}(t,\na_{\fz,v})\rr{v^\alpha Q(t)}}_{\lta L^2_{z,v}}\right)\nnr{\fA_1(t,\na_{\fz,v})\rr{v^\alpha f(t)}}_{\lta L^2_{z,v}}.
		\ea\]
	
		\noindent\textit{Step 3. Closing the backward estimates.}
		Set
		\[
		X_T:=\nnr{\<t\>^3\fA(t,\na_{\fz,\fv})\rho(t)}_{L^2_t([T,\mathsf{T}_k];L^2_x)},\qquad
		Y_T:=\sup_{T\leq t\leq\mathsf{T}_k}\nnr{\fA_1(t,\na_{\fz,v})\rr{v^\alpha f(t)}}_{\lta L^2_{z,v}}.
		\]
		The initial-data and source terms are bounded by $\mathsf K+Z(T)$. The estimate above, the analogous bound for the interaction with $\underline\rho$, and the backward energy inequality give
		\[
		X_T\leq C_\theta\bigl(\mathsf K+Z(T)\bigr)+\eps_0^{0.9}Y_T,
		\qquad
		Y_T\leq \mathsf K+Z(T)+\bigl(C_\theta+\eps_0^{0.9}\bigr)X_T.
		\]
		For $\eps_0$ sufficiently small, a continuity argument absorbs the term proportional to $\eps_0^{0.9}Y_T$ and yields $Y_T\lesssim_\theta\bigl(\mathsf K+Z(T)\bigr)$, uniformly for $T\in[\mathfrak{T}_k,\mathsf{T}_k]$. This proves part (ii).
	\end{proof}
	
	Next we establish a frequency localization estimate on time interval $[10,\mathfrak{T}_k]$.
	\subsection{A short-time estimate}
	The next estimate is centered at the zero frequency. Unlike Proposition \ref{longenergyl}, it is used after the relevant wave packets have been summed, so no localization center $(k',A')$ is needed.
	\begin{prop}\label{shortenergy}
		Let $k,A_*,\mathsf{A}_k, \mathsf{T}_k, \mathfrak{T}_k,\sss$ be the parameters defined in \eqref{DefTk} and Proposition \ref{propudlf}. Assume that $0\leq\mathsf{r}\leq1$, $0\leq\mathtt{r}\leq A_*^{-1}$, and $\mathsf{K}\geq0$.
		For $10\leq t\leq \mathfrak{T}_k$, if $f(\mathfrak{T}_k)$ and $Q(t)$ satisfy
		\[\ba&\nnr{\<\na_{\fz,\fv}\>^3e^{\f{\mathsf{r}}{4}\rr{1-\f{1}{\mathfrak{T}_k}}\<\mathtt{r}\na_{\fz,\fv}\>^\sss+\<\na_{\mrv}\>^\sss} \rr{v^\alpha f(\mathfrak{T}_k)}}_{\lta L^2_{z,v}}\leq \mathsf{K},\\
		&\int_t^{\mathfrak{T}_k}\nnr{\<\na_{\fz,\fv}\>^3e^{\f{\mathsf{r}}{4}\rr{1-\f{1}{\tau}}\<\mathtt{r}\na_{\fz,\fv}\>^\sss+\<\na_{\mrv}\>^\sss} \rr{v^\alpha Q(\tau)}}_{\lta L^2_{z,v}}d\tau\leq Z(t),\ea\]
		where $Z(t)\geq0$ is a decreasing function, then the solution of $\partial_t f=\mathcal Lf+Q$ on $[10,\mathfrak{T}_k]$, with terminal value $f(\mathfrak{T}_k)$, satisfies
		\[\nnr{\<\na_{\fz,\fv}\>^3e^{\f{\mathsf{r}}{4}\rr{1-\f{1}{t}}\<\mathtt{r}\na_{\fz,\fv}\>^\sss+\<\na_{\mrv}\>^\sss} \rr{v^\alpha f(t)}}_{\lta L^2_{z,v}}\leq e^{C_\theta (\mathfrak{T}_k-t)+\eps_0^{0.9}(\mathfrak{T}_k^2-t^2)} \rr{\mathsf{K}+Z(t)}.\]
	\end{prop}
	\begin{proof} 
		If we set $\fB(t,\mu,\xi):=\<\upmu,\upxi\>^3e^{\f{\mathsf{r}}{4}\rr{1-\f{1}{t}}\<\mathtt{r}|\upmu,\upxi|\>^\sss+\<\mrxi\>^\sss}$, then
       $\fB(t,\mu,\xi)=\fB(t,\upmu\bfe,\xi)$. We shall also use the notation  $\sF(\fB(t,\na_{\fz,\fv})\rho)(t,\upmu\bfe):=\fB(t,\upmu,\upmu t\bfe) \wh{\rho}(t,\upmu\bfe)$   for simplicity. By basic energy method, we may derive that 
\begin{equation}\label{fBestimate}
	-\f12\f{d}{dt}\nnr{\fB(t,\na_{\fz,v})\rr{v^\alpha f(t)}}_{\lta L^2_{z,v}}^2
		+\f{\mathsf{r}}{4t^2}\nnr{\<\mathtt{r}\na_{\fz,\fv}\>^{\f\sss2}\fB(t,\na_{\fz,v})\rr{v^\alpha f(t)}}_{\lta L^2_{z,v}}^2=\sum_{i=1}^{6}\mathscr{B}_i,
\end{equation}
where 
		\[\ba &\mathscr{B}_1=\sum_{|\alpha|\leq2}\sum_{\upmu\neq0}\int\wh{\rho}(t,\upmu\bfe)\f1\upmu\pa_{\xi}^\alpha\rr{ (\upxi-\upmu t)\wh{f^0}(\xi-\upmu t\bfe)}\fB(t,\upmu\bfe,\xi)
		\ov{\sF\rr{\fB(t,\na_{\fz,v})(v^\alpha f)}(t,\upmu\bfe,\xi)}d\xi, \\
		&\mathscr{B}_2=\sum_{|\alpha|\leq2}\sum_{\substack{\upmu\in\Z\\\upnu\neq0}}\int\wh{\rho}(t,\upnu\bfe)\f{\upxi-\upmu t}\upnu\pa_{\xi}^\alpha \wh{\udl{f}}(t,(\upmu-\upnu)\bfe,\xi-\upnu t\bfe)\fB(t,\upmu\bfe,\xi) \ov{\sF\rr{\fB(t,\na_{\fz,v})(v^\alpha f)}(t,\upmu\bfe,\xi)}d\xi,\\
		&\mathscr{B}_3=\sum_{|\alpha|\leq2}\sum_{\substack{\upmu\in\Z\\\upnu\neq0}}\int\wh{\udl{\rho}}(t,\upnu\bfe)\f{\upxi-\upmu t}\upnu\pa_{\xi}^\alpha \wh{f}(t,(\upmu-\upnu)\bfe,\xi-\upnu t\bfe) \rr{\fB(t,\upmu\bfe,\xi)-\fB(t,(\upmu-\upnu)\bfe,\xi-\upnu t\bfe)}\\
		&\qquad\qquad\times\ov{\sF\rr{\fB(t,\na_{\fz,v})(v^\alpha f)}(t,\upmu\bfe,\xi)}d\xi,\\
		 &\mathscr{B}_4=\sum_{\substack{|\alpha|\leq2\\\alpha_1>0}}C_\alpha\sum_{\substack{\upmu\in\Z\\\upnu\neq0}}\int\wh{\rho}(t,\upnu\bfe)\f{1}\upnu\pa^{\alpha-\bfe}_{\xi}\wh{\udl{f}}(t,(\upmu-\upnu)\bfe,\xi-\upnu t\bfe)\fB(t,\upmu\bfe,\xi)\ov{\sF\rr{\fB(t,\na_{\fz,v})(v^\alpha f)}(t,\upmu\bfe,\xi)}d\xi,\\
		&\mathscr{B}_5=\sum_{\substack{|\alpha|\leq2\\\alpha_1>0}}C_\alpha\sum_{\substack{\upmu\in\Z\\\upnu\neq0}}\int\wh{\udl{\rho}}(t,\upnu\bfe)\f1\upnu\pa^{\alpha-\bfe}_\xi\wh{f}(t,(\upmu-\upnu)\bfe,\xi-\upnu t\bfe)\fB(t,\upmu\bfe,\xi)\ov{\sF\rr{\fB(t,\na_{\fz,v})(v^\alpha f)}(t,\upmu\bfe,\xi)}d\xi,\\
		&\mathscr{B}_6=-\sum_{|\alpha|\leq2}\sum_\upmu\int\sF\rr{\fB(t,\na_{\fz,v})v^\alpha Q}(t,\upmu\bfe,\xi)\ov{\sF\rr{\fB(t,\na_{\fz,v})(v^\alpha f)}(t,\upmu\bfe,\xi)}d\xi.
		\ea\]

		The polynomial and localization factors in $\fB$ satisfy the same submultiplicative and commutator bounds as the corresponding factors in $\fA$; omitting the Gevrey factor only weakens the weight. Moreover, each $\mathscr{B}_i$ has the same algebraic form as $\mathscr{A}_i$ in the proof of Proposition \ref{longenergyl}. Applying \eqref{claimflen3} and \eqref{commusymbol1} with localization center $(0,0)$ therefore gives
		\[\ba
		&|\mathscr{B}_1|\ls \nnr{\fB(t,\na_{\fz,\fv})\rho(t)}_{L^2_{x}}\nnr{\<\na_{v}\>^5e^{\<\na_{v}\>^\sss}\rr{v^\alpha f^0}}_{\lta L^2_{v}}\nnr{\fB(t,\na_{\fz,v})\rr{v^\alpha f(t)}}_{\lta L^2_{z,v}}\\
		&\qquad\ls C_\theta\nnr{\fB(t,\na_{\fz,v})\rr{v^\alpha f(t)}}_{\lta L^2_{z,v}}^2,\\ 
		&|\mathscr{B}_2|\leq \nnr{\fB(t,\na_{\fz,\fv})\rho(t)}_{L^2_{x}}t\nnr{\<\na_{\fz,\fv}\>^5e^{\f14\<\f{\na_{\fz,\fv}}{A_*}\>^\sss+\<\na_{\mrv}\>^\sss}\rr{v^\alpha\udl{f}(t)}}_{\lta L^2_{z,v}}\nnr{\fB(t,\na_{\fz,v})\rr{v^\alpha f(t)}}_{\lta L^2_{z,v}}\\
		&\qquad\ls\eps_0t\nnr{\fB(t,\na_{\fz,v})\rr{v^\alpha f(t)}}_{\lta L^2_{z,v}}^2,\\ 
		&|\mathscr{B}_3|\ls\mathsf{r}\nnr{\<\na_{\fz,\fv}\>^6e^{\f14\<\f{\na_{\fz,\fv}}{A_*}\>^\sss}\udl{\rho}(t)}_{L^2_{x}}
		\nnr{\<\mathtt{r}\na_{\fz,\fv}\>^{\f\sss2}\fB(t,\na_{\fz,v})(v^\alpha f(t))}_{\lta L^2_{z,v}}^2\\
		&\qquad+\nnr{\<\na_{\fz,\fv}\>^6e^{\f14\<\f{\na_{\fz,\fv}}{A_*}\>^\sss}\udl{\rho}(t)}_{L^2_{x}}
		\nnr{\fB(t,\na_{\fz,v})(v^\alpha f(t))}_{\lta L^2_{z,v}}^2\\
		&\qquad\ls\eps_0\<t\>^{-8}\mathsf{r}\nnr{\<\mathtt{r}\na_{\fz,\fv}\>^{\f\sss2}\fB(t,\na_{\fz,v})(v^\alpha f(t))}_{\lta L^2_{z,v}}^2
		+\eps_0\<t\>^{-8}\nnr{\fB(t,\na_{\fz,v})(v^\alpha f(t))}_{\lta L^2_{z,v}}^2,\\ 
		&|\mathscr{B}_4|+|\mathscr{B}_5|\ls\nnr{\fB(t,\na_{\fz,\fv})\rho(t)}_{L^2_x} \nnr{\<\na_{\fz,\fv}\>^{5}e^{\f14\<\f{\na_{\fz,\fv}}{A_*}\>^\sss+\<\na_{\mrv}\>^\sss}\rr{v^\alpha\udl{f}(t)}}_{\lta L^2_{z,v}}\nnr{\fB(t,\na_{\fz,v})\rr{v^\alpha f(t)}}_{\lta L^2_{z,v}}\\
		&\qquad+\nnr{\<\na_{\fz,\fv}\>^6e^{\f14\<\f{\na_{\fz,\fv}}{A_*}\>^\sss}\udl{\rho}(t)}_{L^2_x}
		\nnr{\fB(t,\na_{\fz,v})\rr{v^\alpha f(t)}}^2_{\lta L^2_{z,v}}\\
		&\qquad\ls\eps_0\nnr{\fB(t,\na_{\fz,v})\rr{v^\alpha f(t)}}_{\lta L^2_{z,v}}^2,\\
		&|\mathscr{B}_6|\leq\nnr{\fB(t,\na_{\fz,v})\rr{v^\alpha f(t)}}_{\lta L^2_{z,v}}\nnr{\fB(t,\na_{\fz,v})\rr{v^\alpha Q(t)}}_{\lta L^2_{z,v}}\ea
		\]
     
      In the estimates for $\mathscr{B}_1$ and $\mathscr{B}_2$ we use Sobolev embedding in Fourier space. The first term in the bound for $\mathscr{B}_3$ is absorbed by the damping term on the left-hand side of \eqref{fBestimate}.

       Combining these estimates gives
		\[-\f{d}{dt}\nnr{\fB(t,\na_{\fz,v})\rr{v^\alpha f(t)}}_{\lta L^2_{z,v}}\leq\rr{C_\theta+\eps_0^{0.9}t}\nnr{\fB(t,\na_{\fz,v})\rr{v^\alpha f(t)}}_{\lta L^2_{z,v}}+\nnr{\fB(t,\na_{\fz,v})\rr{v^\alpha Q(t)}}_{\lta L^2_{z,v}}.\]
		Integrating this differential inequality backward from $\mathfrak{T}_k$ gives
		\begin{equation}\label{shotener}
			\begin{aligned}
				e^{C_\theta t+\eps_0^{0.9}t^2}\nnr{\fB(t,\na_{\fz,v})\rr{v^\alpha f(t)}}_{\lta L^2_{z,v}}&\leq e^{C_\theta \mathfrak{T}_k+\eps_0^{0.9}\mathfrak{T}_k^2}\nnr{\fB(\mathfrak{T}_k,\na_{\fz,v})\rr{v^\alpha f(\mathfrak{T}_k)}}_{\lta L^2_{z,v}} \\
				&\quad +\int^{\mathfrak{T}_k}_t e^{C_\theta \tau+\eps_0^{0.9}\tau^2}\nnr{\fB(\tau,\na_{\fz,v})\rr{v^\alpha Q(\tau)}}_{\lta L^2_{z,v}}d\tau.
			\end{aligned}
		\end{equation}
		This is the desired estimate.
	\end{proof}
	
	\subsection{Consequences for the Taylor coefficients}
	We apply Propositions \ref{longenergyl} and \ref{shortenergy} to the Taylor coefficients introduced in \eqref{fnq}, obtaining the Fourier localization needed for the analysis in Sections 3--5. Proposition \ref{frelocalfn} then controls the recombined coefficients on the early-time interval $[10,\mathfrak T_k]$.
	
	\begin{prop}\label{frelocalfnq}  Let $k, A_*,\mathsf{A}_k, \mathsf{T}_k, \mathfrak{T}_k,\sss$ be the parameters defined in \eqref{DefTk} and Proposition \ref{propudlf}.
		Suppose that $1\leq n\leq \mathsf{A}_k^{\mathtt{N}_1}$ and $\mathsf{s}_n:=\f14(1-\f14n\mathsf{A}_k^{-\mathtt{N}_1})$. Then the solution $f^{(n,q)}$ to \eqref{fonemp} or \eqref{fnq} satisfies the following estimates:

		\noindent\textit{(i) Forward evolution.} For $t\geq \mathsf{T}_k$,
		\begin{equation}\label{fnqlocfor}
			\nnr{\<\na_{\fz,\fv}\>^3e^{(\f{\mathsf{T}_k}{t})^{0.01}\<\na_{\fz,\fv}\>^{0.34}}e^{ \mathsf{s}_n(\f34+\f{1}{t})\<\f{\na_\fz-qk,\na_\fv-q\mathsf{A}_k}{k^\f13}\>^\sss+\<\na_{\mrv}\>^\sss}\rr{v^\alpha f^{(n,q)}(t)}}_{\lta L^2_{z,v}}\leq e^{5(2n-1)\mathsf{A}_k^{0.34}}.
		\end{equation}
		
		\noindent\textit{(ii) Backward evolution.} For $\mathfrak{T}_k\leq t\leq \mathsf{T}_k$,
		\begin{equation}\label{fnqlocbac}
			\nnr{\<\na_{\fz,\fv}\>^3e^{\f{\log t}{\log \mathsf{T}_k}\<\na_{\fz,\fv}\>^{0.34}}e^{ \mathsf{s}_n(1-\f{1}{t})\<\f{\na_\fz-qk,\na_\fv-q\mathsf{A}_k}{k^\f13}\>^\sss+\<\na_{\mrv}\>^\sss}\rr{v^\alpha f^{(n,q)}(t)}}_{\lta L^2_{z,v}}\leq e^{5(2n-1)\mathsf{A}_k^{0.34}}.
		\end{equation}
	\end{prop}
	\begin{proof}
		We first provide the proof for the case $t\geq \mathsf{T}_k$. The backward interval $\mathfrak{T}_k\leq t\leq \mathsf{T}_k$ is obtained from part (ii) of Proposition \ref{longenergyl}, with the backward multiplier displayed there. Thanks to \eqref{qmq}, we only need to consider the case $q\geq0$.

		We introduce the following symbol \[\fC_a(t,n,q,\mu,\xi):=\<\upmu,\upxi\>^{a}e^{(\f{\mathsf{T}_k}{t})^{0.01}\<\upmu,\upxi\>^{0.34}}e^{ \mathsf{s}_n(\f34+\f{1}{t})\<\f{\upmu-qk,\upxi-q\mathsf{A}_k}{k^\f13}\>^\sss}e^{\<\mrxi\>^\sss}.\]
		Again we use the notation that $\sF(\fC_a(t,n,q,\na_{\fz,\fv})\rho)(t,\upmu\bfe):=\fC_a(t,n,q,\upmu\bfe,\upmu t\bfe)\wh{\rho}(t,\upmu\bfe).$

	\noindent\textit{Base case.} For $f^{(1,1)}$, one has
		\[\ba&\nnr{\fC_3(t,1,1,\na_{\fz,v})\rr{v^{\alpha}f^{(1,1)}(t)}}_{\lta L^2_{z,v}}\ls_{\theta} \nnr{\fC_6(\mathsf{T}_k,1,1,\na_{\fz,v})\rr{v^{\alpha}f^{(1,1)}(\mathsf{T}_k)}}_{\lta L^2_{z,v}}\\
		\leq&\f{\eps_k}{2k^\f13}\nnr{\<k,\upxi\>^{6}e^{\<k,\upxi\>^{0.34}}e^{ \f14\<\f{\upxi-\mathsf{A}_k}{k^\f13}\>^\sss}\pa^{\alpha_1}_{\upxi}\wh{\psi}\rr{\f{\upxi-\mathsf{A}_k}{k^\f13}} }_{\ell^2(|\alpha_1|\leq2)L^2_{\upxi}}\nnr{e^{\<\mrxi\>^\sss}\pa_{\mrxi}^{\alpha_2}\wh{\phi}(\mrxi)}_{\ell^2(|\alpha_2|\leq2)L^2_{\mrxi}}\\
		\ls&\f{\eps_k}{2k^\f13}\nnr{\<k,\upxi\>^{6}e^{\<k,\upxi\>^{0.34}}e^{-9\nr{\f{\upxi-\mathsf{A}_k}{k^\f13}}^\sss}}_{L^2_{\upxi}}\ls\f{\eps_k}{2k^\f13}\<k,\mathsf{A}_k\>^{6}e^{\<k\>^{0.34}+\<\mathsf{A}_k\>^{0.34}}\nnr{e^{-\nr{\f{\upxi-\mathsf{A}_k}{k^\f13}}^\sss}}_{L^2_{\upxi}}\leq e^{2\mathsf{A}_k^{0.34}}.\ea\]
    In the first inequality we use Proposition \ref{longenergyl} with $Q=0$. In the second inequality we use \eqref{fonemp}. In the third inequality we use the property of $\psi$ and $\phi$. In the fourth inequality we use the estimates $\<k,\upxi\>^{0.34}\leq\<k\>^{0.34}+\<\upxi\>^{0.34}$ and $\<\upxi\>^{0.34}\leq\<\mathsf{A}_k\>^{0.34}+\nr{\f{\upxi-\mathsf{A}_k}{k^\f13}}^\sss$; the latter is proved in the same manner as the first inequality in \eqref{claimflen1}.
    
	\noindent\textit{Inductive step.} Assume that \eqref{fnqlocfor} holds for every $n'\leq n-1$. For $(n_1,q_1)$ and $(n_2,q_2)$ satisfying $(n_1,q_1)+(n_2,q_2)=(n,q)$, \eqref{xyaineq} gives
		\[\ba
		&\nnr{\fC_6(\tau,n,q,\na_{\fz,v})\rr{v^{\alpha}\cN(f^{(n_1,q_1)},f^{(n_2,q_2)})}(\tau)}_{\lta L^2_{z,v}}\\
		&\leq\left\|\sum_\upnu\<\upmu,\upxi\>^{6}e^{-\f1{32}\mathsf{A}_k^{-\mathtt{N}_1}\<\f{\upmu-qk,\upxi-q\mathsf{A}_k}{k^\f13}\>^\sss}\right.\\
		&\quad\left.\times\nr{\sF\rr{\fC_0(\tau,n-1,q_1,\na_{\fz,\fv})\rho^{(n_1,q_1)}}(\tau,\upnu\bfe)}\<\tau\>\right. \\
		&\qquad\left.\nr{\sF\rr{\fC_1(\tau,n-1,q_2,\na_{\fz,v})\rr{v^{\alpha}f^{(n_2,q_2)}} }(\tau,(\upmu-\upnu)\bfe,\xi-\upnu\tau\bfe)}\right\|_{\lta L^2_{\mu,\xi}}\\
		&\leq \mathsf{A}_k^{7\mathtt{N}_1}\nnr{\fC_3(\tau,n_1,q_1,\na_{\fz,\fv})\rho^{(n_1,q_1)}(\tau)}_{L^2_x}\nnr{\fC_3(\tau,n_2,q_2,\na_{\fz,v})\rr{v^{\alpha}f^{(n_2,q_2)}}}_{\lta L^2_{z,v}}\<\tau\>^{-2}\\
		&\ls \mathsf{A}_k^{7\mathtt{N}_1}\nnr{\fC_3(\tau,n_1,q_1,\na_{\fz,v})\rr{v^{\alpha}f^{(n_1,q_1)}}}_{\lta L^2_{z,v}}\nnr{\fC_3(\tau,n_2,q_2,\na_{\fz,v})\rr{v^{\alpha}f^{(n_2,q_2)}}}_{\lta L^2_{z,v}}\<\tau\>^{-2}.
		\ea\]
		We remark that in the second inequality we use 
		\[
			\begin{aligned}
				&\<\upmu,\upxi\>^{6}e^{-\f1{32}\mathsf{A}_k^{-\mathtt{N}_1}\<\f{\upmu-qk,\upxi-q\mathsf{A}_k}{k^\f13}\>^\sss} \\
				&\ls\<\upmu-qk,\upxi-q\mathsf{A}_k\>^{6}e^{-\f1{32}\mathsf{A}_k^{-\mathtt{N}_1}\<\f{\upmu-qk,\upxi-q\mathsf{A}_k}{k^\f13}\>^\sss}+\<qk,q\mathsf{A}_k\>^6 \\
				&\ls\rr{\mathsf{A}_k^{\f{\mathtt{N}_1}{s}}k^\f13}^6+\rr{n\mathsf{A}_k}^6\ls\mathsf{A}_k^{6.5\mathtt{N}_1}.
			\end{aligned}
		\]
		Thus 
		\[\sum_{\substack{n_1+n_2=n\\ q_1+q_2=q}}\int_{\mathsf{T}_k}^t\nnr{\fC_6(\tau,n,q,\na_{\fz,v})\rr{v^{\alpha}\cN(f^{(n_1,q_1)},f^{(n_2,q_2)})}(\tau)}_{\lta L^2_{z,v}}d\tau\leq n^2\mathsf{A}_k^{8\mathtt{N}_1}e^{5(2n-2)\mathsf{A}_k^{0.34} }\leq e^{5(2n-1.1)\mathsf{A}_k^{0.34} }.\]
		Proposition \ref{longenergyl} now gives the desired estimate for $f^{(n,q)}$ and completes the induction.
	\end{proof}
	
	\begin{prop}\label{frelocalfn}  Let $k, A_*,\mathsf{A}_k, \mathsf{T}_k, \mathfrak{T}_k,\sss$ be the parameters defined in \eqref{DefTk} and Proposition \ref{propudlf}.
		Suppose that $1\leq n\leq \mathsf{A}_k^{\mathtt{N}_1}$ and $\mathsf{s}_n:=\f14(1-\f14n\mathsf{A}_k^{-\mathtt{N}_1})$. Then, for $10\leq t\leq \mathfrak{T}_k$, the solution $f^{(n)}$ to \eqref{fn} satisfies
		\begin{equation}\label{fnlocbac}
			\nnr{\<\na_{\fz,\fv}\>^3e^{ \mathsf{s}_n(1-\f{1}{t})\<\f{\na_{\fz,\fv}}{k^\f13}\>^\sss+\<\na_{\mrv}\>^\sss}\rr{v^{\alpha}f^{(n)}(t)}}_{\lta L^2_{z,v}}\leq e^{0.26(2n-1)\nr{\f{\mathsf{A}_k}{k^\f13}}^\sss}.
		\end{equation}
	\end{prop}
	\begin{proof}
		We set $\fD_a(t,n,\mu,\xi):=\<\upmu,\upxi\>^{a}e^{ \mathsf{s}_n(1-\f{1}{t})\<\f{\upmu,\upxi}{k^\f13}\>^\sss+\<\mrxi\>^\sss}$. By \eqref{fnqlocbac} it holds that
		\[\ba&\nnr{\fD_3(\mathfrak{T}_k,n,\na_{\fz,v})\rr{v^{\alpha}f^{(n)}(\mathfrak{T}_k)}}_{\lta L^2_{z,v}}\leq\sum_{\substack{|q|\leq n\\2\mid n-q}}\nnr{\<\na_{\fz,\fv}\>^3e^{\mathsf{s}_n(1-\f{1}{\mathfrak{T}_k})\<\f{\na_\fz-qk,\na_\fv-q\mathsf{A}_k}{k^\f13}\>^\sss+\<\na_{\mrv}\>^\sss}\rr{v^{\alpha}f^{(n,q)}(\mathfrak{T}_k)}}_{\lta L^2_{z,v}}\\&\times e^{\f14\<\f{qk,q\mathsf{A}_k}{k^\f13}\>^\sss}
		\leq ne^{\f14\<\f{nk,n\mathsf{A}_k}{k^\f13}\>^\sss}e^{5(2n-1)\mathsf{A}_k^{0.34}}\leq e^{0.255n\nr{\f{\mathsf{A}_k}{k^\f13}}^\sss}.
		\ea\]
		For $n=1$, by Proposition \ref{shortenergy}, we have 
		\[\nnr{\fD_3(t,1,\na_{\fz,v})\rr{v^{\alpha}f^{(1)}(t)}}_{\lta L^2_{z,v}}\leq e^{\mathfrak{T}_k^2}\nnr{\fD_3(\mathfrak{T}_k,1,\na_{\fz,v})\rr{v^{\alpha}f^{(1)}(\mathfrak{T}_k)}}_{\lta L^2_{z,v}}\leq e^{0.26\nr{\f{\mathsf{A}_k}{k^\f13}}^\sss}.\]
		We assume \eqref{fnlocbac} holds for $n'\leq n-1$. In what follows we prove that \eqref{fnlocbac} holds for $n'=n$. For $n_1$, $n_2$ satisfying $n=n_1+n_2$, it holds that
		\[\ba
		&\nnr{\fD_3(\tau,n,\na_{\fz,v})\rr{v^{\alpha}\cN(f^{(n_1)},f^{(n_2)})}(\tau)}_{\lta L^2_{z,v}}\leq\left\|\sum_\upnu\<\upmu,\upxi\>^{3}e^{-\f1{32}\mathsf{A}_k^{-\mathtt{N}_1}\<\f{\upmu,\upxi}{k^\f13}\>^\sss}\right.\\&\left.\times\nr{\sF\rr{\fD_1(\tau,n-1,\na_{\fz,\fv})\rho^{(n_1)}}(\tau,\upnu\bfe)}\nr{\sF\rr{\fD_1(\tau,n-1,\na_{\fz,v})\rr{v^{\alpha}f^{(n_2)}} }(\tau,(\upmu-\upnu)\bfe,\xi-\upnu\tau\bfe)}\right\|_{\lta L^2_{\mu,\xi}}\\
		&\leq \mathsf{A}_k^{4\mathtt{N}_1}\nnr{\fD_3(\tau,n_1,\na_{\fz,v})\rr{v^{\alpha}f^{(n_1)}}}_{\lta L^2_{z,v}}\nnr{\fD_3(\tau,n_2,\na_{\fz,v})\rr{v^{\alpha}f^{(n_2)}}}_{\lta L^2_{z,v}}\<\tau\>^{-2}.
		\ea\]
		Thus
		\[\ba&\sum_{n_1+n_2=n}\int_t^{\mathfrak{T}_k}\nnr{\fD_3(\tau,n,\na_{\fz,v})\rr{v^{\alpha}\cN(f^{(n_1)},f^{(n_2)})}(\tau)}_{\lta L^2_{z,v}}d\tau\\
		&\leq \sum_{n_1+n_2=n}\mathsf{A}_k^{4\mathtt{N}_1}e^{0.26(2\rr{n_1+n_2}-2)\nr{\f{\mathsf{A}_k}{k^\f13}}^\sss}\leq2n\mathsf{A}_k^{4\mathtt{N}_1}e^{0.26(2n-1.9)\nr{\f{\mathsf{A}_k}{k^\f13}}^\sss}\leq e^{0.26(2n-1.8)\nr{\f{\mathsf{A}_k}{k^\f13}}^\sss}.\ea\]
		  Proposition \ref{shortenergy} therefore gives
		\[\nnr{\fD_3(t,n,\na_{\fz,v})\rr{v^{\alpha}f^{(n)}(t)}}_{\lta L^2_{z,v}}\leq e^{\mathfrak{T}_k^2}e^{0.26(2n-1.8)\nr{\f{\mathsf{A}_k}{k^\f13}}^\sss}\leq e^{0.26(2n-1)\nr{\f{\mathsf{A}_k}{k^\f13}}^\sss}.\]
		This completes the induction.
	\end{proof}

	\section{Homogeneous linearized dynamics}

This section analyzes the source-free linearized equation \eqref{fone} around the spatially inhomogeneous background constructed in Proposition \ref{propudlf}:
\begin{equation}\label{Defdulfinfty}
	\left\{\begin{aligned}
		&\pa_tf=\cL f=\cN(f,f^0+\udl{f_\infty}+ f_{*,err})+\cN(\udl{f_\infty}+ f_{*,err},f), \\
		&\udl{f_\infty}(z,v):=\eps_0\vf(v)\cos \fz+f_{*,\infty}(z,v), \\
		&f|_{t=\mathsf{T}_k}=\frac{\eps_k}{2} \psi(k^{\f13} \fv)\phi(\mrv) e^{ i(k \fz+\mathsf{A}_k \fv)},\quad z=(\fz,\mathring{z})\in \T^3, v=(\fv,\mrv)\in \R^3.
	\end{aligned}\right.
\end{equation}
We study the evolution of the prescribed high-frequency perturbation around the time-dependent background. Our goal is to establish matching upper and lower bounds for the amplification along the echo chain, while controlling the backward evolution. The analysis combines the Fourier localization estimates from Section 2 with precise estimates of resonant interactions using the Duhamel principle.

\subsection{Forward echo cascade on $[\mathsf{T}_k,\infty)$}
The plane-wave packet prescribed at $\mathsf{T}_k=t_k$ is centered near the frequency $(k,\mathsf{A}_k)$. On $[\mathsf{T}_k,t_0]$, the resonance relation $\upmu t\approx\mathsf{A}_k$ is met successively by the modes $\upmu=k,k-1,\ldots,1$, so this interval contains the complete forward echo chain. Our analysis proceeds in four stages. First, Proposition \ref{prop:forward-recurrence} isolates the single resonant mode on each interval $[t_\upmu,t_{\upmu-1}]$ and reduces its effect to an approximate recurrence, with the Penrose condition entering through the Volterra resolvent. Second, Proposition \ref{prop:essential-support} propagates the effective support of the reduced amplitude $\beta(t_\upmu,\upmu,\fv)$; together with frequency localization near $\upxi=\mathsf{A}_k$, this ensures that the velocity derivative has the size required by the leading transfer. Third, Lemma \ref{lem:J-bounds} and Proposition \ref{propertybetatmumu} estimate the recurrence coefficients from above and below and iterate them along the chain, yielding the amplification factor $\eps_0\mathsf{A}_k/\upmu^3$ at each step while keeping off-chain interactions negligible. Finally, Proposition \ref{prop:forward-profile} transfers these estimates from $\beta$ back to the full profile $f$, proves the desired lower bound, and shows that for $t\geq t_0=2\mathsf{A}_k$ no further resonance occurs and the profile changes only by an exponentially small amount. We begin with the recurrence on a single critical interval.

\begin{prop}\label{prop:forward-recurrence} Let $k, A_*,\mathsf{A}_k, \mathsf{T}_k, \mathfrak{T}_k,\sss$ be the parameters defined in \eqref{DefTk} and Proposition \ref{propudlf}.
		Recall that $t_\upmu=\f{\mathsf{A}_k}{\upmu+\f12}$ with $0\leq\upmu\leq k$. For $1\leq\upmu\leq k$, the following approximate recurrence relation holds on $[t_\upmu,t_{\upmu-1}]$:
		\begin{equation}\label{aprefor}
			\sF_zf(t_{\upmu-1},\upsigma\bfe,v)=\sF_zf(t_{\upmu},\upsigma\bfe,v)+\cD_{\upsigma,\upmu}\beta(t_{\upmu},\upmu,\fv)+\sR_{\upsigma,\upmu}(v).
		\end{equation}
More precisely:
 \begin{itemize}
	\item For $\upsigma\in\Z$ and $\upmu\in\Z\backslash\{0\}$, the transfer operator $\cD_{\upsigma,\upmu}$ is defined by
	\begin{equation}\label{DefcDsigmamu}
		\cD_{\upsigma,\upmu}f(\fv):=2\pi i\rr{\f{1}{\upmu^2}\pa_\fv\sF_z(f^0+\udl{f_\infty})((\upsigma-\upmu)\bfe,v)J(\upmu,\fv)f(\fv)-\f{\upsigma-\upmu}{\upmu^3}\sF_z\udl{f_\infty}((\upsigma-\upmu)\bfe,v)\pa_\fv\rr{J(\upmu,\fv)f(\fv)}}.
	\end{equation}
   where the resolvent factor is
	\begin{equation}\label{DefJmuv}
		J(\upmu,\fv):=\f{\f{\upmu}{|\upmu|}}{1+\int_0^\infty t\rr{\wh{f^0}(\upmu t\bfe)+\wh{\udl{f_\infty}}(0,\upmu t\bfe)}e^{i\upmu t\fv}dt }.
	\end{equation}
\item The reduced amplitude $\beta$ and the corresponding background profiles are
	\begin{equation}\label{Defbetatmuv}
		\begin{aligned}
			\beta(t,\upmu,\fv)&:=\int_{\R^2}\sF_zf(t,\upmu\bfe,v)d\mrv, \qquad \beta^0(\fv):=\int_{\R^2}f^0(v)d\mrv, \\
			\udl{\beta_\infty}(\upmu,\fv)&:=\int_{\R^2}\sF_z\udl{f_\infty}(\upmu\bfe,v)d\mrv, \qquad \beta_{*,\infty}(\upmu,\fv):=\int_{\R^2}\sF_zf_{*,\infty}(\upmu\bfe,v)d\mrv.
		\end{aligned}
	\end{equation}
	\item The error term $\sR_{\upsigma,\upmu}(v)$ satisfies
		\begin{equation}\label{sRsigmamu}
			\nnr{e^{0.12\<\f{\na_\fv-\mathsf{A}_k}{k^\f13}\>^\sss+\<\na_{\mrv}\>^\sss}\rr{v^{\alpha}\sR_{\upsigma,\upmu}}}_{\lta L^2_v}\leq e^{-97\mathsf{A}_k^{0.34}}.
		\end{equation}
	\item On the whole critical interval $[t_\upmu,t_{\upmu-1}]$, the profile variation satisfies
		\begin{equation}\label{tmutupper}
			\nnr{e^{\<\na_{\mrv}\>^\sss}v^{\alpha}\rr{\sF_zf(t,\upsigma\bfe,v)-\sF_zf(t_{\upmu},\upsigma\bfe,v)}}_{\lta L^2_v}\leq\f{\mathsf{A}_k}{\upmu^3}\nnr{\beta(t_{\upmu},\upmu,\fv)}_{L^2_\fv}\<\upsigma-\upmu\>^{-8}+e^{-96\mathsf{A}_k^{0.34}}.
		\end{equation}
 \end{itemize}
\end{prop}
	\begin{proof}
	We divide the proof into four steps: isolate the resonant density mode, control its Volterra resolvent, extract the transfer operator, and estimate the profile throughout the critical interval.
		\smallskip

		\noindent\textit{Step 1. Isolating the resonant mode.}
		On $[t_\upmu,t_{\upmu-1}]$ we prove
		\begin{equation}\label{rerho}
			\wh{\rho}(t,\upnu\bfe)=\wh{\udl{\sR}_1}(t,\upnu\bfe),\quad \text{for}\ \upnu\neq\upmu,
		\end{equation}
		\begin{equation}\label{rhomu}
			\wh{\rho}(t,\upmu\bfe)=\wh{f}(t_\upmu,\upmu\bfe,\upmu t\bfe)-\int_{t_\upmu}^t\wh{\rho}(\tau,\upmu\bfe)(t-\tau)\rr{\wh{f^0}(\upmu (t-\tau)\bfe)+\wh{\udl{f_\infty}}(0,\upmu (t-\tau)\bfe)}d\tau+\wh{\udl{\sR}_2}(t,\upmu\bfe),
		\end{equation}
		\begin{equation}\label{fref}
			\wh{f}(t,\upsigma\bfe,\xi)=\wh{f}(t_\upmu,\upsigma\bfe,\xi)-\int_{t_\upmu}^t\wh{\rho}(\tau,\upmu\bfe)\f{\upxi-\upsigma \tau}{\upmu}\sF\rr{ f^0+\udl{f_\infty} }((\upsigma-\upmu)\bfe,\xi-\upmu\tau\bfe) d\tau+\wh{\udl{\sR}_3}(t,\upsigma\bfe,\xi),
		\end{equation}
		where $\udl{\sR}_1$, $\udl{\sR}_2$, and $\udl{\sR}_3$ satisfy
		\begin{equation}\label{udlsr12}
			\nr{\wh{\udl{\sR}_1}(t,\upnu\bfe)}\leq e^{-99\mathsf{A}_k^{0.34}}e^{-0.124\<\f{\upnu t-\mathsf{A}_k}{k^\f13}\>^\sss},\, \nr{\wh{\udl{\sR}_2}(t,\upmu\bfe)}\leq e^{-99\mathsf{A}_k^{0.34}}e^{-0.124\<\f{\upmu t-\mathsf{A}_k}{k^\f13}\>^\sss},
		\end{equation}
		\begin{equation}\label{udlsr3}
			\nnr{e^{0.124\<\f{\na_\fv-\mathsf{A}_k}{k^\f13}\>^\sss+\<\na_{\mrv}\>^\sss}\rr{v^{\alpha}\udl{\sR}_3(t,\upsigma\bfe,v)}}_{\lta L^2_v}\leq e^{-99\mathsf{A}_k^{0.34}}.
		\end{equation}
		
		\noindent\textit{Off-resonant density modes.}
		Since $f$ solves \eqref{fone}, Proposition \ref{frelocalfnq} gives
		\begin{equation}\label{zsd}
			\nr{\wh{f}(t,\upsigma\bfe,\upxi\bfe)}\leq e^{5\mathsf{A}_k^{0.34}}e^{-0.125\<\f{\upsigma-k,\upxi-\mathsf{A}_k}{k^\f13}\>^\sss},\quad\nr{\wh{\rho}(t,\upnu\bfe)}\leq e^{5\mathsf{A}_k^{0.34}}e^{-0.125\<\f{\upnu-k,\upnu t-\mathsf{A}_k}{k^\f13}\>^\sss}.
		\end{equation}
		If $\upnu\neq\upmu$, then $|\upnu t-\mathsf{A}_k|\ge k^{1/3}\mathsf{A}_k^{0.35}$. Indeed, for $\upmu\ge2$,
		\[
		|\upnu t-\mathsf{A}_k|\geq |\upnu-\upmu|t-|\upmu t-\mathsf{A}_k|
		\geq t_\upmu-\f{\f12}{\upmu-\f12}\mathsf{A}_k
		=\f{\f12\upmu-\f34}{\upmu^2-\f14}\mathsf{A}_k.
		\]
		For $\upmu=1$, the modes $\upnu\leq0$ are immediate, while $\upnu\geq2$ gives $|\upnu t-\mathsf{A}_k|\geq\mathsf{A}_k/3$. Inserting this separation into \eqref{zsd} yields
		\begin{equation}\label{uprhohatnut}
			\nr{\wh{\rho}(t,\upnu\bfe)}\leq e^{-100\mathsf{A}_k^{0.34}}e^{-0.1249\<\f{\upnu t-\mathsf{A}_k}{k^\f13}\>^\sss}.
		\end{equation}

		\noindent\textit{The resonant density mode.}
		Separating the $\upnu=\upmu$ term in \eqref{rhohateq} gives
		\[\ba
		\wh{\udl{\sR}_2}(t,\upmu\bfe)=&-\sum_{\upnu\neq\upmu,0}\int_{t_\upmu}^t\wh{\rho}(\tau,\upnu\bfe)\f{\upmu(t-\tau)}{\upnu}\wh{\udl{f}}(\tau,(\upmu-\upnu)\bfe,(\upmu t-\upnu\tau)\bfe)d\tau	-\int_{t_\upmu}^t\wh{\rho}(\tau,\upmu\bfe)(t-\tau)\wh{f_{*,err}}(\tau,0,\upmu(t-\tau)\bfe)d\tau
		\\&-\sum_{\upnu\neq0}\int_{t_\upmu}^t\wh{\udl{\rho}}(\tau,\upnu\bfe)\f{\upmu(t-\tau)}{\upnu}\wh{f}(\tau,(\upmu-\upnu)\bfe,(\upmu t-\upnu\tau)\bfe)d\tau:=\sum_{i=1}^3\mathscr{C}_i.
		\ea\]
		We estimate $\mathscr C_i$, $i=1,2,3$, separately. By \eqref{zsd} and \eqref{uprhohatnut},
		\[\ba
		|\mathscr{C}_1|&\ls\sum_{\upnu\neq\upmu,0}\int_{t_\upmu}^te^{-100\mathsf{A}_k^{0.34}}e^{-0.1249\<\f{\upnu\tau-\mathsf{A}_k}{k^\f13}\>^\sss}\f{\tau}{|\upnu|}\nr{\wh{\<\na_{\fz,\fv}\>\udl{f}}(\tau,(\upmu-\upnu)\bfe,(\upmu t-\upnu\tau)\bfe}d\tau\\
		&\leq\sum_{\upnu\neq\upmu,0}\int_{t_\upmu}^te^{-99.5\mathsf{A}_k^{0.34}}e^{-0.1248\<\f{\upnu \tau-\mathsf{A}_k}{k^\f13}\>^\sss}\nr{\wh{\<\na_{\fz,\fv}\>\udl{f}}(\tau,(\upmu-\upnu)\bfe,(\upmu t-\upnu\tau)\bfe}d\tau\\
		&\leq e^{-99.5\mathsf{A}_k^{0.34}}e^{-0.1248\<\f{\upmu t-\mathsf{A}_k}{k^\f13}\>^\sss}\sup_{t_\upmu\le\tau\le t_{\upmu-1} }\nnr{\<\na_{\fz,\fv}\>^4e^{\f14\<\f{\na_{\fz,\fv}}{A_*}\>^\sss}\rr{v^{\alpha}\udl{f}(\tau)}}_{\lta L^2_{z,v}};\\
		|\mathscr{C}_2|&\leq\int_{t_\upmu}^t e^{5\mathsf{A}_k^{0.34}}e^{-0.125\<\f{\upmu t-\mathsf{A}_k}{k^\f13}\>^\sss}\<\upmu(t-\tau)\>^{-2}d\tau\ \sup_{t_\upmu\le\tau\le t_{\upmu-1}}\nnr{\<\na_{\fz,\fv}\>^3e^{\f18\<\f{\na_{\fz,\fv}}{A_*}\>^\sss}\rr{v^{\alpha}f_{*,err}(\tau)}}_{\lta L^2_{z,v}}\\
		&\leq e^{-\f12\<\f{t_\upmu}{A_*}\>^\sss}e^{5\mathsf{A}_k^{0.34}}e^{-0.125\<\f{\upmu t-\mathsf{A}_k}{k^\f13}\>^\sss}\leq e^{-100\mathsf{A}_k^{0.34}}e^{-0.125\<\f{\upmu t-\mathsf{A}_k}{k^\f13}\>^\sss}.\\
		|\mathscr{C}_3|&\ls\eps_0\sum_{\upnu\neq0}\int_{t_\upmu}^te^{-\f12\<\f{t_\upmu}{A_*}\>^\sss}e^{-\f12\<\f{\upnu,\upnu\tau}{A_*}\>^\sss}\<\upnu,\upnu\tau\>^{-10}\tau e^{6\mathsf{A}_k^{0.34}}e^{-0.1248\<\f{\upmu t-\upnu\tau-\mathsf{A}_k}{k^\f13}\>^\sss}d\tau\\
		&\leq\eps_0 e^{-\f12\<\f{t_\upmu}{A_*}\>^\sss}e^{6\mathsf{A}_k^{0.34}}e^{-0.1248\<\f{\upmu t-\mathsf{A}_k}{k^\f13}\>^\sss}\leq e^{-100\mathsf{A}_k^{0.34}}e^{-0.124\<\f{\upmu t-\mathsf{A}_k}{k^\f13}\>^\sss}.
		\ea\]
		We remark that for the estimate of $\mathscr{C}_3$ we use $|\upmu(t-\tau)|\leq\<\upmu-\upnu,\upmu t-\upnu\tau\>\tau\ls\tau\mathsf{A}_ke^{10^{-4}\<\f{\upmu-\upnu-k,\upmu t-\upnu\tau-\mathsf{A}_k}{k^\f13}\>^\sss}$.
		These bounds prove the estimate for $\udl{\sR}_2$ in \eqref{udlsr12}.

		\noindent\textit{The profile equation.}
		By \eqref{whfeq}, it remains to show that
		\[\ba
		& \wh{\udl{\sR}_3}(t,\upsigma\bfe,\xi)=-\sum_{\upnu\neq\upmu,0}\int_{t_\upmu}^t\wh{\rho}(\tau,\upnu\bfe)\f{\upxi-\upsigma \tau}{\upnu}\sF\rr{\udl{f}+f^0}(\tau,(\upsigma-\upnu)\bfe,\xi-\upnu\tau\bfe)d\tau		\\&-\int_{t_\upmu}^t\wh{\rho}(\tau,\upmu\bfe)\f{\upxi-\upsigma \tau}{\upmu}\wh{f_{*,err}}(\tau,(\upsigma-\upmu)\bfe,\xi-\upmu\tau\bfe)d\tau
	    -\sum_{\upnu\neq0}\int_{t_\upmu}^t\wh{\udl{\rho}}(\tau,\upnu\bfe)\f{\upxi-\upsigma \tau}{\upnu}\wh{f}(\tau,(\upsigma-\upnu)\bfe,\xi-\upnu\tau\bfe)d\tau
		\ea\]
		satisfies \eqref{udlsr3}. The three terms are controlled by \eqref{uprhohatnut}, the estimate for the resonant density mode above, and \eqref{udlrho}, respectively:
		\[\ba
		&\nnr{e^{0.124\<\f{\upxi-\mathsf{A}_k}{k^\f13}\>^\sss+\<\mrxi\>^\sss}\pa_\xi^\alpha\wh{\udl{\sR}_3}(t,\upsigma\bfe,\xi)}_{\lta L^2_\xi}\ls\sum_{\upnu\neq\upmu,0}\int_{t_\upmu}^te^{0.124\<\f{\upnu\tau-\mathsf{A}_k}{k^\f13}\>^\sss}\nr{\wh{\rho}(\tau,\upnu\bfe)}\tau\left\|e^{\f18\<\f{\upxi}{A_*}\>^\sss+\<\mrxi\>^\sss}\sF\rr{\<\na_{\fz,\fv}\>v^\alpha\rr{\udl{f}+f^0}}\right.\\&\left.(\tau,(\upsigma-\upnu)\bfe,\xi)\right\|_{\lta L^2_\xi}d\tau+\int_{t_\upmu}^te^{0.124\<\f{\upmu\tau-\mathsf{A}_k}{k^\f13}\>^\sss}\nr{\wh{\rho}(\tau,\upmu\bfe)}\tau\nnr{e^{\f18\<\f{\upxi}{A_*}\>^\sss+\<\mrxi\>^\sss}\sF\rr{\<\na_{\fz,\fv}\>v^\alpha f_{*,err}}(\tau,(\upsigma-\upmu)\bfe,\xi)}_{\lta L^2_\xi}d\tau\\
		&+\sum_{\upnu\neq0}\int_{t_\upmu}^te^{-\f12\<\f{t_\upmu}{A_*}\>^\sss}\<\upnu,\upnu\tau\>^{-10}\nnr{e^{\f18\<\f{\upxi-\mathsf{A}_k}{k^\f13}\>^\sss+\<\mrxi\>^\sss}\sF\rr{\<\na_{\fz,\fv}\>v^\alpha f}(\tau,(\upsigma-\upnu)\bfe,\xi)}_{\lta L^2_\xi}d\tau\leq e^{-99.5\mathsf{A}_k^{0.34}}.
		\ea\]
		Thus we get \eqref{fref} and \eqref{udlsr3}.
		\smallskip

		\noindent\textit{Step 2. The Volterra resolvent.}
		We estimate the resolvent $R(t,\upmu\bfe)$ associated with the kernel $t\sF\rr{f^0+\udl{f_{\infty}}}(0,\upmu t\bfe)\vv{1}_{t>0}$ for $\upmu\neq0$:
		\begin{equation}\label{DefR}
			R(t,\upmu\bfe)=t\sF\rr{f^0+f_{*,\infty}}(0,\upmu t\bfe)\vv{1}_{t>0}-\int_0^tR(\tau,\upmu\bfe)(t-\tau)\sF\rr{f^0+f_{*,\infty}}(0,\upmu(t-\tau)\bfe)d\tau.
		\end{equation}
		We remark that $\wh{\udl{f_{\infty}}}(0,\xi)=\wh{f_{*,\infty}}(0,\xi)$. To simplify notation we set
		\[K^0_{\upmu}(t)=t\wh{f^0}(\upmu t\bfe)\vv{1}_{t>0},\quad \udl{K}_{\upmu}(t)=t\wh{f_{*,\infty}}(0,\upmu t\bfe)\vv{1}_{t>0}, \quad R_{\upmu}(t)=R(t,\upmu\bfe),\quad R^0_{\upmu}(t)=R^0(t,\upmu\bfe),\]
		where $R^0(t,\upmu\bfe)$ is defined in \eqref{Rz} and satisfies the estimate \eqref{eRz}. Let $*$ stand for the convolution in $t$. Then $R_{\upmu}(t)$ obeys
		\[R_{\upmu}=K^0_{\upmu}+\udl{K}_{\upmu}-R_{\upmu}*(K^0_{\upmu}+\udl{K}_{\upmu})=\left\{K^0_{\upmu}+\udl{K}_{\upmu}-R_{\upmu}*\udl{K}_{\upmu}\right\}-R_{\upmu}*K^0_{\upmu}.\]
		Recalling the definition of $R^0(t,\upmu\bfe)$(see \eqref{Rz}), it yields
		\[\ba &R_{\upmu}=\left\{K^0_{\upmu}+\udl{K}_{\upmu}-R_{\upmu}*\udl{K}_{\upmu}\right\}-\left\{K^0_{\upmu}+\udl{K}_{\upmu}-R_{\upmu}*\udl{K}_{\upmu}\right\}*R^0_{\upmu}\\&\quad=(R^0_{\upmu}+\udl{K}_{\upmu}-\udl{K}_{\upmu}*R^0_{\upmu})-\rr{\udl{K}_{\upmu}-\udl{K}_{\upmu}*R^0_{\upmu}}*R_{\upmu}.\ea\]
		If we define $\udl{R}_{\upmu}(t)$ as the solution to
		\begin{equation}\label{Rstar}
			\udl{R}_{\upmu}=\rr{\udl{K}_{\upmu}-\udl{K}_{\upmu}*R^0_{\upmu}}-\rr{\udl{K}_{\upmu}-\udl{K}_{\upmu}*R^0_{\upmu}}*\udl{R}_{\upmu},
		\end{equation}
		then $R_{\upmu}$ satisfies
		\[R_{\upmu}=\left\{R^0_{\upmu}+\udl{K}_{\upmu}-\udl{K}_{\upmu}*R^0_{\upmu}\right\}-\left\{R^0_{\upmu}+\udl{K}_{\upmu}-\udl{K}_{\upmu}*R^0_{\upmu}\right\}*\udl{R}_{\upmu}=R^0_{\upmu}+\udl{R}_{\upmu}-R^0_{\upmu}*\udl{R}_{\upmu}.\]
		
		Now we estimate $\udl{R}_{\upmu}(t)$. By \eqref{udlf} and \eqref{eRz} we have
		\[\nr{\udl{K}_{\upmu}-\udl{K}_{\upmu}*R^0_{\upmu}}\ls_\theta\eps_0\eps_*\f{1}{|\upmu|}\<\upmu t\>^{-8}e^{-\<\f{\upmu t}{A_*}\>^\sss}+\int_0^t\eps_0\eps_*\f{1}{|\upmu|}\<\upmu \tau\>^{-8}e^{-\<\f{\upmu \tau}{A_*}\>^\sss}e^{-\f12|\upmu(t-\tau)|}d\tau\ls\eps_0\eps_*\f{1}{|\upmu|}\<\upmu t\>^{-8}e^{-\<\f{\upmu t}{A_*}\>^\sss}.\]		
		Applying a standard contraction mapping argument to the equation \eqref{Rstar} gives
		\[\nr{\udl{R}_{\upmu}(t)}\ls_\theta \eps_0\eps_*\f{1}{|\upmu|}\<\upmu t\>^{-8}e^{-\<\f{\upmu t}{A_*}\>^\sss}.\]
		Together with \eqref{eRz} this yields the following estimate for $R(t,\upmu\bfe)$
		\begin{equation}\label{eR}
			\nr{R(t,\upmu\bfe)}\ls_\theta \f{1}{|\upmu|}\rr{e^{-\f12|\upmu t|}+\eps_0\eps_*\<\upmu t\>^{-8}e^{-\<\f{\upmu t}{A_*}\>^\sss}}.
		\end{equation}
		
		\noindent\textit{Step 3. Extracting the recurrence.}
		By \eqref{DefR}, we rewrite \eqref{rhomu} as
		\begin{equation}\label{rhotmu}
			\wh{\rho}(t,\upmu\bfe)=\wh{f}(t_\upmu,\upmu\bfe,\upmu t\bfe)-\int_\R\wh{f}(t_\upmu,\upmu\bfe,\upmu \tau\bfe)R(t-\tau,\upmu\bfe)\vv{1}_{t-\tau>0}d\tau+\wh{\udl{\sR}_{21}}(t,\upmu\bfe).
		\end{equation}
		Here $\wh{\udl{\sR}_{21}}(t,\upmu\bfe)=\wh{\udl{\sR}_{2}}(t,\upmu\bfe)-\int_{t_\upmu}^t \wh{\udl{\sR}_{2}}(\tau,\upmu\bfe)R(t-\tau,\upmu\bfe)d\tau+\int_\R\wh{f}(t_\upmu,\upmu\bfe,\upmu \tau\bfe)R(t-\tau,\upmu\bfe)\vv{1}_{\tau\le t_\upmu}d\tau$. The convolution bound \eqref{eR}, together with \eqref{udlsr12} and \eqref{zsd}, shows that both added terms satisfy the same weighted estimate as $\udl{\sR}_2$. Hence $\udl{\sR}_{21}$ satisfies \eqref{udlsr12}. Combining \eqref{fref} and \eqref{rhotmu}, we obtain
		\[\ba
		&\wh{f}(t_{\upmu-1},\upsigma\bfe,\xi)=\wh{f}(t_\upmu,\upsigma\bfe,\xi)-\int_{t_\upmu}^{t_{\upmu-1}}\wh{\rho}(\tau,\upmu\bfe)\f{\upxi-\upsigma \tau}{\upmu}\sF\rr{ f^0+\udl{f_\infty} }((\upsigma-\upmu)\bfe,\xi-\upmu\tau\bfe) d\tau+\wh{\udl{\sR}_3}(t_{\upmu-1},\upsigma\bfe,\xi)
		\\&=\wh{f}(t_\upmu,\upsigma\bfe,\xi)-\int_{t_\upmu}^{t_{\upmu-1}}\int_\R\wh{f}(t_\upmu,\upmu\bfe,\upmu s\bfe)\rr{\delta_0(\tau-s)-R(\tau-s,\upmu\bfe)\vv{1}_{\tau-s>0}}\f{\upxi-\upsigma \tau}{\upmu}\\
		&\times
		\sF\rr{ f^0+\udl{f_\infty} }((\upsigma-\upmu)\bfe,\xi-\upmu\tau\bfe) dsd\tau+\wh{\udl{\sR}_{31}}(t_{\upmu-1},\upsigma\bfe,\xi)
		=\wh{f}(t_\upmu,\upsigma\bfe,\xi)
		-\int_\R\int_\R\wh{f}(t_\upmu,\upmu\bfe,\upmu s\bfe)(\delta_0(\tau-s)\\&-R(\tau-s,\upmu\bfe)\vv{1}_{\tau-s>0})\f{\upxi-\upsigma \tau}{\upmu} 
	 \sF\rr{ f^0+\udl{f_\infty} }((\upsigma-\upmu)\bfe,\xi-\upmu\tau\bfe) dsd\tau+\wh{\sR_{\upsigma,\upmu}}(t_{\upmu-1},\upsigma\bfe,\xi)
		\\
		&=\wh{f}(t_\upmu,\upsigma\bfe,\xi)-\int_\R\int_\R\wh{f}(t_\upmu,\upmu\bfe,\upeta\bfe)\rr{\delta_0(\upzeta-\upeta)\f{\upmu}{|\upmu|}-\f1\upmu R(\f{\upzeta-\upeta}{\upmu},\upmu\bfe)\vv{1}_{\f{\upzeta-\upeta}{\upmu}>0}}\big(\f{\upxi-\upzeta}{\upmu^2} -\f{\upsigma-\upmu}{\upmu^3}\upzeta\big)\\&\qquad\times
		\sF\rr{ f^0+\udl{f_\infty} }((\upsigma-\upmu)\bfe,\xi-\upzeta\bfe) d\upeta d\upzeta+\wh{\sR_{\upsigma,\upmu}}(t_{\upmu-1},\upsigma\bfe,\xi),
		\ea\]
		where
		\[\ba&\wh{\udl{\sR}_{31}}(t_{\upmu-1},\upsigma\bfe,\xi):=\wh{\udl{\sR}_{3}}(t_{\upmu-1},\upsigma\bfe,\xi)- \int_{t_\upmu}^{t_{\upmu-1}}\wh{\udl{\sR}_{21}}(\tau,\upmu\bfe)\f{\upxi-\upsigma \tau}{\upmu}\sF\rr{ f^0+\udl{f_\infty} }((\upsigma-\upmu)\bfe,\xi-\upmu\tau\bfe) d\tau;\\&\wh{\sR_{\upsigma,\upmu}}(t_{\upmu-1},\upsigma\bfe,\xi):=\wh{\udl{\sR}_{31}}(t_{\upmu-1},\upsigma\bfe,\xi)+ \int_\R\int_\R\wh{f}(t_\upmu,\upmu\bfe,\upmu s\bfe)(\delta_0(\tau-s)-R(\tau-s,\upmu\bfe)\vv{1}_{\tau-s>0})\f{\upxi-\upsigma \tau}{\upmu}\\&\qquad\qquad\qquad\times \sF\rr{ f^0+\udl{f_\infty} }((\upsigma-\upmu)\bfe,\xi-\upmu\tau\bfe)1_{\tau\in (-\infty,t_{\upmu})\cup (t_{\upmu-1},+\infty)} dsd\tau.\\
		\ea\]
		Taking the inverse Fourier transform and using \eqref{DefR} gives \eqref{aprefor}. The bounds \eqref{udlsr12}, \eqref{udlsr3}, and \eqref{zsd} control the three components of $\sR_{\upsigma,\upmu}$ and yield \eqref{sRsigmamu}.
		
		\noindent\textit{Step 4. Control inside the critical interval.}
		By \eqref{rhotmu} and \eqref{eR}, we first have
		\[\nnr{\wh{\rho}(t,\upmu\bfe)}_{L^2[t_\upmu,t_{\upmu-1}]}\ls_\theta \f{1}{\upmu^\f12}\nnr{\wh{f}(t_\upmu,\upmu\bfe,\upzeta\bfe)}_{L^2_\upzeta}+e^{-97\mathsf{A}_k^{0.34}}.\]
	    By \eqref{fref}, it remains to estimate the resonant integral. The Cauchy--Schwarz inequality gives
		\[\ba
		&\nr{\int_{t_\upmu}^t\wh{\rho}(\tau,\upmu\bfe)\f{\upxi-\upsigma \tau}{\upmu}e^{\<\mrxi\>^\sss}\pa_\xi^\alpha\wh{\udl{f_\infty}}((\upsigma-\upmu)\bfe,\xi-\upmu\tau\bfe) d\tau}^2\leq\rr{\f{3\mathsf{A}_k}{\upmu^2}}^2\f1\upmu\int_\R\nr{\sF\rr{e^{\<\na_{\mrv}\>^\sss}\<\na_{\fz,\fv}\>v^\alpha\udl{f_\infty}}((\upsigma-\upmu)\bfe,\xi)}d\upxi\\&\times
		\int_{t_\upmu}^t\nr{\wh{\rho}(\tau,\upmu\bfe)}^2\nr{\sF\rr{e^{\<\na_{\mrv}\>^\sss}\<\na_{\fz,\fv}\>v^\alpha\udl{f_\infty}}((\upsigma-\upmu)\bfe,\xi-\upmu\tau\bfe)} d\tau,
		\ea\]
		which implies that 
		\[\ba
		&\nnr{\int_{t_\upmu}^t\wh{\rho}(\tau,\upmu\bfe)\f{\upxi-\upsigma \tau}{\upmu}e^{\<\mrxi\>^\sss}\pa_\xi^\alpha\wh{\udl{f_\infty} }((\upsigma-\upmu)\bfe,\xi-\upmu\tau\bfe) d\tau}_{L^2_\xi}^2
		\leq\rr{\f{3\mathsf{A}_k}{\upmu^2}}^2\f1\upmu\nnr{\sF\rr{e^{\<\na_{\mrv}\>^\sss}\<\na_{\fz,\fv}\>v^\alpha\udl{f_\infty}}((\upsigma-\upmu)\bfe,\xi)}_{L^2_{\mrxi}L^1_{\upxi}}\\&\times\int_{t_\upmu}^t\nr{\wh{\rho}(\tau,\upmu\bfe)}^2\nnr{\sF\rr{e^{\<\na_{\mrv}\>^\sss}\<\na_{\fz,\fv}\>v^\alpha\udl{f_\infty}}((\upsigma-\upmu)\bfe,\xi-\upmu\tau\bfe)}_{L^1_\upxi L^2_{\mrxi}} d\tau\\
		\leq&\rr{\f{3\mathsf{A}_k}{\upmu^2}}^2\f1\upmu\nnr{\wh{\rho}(t,\upmu\bfe)}^2_{L^2[t_\upmu,t_{\upmu-1}]}\nnr{\sF\rr{\<e^{\<\na_{\mrv}\>^\sss}\na_{\fz,\fv}\>v^\alpha\udl{f_\infty}}((\upsigma-\upmu)\bfe,\xi)}^2_{L^1_{\upxi}L^2_{\mrxi}}.
		\ea\]
		Now we are led to that
		\[\ba
		&\nnr{e^{\<\mrxi\>^\sss}\rr{\pa_\xi^\alpha\wh{f}(t,\upsigma\bfe,\xi)-\pa_\xi^\alpha\wh{f}(t_\upmu,\upsigma\bfe,\xi)}}_{\lta L^2_\xi}\\
		&\quad\ls_\theta\f{\mathsf{A}_k}{\upmu^3}\nnr{\wh{f}(t_\upmu,\upmu\bfe,\upzeta\bfe)}_{L^2_\upzeta}
		\nnr{\sF\rr{e^{\<\na_{\mrv}\>^\sss}\<\na_{\fz,\fv}\>^3(v^\alpha\udl{f_\infty})}((\upsigma-\upmu)\bfe,\xi)}_{\lta L^2_{\upxi}L^2_{\mrxi}}\\
		&\qquad+\f{1}{\upmu^2}\nnr{\wh{f}(t_\upmu,\upmu\bfe,\upzeta\bfe)}_{L^2_\upzeta}
		\nnr{e^{\<\na_{\mrv}\>^\sss}\<\na_\fv\>^3v^\alpha f^0(v)}_{\lta L^2_v}\vv{1}_{\upsigma=\upmu}
		+e^{-97\mathsf{A}_k^{0.34}}\\
		&\quad\leq\f{\eps_0^{0.9}\mathsf{A}_k}{\upmu^3}\nnr{\beta(t_{\upmu},\upmu,\fv)}_{L^2_\fv}\<\upsigma-\upmu\>^{-8}
		+e^{-97\mathsf{A}_k^{0.34}}.\ea\]
		This proves \eqref{tmutupper} and completes the proof.
	\end{proof}

    Integrating \eqref{DefcDsigmamu} in $\mrv\in\R^2$ and using \eqref{Defbetatmuv}, we obtain the reduced transfer operator
    \begin{equation}\label{DefcDsigmamu2}
    	\fD_{\upsigma,\upmu}f(\fv)=2\pi i\rr{\f{1}{\upmu^2}\pa_\fv\rr{\beta^0+\udl{\beta_\infty}}(\upsigma-\upmu,\fv)J(\upmu,\fv)f(\fv)-\f{\upsigma-\upmu}{\upmu^3}\udl{\beta_\infty}(\upsigma-\upmu,\fv)\pa_\fv\rr{J(\upmu,\fv)f(\fv)}}.
    \end{equation}
   The trace theorem gives
	\begin{equation}\label{traceLinftyfv}
		\nnr{\udl{\beta_\infty}(\upsigma-\upmu,\fv)}_{L^\infty_{\fv}}\ls\nnr{\<\na_{\fv}\>v^\alpha\sF_z\udl{f_\infty}((\upsigma-\upmu)\bfe,v)}_{\lta L^2_v}.
	\end{equation}

	We next estimate the size of the cascade amplitudes $\beta(t_\upmu,\upmu,\fv)$. The leading term in \eqref{DefcDsigmamu2} contains \(\partial_{\mathrm v}\beta\). Thus, to obtain two-sided transfer estimates, we need both Fourier localization near \(\mathsf \xi=\mathsf A_k\) and control of the effective support of \(\beta\) in \(\mathrm v\).
	 To propagate this effective support along the cascade, define for $1\leq\upmu\leq k$
	\begin{equation}\label{chimufv}
		\chi_{\upmu}(\fv):=\f{1}{\int\psi(\mathfrak{u})d\mathfrak{u}}\int_{|\fv-\mathfrak{u}|\leq(2k-\upmu+\f12)k^{-\f43}}4k^\f43\psi(4k^\f43\mathfrak{u})d\mathfrak{u}.
	\end{equation}
	A direct Fourier transform gives
	\begin{equation}\label{chifour}
		\wh{\chi_{\upmu}}(\upxi)=\f{1}{\int\psi(\mathfrak{u})d\mathfrak{u}}\wh{\psi}(\f{\upxi}{4k^\f43})\f2\upxi\sin\rr{\rr{2k-\upmu+\f12}k^{-\f43}\upxi}.
	\end{equation}
	\begin{prop}\label{prop:essential-support} Let $k, A_*,\mathsf{A}_k, \mathsf{T}_k, \mathfrak{T}_k,\sss$ be the parameters defined in \eqref{DefTk} and Proposition \ref{propudlf}.
		For $1\leq\upmu\leq k$, the function $\beta(t_\upmu,\upmu,\fv)$ defined in \eqref{Defbetatmuv} satisfies
		\begin{equation}\label{esssup}
			\nnr{\beta(t_\upmu,\upmu,\fv)\rr{1-\chi_{\upmu}(\fv)}}_{L^2_\fv}\leq \vv{1}_{\upmu\leq k-1}\rr{\prod_{\upsigma=\upmu+1}^{k}\f{10C_\theta k^3}{\upsigma^3}}e^{-96\mathsf{A}_k^{0.34}};
		\end{equation}
		\begin{equation}\label{esssupp}
			\nnr{\pa_\fv\beta(t_\upmu,\upmu,\fv)\rr{1-\chi_{\upmu-1}(\fv)}}_{L^2_\fv}\leq e^{-95\mathsf{A}_k^{0.34}}.
		\end{equation}
	\end{prop}
	\begin{proof}
	We prove \eqref{esssup} by descending induction. By \eqref{Defofpsi},
	\[
	\chi_\upmu(\fv)=1\quad\text{if }|\fv|\leq(2k-\upmu+\tfrac14)k^{-4/3},
	\qquad
	\chi_\upmu(\fv)=0\quad\text{if }|\fv|\geq(2k-\upmu+\tfrac34)k^{-4/3}.
	\]
	At $\upmu=k$, the prescribed packet satisfies
	$\beta(t_k,k,\fv)=\beta(\mathsf{T}_k,k,\fv)$ and is supported in $\{|\fv|\leq k^{-1/3}\}$. Hence
	\[ \beta(t_k,k,\fv)\bigl(1-\chi_k(\fv)\bigr)=0, \]
	which proves the base case.

	Assume that \eqref{esssup} holds for $\upmu\leq\upnu\leq k$. We prove it at $\upmu-1$. Integrating \eqref{aprefor} in $\mrv$ and using \eqref{DefcDsigmamu2} together with $\beta(\mathsf{T}_k,\upnu,\fv)=0$ for $\upnu\neq k$, we obtain
		\[\beta(t_{\upmu-1},\upmu-1,\fv)=\sum_{\upnu=\upmu}^{k}\fD_{\upmu-1,\upnu}\beta(t_\upnu,\upnu,\fv)+\sR_{\upmu-1}(\fv),\quad\mbox{where}\quad\nr{\wh{\sR_{\upmu-1}}(\upxi)}\leq e^{-96.9\mathsf{A}_k^{0.34}}e^{-0.12\<\f{\upxi-\mathsf{A}_k}{k^\f13}\>^\sss}.\]
		We split the differential operator $\fD_{\upmu-1,\upnu}$ into two parts:
		\[
		\big[\fD_{\upmu-1,\upnu}\beta(t_\upnu,\upnu,\fv)\big]\rr{1-\chi_{\upmu-1}(\fv)}=\mathscr{D}_1+\mathscr{D}_2,\]
where \[\ba &\mathscr{D}_1:
=2\pi i\rr{\f1{\upnu^2}\pa_\fv\udl{\beta_\infty}(\upmu-1-\upnu,\fv)J(\upnu,\fv)-\f{\upmu-1-\upnu}{\upnu^3}\udl{\beta_\infty}(\upmu-1-\upnu,\fv)\pa_\fv J(\upnu,\fv)} 
	\beta(t_\upnu,\upnu,\fv)\rr{1-\chi_{\upmu-1}(\fv)},\\
		&\mathscr{D}_2:=-2\pi i\f{\upmu-\upnu-1}{\upnu^3}\udl{\beta_\infty}(\upmu-1-\upnu,\fv)J(\upnu,\fv)\pa_\fv\beta(t_\upnu,\upnu,\fv)\rr{1-\chi_{\upmu-1}(\fv)}.
		\ea\]
		
		For $\mathscr{D}_1$, it holds that
		\[\ba
		\nnr{\mathscr{D}_1}_{L^2_\fv}\leq\f{\eps_0 C_\theta}{\upnu^2}\<\upmu-1-\upnu\>^{-8}\nnr{\beta(t_\upnu,\upnu,\fv)\rr{1-\chi_{\upmu-1}(\fv)}}_{L^2_\fv}
		\leq\f{\eps_0 C_\theta}{\upnu^2}\<\upmu-1-\upnu\>^{-8}\nnr{\beta(t_\upnu,\upnu,\fv)\rr{1-\chi_{\upnu}(\fv)}}_{L^2_\fv}.
		\ea\]
		For $\mathscr{D}_2$, we first estimate $\nnr{\pa_\fv\rr{\beta(t_\upnu,\upnu,\fv)\rr{1-\chi_{\upmu-1}(\fv)}}}_{L^2_\fv}$.
		By \eqref{Defofpsi}, \eqref{zsd}, and \eqref{chifour},
		\begin{equation}\label{beta1chimu}
			\sF_\fv\rr{\beta(t_\upnu,\upnu,\fv)}=\wh{f}(t_\upnu,\upnu\bfe,\upxi\bfe),\quad  \nr{\sF_\fv\rr{\beta(t_\upnu,\upnu,\fv)\rr{1-\chi_{\upmu-1}(\fv)}}(\upxi) }\leq e^{6\mathsf{A}_k^{0.34}}e^{-\f14\<\f{\upxi-\mathsf{A}_k}{k^\f43}\>^\sss},
		\end{equation}
		which implies
		\[\ba
		&\nnr{\pa_\fv\rr{\beta(t_\upnu,\upnu,\fv)\rr{1-\chi_{\upmu-1}(\fv)}}}_{L^2_\fv}\leq\f1{(2\pi)^\f12}\nnr{\upxi\sF_\fv\rr{\beta(t_\upnu,\upnu,\fv)\rr{1-\chi_{\upmu-1}(\fv)}}(\upxi)\vv{1}_{|\upxi-\mathsf{A}_k|\leq \mathsf{A}_k^{0.35}k^\f43}}_{L^2_\upxi}\\
		+&\f1{(2\pi)^\f12}\nnr{\upxi\sF_\fv\rr{\beta(t_\upnu,\upnu,\fv)\rr{1-\chi_{\upmu-1}(\fv)}}(\upxi)\vv{1}_{|\upxi-\mathsf{A}_k|\geq \mathsf{A}_k^{0.35}k^\f43}}_{L^2_\upxi}\leq2\mathsf{A}_k\nnr{\beta(t_\upnu,\upnu,\fv)\rr{1-\chi_{\upmu-1}(\fv)}}_{L^2_\fv}+e^{-100\mathsf{A}_k^{0.34}}.
		\ea\]
		Now we can estimate $\mathscr{D}_2$ as follows
		\[\ba
		&\nnr{\mathscr{D}_2}_{L^2_\fv}\leq\f{\eps_0 C_\theta}{\upnu^3}\<\upmu-1-\upnu\>^{-8}\rr{\nnr{\pa_\fv\rr{\beta(t_\upnu,\upnu,\fv)\rr{1-\chi_{\upmu-1}(\fv)}}}_{L^2_\fv}+\nnr{\rr{\beta(t_\upnu,\upnu,\fv)\pa_\fv\chi_{\upmu-1}(\fv)}}_{L^2_\fv}}\\
		\leq&\f{\eps_0 C_\theta}{\upnu^3}\<\upmu-1-\upnu\>^{-8}\rr{2\mathsf{A}_k\nnr{\beta(t_\upnu,\upnu,\fv)\rr{1-\chi_{\upnu}(\fv)}}_{L^2_\fv}+k^\f53\nnr{\beta(t_\upnu,\upnu,\fv)\rr{1-\chi_{\upnu}(\fv)}}_{L^2_\fv}}+e^{-100\mathsf{A}_k^{0.34}}\\
		\leq&\f{3\eps_0 C_\theta \mathsf{A}_k}{\upnu^3}\<\upmu-1-\upnu\>^{-8}\nnr{\beta(t_\upnu,\upnu,\fv)\rr{1-\chi_{\upnu}(\fv)}}_{L^2_\fv}+e^{-100\mathsf{A}_k^{0.34}}.
		\ea\]
		Combining all the estimates, we arrive at
		\[\ba
		\nnr{\beta(t_{\upmu-1},\upmu-1,\fv)\rr{1-\chi_{\upmu-1}(\fv)}}_{L^2_\fv}\leq\sum_{\upnu=\upmu}^{k-1}\nnr{\fD_{\upmu-1,\upnu}\beta(t_\upnu,\upnu,\fv)\rr{1-\chi_{\upmu-1}(\fv)}}+e^{-96.8\mathsf{A}_k^{0.34}}\\
		\leq\sum_{\upnu=\upmu}^{k-1}\f{4C_\theta k^3}{\upnu^3}\<\upmu-1-\upnu\>^{-8}\rr{\prod_{\upsigma=\upmu+1}^{k}\f{10C_\theta k^3}{\upsigma^3}}e^{-96\mathsf{A}_k^{0.34}}+e^{-96.8\mathsf{A}_k^{0.34}}\leq\rr{\prod_{\upsigma=\upmu}^{k}\f{10C_\theta k^3}{\upsigma^3}}e^{-96\mathsf{A}_k^{0.34}}.
		\ea\]

		To prove \eqref{esssupp}, we use the same near-frequency/tail decomposition:
		\[\ba
		&\nnr{\pa_\fv\beta(t_\upmu,\upmu,\fv)\rr{1-\chi_{\upmu-1}(\fv)}}_{L^2_\fv}\leq\nnr{\pa_\fv\rr{\beta(t_\upmu,\upmu,\fv)\rr{1-\chi_{\upmu-1}(\fv)}}}_{L^2_\fv}+\nnr{\beta(t_\upmu,\upmu,\fv)\pa_\fv\chi_{\upmu-1}(\fv)}_{L^2_\fv}\\
		&\leq2\mathsf{A}_k\nnr{\beta(t_\upmu,\upmu,\fv)\rr{1-\chi_{\upmu-1}(\fv)}}_{L^2_\fv}+k^\f53\nnr{\beta(t_\upmu,\upmu,\fv)\rr{1-\chi_{\upmu}(\fv)}}_{L^2_\fv}\leq e^{-95\mathsf{A}_k^{0.34}}.
		\ea\]
		The last inequality uses \eqref{beta1chimu}. This completes the induction and the proof.
	\end{proof}
	
	The preceding proposition gives the velocity-space input for the leading echo transfer: on the effective support of $\beta$, the Gaussian coefficient in the background low mode is close to one. Meanwhile, frequency localization near $\upxi=\mathsf{A}_k$ implies $\|\pa_\fv\beta\|_2\approx\mathsf{A}_k\|\beta\|_2$. It remains to control the resolvent factor $J$, after which these two facts produce a two-sided recurrence for the cascade amplitudes.
	\begin{equation}\label{etheta}
		\etmu:=\left\{\begin{aligned}
			&1-2\f{\theta^{-2}}{\upmu^2},\quad&\upmu\geq2\theta^{-1}, \\
			&\f12\theta^2,\quad&\upmu\leq2\theta^{-1},
		\end{aligned}\right.\quad\Etmu:=\left\{\ba&1+2\f{\theta^{-2}}{\upmu^2},\quad&\upmu\geq2\theta^{-1},\\
	&2\theta^{-1},\quad&\upmu\leq2\theta^{-1},\ea\right.
	\end{equation}
	\begin{lem}\label{lem:J-bounds}
	For every $1\leq\upmu\leq k$ and $\fv\in\R$, the function $J$ defined in \eqref{DefJmuv} satisfies
	\begin{equation}\label{luofJmufv}
		\etmu\leq |J(\upmu,\fv)|\leq\Etmu.
	\end{equation}
	\end{lem}
	\begin{proof}
	For $\upmu\geq2\theta^{-1}$, we first notice that
	\[\big|\int_0^\infty t\rr{\wh{f^0}(\upmu t\bfe)+\wh{\udl{f_\infty}}(0,\upmu t\bfe)}e^{i\upmu t\fv}dt\big|\le \int_0^\infty t(\theta^{-2}e^{-2|\upmu|t}+\eps_0^{0.99}\<|\upmu|t\>^{-20})dt\le \f{0.26\theta^{-2}}{|\upmu|^2}. \]
Thus
\[1-2\f{\theta^{-2}}{\upmu^2}\le \f{1}{1+\f{0.26\theta^{-2}}{|\upmu|^2}}\le |J(\upmu,\fv)|\leq\f{1}{1-\f{0.26\theta^{-2}}{|\upmu|^2}}\le 1+2\f{\theta^{-2}}{\upmu^2}.\]
	For $1\leq\upmu\leq2\theta^{-1}$, the direct lower bound gives
	\[|J(\upmu,\fv)|\ge \f{1}{1+0.16\theta^{-2}}\ge\f12\theta^2.\]
	The Penrose condition \eqref{Penrose}, together with the smallness of the zero spatial mode of $\udl{f_\infty}-\eps_0\vf\cos\fz$, gives
	\[|J(\upmu,\fv)|\le \f{1}{\theta-\eps_0^{\f12}}\le 2\theta^{-1}.\]
	This proves \eqref{luofJmufv}.
	\end{proof}

	\begin{prop}\label{propertybetatmumu}
		For $1\leq\upmu\leq k$, it holds that
		\begin{equation}\label{mumu1}
			\f{(1+\eps_0^6)k^3}{\upmu^3}\rr{1-10^2\theta^{-4}\eps_{*}}\etmu\leq\f{\nnr{\beta(t_{\upmu-1},\upmu-1,\fv)}_{L^2_\fv}}{\nnr{\beta(t_\upmu,\upmu,\fv)}_{L^2_\fv}}\leq\f{(1+\eps_0^6)k^3}{\upmu^3}\rr{1+10^2\theta^{-4}\eps_{*}}\Etmu,
		\end{equation}
		\begin{equation}\label{betatkk}
			\nnr{\beta(t_{k},k,\fv)}_{L^2_\fv}=\nnr{\beta(\mathsf{T}_k,k,\fv)}_{L^2_\fv}\sim \eps_k k^{-\f16}.
		\end{equation}
		In particular, for sufficiently large $k$, the sequence $\{\nnr{\beta(t_{\upmu},\upmu,\fv)}_{L^2_\fv}\}_{0\le\upmu\le k}$ is decreasing as $\upmu$ increases.
	\end{prop}
	\begin{proof}
		The proof has three steps: estimate the velocity derivative of each reduced amplitude, obtain upper and lower bounds for a single transfer, and then iterate the recurrence.

		\noindent\textit{Step 1. Size of the velocity derivative.}
		Frequency localization gives
		\begin{equation}\label{pafvbetaupper}
			\begin{aligned}
				\nnr{\pa_\fv\beta(t_\upnu,\upnu,\fv)}_{L^2_\fv}&\leq\f{1}{(2\pi)^\f12}\nnr{\upxi\wh{\beta}(t_\upnu,\upnu,\upxi)\vv{1}_{|\upxi-\mathsf{A}_k|\leq \mathsf{A}_k^{0.35}k^\f13}}_{L^2_\upxi} \\
				&\quad +\f{1}{(2\pi)^\f12}\nnr{\upxi\wh{\beta}(t_\upnu,\upnu,\upxi)\vv{1}_{|\upxi-\mathsf{A}_k|\geq \mathsf{A}_k^{0.35}k^\f13}}_{L^2_\upxi} \\
				&\leq\rr{\mathsf{A}_k+\mathsf{A}_k^{0.35}k^\f13}\f{1}{(2\pi)^\f12}\nnr{\wh{\beta}(t_\upnu,\upnu,\upxi)}_{L^2_\upxi}+e^{-100\mathsf{A}_k^{0.34}} \\
				&=\rr{\mathsf{A}_k+\mathsf{A}_k^{0.35}k^\f13}\nnr{\beta(t_\upnu,\upnu,\fv)}_{L^2_\fv}+e^{-100\mathsf{A}_k^{0.34}},
			\end{aligned}
		\end{equation}
		and
		\begin{equation}\label{pafvbetalower}
			\begin{aligned}
				\nnr{\pa_\fv\beta(t_\upnu,\upnu,\fv)}_{L^2_\fv}&\geq\f{1}{(2\pi)^\f12}\nnr{\upxi\wh{\beta}(t_\upnu,\upnu,\upxi)\vv{1}_{|\upxi-\mathsf{A}_k|\leq \mathsf{A}_k^{0.35}k^\f13}}_{L^2_\upxi} \\
				&\quad -\f{1}{(2\pi)^\f12}\nnr{\upxi\wh{\beta}(t_\upnu,\upnu,\upxi)\vv{1}_{|\upxi-\mathsf{A}_k|\geq \mathsf{A}_k^{0.35}k^\f13}}_{L^2_\upxi} \\
				&\geq\rr{\mathsf{A}_k-\mathsf{A}_k^{0.35}k^\f13}\f{1}{(2\pi)^\f12}\rr{\nnr{\wh{\beta}(t_\upnu,\upnu,\upxi)}_{L^2_\upxi}-\nnr{\wh{\beta}(t_\upnu,\upnu,\upxi)\vv{1}_{|\upxi-\mathsf{A}_k|\geq \mathsf{A}_k^{0.35}k^\f13}}_{L^2_\upxi}} \\
				&\quad -\f{1}{(2\pi)^\f12}\nnr{\upxi\wh{\beta}(t_\upnu,\upnu,\upxi)\vv{1}_{|\upxi-\mathsf{A}_k|\geq \mathsf{A}_k^{0.35}k^\f13}}_{L^2_\upxi} \\
				&\geq\rr{\mathsf{A}_k-\mathsf{A}_k^{0.35}k^\f13}\nnr{\beta(t_\upnu,\upnu,\fv)}_{L^2_\fv}-e^{-100\mathsf{A}_k^{0.34}}.
			\end{aligned}
		\end{equation}

		\noindent\textit{Step 2. Bounds for one transfer.}
		We claim that, for $1\leq\upnu\leq k$,
		\begin{equation}\label{Dmnup}
			\nnr{\fD_{\upsigma,\upnu}\beta(t_\upnu,\upnu,\fv)}_{L^2_\fv}\leq e^{-99\mathsf{A}_k^{0.34}}+\nnr{\beta(t_\upnu,\upnu,\fv)}_{L^2_\fv}\Etnu\left\{\begin{aligned}
				&\f{\eps_0 \mathsf{A}_k}{\upnu^3}\rr{1+20\eps_*},\quad&\upsigma=\upnu\pm1, \\
				&\f{50\eps_0 \mathsf{A}_k}{\upnu^3}\eps_*\<\upsigma-\upnu\>^{-8},\quad&\upsigma\neq\upnu\pm1,\upnu, \\
				&\f{10^2\theta^{-2}}{\upnu^2},\quad&\upsigma=\upnu.
			\end{aligned}\right.
		\end{equation}

		For $1\leq\upmu\leq k$, the leading transfer also satisfies the lower bound
		\begin{equation}\label{Dmulow}
			\nnr{\fD_{\upmu-1,\upmu}\beta(t_\upmu,\upmu,\fv)}_{L^2_\fv}\geq\f{\eps_0 \mathsf{A}_k}{\upmu^3}\rr{1-50\theta^{-4}\eps_*}\nnr{\beta(t_\upmu,\upmu,\fv)}_{L^2_\fv}\etmu-e^{-94\mathsf{A}_k^{0.34}}.
		\end{equation}
		We now prove the transfer bounds. Suppose first that $\upsigma\neq\upnu$. From \eqref{DefcDsigmamu2},
		\[\ba
		&\nnr{\fD_{\upsigma,\upnu}\beta(t_\upnu,\upnu,\fv)}_{L^2_\fv}\leq2\pi\rr{\f{1}{\upnu^2}\nnr{\pa_\fv\udl{\beta_\infty}(\upsigma-\upnu,\fv)J(\upnu,\fv)}_{L^\infty_\fv}+\f{1}{\upnu^3}\nnr{(\upsigma-\upnu)\udl{\beta_\infty}(\upsigma-\upnu,\fv)\pa_\fv J(\upnu,\fv)}_{L^\infty_\fv}}\\
		&\times\nnr{\beta(t_\upnu,\upnu,\fv)}_{L^2_\fv}
		+\f{2\pi}{\upnu^3}\nnr{(\upsigma-\upnu)\udl{\beta_\infty}(\upsigma-\upnu,\fv)J(\upnu,\fv)}_{L^\infty_\fv}\nnr{\pa_\fv\beta(t_\upnu,\upnu,\fv)}_{L^2_\fv}\\
		\leq&\f{C_\theta}{\upnu^2}\eps_0\<\upsigma-\upnu\>^{-8}\nnr{\beta(t_\upnu,\upnu,\fv)}_{L^2_\fv}+e^{-100\mathsf{A}_k^{0.34}}\\
		&+2\pi\f{\mathsf{A}_k+\mathsf{A}_k^{0.35}k^\f13}{\upnu^3}\nnr{\beta(t_\upnu,\upnu,\fv)}_{L^2_\fv}\Etnu\left\{\ba
		&\nnr{\eps_0\f{1}{2\pi}e^{-\f{\fv^2}{2}}}_{L^\infty_\fv}+\nnr{\beta_{*,\infty}(\upsigma-\upnu,\fv)}_{L^\infty_\fv},\quad&\upsigma=\upnu\pm1,\\
		&\nnr{(\upsigma-\upnu)\beta_{*,\infty}(\upsigma-\upnu,\fv)}_{L^\infty_\fv},\quad&\upsigma\neq\upnu\pm1,\upnu,
		\ea\right.\\
		\leq&\f{C_\theta}{\upnu^2}\eps_0\<\upsigma-\upnu\>^{-8}\nnr{\beta(t_\upnu,\upnu,\fv)}_{L^2_\fv}+e^{-100\mathsf{A}_k^{0.34}}\\
		&+\f{\eps_0\mathsf{A}_k}{\upnu^3}\nnr{\beta(t_\upnu,\upnu,\fv)}_{L^2_\fv}\Etnu\left\{\ba
		&(1+19\eps_*)\rr{1+\mathsf{A}_k^{-0.65}k^\f13} ,\quad&\upsigma=\upnu\pm1,\\
		&47\pi\eps_{*}\<\upsigma-\upnu\>^{-8},\quad&\upsigma\neq\upnu\pm1,\upnu,
		\ea\right.\\
		\leq&e^{-100\mathsf{A}_k^{0.34}}+\left\{\ba
		&\f{\eps_0\mathsf{A}_k}{\upnu^3}(1+20\eps_*) \nnr{\beta(t_\upnu,\upnu,\fv)}_{L^2_\fv}\Etnu,\quad&\upsigma=\upnu\pm1,\\
		&\f{50\eps_0\mathsf{A}_k}{\upnu^3}\eps_{*}\<\upsigma-\upnu\>^{-8}\nnr{\beta(t_\upnu,\upnu,\fv)}_{L^2_\fv}\Etnu,\quad&\upsigma\neq\upnu\pm1,\upnu.
		\ea\right.
		\ea\]
		Here we use $|\pa_\fv J(\upnu,\fv)|\leq C_\theta|J(\upnu,\fv)|$.
		For $\upsigma=\upnu$,
		\[\nnr{\fD_{\upnu,\upnu}\beta(t_\upnu,\upnu,\fv)}_{L^2_\fv}\leq\f{2\pi}{\upnu^2}\nnr{\pa_\fv\rr{\beta^0+\beta_{*,\infty}}(0,\fv)J(\upnu,\fv)}_{L^\infty_\fv}\nnr{\beta(t_\upnu,\upnu,\fv)}_{L^2_\fv}\leq\f{10\theta^{-2}}{\upnu^2}\nnr{\beta(t_\upnu,\upnu,\fv)}_{L^2_\fv}\Etnu.\]
		
		For the lower bound \eqref{Dmulow}, \eqref{DefcDsigmamu2} gives
		\[\ba&\fD_{\upmu-1,\upmu}\beta(t_\upmu,\upmu,\fv)=\f{ i}{\upmu^3}\eps_0e^{-\f{\fv^2}{2}}J(\upmu,\fv)\pa_\fv\beta(t_\upmu,\upmu,\fv)+2\pi i\left(\f1{\upmu^2}\rr{-\f{\eps_0}{2\pi}\fv e^{-\f{\fv^2}{2}}+\pa_\fv\beta_{*,\infty}(-1,\fv)}J(\upmu,\fv)\right.\\&\left.+\f{1}{\upmu^3}\rr{\f{\eps_0}{2\pi}e^{-\f{\fv^2}{2}}+\beta_{*,\infty}(-1,\fv)}\pa_\fv J(\upmu,\fv)\right)\beta(t_\upmu,\upmu,\fv)+\f{2\pi i}{\upmu^3}\beta_{*,\infty}(-1,\fv)J(\upmu,\fv)\pa_\fv\beta(t_\upmu,\upmu,\fv)=\sum_{i=1}^3\mathscr{E}_i.
		\ea\]
		By \eqref{pafvbetalower}, we first provide a lower bound of $\mathscr{E}_1$ as follows
		\[\ba
		&\nnr{\mathscr{E}_1}_{L^2_\fv}\geq\nnr{\f{ i}{\upmu^3}\eps_0e^{-\f{\fv^2}{2}}J(\upmu,\fv)\pa_\fv\beta(t_\upmu,\upmu,\fv)\chi_{\upmu-1}(\fv)}_{L^2_\fv}
		\geq\rr{1-\f52k^{-\f23}}\f{\eps_0}{\upmu^3}\etmu\nnr{\pa_\fv\beta(t_\upmu,\upmu,\fv)\chi_{\upmu-1}(\fv)}_{L^2_\fv}\\
		&\geq\rr{1-\f52k^{-\f23}}\f{\eps_0}{\upmu^3}\etmu\rr{\nnr{\pa_\fv\beta(t_\upmu,\upmu,\fv)}_{L^2_\fv}-\nnr{\pa_\fv\beta(t_\upmu,\upmu,\fv)\rr{1-\chi_{\upmu-1}(\fv)}}_{L^2_\fv}}\\
		&\geq\f{\eps_0\mathsf{A}_k}{\upmu^3}\rr{1-\f52k^{-\f23}}\rr{1-\mathsf{A}_k^{-0.65}k^\f13}\etmu\nnr{\beta(t_\upmu,\upmu,\fv)}_{L^2_\fv}-e^{-95\mathsf{A}_k^{0.34}}.
		\ea\]
		In the second inequality we use $e^{-\fv^2/2}\geq1-\frac52k^{-2/3}$ on the support of $\chi_{\upmu-1}$ and $|J(\upmu,\fv)|\geq\etmu$; the last inequality uses \eqref{esssupp}. Thus
		\[\ba
		&\nnr{\fD_{\upmu-1,\upmu}\beta(t_\upmu,\upmu,\fv)}_{L^2_\fv}\geq\nnr{\mathscr{E}_1}_{L^2_\fv}- \nnr{\mathscr{E}_2+\mathscr{E}_3}_{L^2_\fv}\geq\f{\eps_0\mathsf{A}_k}{\upmu^3}\rr{1-3k^{-\f23}}\nnr{\beta(t_\upmu,\upmu,\fv)}_{L^2_\fv}\etmu\\
		&-e^{-95\mathsf{A}_k^{0.34}}-\f{20\eps_0\eps_*\mathsf{A}_k}{\upmu^3}\nnr{\beta(t_\upmu,\upmu,\fv)}_{L^2_\fv}\Etmu\geq\f{\eps_0 \mathsf{A}_k}{\upmu^3}\rr{1-50\theta^{-4}\eps_*}\nnr{\beta(t_\upmu,\upmu,\fv)}_{L^2_\fv}\etmu-e^{-94\mathsf{A}_k^{0.34}}.
		\ea\]
		
		\noindent\textit{Step 3. Iterating the cascade.}
		By \eqref{Dmnup}, \eqref{Dmulow}, and $\beta(\mathsf{T}_k,\upnu,\fv)=0$ for $\upnu\neq k$, the first transfer satisfies
		\[\ba
		&\nnr{\beta(t_{k-1},k-1,\fv)}_{L^2_\fv}\leq\nnr{\fD_{k-1,k}\beta(t_k,k,\fv)}_{L^2_\fv}+e^{-96\mathsf{A}_k^{0.34}}\\\leq&\f{\eps_0 \mathsf{A}_k}{k^3}\rr{1+20\eps_*}\nnr{\beta(t_k,k,\fv)}_{L^2_\fv}\Etmu
		+e^{-95\mathsf{A}_k^{0.34}}	\leq \f{(1+\eps_0^6)k^3}{k^3}\rr{1+10^2\eps_*}\nnr{\beta(t_k,k,\fv)}_{L^2_\fv}\Etmu
		\ea\] and
		\[\ba
		&\nnr{\beta(t_{k-1},k-1,\fv)}_{L^2_\fv}\geq\nnr{\fD_{k-1,k}\beta(t_k,k,\fv)}_{L^2_\fv}-e^{-96\mathsf{A}_k^{0.34}}\\\geq&\f{\eps_0 \mathsf{A}_k}{k^3}\rr{1-50\theta^{-4}\eps_*}\nnr{\beta(t_k,k,\fv)}_{L^2_\fv}\etmu
		-e^{-93\mathsf{A}_k^{0.34}}	\geq \f{(1+\eps_0^6)k^3}{k^3}\rr{1-10^2\theta^{-4}\eps_*}\nnr{\beta(t_k,k,\fv)}_{L^2_\fv}\etmu.
		\ea\]
		This proves \eqref{mumu1} for $\upmu=k$. Assume it holds for $\upmu+1\leq\upnu\leq k$. The resulting bounds show that $\{\nnr{\beta(t_{\upnu},\upnu,\fv)}_{L^2_\fv}\}_{\upmu\leq\upnu\leq k}$ decreases with $\upnu$ for sufficiently large $k$. To prove \eqref{mumu1} at $\upmu$, we use
		\[\ba
		&\nnr{\beta(t_{\upmu-1},\upmu-1,\fv)}_{L^2_\fv}\leq\sum_{\upnu=\upmu+1}^{k}\nnr{\fD_{\upmu-1,\upnu}\beta(t_\upnu,\upnu,\fv)}_{L^2_\fv}+\nnr{\fD_{\upmu-1,\upmu}\beta(t_\upmu,\upmu,\fv)}_{L^2_\fv}+e^{-96\mathsf{A}_k^{0.34}}\\
		\leq&\sum_{\upnu=\upmu+1}^{k}\f{50\eps_0 \mathsf{A}_k}{\upnu^3}\eps_*\<\upmu-1-\upnu\>^{-8}\nnr{\beta(t_\upnu,\upnu,\fv)}_{L^2_\fv}\Etnu+\f{\eps_0 \mathsf{A}_k}{\upmu^3}\rr{1+20\eps_*}\nnr{\beta(t_\upmu,\upmu,\fv)}_{L^2_\fv}\Etmu+e^{-95\mathsf{A}_k^{0.34}}\\
		\leq&\sum_{\upnu=\upmu+1}^{k}\f{50\eps_0 \mathsf{A}_k}{\upmu^3}\eps_*\<\upmu-1-\upnu\>^{-8}\nnr{\beta(t_\upmu,\upmu,\fv)}_{L^2_\fv}\Etmu+\f{\eps_0 \mathsf{A}_k}{\upmu^3}\rr{1+20\eps_*}\nnr{\beta(t_\upmu,\upmu,\fv)}_{L^2_\fv}\Etmu+e^{-95\mathsf{A}_k^{0.34}}\\
		\leq&\f{(1+\eps_0^6)k^3}{\upmu^3}(1+10^2\eps_{*})\nnr{\beta(t_\upmu,\upmu,\fv)}_{L^2_\fv}\Etmu,
		\ea\]
		and
		\[\ba
		&\nnr{\beta(t_{\upmu-1},\upmu-1,\fv)}_{L^2_\fv}\geq\nnr{\fD_{\upmu-1,\upmu}\beta(t_\upmu,\upmu,\fv)}_{L^2_\fv}-\sum_{\upnu=\upmu+1}^{k}\nnr{\fD_{\upmu-1,\upnu}\beta(t_\upnu,\upnu,\fv)}_{L^2_\fv}-e^{-96\mathsf{A}_k^{0.34}}\\
		\geq&\f{\eps_0 \mathsf{A}_k}{\upmu^3}\rr{1-50\theta^{-4}\eps_*}\nnr{\beta(t_\upmu,\upmu,\fv)}_{L^2_\fv}\etmu-\sum_{\upnu=\upmu+1}^{k}\f{50\eps_0 \mathsf{A}_k}{\upnu^3}\eps_*\<\upmu-1-\upnu\>^{-8}\nnr{\beta(t_\upnu,\upnu,\fv)}_{L^2_\fv}\Etnu-e^{-93\mathsf{A}_k^{0.34}}\\
		\geq&\f{(1+\eps_0^6)k^3}{\upmu^3}\rr{1-10^2\theta^{-4}\eps_{*}}\nnr{\beta(t_\upmu,\upmu,\fv)}_{L^2_\fv}\etmu.
		\ea\]
		This completes the induction and proves the proposition.
	\end{proof}
	
	We now transfer the estimates for the reduced amplitude $\beta$ back to the full profile. The result gives uniform control during the cascade, a lower bound at each completed transfer, and stabilization after the last critical time.
	\begin{prop}\label{prop:forward-profile}  Recall that for $0\leq\upmu\leq k$, $t_\upmu=\f{\mathsf{A}_k}{\upmu+\f12}$. Then we have
	\begin{itemize}
	\item[(i).] \textnormal{(Control during the cascade.)} For any $\upsigma\in\Z$ and $t\in[t_\upmu,t_{\upmu-1}]$ with $1\leq\upmu\leq k$,
\begin{equation}\label{uppsigma}
	\nnr{e^{\<\na_{\mrv}\>^\sss}\rr{v^{\alpha}\sF_zf(t,\upsigma\bfe,v)}}_{\lta L^2_v}\leq 5\mathsf{A}_k\nnr{\beta(t_\upmu,\upmu,\fv)}_{L^2_\fv}.
\end{equation}
			 \item[(ii).] \textnormal{(Post-resonance stabilization.)} For any $\upsigma\in\Z$ and $t\geq t_0$,
		\begin{equation}\label{tztinfty}
			\nnr{e^{\<\na_{\mrv}\>^\sss}v^{\alpha}\rr{\sF_zf(t,\upsigma\bfe,v)-\sF_zf(t_0,\upsigma\bfe,v)}}_{\lta L^2_v}\leq e^{-98\mathsf{A}_k^{0.34}}.
		\end{equation}
			 \item[(iii).] \textnormal{(Lower bound after one transfer.)} For any $1\leq\upmu\leq k$,
		\begin{equation}\label{ecq9}
			\nnr{\sF_zf(t_{\upmu-1},(\upmu-1)\bfe,v)}_{L^2_v}\geq\f{1}{2\pi^\f12}\rr{1-300\theta^{-4}\eps_*}\f{\etmu}{\Etmu}\nnr{\beta(t_{\upmu-1},\upmu-1,\fv)}_{L^2_\fv}.
		\end{equation}
			\end{itemize}
	\end{prop}
	\begin{proof}
		\noindent\textit{Control during the cascade.}
		To prove \eqref{uppsigma}, we first claim that
		\begin{equation}\label{DmnupHR}
			\nnr{e^{\<\na_{\mrv}\>^\sss}v^\alpha\cD_{\upsigma,\upnu}\beta(t_\upnu,\upnu,\fv)}_{\lta L^2_v}\ls e^{-99\mathsf{A}_k^{0.34}}+\nnr{\beta(t_\upnu,\upnu,\fv)}_{L^2_\fv}\Etnu\left\{\begin{aligned}
				&\f{\eps_0 \mathsf{A}_k}{\upnu^3},\quad&\upsigma=\upnu\pm1, \\
				&\f{\eps_0\eps_* \mathsf{A}_k}{\upnu^3}\<\upsigma-\upnu\>^{-8},\quad&\upsigma\neq\upnu\pm1,\upnu, \\
				&\f{\theta^{-2}}{\upnu^2},\quad&\upsigma=\upnu.
			\end{aligned}\right.
		\end{equation}
		This is the full-profile counterpart of \eqref{Dmnup}. For example, when $\upsigma\neq\upnu\pm1,\upnu$, equations \eqref{DefcDsigmamu}, \eqref{pafvbetaupper}, and \eqref{udlf} give
		\[\ba
		&\nnr{e^{\<\na_{\mrv}\>^\sss}v^\alpha\cD_{\upsigma,\upnu}\beta(t_\upnu,\upnu,\fv)}_{\lta L^2_v}\leq
		2\pi\left(\f{1}{\upnu^2}\nnr{e^{\<\na_{\mrv}\>^\sss}v^\alpha\pa_\fv\sF_z\udl{f_\infty}((\upsigma-\upnu)\bfe,v)}_{\lta L^{\infty}_{\fv} L^2_{\mrv}}\nnr{J(\upnu,\fv)}_{L^\infty_\fv}\right.\\
		&\left.+\f{1}{\upnu^3}\nnr{(\upsigma-\upnu)e^{\<\na_{\mrv}\>^\sss}v^\alpha\sF_z\udl{f_\infty}((\upsigma-\upnu)\bfe,v)}_{\lta L^{\infty}_{\fv} L^2_{\mrv}}\nnr{\pa_\fv J(\upnu,\fv)}_{L^\infty_\fv}\right)\nnr{\beta(t_\upnu,\upnu,\fv)}_{L^2_\fv}\\
		&+\f{2\pi}{\upnu^3}\nnr{(\upsigma-\upnu)e^{\<\na_{\mrv}\>^\sss}v^\alpha\sF_z\udl{f_\infty}((\upsigma-\upnu)\bfe,v)}_{\lta L^{\infty}_{\fv} L^2_{\mrv}}\nnr{J(\upnu,\fv)}_{L^\infty_\fv}\nnr{\pa_\fv\beta(t_\upnu,\upnu,\fv)}_{L^2_\fv}.
		\ea\]
		Now from this together with \eqref{aprefor}, \eqref{tmutupper} and Proposition \ref{propertybetatmumu}, we get that
		\[\ba
		&\nnr{e^{\<\na_{\mrv}\>^\sss}v^{\alpha}\sF_zf(t,\upsigma\bfe,v)}_{\lta L^2_v}\leq\nnr{e^{\<\na_{\mrv}\>^\sss}v^{\alpha}\sF_zf(\mathsf{T}_k,\upsigma\bfe,v)}_{\lta L^2_v}+\sum_{\upnu=\upmu+1}^{k}\nnr{e^{\<\na_{\mrv}\>^\sss}v^{\alpha}\cD_{\upsigma,\upnu}\beta(t_\upnu,\upnu,\fv)}_{\lta L^2_v}\\
		&+\nnr{e^{\<\na_{\mrv}\>^\sss}\rr{v^{\alpha}\sF_zf(t,\upsigma\bfe,v)-v^{\alpha}\sF_zf(t_\upmu,\upsigma\bfe,v)}}_{\lta L^2_v}+e^{-97\mathsf{A}_k^{0.34}}
		\leq \nnr{e^{\<\na_{\mrv}\>^\sss}v^{\alpha}\sF_zf(\mathsf{T}_k,\upsigma\bfe,v)}_{\lta L^2_v}\\
		&+\sum_{\upnu=\upmu+1}^{k}\f{\mathsf{A}_k}{\upnu^3}\<\upsigma-\upnu\>^{-8}\nnr{\beta(t_\upnu,\upnu,\fv)}_{L^2_\fv}+\f{\mathsf{A}_k}{\upmu^3}\nnr{\beta(t_\upmu,\upmu,\fv)}_{L^2_\fv}+e^{-90\mathsf{A}_k^{0.34}}
		\leq5\mathsf{A}_k\nnr{\beta(t_\upmu,\upmu,\fv)}_{L^2_\fv}.
		\ea\]

		\noindent\textit{Lower bound after one transfer.}
		By \eqref{aprefor}, \eqref{Dmulow}, \eqref{Dmnup}, and $\sF_zf(\mathsf{T}_k,(\upmu-1)\bfe,v)=0$,
		\[\ba
		&\nnr{\sF_zf(t_{\upmu-1},(\upmu-1)\bfe,v)}_{L^2_v}\geq\nnr{\cD_{\upmu-1,\upmu}\beta(t_\upmu,\upmu,\fv)}_{L^2_v}-\sum_{\upnu=\upmu+1}^{k}\nnr{\cD_{\upmu-1,\upnu}\beta(t_\upnu,\upnu,\fv)}_{L^2_v}-e^{-96\mathsf{A}_k^{0.34}}\\
		&\geq\f{\eps_0\mathsf{A}_k}{2\pi^\f12\upmu^3}\rr{1-10^2\theta^{-4}\eps_*}\etmu\nnr{\beta(t_\upmu,\upmu,\fv)}_{L^2_\fv}.
		\ea\]
		Combining this estimate with \eqref{mumu1} proves \eqref{ecq9}.
	 
		\noindent\textit{Post-resonance stabilization.}
		For $t\geq t_0$, every spatial mode satisfies $|\upnu t-\mathsf{A}_k|\geq\mathsf{A}_k$. Hence \eqref{zsd} implies
		\[\nr{\wh{\rho}(t,\upnu\bfe)}\leq e^{-100\mathsf{A}_k^{0.34}}e^{-0.1249\<\f{\upnu t-\mathsf{A}_k}{k^\f13}\>^\sss}.\]
		From this together with \eqref{zsd}, \eqref{udlf} and \eqref{udlrho}, \eqref{whfeq} with $Q=0$ yields that
		\[\ba
		\nnr{e^{0.124\<\f{\na_\fv-\mathsf{A}_k}{k^\f13}\>^\sss+\<\na_{\mrv}\>^\sss}v^{\alpha}\rr{\sF_zf(t,\upsigma\bfe,v)-\sF_zf(t_0,\upsigma\bfe,v)}}_{\lta L^2_v}\leq e^{-99\mathsf{A}_k^{0.34}}.
		\ea\]
		This proves \eqref{tztinfty} and completes the proof.
	\end{proof}
	
	\subsection{Backward stability on $[\mathfrak{T}_k,\mathsf{T}_k]$}
	Set
		\[M_k:=\left\lfloor\mathsf{A}_k^{0.35}k^{\f13}\right\rfloor,
		\qquad \mathsf{T}_{1,k}:=\f{\mathsf{A}_k}{M_k+\f12}.\]
		The splitting at $\mathsf{T}_{1,k}$ separates two different backward dynamics. On $[\mathsf{T}_{1,k},\mathsf{T}_k]$, the resonance relation $\upnu t\approx\mathsf{A}_k$ can be met by the discrete modes $k+1\leq\upnu\leq M_k$. Their critical times $t_\upnu$ must therefore be resolved one by one through a backward recurrence. On the earlier interval $[\mathfrak{T}_k,\mathsf{T}_{1,k}]$, no Fourier mode can remain close to both the spatial center $k$ and the velocity-frequency center $\mathsf{A}_k$: if $\upnu\geq M_k$, then $|\upnu-k|/k^{1/3}\gtrsim\mathsf{A}_k^{0.35}$, while if $\upnu\leq M_k$, then $|\upnu t-\mathsf{A}_k|\gtrsim\mathsf{A}_k^{0.35}k^{1/3}$. Frequency localization consequently makes the density exponentially small there, so a global perturbative estimate replaces the mode-by-mode recurrence. We first analyze the recurrence on $[\mathsf{T}_{1,k},\mathsf{T}_k]$.

	\begin{prop}\label{prop:backward-recurrence}  Let $k, A_*,\mathsf{A}_k, \mathsf{T}_k, \mathfrak{T}_k,\sss$ be the parameters defined in \eqref{DefTk} and Proposition \ref{propudlf}.
		Let $t_\upmu=\f{\mathsf{A}_k}{\upmu+\f12}$ with $k+1\leq\upmu\leq M_k$. Then the following approximate recurrence relation holds on $[t_\upmu,t_{\upmu-1}]$:
		\begin{equation}\label{aprefor2}
			\sF_zf(t_{\upmu},\upsigma\bfe,v)=\sF_zf(t_{\upmu-1},\upsigma\bfe,v)+\cE_{\upsigma,\upmu}\beta(t_{\upmu-1},\upmu,\fv)+\sR_{\upsigma,\upmu}(v).
		\end{equation}
More precisely:
 \begin{itemize}
 	\item For $\upsigma\in\Z$ and $\upmu\in\Z\backslash\{0\}$, the backward transfer operator $\cE_{\upsigma,\upmu}$ is defined by
	\begin{equation}\label{DefcEsigmamu}
		\cE_{\upsigma,\upmu}f(v)=2\pi i\rr{\f{1}{\upmu^2}\pa_\fv\sF_z(f^0+\udl{f_\infty})((\upsigma-\upmu)\bfe,v)J(-\upmu,\fv)f(\fv)-\f{\upsigma-\upmu}{\upmu^3}\sF_z\udl{f_\infty}((\upsigma-\upmu)\bfe,v)\pa_\fv\rr{J(-\upmu,\fv)f(\fv)}}.\,\,
	\end{equation}
   Here $J$ is defined in \eqref{DefJmuv}, and $\beta$ is defined in \eqref{Defbetatmuv}.
	\item The error term $\sR_{\upsigma,\upmu}(v)$ satisfies
	\begin{equation}\label{sEsigmamu}
		\nnr{e^{0.12\<\f{\na_\fv-\mathsf{A}_k}{k^\f13}\>^\sss+\<\na_{\mrv}\>^\sss}\rr{v^{\alpha}\sR_{\upsigma,\upmu}}}_{\lta L^2_v}\leq e^{-97.5\mathsf{A}_k^{0.34}}.
	\end{equation}
	\item On the whole critical interval $[t_\upmu,t_{\upmu-1}]$, the profile variation satisfies
		\begin{equation}\label{tmutupper2}
			\nnr{e^{\<\na_{\mrv}\>^\sss}v^{\alpha}\rr{\sF_zf(t,\upsigma\bfe,v)-\sF_zf(t_{\upmu-1},\upsigma\bfe,v)}}_{\lta L^2_v}\leq\f{\mathsf{A}_k}{\upmu^3}\nnr{\beta(t_{\upmu-1},\upmu,\fv)}_{L^2_\fv}\<\upsigma-\upmu\>^{-8}+e^{-96.9\mathsf{A}_k^{0.34}}.
		\end{equation}
	 \end{itemize}
	\end{prop}
	\begin{proof}
	The proof uses the argument of Proposition \ref{prop:forward-recurrence} with the time orientation reversed. We record the three changes needed below.
	
	\noindent\textit{Backward Step 1. Frequency-space approximation.}
	The identities corresponding to \eqref{rerho}--\eqref{fref} are
		\[\ba
		&\wh{\rho}(t,\upnu\bfe)= \wh{\udl{\sR}_1}(t,\upnu\bfe),\quad \text{for}\ \upnu\neq\upmu,\\
		&\wh{\rho}(t,\upmu\bfe)=\wh{f}(t_{\upmu-1},\upmu\bfe,\upmu t\bfe)+\int_t^{t_{\upmu-1}}\wh{\rho}(\tau,\upmu\bfe)(t-\tau)\rr{\wh{f^0}(\upmu (t-\tau)\bfe)+\wh{\udl{f_\infty}}(0,\upmu (t-\tau)\bfe)}d\tau+\wh{\udl{\sR}_2}(t,\upmu\bfe),\\
		&\wh{f}(t,\upsigma\bfe,\xi)=\wh{f}(t_{\upmu-1},\upsigma\bfe,\xi)+\int^{t_{\upmu-1}}_t\wh{\rho}(\tau,\upmu\bfe)\f{\upxi-\upsigma \tau}{\upmu}\sF\rr{ f^0+\udl{f_\infty} }((\upsigma-\upmu)\bfe,\xi-\upmu\tau\bfe) d\tau+\wh{\udl{\sR}_3}(t,\upsigma\bfe,\xi).
		\ea\]

	\noindent\textit{Backward Step 2. Volterra inversion.}
	Using \eqref{DefR} and \eqref{eR}, the density equation becomes
		\[\wh{\rho}(t,\upmu\bfe)=\wh{f}(t_{\upmu-1},\upmu\bfe,\upmu t\bfe)-\int_\R\wh{f}(t_{\upmu-1},\upmu\bfe,\upmu \tau\bfe)R(-(t-\tau),-\upmu\bfe)\vv{1}_{t-\tau<0}d\tau+\wh{\udl{\sR}_{21}}(t,\upmu\bfe).\]
		
	\noindent\textit{Backward Step 3. Recurrence and errors.}
	Combining the preceding identities gives
		\[\ba
		&\wh{f}(t_{\upmu},\upsigma\bfe,\xi)
		=\wh{f}(t_{\upmu-1},\upsigma\bfe,\xi)+\wh{\sR_{\upsigma,\upmu}}(t_\upmu,\upsigma\bfe,\xi)+\int_\R\int_\R\wh{f}(t_{\upmu-1},\upmu\bfe,\upeta\bfe)\rr{\delta_0(\upzeta-\upeta)\f{\upmu}{|\upmu|}-\f1\upmu R(-\f{\upzeta-\upeta}{\upmu},-\upmu\bfe)\vv{1}_{\f{\upzeta-\upeta}{\upmu}<0}}\\
		&\left(\f{\upxi-\upzeta}{\upmu^2}-\f{\upsigma-\upmu}{\upmu^3}\upzeta\right)
		\sF\rr{ f^0+\udl{f_\infty} }((\upsigma-\upmu)\bfe,\xi-\upzeta\bfe) d\upeta d\upzeta.
		\ea\]
		Then \eqref{aprefor2} follows by performing the inverse Fourier transform. Tracking the decay margins in the same estimates used for \eqref{sRsigmamu} and \eqref{tmutupper} gives the strengthened bounds \eqref{sEsigmamu} and \eqref{tmutupper2}.
	\end{proof}
	
	Next we estimate the backward recurrence coefficients on $[\mathsf{T}_{1,k},\mathsf{T}_k]$.
		 \begin{prop}\label{propertybeta1} Recall that $t_\upmu=\f{\mathsf{A}_k}{\upmu+\f12}$ with $k+1\leq\upmu\leq M_k$. Then:
	 \begin{itemize}
	 	\item[(i).] Integrating \eqref{aprefor2} in $\mrv\in\R^2$ gives
		\begin{equation}\label{betaindu}
			\beta(t_{\upmu},\upmu+1,\fv)=\sum_{\upnu=k+2}^{\upmu}\fE_{\upmu+1,\upnu}\beta(t_{\upnu-1},\upnu,\fv)+\sR_{\upmu}(\fv),\quad\mbox{where}\quad\nr{\wh{\sR_{\upmu}}(\upxi)}\leq e^{-96.9\mathsf{A}_k^{0.34}}e^{-0.12\<\f{\upxi-\mathsf{A}_k}{k^\f13}\>^\sss},
		\end{equation}
\begin{equation}\label{DeffEsigmamu}
	\fE_{\upsigma,\upmu}f(\fv):=2\pi i\rr{\f{1}{\upmu^2}\pa_\fv\rr{\beta^0+\udl{\beta_\infty}}(\upsigma-\upmu,\fv)J(-\upmu,\fv)f(\fv)-\f{\upsigma-\upmu}{\upmu^3}\udl{\beta_\infty}(\upsigma-\upmu,\fv)\pa_\fv\rr{J(-\upmu,\fv)f(\fv)}}.
\end{equation}
\item[(ii).] The estimates \eqref{Dmnup} and \eqref{DmnupHR} remain valid after replacing $(\fD,\cD,J(\upnu,\cdot),\beta(t_\upnu,\upnu))$ by $(\fE,\cE,J(-\upnu,\cdot),\beta(t_{\upnu-1},\upnu))$, respectively.
		 \item[(iii).] For every $k+1\leq\upmu\leq M_k$,
			\begin{equation}\label{zpm}
				\nnr{\beta(t_{\upmu-1},\upmu,\fv)}_{L^2_\fv}\leq e^{-96.8\mathsf{A}_k^{0.34}}.
			\end{equation}
	 \end{itemize}
		\end{prop}

	\begin{proof}
	Part (i) follows from \eqref{aprefor2}, \eqref{DeffEsigmamu}, and $\beta(t_k,\upmu,\fv)=0$ for $\upmu\neq k$. For part (ii), Lemma \ref{lem:J-bounds} and the corresponding derivative bound for $J$ are unchanged when $\upnu$ is replaced by $-\upnu$. The localization estimates for $\beta(t_{\upnu-1},\upnu)$ are also identical. Substitution into \eqref{DeffEsigmamu} and \eqref{DefcEsigmamu} therefore reproduces the proofs of \eqref{Dmnup} and \eqref{DmnupHR} with the stated replacements.

	     \noindent\textit{The endpoint $\upmu=k+1$.}
	     Here $t_{\upmu-1}=t_k=\mathsf{T}_k$ and
	     \[
	     \beta(t_k,k+1,\fv)=0,
	     \]
	     because the prescribed packet at $\mathsf{T}_k$ is supported only at spatial modes $\pm k$. Hence \eqref{zpm} holds for $\upmu=k+1$. Moreover, taking $\upsigma=k+1$ and $\upmu=k+1$ in \eqref{tmutupper2}, and using the same vanishing, gives for every $t\in[t_{k+1},t_k]$
	     \[
	     \nnr{\beta(t,k+1,\fv)}_{L^2_\fv}
	     \leq C_\theta e^{-96.9\mathsf{A}_k^{0.34}}
	     \leq e^{-96.8\mathsf{A}_k^{0.34}}.
	     \]
	     This also controls the first backward critical interval. It remains to prove \eqref{zpm} for $\upmu\geq k+2$.

	     \noindent\textit{Moderate backward modes.}
	     For $k+2\leq\upmu\leq\lfloor\f54k\rfloor$, we prove
		\begin{equation}\label{indu1}
			\nnr{\beta(t_{\upmu-1},\upmu,\fv)}_{L^2_\fv}\leq (1+2\eps_0^6)^{\upmu-k-2}e^{-96.9\mathsf{A}_k^{0.34}}.
		\end{equation}
		For $\upmu=k+2$, this follows from \eqref{betaindu} with an empty interaction sum. Assume \eqref{indu1} holds for $k+2\leq\upnu\leq\upmu$. Applying \eqref{betaindu} at the next mode gives
		\[\ba
		&\nnr{\beta(t_{\upmu},\upmu+1,\fv)}_{L^2_\fv}\leq\sum_{\upnu=k+2}^{\upmu-1}\nnr{\fE_{\upmu+1,\upnu}\beta(t_{\upnu-1},\upnu,\fv)}_{L^2_\fv}+\nnr{\fE_{\upmu+1,\upmu}\beta(t_{\upmu-1},\upmu,\fv)}_{L^2_\fv}+e^{-96.95\mathsf{A}_k^{0.34}}\\
		\leq& \sum_{\upnu=k+2}^{\upmu-1}\f{50(1+\eps_0^6)k^3}{\upnu^3}\eps_*\<\upmu+1-\upnu\>^{-8}\nnr{\beta(t_{\upnu-1},\upnu,\fv)}_{L^2_\fv}\Etnu+\f{(1+\eps_0^6)k^3}{\upmu^3}(1+20\eps_*)\nnr{\beta(t_{\upmu-1},\upmu,\fv)}_{L^2_\fv}\Etmu\\
		\leq&(1+\eps_0^6)(1+10^2\eps_*)(1+2\eps_0^6)^{\upmu-k-2}e^{-96.9\mathsf{A}_k^{0.34}}\leq(1+2\eps_0^6)^{\upmu-k-1}e^{-96.9\mathsf{A}_k^{0.34}}.
		\ea\]
		This proves \eqref{indu1}.

		\noindent\textit{Large backward modes.}
		For $\lfloor\f54k\rfloor\leq\upmu\leq M_k$, we claim that
		\begin{equation}\label{indu2}
			\nnr{\beta(t_{\upmu-1},\upmu,\fv)}_{L^2_\fv}\leq e^{-96.8\mathsf{A}_k^{0.34}}.
		\end{equation}
		The initial case follows from \eqref{indu1}. Assume \eqref{indu2} holds for $\lfloor\f54k\rfloor\leq\upnu\leq\upmu$. Equation \eqref{betaindu} gives
		\[\ba
		&\nnr{\beta(t_{\upmu},\upmu+1,\fv)}_{L^2_\fv}\leq\sum_{\upnu=k+2}^{\upmu-1}\nnr{\fE_{\upmu+1,\upnu}\beta(t_{\upnu-1},\upnu,\fv)}_{L^2_\fv}+\nnr{\fE_{\upmu+1,\upmu}\beta(t_{\upmu-1},\upmu,\fv)}_{L^2_\fv}+e^{-96.95\mathsf{A}_k^{0.34}}\\
		\leq&\sum_{\upnu=k+2}^{\upmu-1}\f{50(1+\eps_0^6)k^3}{\upnu^3}\eps_*\<\upmu+1-\upnu\>^{-8}\Etnu e^{-96.8\mathsf{A}_k^{0.34}}+\f{(1+\eps_0^6)k^3}{\upmu^3}(1+20\eps_*)\Etmu e^{-96.8\mathsf{A}_k^{0.34}}+e^{-96.95\mathsf{A}_k^{0.34}}\\
		\leq& (10^2\eps_*+0.8)e^{-96.8\mathsf{A}_k^{0.34}}\leq e^{-96.8\mathsf{A}_k^{0.34}}.
		\ea\]
		Here we use $\frac{(1+\eps_0^6)k^3}{\upmu^3}(1+30\eps_*)\Etmu\leq0.8$ for $\upmu\geq\lfloor5k/4\rfloor$. This proves \eqref{indu2}, hence \eqref{zpm}.

	\end{proof}

	We now convert the backward recurrence bounds into estimates for the full profile on $[\mathfrak{T}_k,\mathsf{T}_k]$.
	
	\begin{prop}\label{prop:backward-stability} For every $\upsigma\in\Z$, the following estimates hold:
	\begin{itemize} 
		\item[(i).] If $t\in[\mathsf{T}_{1,k},\mathsf{T}_k]$, then
	\begin{equation}\label{ttwot}
		\nnr{e^{\<\na_{\mrv}\>^\sss}v^{\alpha}\rr{\sF_zf(t,\upsigma\bfe,v)-\sF_zf(\mathsf{T}_k,\upsigma\bfe,v)}}_{\lta L^2_v}\leq e^{-95\mathsf{A}_k^{0.34}}.
	\end{equation}
	\item[(ii).] If $t\in[\mathfrak{T}_k,\mathsf{T}_{1,k}]$, then \begin{equation}\label{t1t2}
		\nnr{e^{\<\na_{\mrv}\>^\sss}v^{\alpha}\rr{\sF_zf(t,\upsigma\bfe,v)-\sF_zf(\mathsf{T}_{1,k},\upsigma\bfe,v)}}_{\lta L^2_v}\leq e^{-\f12\<\f{\mathfrak{T}_k}{A_*}\>^\sss}\eps_k k^{-\f16}.
	\end{equation}
	\end{itemize}
	\end{prop} 
	\begin{proof}
	\noindent\textit{The interval $[\mathsf{T}_{1,k},\mathsf{T}_k]$.}
	For $t$ in this interval, choose $k+1\leq\upmu\leq M_k$ such that $t_\upmu\leq t\leq t_{\upmu-1}$. By \eqref{aprefor2} and Proposition \ref{propertybeta1},
		\[\ba
		&\nnr{e^{\<\na_{\mrv}\>^\sss}v^{\alpha}\rr{\sF_zf(t,\upsigma\bfe,v)-\sF_zf(\mathsf{T}_k,\upsigma\bfe,v)}}_{\lta L^2_v}\leq\nnr{e^{\<\na_{\mrv}\>^\sss}v^{\alpha}\rr{\sF_zf(t,\upsigma\bfe,v)-\sF_zf(t_{\upmu-1},\upsigma\bfe,v)}}_{\lta L^2_v}\\
		+&\sum_{\upnu=k+1}^{\upmu-1}\nnr{e^{\<\na_{\mrv}\>^\sss}v^{\alpha}\cE_{\upsigma,\upnu}\beta(t_{\upnu-1},\upnu,\fv)}_{\lta L^2_v}+e^{-96.9\mathsf{A}_k^{0.34}}
		\leq\sum_{\upnu=k+1}^{\upmu-1}\f{\mathsf{A}_k}{\upnu^3}\<\upsigma-\upnu\>^{-8}\nnr{\beta(t_{\upnu-1},\upnu,\fv)}_{L^2_\fv}\\
		+&\f{\mathsf{A}_k}{\upmu^3}\nnr{\beta(t_{\upmu-1},\upmu,\fv)}_{L^2_\fv}+e^{-96.8\mathsf{A}_k^{0.34}}\leq e^{-95\mathsf{A}_k^{0.34}}.
		\ea\]

		\noindent\textit{The interval $[\mathfrak{T}_k,\mathsf{T}_{1,k}]$.}
		We first claim that
		\begin{equation}\label{appttwot}
			\begin{aligned}
				&\wh{\rho}(t,\upnu\bfe)=\wh{\sR_1}(t,\upnu\bfe)\quad\text{for any}\ \upnu\in\Z, \\
				&\wh{f}(t,\upsigma\bfe,\xi)-\wh{f}(\mathsf{T}_{1,k},\upsigma\bfe,\xi)=\sum_{\upnu\neq0}\int_t^{\mathsf{T}_{1,k}}\wh{\udl{\rho}}(\tau,\upnu\bfe)\f{\upxi-\upsigma \tau}{\upnu}\wh{f}(\tau,(\upsigma-\upnu)\bfe,\xi-\upnu\tau\bfe) d\tau+\wh{\sR_2}(t,\upsigma\bfe,\xi),
			\end{aligned}
		\end{equation}
		with the error terms satisfying
		\[\nr{\wh{\sR_1}(t,\upnu\bfe)}\leq e^{-99\mathsf{A}_k^{0.34}}e^{-0.124\<\f{\upnu t-\mathsf{A}_k}{k^\f13}\>^\sss},\quad\nnr{e^{0.124\<\f{\na_\fv-\mathsf{A}_k}{k^\f13}\>^\sss+\<\na_{\mrv}\>^\sss}\rr{v^{\alpha}\sR_2(t,\upsigma\bfe,v)}}_{\lta L^2_v}\leq e^{-99\mathsf{A}_k^{0.34}}.\]
		Indeed, \eqref{zsd} gives $\nr{\wh{\rho}(t,\upnu\bfe)}\leq e^{5\mathsf{A}_k^{0.34}}e^{-0.125\<\f{\upnu-k,\upnu t-\mathsf{A}_k}{k^\f13}\>^\sss}$. If $\upnu\geq M_k$, then $|\upnu-k|/k^{1/3}\gtrsim\mathsf{A}_k^{0.35}$; if $\upnu\leq M_k$, then $|\upnu t-\mathsf{A}_k|\geq\mathsf{A}_k^{0.35}k^{1/3}$ on the interval under consideration. Thus every mode is separated from at least one localization center, which proves the estimate for $\sR_1$. Substituting this density bound into the Duhamel formula and repeating the estimates leading to \eqref{fref} gives the stated bound for $\sR_2$.
		  
		By Proposition \ref{frelocalfnq}, \eqref{udlrho}, and \eqref{appttwot},
		\[\ba
		&\nnr{e^{\<\mrxi\>^\sss}\rr{\pa_\xi^\alpha\wh{f}(t,\upsigma\bfe,\xi)-\pa_\xi^\alpha\wh{f}(\mathsf{T}_{1,k},\upsigma\bfe,\xi)}}_{\lta L^2_\xi}\ls \bigg\|\sum_{|\upsigma-\upnu|\leq M_k}\int_t^{\mathsf{T}_{1,k}}|\wh{\udl{\rho}}(\tau,\upnu\bfe)|\tau\bigg|\sF\rr{e^{\<\na_{\mrv}\>^\sss}\<\na_{\fz,\fv}\>v^\alpha f}\\
		&(\tau,(\upsigma-\upnu)\bfe,\xi-\upnu\tau\bfe)\bigg|\vv{1}_{|\upxi-\upnu\tau-\mathsf{A}_k|\leq M_k} d\tau\bigg\|_{\lta L^2_\xi}+e^{-98\mathsf{A}_k^{0.34}}\leq3\mathsf{A}_k\sum_{\upnu}\int_t^{\mathsf{T}_{1,k}}\nr{\wh{\udl{\rho}}(\tau,\upnu\bfe)}\tau d\tau\\
	    & \times\sup_{t\leq\tau\leq \mathsf{T}_{1,k}}\nnr{e^{\<\mrxi\>^\sss}\pa_\xi^\alpha\wh{f}(\tau,(\upsigma-\upnu)\bfe,\xi)}_{\lta L^2_\xi}+e^{-98\mathsf{A}_k^{0.34}}\leq e^{-\f34\<\f{\mathfrak{T}_k}{A_*}\>^\sss}\sup_{\substack{t\leq\tau\leq \mathsf{T}_{1,k}\\\upsigma\in\Z}}\nnr{e^{\<\mrxi\>^\sss}\pa_\xi^\alpha\wh{f}(\tau,\upsigma\bfe,\xi)}_{\lta L^2_\xi}+e^{-98\mathsf{A}_k^{0.34}}.
		\ea\]
		Moreover, \eqref{ttwot} gives
		 $\nnr{e^{\<\mrxi\>^\sss}\pa_\xi^\alpha\wh{f}(\mathsf{T}_{1,k},\upsigma\bfe,\xi)}_{\lta L^2_\xi}\leq50\eps_k k^{-\f16}+e^{-95\mathsf{A}_k^{0.34}}$. Therefore
		\[\ba &\nnr{e^{\<\mrxi\>^\sss}\pa_\xi^\alpha\rr{\wh{f}(t,\upsigma\bfe,\xi)-\wh{f}(\mathsf{T}_{1,k},\upsigma\bfe,\xi)}}_{\lta L^2_\xi}\\& \leq e^{-\f34\<\f{\mathfrak{T}_k}{A_*}\>^\sss}\sup_{\substack{\upsigma\in\Z\\t\leq\tau\leq \mathsf{T}_{1,k}}}\nnr{e^{\<\mrxi\>^\sss}\pa_\xi^\alpha\rr{\wh{f}(\tau,\upsigma\bfe,\xi)-\wh{f}(\mathsf{T}_{1,k},\upsigma\bfe,\xi)}}_{\lta L^2_\xi}+e^{-0.7\<\f{\mathfrak{T}_k}{A_*}\>^\sss}\eps_k k^{-\f16}.\ea\]
		A standard continuity argument now yields \eqref{t1t2} and completes the proof.
	\end{proof}

	\subsection{Output of the homogeneous analysis}
	The conclusions of this section used later in the nonlinear construction are collected below.
	\begin{prop}\label{prop:homogeneous-output}
	Let $f$ solve the homogeneous linear problem stated at the beginning of this section. Then:
	\begin{itemize}
		\item[(i)] the forward echo chain satisfies the two-sided amplification estimate \eqref{mumu1}, while the full profile satisfies \eqref{uppsigma}, \eqref{tztinfty}, and the lower bound \eqref{ecq9};
		\item[(ii)] the backward evolution creates no comparable cascade: its recurrence coefficients satisfy \eqref{zpm}, and the solution satisfies the stability estimates \eqref{ttwot} and \eqref{t1t2} on $[\mathfrak{T}_k,\mathsf{T}_k]$.
	\end{itemize}
	\end{prop}
	\begin{proof}
	Part (i) follows from Propositions \ref{propertybetatmumu} and \ref{prop:forward-profile}. Part (ii) follows from Propositions \ref{propertybeta1} and \ref{prop:backward-stability}.
	\end{proof}

	\section{Inhomogeneous linear response}

	We now study the response of the linearized equation to a source localized near one frequency. More precisely, we consider
\begin{equation}\label{inhomequ}
	\left\{\begin{aligned}
		&\pa_tf=\cL f+Q=\cN(f,f^0+\udl{f_\infty}+ f_{*,err})+\cN(\udl{f_\infty}+ f_{*,err},f)+Q, \\
		&\udl{f_\infty}(z,v):=\eps_0\vf(v)\cos \fz+f_{*,\infty}(z,v), \\
		&f|_{t=\mathsf{T}_k}=0,\quad z=(\fz,\mathring{z})\in \T^3, v=(\fv,\mathring{v})\in \R^3.
	\end{aligned}\right.
\end{equation}
	on the time interval $[\mathfrak{T}_k,\infty)$. The solution vanishes at $\mathsf{T}_k$. Thus, \(\int_{\mathsf T_k}^t Q(\tau)\,d\tau\) is the direct contribution of the source.
	We introduce an envelope \(Z(t)\) to control the contribution of the source. Our purpose is to show that if \(Z(t)\) satisfies some proper growth conditions, see \eqref{growth1} and \eqref{growth2}, then the solution \(f\) can be bounded by \(Z(t)\) multiplied by suitable powers of \(|k',A'|\).
	 This estimate will later allow us to control every inhomogeneous term in the Taylor expansion.

	The assumptions on the source are the following.
 \begin{itemize}
	\item[(i).] \emph{One-dimensional frequency localization.} The source satisfies $Q(t,z,v)=Q(t,\fz\bfe,v)$, $\wh{Q}(t,0,0)=0$ and
	\begin{equation}\label{AssumptionQ}
		\mathrm{Supp}\ \wh{Q}(t)\subset \{k'\bfe\}\times[A'-1,A'+1]\times\R^2,\quad \int_{\min\{t,\mathsf{T}_k\}}^{\max\{t,\mathsf{T}_k\}}\nnr{e^{\<\na_{\mrv}\>^{\sss}}(v^\alpha Q(\tau))}_{\lta L^2_v}d\tau\leq Z(t).
	\end{equation}
	\item[(ii).] \emph{Monotone source envelope.} The function $Z(t)$ is decreasing on $[\mathfrak{T}_k,\mathsf{T}_k]$ and increasing on $[\mathsf{T}_k,\infty)$. Thus $Z(t)$ controls the source accumulated from $t$ to $\mathsf{T}_k$ in the backward direction and from $\mathsf{T}_k$ to $t$ in the forward direction.
	\item[(iii).] \emph{Admissible source frequencies.} The pair $(k',A')$ satisfies
	\begin{equation}\label{AssumptionkA}
		\max\bigg\{\exp(A_*),\exp(\f1{\eps_*})\bigg\}\leq|k',A'|\leq \mathsf{A}_k^{\mathtt{N}_1+2}.
	\end{equation}
	By symmetry in the velocity frequency, we may assume $A'\geq0$.
 \end{itemize}  

	The geometry behind the proof is the same as in Section 3. A packet near $(k',A')$ can produce a density response only when $\upmu t\approx A'$ for some spatial mode $\upmu$. If $|k'|\geq A'$, the resonance associated with the source mode lies before the time interval under consideration, and the response is perturbative. If $|k'|\leq A'$, the critical times
	\[
		t'_\upmu=\frac{A'}{\upmu+\frac12},\qquad \upmu\in\mathbb{Z}_{\geq0},
	\]
	organize a possible echo chain. We use the threshold $k'_2=\lfloor A'^{1/3}\rfloor$: for $\upmu\geq k'_2$, the one-step transfer factor $\eps_0A'/\upmu^3$ is perturbative, whereas the modes below $k'_2$ must be followed by a discrete recurrence. Set $k'_0=0$ if $A'\leq\mathsf{T}_k/2$; when $A'>\mathsf{T}_k/2$, let $k'_0$ be determined by
	\[
		t'_{k'_0}\leq\mathsf{T}_k<t'_{k'_0-1}.
	\]
	Similarly, set $k'_1=0$ if $A'\leq\mathfrak{T}_k/2$; when $A'>\mathfrak{T}_k/2$, let $k'_1$ satisfy
	\[
		t'_{k'_1}\leq\mathfrak{T}_k<t'_{k'_1-1}.
	\]
	The indices \(k'_0\) and \(k'_1\) identify the critical intervals containing \(\mathsf T_k\) and \(\mathfrak T_k\), respectively, while \(k'_2\) separates the perturbative high-mode regime from the potentially amplifying low-mode regime.

	We first dispose of the nonresonant frequency configuration.
	\begin{prop}\label{prop:far-frequency-response}  Let $(k',A')$ verify \eqref{AssumptionkA} and $|k'|\geq A'$. Then
		the solution $f$ to \eqref{inhomequ} satisfies
		\begin{equation}\label{kgeqA}
			\nnr{e^{\<\na_{\mrv}\>^{\sss}}v^\alpha\rr{\sF_zf(t,\upsigma\bfe,v)-\int_{\mathsf{T}_k}^t\sF_zQ(\tau,\upsigma\bfe,v)d\tau}}_{\lta L^2_v}\leq\rr{e^{-\f12\<\f{\mathfrak{T}_k}{A_*}\>^\sss}+e^{-95|k'|^{0.34}}}Z(t)
		\end{equation}
		for any $\upsigma\in\Z$	on $[\mathfrak{T}_k,\infty)$.
	\end{prop}
	\begin{proof}
		By Proposition \ref{longenergyl}, it holds that
		\[\nnr{\<\na_{\fz,\fv}\>^3e^{\f18\<\f{\na_\fz-k',\na_\fv-A'}{A_*}\>^\sss+\<\na_{\mrv}\>^{\sss}}\rr{v^\alpha f(t)}}_{\lta L^2_{z,v}}\leq e^{3|k'|^{0.34}}Z(t),\quad \nr{\wh{\rho}(t,\upnu\bfe)}\leq e^{3|k'|^{0.34}-\f18\<\f{\upnu-k',\upnu t-A'}{A_*}\>^\sss}Z(t).\]
		The weighted estimate localizes the solution near $(k',A')$. Since $|k'|\geq A'$, the relation $\upnu t\approx A'$ cannot hold on $[\mathfrak{T}_k,\infty)$ for modes $\upnu$ in this localization region. Consequently, the density and the interaction part of the profile can be placed into remainders:
		\[\ba
		&\wh{\rho}(t,\upnu\bfe)=\wh{\sR}_1(t,\upnu\bfe)\quad\text{for any}\ \upnu\in\Z,\\
		&\wh{f}(t,\upsigma\bfe,\xi)-\int_{\mathsf{T}_k}^t\wh{Q}(\tau,\upsigma\bfe,\xi)d\tau=-\sum_{\upnu\neq0}\int_{\mathsf{T}_k}^{t}\wh{\udl{\rho}}(\tau,\upnu\bfe)\f{\upxi-\upsigma \tau}{\upnu}\wh{f}(\tau,(\upsigma-\upnu)\bfe,\xi-\upnu\tau\bfe) d\tau+\wh{\sR}_2(t,\upsigma\bfe,\xi),
		\ea\]
		with the error terms satisfying
		\[\nr{\wh{\sR}_1(t,\upnu\bfe)}\leq e^{-100|k'|^{0.34}}e^{-0.124\<\f{\upnu-k',\upnu t-A'}{A_*}\>^\sss}Z(t),\quad\nnr{e^{0.124\<\f{\na_\fv-A'}{A_*}\>^\sss+\<\na_{\mrv}\>^{\sss}}\rr{v^\alpha\sR_2(t,\upsigma\bfe,v)}}_{\ltaa L^2_v}\leq e^{-100|k'|^{0.34}}Z(t).\]
		To verify these bounds, split into $|\upnu-k'|\geq A_*|k'|^{0.35}$ and $|\upnu-k'|\leq A_*|k'|^{0.35}$. In the first region, \(|\upnu-k'|\) is sufficiently large to yield exponential decay. In the second region, the same decay follows from the lower bound on \(|\upnu t-A'|\).
		This is the decomposition used in \eqref{appttwot}.

		We explain the closure on the backward interval; the forward interval is identical after reversing the orientation of the source integral. For $t\in[\mathfrak{T}_k,\mathsf{T}_k]$, the preceding approximation yields
		\[\ba
		&\nnr{e^{\<\mrxi\>^\sss}\pa_\xi^\alpha\rr{\wh{f}(t,\upsigma\bfe,\xi)+\int_t^{\mathsf{T}_k}\wh{Q}(\tau,\upsigma\bfe,\xi)d\tau}}_{\ltaa L^2_\xi}
		\leq \sum_{\upnu}\int_t^{\mathsf{T}_k}\nr{\wh{\udl{\rho}}(\tau,\upnu\bfe)}\<\tau\> d\tau\sup_{\substack{t\leq\tau\leq \mathsf{T}_k\\\upsigma\in\Z}}\nnr{e^{\<\mrxi\>^\sss}\<\upsigma,\upxi\>\pa_\xi^\alpha\wh{f}(\tau,\upsigma\bfe,\xi)\right.\\&\left.\vv{1}_{|\upsigma-k',\upxi-A'|\leq|k'|^{0.35}} }_{\ltaa L^2_\xi}+e^{-99|k'|^{0.34}}Z(t)
		\leq e^{-0.9\<\f{\mathfrak{T}_k}{A_*}\>^\sss}3|k'|\sup_{\substack{t\leq\tau\leq \mathsf{T}_k\\\upsigma\in\Z}}\nnr{e^{\<\mrxi\>^\sss}\pa_\xi^\alpha\wh{f}(\tau,\upsigma\bfe,\xi)}_{\ltaa L^2_\xi}+e^{-99|k'|^{0.34}}Z(t)\\
		&\leq e^{-0.8\<\f{\mathfrak{T}_k}{A_*}\>^\sss}\sup_{\substack{t\leq\tau\leq \mathsf{T}_k\\\upsigma\in\Z}}\nnr{e^{\<\mrxi\>^\sss}\pa_\xi^\alpha\rr{\wh{f}(\tau,\upsigma\bfe,\xi)+\int_\tau^{\mathsf{T}_k}\wh{Q}(s,\upsigma\bfe,\xi)ds}}_{\ltaa L^2_\xi}+\rr{e^{-99|k'|^{0.34}}+e^{-0.8\<\f{\mathfrak{T}_k}{A_*}\>^\sss}}Z(t).
		\ea\]
		Here the small coefficient follows from the definition of $\mathfrak{T}_k$ and $\log |k'|\leq3(\mathtt{N}_1+2)\log k$. It can therefore be absorbed by a continuity argument, which proves \eqref{kgeqA} on $[\mathfrak{T}_k,\mathsf{T}_k]$. On $[\mathsf{T}_k,\infty)$, the same estimate is applied to $f-\int_{\mathsf{T}_k}^tQ(\tau)d\tau$ and gives the stated bound.
	\end{proof}
	
	For the remainder of this section, we assume $|k'|\leq A'$. We first estimate the solution forward in time from \(\mathsf{T}_k\) and then backward in time to \(\mathfrak{T}_k\).
	
	\subsection{Forward response on $[\mathsf{T}_k,\infty)$}
	If $A'\leq\mathsf{T}_k/2$, then $t'_0=2A'\leq\mathsf{T}_k$ and there is no future critical time; this case is covered directly by Proposition \ref{4.6}. We therefore begin with $A'>\mathsf{T}_k/2$. There are two possible entry points into the echo chain. If $k'_0\geq k'_2$, the solution first crosses only weak high-mode resonances, which can be controlled simultaneously by an $L^2_t$ estimate. If $k'_0<k'_2$, the time $\mathsf T_k$ lies in the low-mode regime, so the first partial critical interval requires separate treatment.

	We start with the weak high-mode regime $k'_0\geq k'_2$.
	\begin{prop}\label{4.2}
		Assume $A'>\mathsf{T}_k/2$ and $k'_0\geq k'_2$. On $[\mathsf{T}_k,t'_{k'_2-1}]$,
		\begin{equation}\label{ttk2}
			\nnr{e^{\<\na_{\mrv}\>^\sss}v^\alpha\rr{\sF_zf(t,\upsigma\bfe,v)-\int_{\mathsf{T}_k}^t\sF_zQ(\tau,\upsigma\bfe,v)d\tau}}_{\lta L^2_v}\leq \eps_0^{0.9}Z(t)
		\end{equation}
		for any $\upsigma\in\Z$.
	\end{prop}
	\begin{proof}
		The proof has four steps. We first separate the active high modes from exponentially small off-resonant modes. We then invert the diagonal Volterra equation for the density, recover the profile from the density estimate, and finally close the two estimates simultaneously.
		
		\noindent\textit{Step 1. Frequency-space reduction.}
		By Proposition \ref{longenergyl} it holds that
		\[\nnr{\<\na_{\fz,\fv}\>^3e^{\f18\<\f{\na_\fz-k',\na_\fv-A'}{A_*}\>^\sss+\<\na_{\mrv}\>^\sss}\rr{v^{\alpha}f(t)}}_{\lta L^2_{z,v}}\leq e^{3A'^{0.34}}Z(t),\quad \nr{\wh{\rho}(t,\upnu\bfe)}\leq e^{3A'^{0.34}-\f18\<\f{\upnu-k',\upnu t-A'}{A_*}\>^\sss}Z(t).\]
		On $[\mathsf{T}_k,t'_{k'_2-1}]$, only modes $\upnu\geq k'_2$ can satisfy $\upnu t\approx A'$. We therefore retain their density equations and place the modes $\upnu<k'_2$, together with interactions involving $f_{*,err}$ and the unknown density $\udl{\rho}$, into remainders. This gives
		\begin{equation*}
			\wh{\rho}(t,\upnu\bfe)=\wh{\sR_2}(t,\upnu\bfe)\quad\text{for any}\ \upnu<k'_2,
		\end{equation*}
		\begin{equation}\label{rho42}
			\begin{aligned}
				\wh{\rho}(t,\upmu\bfe)&=\int_{\mathsf{T}_k}^t\wh{Q}(\tau,\upmu\bfe,\upmu t\bfe)d\tau-\int_{\mathsf{T}_k}^{t}\wh{\rho}(\tau,\upmu\bfe)(t-\tau)\sF\rr{f^0+\udl{f_\infty}}(0,\upmu(t-\tau)\bfe) d\tau \\
				&\quad -\sum_{\substack{\upnu\geq k'_2\\\upnu\neq \upmu}}\int_{\mathsf{T}_k}^{t}\wh{\rho}(\tau,\upnu\bfe)\f{\upmu(t-\tau)}{\upnu}\wh{\udl{f_\infty}}((\upmu-\upnu)\bfe,(\upmu t-\upnu\tau)\bfe) d\tau
	    +\wh{\sR_1}(t,\upmu\bfe)+\wh{\sR_2}(t,\upmu\bfe)\quad\text{for any}\ \upmu\geq k'_2,
			\end{aligned}
		\end{equation}
		\begin{equation}\label{f41}
			\begin{aligned}
				\wh{f}(t,\upsigma\bfe,\xi)-\int_{\mathsf{T}_k}^t\wh{Q}(\tau,\upsigma\bfe,\xi)d\tau&=-\sum_{\upnu\geq k'_2}\int_{\mathsf{T}_k}^{t}\wh{\rho}(\tau,\upnu\bfe)\f{\upxi-\upsigma \tau}{\upnu}\sF\rr{f^0+\udl{f_\infty}}((\upsigma-\upnu)\bfe,\xi-\upnu\tau\bfe) d\tau \\
				&\quad +\wh{\sR_1}(t,\upsigma\bfe,\xi)+\wh{\sR_2}(t,\upsigma\bfe,\xi),
			\end{aligned}
		\end{equation}
		with the error terms satisfying
		\[\ba
		&\nr{\wh{\sR_2}(t,\upnu\bfe)}\leq e^{-100A'^{0.34}-0.124\<\f{\upnu-k',\upnu t-A'}{A_*}\>^\sss}Z(t),\quad\nnr{e^{0.124\<\f{\na_\fv-A'}{A_*}\>^\sss+\<\na_{\mrv}\>^{\sss}}\rr{v^{\alpha}\sR_2(t,\upsigma\bfe,v)}}_{\ltaa L^2_v}\leq e^{-100A'^{0.34}}Z(t).\\
		&\nnr{\wh{\sR_1}(t,\upmu\bfe)}_{L^2([\mathsf{T}_k,\tilde{T}])}\upmu^\f12\leq e^{-0.9\<\f{\mathsf{T}_k}{A_*}\>^\sss}\rr{\sup_{\upmu\geq k'_2}\nnr{\wh{\rho}(t,\upmu\bfe)}_{L^2([\mathsf{T}_k,\tilde{T}])}+\sup_{\substack{\mathsf{T}_k\leq\tau\leq\tilde{T}\\\upsigma\in\Z}}\nnr{\wh{f}(\tau,\upsigma\bfe,\upxi\bfe)}_{L^2_\upxi}}+e^{-99{A'}^{0.34}}Z(\tilde{T}),\\
		&\nnr{e^{\<\na_{\mrv}\>^{\sss}}v^{\alpha}\sR_1(t,\upsigma\bfe,v)}_{\ltaa L^2_v}\leq e^{-0.9\<\f{\mathsf{T}_k}{A_*}\>^\sss}\rr{\sup_{\upmu\geq k'_2}\nnr{\wh{\rho}(\tau,\upmu\bfe)}_{L^2([\mathsf{T}_k,t])}+\sup_{\substack{\mathsf{T}_k\leq\tau\leq t\\\upsigma\in\Z}}\nnr{e^{\<\mrxi\>^{\sss}}\pa_\xi^\alpha\wh{f}(\tau,\upsigma\bfe,\xi)}_{\ltaa L^2_\xi}}+e^{-99{A'}^{0.34}}Z(t),
		\ea\]
		for any $\mathsf{T}_k\leq\tilde{T}\leq t'_{k'_2-1}$, $\mathsf{T}_k\leq t\leq t'_{k'_2-1}$. In fact we have
		\[\ba
		&\wh{\sR_1}(t,\upmu\bfe)=-\sum_{\upnu\geq k'_2}\int_{\mathsf{T}_k}^{t}\wh{\rho}(\tau,\upnu\bfe)\f{\upmu(t-\tau)}{\upnu}\wh{f_{*,err}}(\tau,(\upmu-\upnu)\bfe,(\upmu t-\upnu\tau)\bfe) d\tau\\	&-\sum_{\upnu\neq0}\int_{\mathsf{T}_k}^t\wh{\udl{\rho}}(\tau,\upnu\bfe)\f{\upmu(t-\tau)}{\upnu}\wh{f}(\tau,(\upmu-\upnu)\bfe,(\upmu t-\upnu\tau)\bfe)d\tau,\\
		&\wh{\sR_2}(t,\upmu\bfe)=-\sum_{\substack{\upnu\leq k'_2-1\\\upnu\neq0}}\int_{\mathsf{T}_k}^t\wh{\rho}(\tau,\upnu\bfe)\f{\upmu(t-\tau)}{\upnu}\wh{\udl{f}}(\tau,(\upmu-\upnu)\bfe,(\upmu t-\upnu\tau)\bfe)d\tau,\quad\text{for any}\ \upmu\geq k'_2.
		\ea\]
		\[\ba
		&\wh{\sR_1}(t,\upsigma\bfe,\xi)=-\sum_{\upnu\geq k'_2}\int_{\mathsf{T}_k}^{t}\wh{\rho}(\tau,\upnu\bfe)\f{\upxi-\upsigma\tau}{\upnu}\wh{f_{*,err}}(\tau,(\upsigma-\upnu)\bfe,\xi-\upnu\tau\bfe) d\tau-\sum_{\upnu\neq0}\int_{\mathsf{T}_k}^t\wh{\udl{\rho}}(\tau,\upnu\bfe)\f{\upxi-\upsigma\tau}{\upnu}\\&\times\wh{f}(\tau,(\upsigma-\upnu)\bfe,\xi-\upnu\tau\bfe)d\tau,\quad
		\wh{\sR_2}(t,\upsigma\bfe,\xi)=-\sum_{\substack{\upnu\leq k'_2-1\\\upnu\neq0}}\int_{\mathsf{T}_k}^t\wh{\rho}(\tau,\upnu\bfe)\f{\upxi-\upsigma\tau}{\upnu}\wh{f^0+\udl{f}}(\tau,(\upsigma-\upnu)\bfe,\xi-\upnu\tau\bfe)d\tau.
		\ea\]
		The $\sR_2$ bounds follow from the same frequency separation used in \eqref{rerho}, \eqref{rhomu}, and \eqref{fref}: a mode below $k'_2$ stays away from the active resonant set on this interval. The terms in $\sR_1$ are small for a different reason. They contain either the exponentially small coefficient $f_{*,err}$ or the background density $\udl{\rho}$, whose decay supplies the factor $e^{-c\<\mathsf{T}_k/A_*\>^\sss}$. More precisely,
		\[\ba
		&\nnr{\wh{\sR_1}(t,\upmu\bfe)}_{L^2([\mathsf{T}_k,\tilde{T}])}\leq\sum_{\upnu\geq k'_2}\int_{\mathsf{T}_k}^{\tilde{T}}\nr{\wh{\rho}(\tau,\upnu\bfe)}\<\tau\>\nnr{\sF\rr{\<\na_{\fz,\fv}\>f_{*,err}}(\tau,(\upmu-\upnu)\bfe,(\upmu t-\upnu\tau)\bfe)}_{L^2_t} d\tau\\
		&+\sum_{\upnu}\int_{\mathsf{T}_k}^{\tilde{T}}\nr{\wh{\udl{\rho}}(\tau,\upnu\bfe)}\<\tau\>\nnr{\sF\rr{\<\na_{\fz,\fv}\>f}(\tau,(\upmu-\upnu)\bfe,(\upmu t-\upnu\tau)\bfe)}_{L^2_t} d\tau\\
		&\leq\sum_{\upnu\geq k'_2}\f{A'^{\f32}}{\upmu^\f12}\nnr{\wh{\rho}(t,\upnu\bfe)}_{L^2([\mathsf{T}_k,\tilde{T}])}\sup_{\mathsf{T}_k\leq\tau\leq \tilde{T}}\nnr{\sF\rr{\<\na_{\fz,\fv}\>f_{*,err}}(\tau,(\upmu-\upnu)\bfe,\upxi\bfe)}_{L^2_\upxi}\\
		&+\f{1}{\upmu^\f12}e^{-0.95\<\f{\mathsf{T}_k}{A_*}\>^\sss}\rr{2A'\sup_{\substack{\mathsf{T}_k\leq\tau\leq \tilde{T}\\\upsigma\in\Z}}\nnr{\wh{f}(\tau,\upsigma\bfe,\upxi\bfe)}_{L^2_\upxi}+e^{-99.9A'^{0.34}}Z(\tilde{T})}\\
		&\leq\f{1}{\upmu^\f12} e^{-0.9\<\f{\mathsf{T}_k}{A_*}\>^\sss}\rr{\sup_{\upmu\geq k'_2}\nnr{\wh{\rho}(t,\upmu\bfe)}_{L^2([\mathsf{T}_k,\tilde{T}])}+\sup_{\substack{\mathsf{T}_k\leq\tau\leq\tilde{T}\\\upsigma\in\Z}}\nnr{\wh{f}(\tau,\upsigma\bfe,\upxi\bfe)}_{L^2_\upxi}}+\f{1}{\upmu^\f12}e^{-99{A'}^{0.34}}Z(\tilde{T}).
		\ea\]
		Similarly we have
		\[\ba
		&\nnr{e^{\<\mrxi\>^{\sss}}\pa_\xi^{\alpha}\wh{\sR_1}(t,\upsigma\bfe,\xi)}_{\ltaa L^2_\xi}\ls\sum_{\upnu\geq k'_2}\int_{\mathsf{T}_k}^{t}\nr{\wh{\rho}(\tau,\upnu\bfe)}\<\tau\>\nnr{\sF\rr{e^{\<\na_{\mrv}\>^{\sss}}\<\na_{\fz,\fv}\>v^\alpha f_{*,err}}(\tau,(\upsigma-\upnu)\bfe,\xi)}_{\ltaa L^2_\xi} d\tau\\
		&+\sum_{\upnu}\int_{\mathsf{T}_k}^{t}\nr{\wh{\udl{\rho}}(\tau,\upnu\bfe)}\<\tau\>\nnr{\sF\rr{e^{\<\na_{\mrv}\>^{\sss}}\<\na_{\fz,\fv}\>v^\alpha f}(\tau,(\upsigma-\upnu)\bfe,\xi)}_{\ltaa L^2_\xi} d\tau\\
		&\leq e^{-0.95\<\f{\mathsf{T}_k}{A_*}\>^\sss}\rr{\sup_{\upmu\geq k'_2}\nnr{\wh{\rho}(\tau,\upmu\bfe)}_{L^2([\mathsf{T}_k,t])}+\sup_{\substack{\mathsf{T}_k\leq\tau\leq t\\\upsigma\in\Z}}\nnr{e^{\<\mrxi\>^{\sss}}\pa_\xi^\alpha\wh{f}(\tau,\upsigma\bfe,\xi)}_{\ltaa L^2_\xi}}+e^{-99.4{A'}^{0.34}}Z(t).
		\ea\]
		
		\noindent\textit{Step 2. Volterra estimate for the active density.}
		Set $\sR=\sR_1+\sR_2$. For $\upmu\geq k'_2$, equation \eqref{rho42} is a diagonal Volterra equation perturbed by transfers between distinct high modes. Inverting its diagonal part gives
		\[\ba
		&\wh{\rho}(t,\upmu\bfe)=H(t,\upmu\bfe)-\int_{\mathsf{T}_k}^tH(\tau,\upmu\bfe)R(t-\tau,\upmu\bfe)d\tau,\\
		&H(t,\upmu\bfe)=\int_{\mathsf{T}_k}^t\wh{Q}(\tau,\upmu\bfe,\upmu t\bfe)d\tau
		-\sum_{\substack{\upnu\geq k'_2\\\upnu\neq \upmu}}\int_{\mathsf{T}_k}^{t}\wh{\rho}(\tau,\upnu\bfe)\f{\upmu(t-\tau)}{\upnu}\wh{\udl{f_\infty}}(\tau,(\upmu-\upnu)\bfe,(\upmu t-\upnu\tau)\bfe) d\tau+\wh{\sR}(t,\upmu\bfe).
		\ea\]
	    By the Schur test and the resolvent bound \eqref{eR},
		\[\nnr{\wh{\rho}(t,\upmu\bfe)}_{L^2_t[\mathsf{T}_k,\tilde{T}]}\ls_\theta \nnr{H(t,\upmu\bfe)}_{L^2_t[\mathsf{T}_k,\tilde{T}]}.\]
		Thus it remains to estimate $H$. Denote its source term and off-diagonal transfer term by $\mathscr{I}_1$ and $\mathscr{I}_2$, respectively. The source localization gives
		\[\nnr{\mathscr{I}_1}_{L^2_t[\mathsf{T}_k,\tilde{T}]}\leq\int_{\mathsf{T}_k}^{\tilde{T}}\nnr{\wh{Q}(\tau,\upmu\bfe,\upmu t\bfe)}_{L^2_t}d\tau\leq \f{(2\pi)^\f12}{\upmu^\f12}Z(\tilde{T}).\]
		For $\mathscr{I}_2$, the localization $|\upnu\tau-A'|\leq A'^{0.35}$ restricts the interaction to a short resonant strip. The Schur test then converts the convolution into the small factor $\eps_0A'/{k'_2}^3$:
		\[\ba
		&\nnr{\mathscr{I}_2}_{L^2_t[\mathsf{T}_k,\tilde{T}]}\leq\nnr{\sum_{\substack{\upnu\geq k'_2\\\upnu\neq \upmu}}\int_{\mathsf{T}_k}^{t}\nr{\wh{\rho}(\tau,\upnu\bfe)}\vv{1}_{|\upnu\tau-A'|\leq A'^{0.35}}\f{\<\tau\>}{\upnu}\nr{\wh{\<\na_{\fz,\fv}\>\udl{f_\infty}}((\upmu-\upnu)\bfe,(\upmu t-\upnu\tau)\bfe)} d\tau}_{L^2_t[\mathsf{T}_k,\tilde{T}]}+\f{1}{\upmu^\f12}e^{-99A'^{0.34}}Z(\tilde{T})\\
		&\leq \f{1}{\upmu^\f12}e^{-99A'^{0.34}}Z(\tilde{T})+\f{2A'}{{k'_2}^{2.5}}\sup_{\upnu\geq k'_2}\rr{\nnr{\wh{\rho}(t,\upnu\bfe)}_{L^2_t[\mathsf{T}_k,\tilde{T}]}\upnu^\f12}\sup_{\mathsf{T}_k\leq\tau\leq \tilde{T}}\rr{\sum_{\upnu\geq k'_2}\int_{\mathsf{T}_k}^{\tilde{T}}\nr{\wh{\<\na_{\fz,\fv}\>\udl{f_\infty}}((\upmu-\upnu)\bfe,(\upmu t-\upnu \tau)\bfe)}dt}^\f12\\
		&\times\sup_{\mathsf{T}_k\leq t\leq \tilde{T}}\rr{\sum_{\upnu\geq k'_2}\int_{\mathsf{T}_k}^{\tilde{T}}\nr{\wh{\<\na_{\fz,\fv}\>\udl{f_\infty}}((\upmu-\upnu)\bfe,(\upmu t-\upnu \tau)\bfe)}d\tau}^\f12
		\leq\f{\eps_0^{0.99}}{\upmu^\f12}\sup_{\upnu\geq k'_2}\rr{\nnr{\wh{\rho}(t,\upnu\bfe)}_{L^2_t[\mathsf{T}_k,\tilde{T}]}\upnu^\f12}+\f{1}{\upmu^\f12}e^{-99A'^{0.34}}Z(\tilde{T}).
		\ea\]
		Here we used $A'/{k'_2}^3\leq2$. The additional factor $\eps_0$ comes from the background wave, so the off-diagonal part is perturbative. Combining the estimates and absorbing it yields
		\begin{equation}\label{rho41}
			\sup_{\upmu\geq k'_2}\nnr{\wh{\rho}(t,\upmu\bfe)}_{L^2_t[\mathsf{T}_k,\tilde{T}]}\upmu^\f12\leq C_\theta(2\pi)^\f12Z(\tilde{T})+e^{-0.85\<\f{\mathsf{T}_k}{A_*}\>^\sss}\sup_{\substack{\mathsf{T}_k\leq \tau\leq\tilde{T}\\\upsigma\in\Z}}\nnr{\wh{f}(\tau,\upsigma\bfe,\xi\bfe)}_{L^2_\upxi}.
		\end{equation}
		
		\noindent\textit{Step 3. Recovery of the profile.}
		We now insert \eqref{rho41} into the profile equation \eqref{f41}. The same time-frequency convolution estimate as in Step 4 of Proposition \ref{prop:forward-recurrence} gives
		\[\ba
		&\nnr{\sum_{\upnu\geq k'_2}\int_{\mathsf{T}_k}^t\wh{\rho}(\tau,\upnu\bfe)\f{\upxi-\upsigma \tau}{\upnu}e^{\<\mrxi\>^\sss}\pa_\xi^\alpha\wh{\udl{f_\infty}}((\upsigma-\upnu)\bfe,\xi-\upnu\tau\bfe) d\tau}_{\ltaa L^2_\xi}\\
		&\leq\nnr{\sum_{\upnu\geq k'_2}\int_{\mathsf{T}_k}^t\nr{\wh{\rho}(\tau,\upnu\bfe)}\vv{1}_{|\upnu\tau-A'|\leq A'^{0.35}}\f{\<\tau\>}{\upnu}\nr{\sF\rr{e^{\<\na_{\mrv}\>^\sss}\<\na_{\fz,\fv}\>v^\alpha\udl{f_\infty}}((\upsigma-\upnu)\bfe,\xi-\upnu\tau\bfe)} d\tau}_{\ltaa L^2_\xi}+e^{-99A'^{0.34}}Z(t)\\
		&\leq\f{3A'}{{k'_2}^3}\sup_{\upnu\geq k'_2}\rr{\nnr{\wh{\rho}(t,\upnu\bfe)}_{L^2_t[\mathsf{T}_k,\tilde{T}]}\upnu^\f12}\rr{\sum_{\upnu}\nnr{\sF\rr{e^{\<\na_{\mrv}\>^\sss}\<\na_{\fz,\fv}\>v^\alpha\udl{f_\infty}}(\upnu\bfe,\xi)}_{\ltaa L^1_{\upxi}L^2_{\mrxi}}}+e^{-99A'^{0.34}}Z(t)\\
		&\leq\eps_0^{0.99}\sup_{\upnu\geq k'_2}\rr{\nnr{\wh{\rho}(t,\upnu\bfe)}_{L^2_t[\mathsf{T}_k,\tilde{T}]}\upnu^\f12}+e^{-99A'^{0.34}}Z(t).
		\ea\]
		Again $A'/{k'_2}^{3}\leq2$ makes the high-mode transfer perturbative. Adding the remainders from Step 1 gives
		\begin{equation}\label{f42}
			\begin{aligned}
				\nnr{e^{\<\mrxi\>^\sss}\pa_\xi^\alpha\rr{\wh{f}(t,\upsigma\bfe,\xi)-\int_{\mathsf{T}_k}^t\wh{Q}(\tau,\upsigma\bfe,\xi)d\tau}}_{\lta L^2_\xi}&\leq \eps_0^{0.97}\sup_{\upnu\geq k'_2}\rr{\nnr{\wh{\rho}(\tau,\upnu\bfe)}_{L^2_\tau[\mathsf{T}_k,t]}\upnu^\f12} \\
				&\quad +e^{-0.9\<\f{\mathsf{T}_k}{A_*}\>^\sss}\sup_{\substack{\mathsf{T}_k\leq \tau\leq t\\\upsigma\in\Z}}\nnr{e^{\<\mrxi\>^\sss}\pa_\xi^\alpha\wh{f}(\tau,\upsigma\bfe,\xi)}_{\ltaa L^2_\xi}+e^{-98.9A'^{0.34}}Z(t).
			\end{aligned}
		\end{equation}
		
		\noindent\textit{Step 4. Coupled closure.}
		The right-hand side of \eqref{rho41} contains the profile norm, while \eqref{f42} contains the density norm. Their coefficients are small by the choice of $\mathfrak{T}_k$, $k'_2$, and $\eps_0$. A continuity argument therefore closes the two bounds simultaneously and gives, for $\mathsf{T}_k\leq t\leq t'_{k'_2-1}$,
		\[\sup_{\upmu\geq k'_2}\nnr{\wh{\rho}(\tau,\upmu\bfe)}_{L^2_t[\mathsf{T}_k,t]}\upmu^\f12\leq 3C_\theta Z(t),\quad\nnr{e^{\<\mrxi\>^\sss}\pa_\xi^\alpha\rr{\wh{f}(t,\upsigma\bfe,\xi)-\int_{\mathsf{T}_k}^t\wh{Q}(\tau,\upsigma\bfe,\xi)d\tau}}_{\ltaa L^2_\xi}\leq \eps_0^{0.95}Z(t).\]
	\end{proof}
	
	We next treat $k'_0<k'_2$. In this configuration $\mathsf{T}_k$ already belongs to the low-mode part of the chain. Before the first complete critical interval begins, only the mode $k'_0$ can resonate. The following proposition isolates this single response and records the spatial modes to which it can transfer: the main background wave creates $k'_0\pm1$, while $f_{*,\infty}$ contributes the decaying tail in $\upsigma-k'_0$.
	
	\begin{prop}\label{4.3}
		Assume $A'>\mathsf{T}_k/2$ and $k'_0<k'_2$. On $[\mathsf{T}_k,t'_{k'_0-1}]$,
		\begin{equation}\label{ttk0}
			\begin{aligned}
				&\nnr{e^{\<\na_{\mrv}\>^\sss} v^{\alpha}\rr{\sF_zf(t,\upsigma\bfe,v)-\int_{\mathsf{T}_k}^t\sF_zQ(\tau,\upsigma\bfe,v)d\tau}}_{\lta L^2_v} \\
				&\leq \f12A'Z(t)\rr{\vv{1}_{\upsigma=k'_0\pm1}+\eps_*\<\upsigma-k'_0\>^{-8}}+\rr{e^{-98A'^{0.34}}+e^{-0.8\<\f{\mathsf{T}_k}{A_*}\>^\sss} }Z(t).
			\end{aligned}
		\end{equation}
		for any $\upsigma\in\Z$.
	\end{prop}
	\begin{proof}
		We split the proof into three steps: reduction to the single active density, Volterra control of that density and the induced profile, and a coupled continuity argument.
		
		\noindent\textit{Step 1. Reduction to the first active mode.}
		On $[\mathsf{T}_k,t'_{k'_0-1}]$, due to the Fourier localization of the source, only the mode \(k'_0\) can be near resonance. Hence
		\[\ba
		&\wh{\rho}(t,\upnu\bfe)=\wh{\sR_2}(t,\upnu\bfe),\quad \text{for}\ \upnu\neq k'_0,\\
		&\wh{\rho}(t,k'_0\bfe)=\int_{\mathsf{T}_k}^t\wh{Q}(\tau,k'_0\bfe,k'_0t\bfe)d\tau-\int_{\mathsf{T}_k}^t\wh{\rho}(\tau,k'_0\bfe)(t-\tau)\sF\rr{ f^0+\udl{f_\infty} }(0,k'_0 (t-\tau)\bfe)d\tau+\wh{\sR}(t,k'_0\bfe),\\
		&\wh{f}(t,\upsigma\bfe,\xi)-\int_{\mathsf{T}_k}^t\wh{Q}(\tau,\upsigma\bfe,\xi)d\tau=-\int_{\mathsf{T}_k}^t\wh{\rho}(\tau,k'_0\bfe)\f{\upxi-\upsigma \tau}{k'_0}\sF\rr{ f^0+\udl{f_\infty} }((\upsigma-k'_0)\bfe,\xi-k'_0\tau\bfe) d\tau+\wh{\sR}(t,\upsigma\bfe,\xi),
		\ea\]
		with the error terms satisfying
		\[\nr{\wh{\sR_2}(t,\upnu\bfe)}\leq e^{-99A'^{0.34}}e^{-0.124\<\f{\upnu t-A'}{A_*}\>^\sss}Z(t)\]
		and
		\[\ba
		&\nnr{\wh{\sR}(t,k'_0\bfe)}_{L^2_t([\mathsf{T}_k,\tilde{T}])}{k'_0}^\f12\leq e^{-0.9\<\f{\mathsf{T}_k}{A_*}\>^\sss}\rr{\nnr{\wh{\rho}(t,k'_0\bfe)}_{L^2_t([\mathsf{T}_k,\tilde{T}])}+\sup_{\substack{\mathsf{T}_k\leq\tau\leq\tilde{T}\\\upsigma\in\Z}}\nnr{\wh{f}(\tau,\upsigma\bfe,\upxi\bfe)}_{L^2_\upxi}}+e^{-99A'^{0.34}}Z(\tilde{T}),\\
		&\nnr{e^{\<\mrxi\>^\sss}\pa_{\xi}^\alpha\wh{\sR}(t,\upsigma\bfe,\xi)}_{L^2_\xi}\leq e^{-0.9\<\f{\mathsf{T}_k}{A_*}\>^\sss}\rr{\nnr{\wh{\rho}(\tau,k'_0\bfe)}_{L^2_\tau([\mathsf{T}_k,t])}+\sup_{\substack{\mathsf{T}_k\leq\tau\leq t\\\upsigma\in\Z}}\nnr{e^{\<\mrxi\>^\sss}\pa_{\xi}^\alpha\wh{f}(\tau,\upsigma\bfe,\xi)}_{\ltaa L^2_\xi}}+e^{-99A'^{0.34}}Z(t)
		\ea\]
		for any $\mathsf{T}_k\leq\tilde{T}\leq t'_{k'_0-1}$ and $\mathsf{T}_k\leq t\leq t'_{k'_0-1}$. These estimates follow from the decomposition in Step 1 of Proposition \ref{4.2}, with the active set $\{\upnu\geq k'_2\}$ replaced by the single mode $k'_0$.
		
		\noindent\textit{Step 2. Density and profile estimates.}
		We use $\sR$ for remainders satisfying the bounds from Step 1; its value may change from line to line. Inverting the scalar Volterra equation for the only active density gives
		\[\wh{\rho}(t,k'_0\bfe)=\int_{\mathsf{T}_k}^t\int_{\mathsf{T}_k}^\tau\wh{Q}(s,k'_0\bfe,k'_0\tau\bfe)ds\rr{\delta_0(t-\tau)-R(t-\tau,k'_0\bfe)}d\tau+ \wh{\sR}(t,k'_0\bfe).\]
		Hence the resolvent estimate \eqref{eR} gives
		\begin{equation}\label{rho44}
			\nnr{\wh{\rho}(t,k'_0\bfe)}_{L^2_t([\mathsf{T}_k,\tilde{T}])}{k'_0}^\f12\leq C_\theta Z(\tilde{T})+e^{-0.89\<\f{\mathsf{T}_k}{A_*}\>^\sss}\sup_{\substack{\mathsf{T}_k\leq\tau\leq\tilde{T}\\\upsigma\in\Z}}\nnr{\wh{f}(\tau,\upsigma\bfe,\upxi\bfe)}_{L^2_\upxi}+e^{-98.9A'^{0.34}}Z(\tilde{T}).
		\end{equation}
		The induced profile is estimated as in Step 3 of Proposition \ref{4.2}. The Fourier support of $\udl{f_\infty}$ separates the leading transfers to $k'_0\pm1$ from the small spatial tail:
		\[\ba
		&\nnr{\int_{\mathsf{T}_k}^t\wh{\rho}(\tau,k'_0\bfe)\f{\upxi-\upsigma \tau}{k'_0}e^{\<\mrxi\>^\sss}\pa_{\xi}^\alpha\wh{\udl{f_\infty} }((\upsigma-k'_0)\bfe,\xi-k'_0\tau\bfe) d\tau}_{L^2_\xi}\leq2A'\nnr{\wh{\rho}(\tau,k'_0\bfe)}_{L^2_t([\mathsf{T}_k,t])}\\
		&\times\nnr{\sF\rr{e^{\<\na_{\mrv}\>^\sss}\<\na_{\fz,\fv}\>v^\alpha\udl{f_\infty}}((\upsigma-k'_0)\bfe,\xi)}_{L^1_{\upxi} L^2_{\mrxi}}
		\leq\eps_0^\f12A'\nnr{\wh{\rho}(\tau,k'_0\bfe)}_{L^2_t([\mathsf{T}_k,t])}\rr{\vv{1}_{\upsigma=k'_0\pm1}+\eps_*\<\upsigma-k'_0\>^{-8}}.
		\ea\]
		Thus
		\begin{equation}\label{f43}
			\begin{aligned}
				&\nnr{e^{\<\mrxi\>^\sss}\pa_{\xi}^\alpha\rr{\wh{f}(t,\upsigma\bfe,\xi)-\int_{\mathsf{T}_k}^t\wh{Q}(\tau,\upsigma\bfe,\xi)d\tau}}_{\ltaa L^2_\xi} \\
				&\quad\leq\eps_0^\f12A'\nnr{\wh{\rho}(\tau,k'_0\bfe)}_{L^2_t([\mathsf{T}_k,t])}\rr{\vv{1}_{\upsigma=k'_0\pm1}+\eps_*\<\upsigma-k'_0\>^{-8}} \\
				&\qquad +e^{-0.9\<\f{\mathsf{T}_k}{A_*}\>^\sss}\Bigl(\nnr{\wh{\rho}(\tau,k'_0\bfe)}_{L^2_\tau([\mathsf{T}_k,t])} \\
				&\qquad\qquad +\sup_{\substack{\mathsf{T}_k\leq\tau\leq t\\\upsigma\in\Z}}\nnr{e^{\<\mrxi\>^\sss}\pa_{\xi}^\alpha \wh{f}(\tau,\upsigma\bfe,\xi)}_{\ltaa L^2_\xi}\Bigr)+e^{-99A'^{0.34}}Z(t).
			\end{aligned}
		\end{equation}
		
		\noindent\textit{Step 3. Closure on the partial critical interval.}
		The density bound \eqref{rho44} and profile bound \eqref{f43} form the same coupled system as in Proposition \ref{4.2}. Absorbing the exponentially small feedback gives, for $\mathsf{T}_k\leq t\leq t'_{k'_0-1}$,
		\[\ba
		&\nnr{\wh{\rho}(\tau,k'_0\bfe)}_{L^2_t([\mathsf{T}_k,t])}{k'_0}^\f12\leq2C_\theta Z(t),\\ &\nnr{e^{\<\mrxi\>^\sss}\pa_{\xi}^\alpha\rr{\wh{f}(t,\upsigma\bfe,\xi)-\int_{\mathsf{T}_k}^t\wh{Q}(\tau,\upsigma\bfe,\xi)d\tau}}_{\ltaa L^2_\xi}\leq\f14A'Z(t)\rr{\vv{1}_{\upsigma=k'_0\pm1}+\eps_*\<\upsigma-k'_0\>^{-8}}\\
		&\qquad+\rr{e^{-98A'^{0.34}}+e^{-0.8\<\f{\mathsf{T}_k}{A_*}\>^\sss}}Z(t).
		\ea\]
	\end{proof}
	
	After the entry interval, the remaining critical intervals are complete and can be treated one at a time. Set
	\[\gamma(t,\upmu,\fv)=\int_{\R^2}\sF_z f(t,\upmu\bfe,v)d\mrv-\int_{\mathsf{T}_k}^{t}\int_{\R^2}\sF_z Q(\tau,\upmu\bfe,v)d\mrv d\tau. \]
	Thus $\gamma$ is the density amplitude after subtracting the response produced directly by $Q$. The next proposition establishes an approximate recurrence relation on each critical interval, with the resonant transfer described by the same one-step echo operator as in Section 3.
	
	\begin{prop}\label{4.4}
		For any $1\leq\upmu\leq\min\{k'_0-1,k'_2-1\}$, we have the following approximate recurrence relation on $[t'_\upmu,t'_{\upmu-1}]$
		\begin{equation}\label{rec41}
			\begin{aligned}
				\sF_zf(t'_{\upmu-1},\upsigma\bfe,v)&=\sF_zf(t'_{\upmu},\upsigma\bfe,v)+\int_{t'_\upmu}^{t'_{\upmu-1}}\sF_z Q(\tau,\upsigma\bfe,v)d\tau \\
				&\quad +\cD_{\upsigma,\upmu}\rr{\gamma(t'_{\upmu},\upmu,\fv)+\sF_\upzeta^{-1}\rr{\int_{\mathsf{T}_k}^{\f\upzeta\upmu} \wh{Q}(\tau,\upmu\bfe,\upzeta\bfe)d\tau}}+\sR_{\upsigma,\upmu}(v)
			\end{aligned}
		\end{equation}
		with the error terms satisfying
		\[\ba
		&\nnr{e^{\<\na_{\mrv}\>^\sss}v^{\alpha}\sR_{\upsigma,\upmu}(v)}_{\lta L^2_v}\leq e^{-0.8\<\f{\mathsf{T}_k}{A_*}\>^\sss}\sup_{\upsigma\in\Z}\nnr{e^{\<\na_{\mrv}\>^\sss}v^\alpha\rr{\sF_z f(t'_{\upmu},\upsigma\bfe,v)-\int_{\mathsf{T}_k}^{t'_\upmu}\sF_z Q(\tau,\upsigma\bfe,v)d\tau}}_{\lta L^2_v}\\
		&+\rr{e^{-0.8\<\f{\mathsf{T}_k}{A_*}\>^\sss}+e^{-98A'^{0.34}}}Z(t'_{\upmu-1})\ea\]
		and
		\[\nnr{\int_{\R^2}\sR_{\upsigma,\upmu}(v)d\mrv}_{ L^2_\fv}\leq e^{-0.8\<\f{\mathsf{T}_k}{A_*}\>^\sss}\sup_{\upsigma\in\Z}\nnr{\gamma(t'_{\upmu},\upsigma,\fv)}_{ L^2_\fv}+\rr{e^{-0.8\<\f{\mathsf{T}_k}{A_*}\>^\sss}+e^{-98A'^{0.34}}}Z(t'_{\upmu-1}).\]
		Moreover the following estimate holds on $[t'_\upmu,t'_{\upmu-1}]$:
		\begin{equation}\label{fjeif}
			\begin{aligned}
				&\nnr{e^{\<\na_{\mrv}\>^\sss}v^{\alpha}\rr{\sF_zf(t,\upsigma\bfe,v)-\sF_zf(t'_{\upmu},\upsigma\bfe,v)-\int_{t'_\upmu}^{t} \sF_z Q(\tau,\upsigma\bfe,v)d\tau}}_{\lta L^2_v} \\
				&\quad\leq \eps_0^{0.8}A'\rr{\nnr{\gamma(t'_{\upmu},\upmu,\fv)}_{L^2_\fv}+Z(t)} \\
				&\qquad + e^{-0.8\<\f{\mathsf{T}_k}{A_*}\>^\sss}\sup_{\upsigma\in\Z}\nnr{e^{\<\na_{\mrv}\>^\sss}v^\alpha\rr{\sF_z f(t'_{\upmu},\upsigma\bfe,v)-\int_{\mathsf{T}_k}^{t'_\upmu}\sF_z Q(\tau,\upsigma\bfe,v)d\tau}}_{\lta L^2_v}.
			\end{aligned}
		\end{equation}
	\end{prop}
	\begin{proof}
		The proof has two steps: isolate the single resonant density and sharpen the error bounds, then integrate the resulting leading term to identify the transfer operator $\cD_{\upsigma,\upmu}$.
		
		\noindent\textit{Step 1. Single-resonance reduction and refined errors.}
		On $[t'_{\upmu},t'_{\upmu-1}]$, only the density mode $\upmu$ is resonant. Therefore
		\[\ba
		&\wh{\rho}(t,\upnu\bfe)=\wh{\sR_2}(t,\upnu\bfe),\quad \text{for}\ \upnu\neq \upmu,\\
		&\wh{\rho}(t,\upmu\bfe)=\wh{f}(t'_\upmu,\upmu\bfe,\upmu t\bfe)+\int_{t'_{\upmu}}^t\wh{Q}(\tau,\upmu\bfe,\upmu t\bfe)d\tau-\int_{t'_{\upmu}}^t\wh{\rho}(\tau,\upmu\bfe)(t-\tau)\sF\rr{ f^0+\udl{f_\infty} }(0,\upmu (t-\tau))d\tau+\wh{\sR}(t,\upmu\bfe),\\
		&\wh{f}(t,\upsigma\bfe,\xi)-\int_{t'_{\upmu}}^t\wh{Q}(\tau,\upsigma\bfe,\xi)d\tau=\wh{f}(t'_\upmu,\upsigma\bfe,\xi)-\int_{t'_{\upmu}}^t\wh{\rho}(\tau,\upmu\bfe)\f{\upxi-\upsigma \tau}{\upmu}\sF\rr{ f^0+\udl{f_\infty} }(\upsigma-\upmu,\xi-\upmu\tau) d\tau+\wh{\sR}(t,\upsigma\bfe,\xi),
		\ea\]
		with the error terms satisfying
		\[\nr{\wh{\sR_2}(t,\upnu\bfe)}\leq e^{-99A'^{0.34}}e^{-0.124\<\f{\upnu t-A'}{A_*}\>^\sss}Z(t)\]
		and
		\begin{equation}\label{wdq}
			\begin{aligned}
				\nnr{\wh{\sR}(t,\upmu\bfe)}_{L^2_t([t'_{\upmu},\tilde{T}])}\upmu^\f12&\leq e^{-0.9\<\f{\mathsf{T}_k}{A_*}\>^\sss}\rr{\nnr{\wh{\rho}(t,\upmu\bfe)}_{L^2_t([t'_{\upmu},\tilde{T}])}+\sup_{\substack{t'_{\upmu}\leq\tau\leq\tilde{T}\\\upsigma\in\Z}}\nnr{\wh{f}(\tau,\upsigma\bfe,\upxi\bfe)}_{L^2_\upxi}}+e^{-99A'^{0.34}}Z(t), \\
				\nnr{\wh{\sR}(t,\upsigma\bfe,\upxi\bfe)}_{L^2_\upxi}&\leq e^{-0.9\<\f{\mathsf{T}_k}{A_*}\>^\sss}\rr{\nnr{\wh{\rho}(\tau,\upmu\bfe)}_{L^2_\tau([t'_{\upmu},t])}+\sup_{\substack{t'_{\upmu}\leq\tau\leq t\\\upsigma\in\Z}}\nnr{\wh{f}(\tau,\upsigma\bfe,\upxi\bfe)}_{L^2_\upxi}}+e^{-99A'^{0.34}}Z(\tilde{T})
			\end{aligned}
		\end{equation}
		\begin{equation}\label{tya}
			\begin{aligned}
				\nnr{e^{\<\mrxi\>^\sss}\pa_{\xi}^\alpha\wh{\sR}(t,\upsigma\bfe,\xi)}_{\ltaa L^2_\xi}&\leq e^{-0.9\<\f{\mathsf{T}_k}{A_*}\>^\sss}\rr{\nnr{\wh{\rho}(\tau,\upmu\bfe)}_{L^2_\tau([t'_{\upmu},t])}+\sup_{\substack{t'_{\upmu}\leq\tau\leq t\\\upsigma\in\Z}}\nnr{e^{\<\mrxi\>^\sss}\pa_{\xi}^\alpha\wh{f}(\tau,\upsigma\bfe,\xi)}_{\ltaa L^2_\xi}} \\
				&\quad +e^{-99A'^{0.34}}Z(t),
			\end{aligned}
		\end{equation}
		for any $t'_{\upmu}\leq \tilde{T}\leq t'_{\upmu-1}$ and $t'_{\upmu}\leq t\leq t'_{\upmu-1}$. The proof is the same frequency separation used in Step 1 of Proposition \ref{4.3}.
		
		For iteration, the error must be expressed in terms of the profile at the left endpoint of the interval, rather than a supremum over the interval. The density satisfies
		\begin{equation}\label{cwa}
			\begin{aligned}
				\wh{\rho}(t,\upmu\bfe)&=\int_{\mathsf{T}_k}^t\rr{\wh{f}(t'_{\upmu},\upmu\bfe,\upmu \tau\bfe)-\int_{\mathsf{T}_k}^{t'_{\upmu}}\wh{Q}(s,\upmu\bfe,\upmu\tau\bfe)ds+\int_{\mathsf{T}_k}^\tau\wh{Q}(s,\upmu\bfe,\upmu\tau\bfe)ds} \\
				&\quad \times\rr{\delta_0(t-\tau)-R(t-\tau,\upmu\bfe)}d\tau+ \wh{\sR}(t,\upmu\bfe).
			\end{aligned}
		\end{equation}
		Thus by \eqref{eR} we get
		\begin{equation}\label{poi}
			\begin{aligned}
				\nnr{\wh{\rho}(\tau,\upmu\bfe)}_{L^2_\tau([t'_{\upmu},t])}\upmu^\f12&\leq C_\theta\rr{\nnr{\wh{f}(t'_{\upmu},\upmu\bfe,\upxi\bfe)-\int_{\mathsf{T}_k}^{t'_{\upmu}}\wh{Q}(s,\upmu\bfe,\upxi\bfe)ds}_{L^2_\upxi}+Z(t)} \\
				&\quad +e^{-0.89\<\f{\mathsf{T}_k}{A_*}\>^\sss}\sup_{\substack{t'_{\upmu}\leq\tau\leq t\\\upsigma\in\Z}}\nnr{\wh{f}(\tau,\upsigma\bfe,\upxi\bfe)}_{L^2_\upxi}.
			\end{aligned}
		\end{equation}
		The profile estimate from Step 3 of Proposition \ref{4.2} gives
		\begin{equation}\label{qwe}
			\begin{aligned}
				\nnr{\wh{f}(t,\upsigma\bfe,\upxi\bfe)-\int_{t'_{\upmu}}^t\wh{Q}(\tau,\upsigma\bfe,\upxi\bfe)d\tau-\wh{f}(t'_\upmu,\upsigma\bfe,\upxi\bfe)}_{L^2_\upxi}&\leq\eps_0^{0.9}A'\nnr{\wh{\rho}(\tau,\upmu\bfe)}_{L^2_\tau([t'_{\upmu},t])}\upmu^\f12 \\
				&\quad +e^{-99A'^{0.34}}Z(t).
			\end{aligned}
		\end{equation}
		A continuity argument applied to \eqref{poi} and \eqref{qwe} now yields, on $[t'_{\upmu},t'_{\upmu-1}]$,
		\[\ba
		&\nnr{\wh{\rho}(\tau,\upmu\bfe)}_{L^2_\tau([t'_{\upmu},t])}\upmu^\f12
		\leq 2C_\theta\rr{\nnr{\wh{f}(t'_{\upmu},\upmu\bfe,\upxi\bfe)-\int_{\mathsf{T}_k}^{t'_{\upmu}}\wh{Q}(s,\upmu\bfe,\upxi\bfe)ds}_{L^2_\upxi}+Z(t)}\\
		&+e^{-0.88\<\f{\mathsf{T}_k}{A_*}\>^\sss}\sup_{\upsigma\in\Z}\nnr{\wh{f}(t'_{\upmu},\upsigma\bfe,\upxi\bfe)-\int_{\mathsf{T}_k}^{t'_\upmu}\wh{Q}(s,\upsigma\bfe,\upxi\bfe)ds}_{L^2_\upxi},\\
		&\sup_{\substack{t'_{\upmu}\leq\tau\leq t\\\upsigma\in\Z}}\nnr{\wh{f}(\tau,\upsigma\bfe,\upxi\bfe)}_{L^2_\upxi}\leq \eps_0^{0.8}A'\rr{\sup_{\upsigma}\nnr{\wh{f}(t'_{\upmu},\upsigma\bfe,\upxi\bfe)-\int_{\mathsf{T}_k}^{t'_\upmu}\wh{Q}(s,\upsigma\bfe,\upxi\bfe)ds}_{L^2_\upxi}+Z(t)}.		
		\ea\]
		Substitution into \eqref{wdq} removes the interval supremum and gives
		\[\nnr{\wh{\sR}(t,\upsigma\bfe,\upxi\bfe)}_{L^2_\upxi}\leq e^{-0.8\<\f{\mathsf{T}_k}{A_*}\>^\sss}\sup_{\upsigma}\nnr{\wh{f}(t'_{\upmu},\upsigma\bfe,\upxi\bfe)-\int_{\mathsf{T}_k}^{t'_\upmu}\wh{Q}(s,\upsigma\bfe,\upxi\bfe)ds}_{L^2_\upxi}+\rr{e^{-0.88\<\f{\mathsf{T}_k}{A_*}\>^\sss}+e^{-98.5A'^{0.34}}}Z(t).\]
		The weighted version and the corresponding control inside the interval are
		\[\ba
		&\nnr{e^{\<\mrxi\>^\sss}\pa_\xi^\alpha\wh{\sR}(t,\upsigma\bfe,\xi)}_{\ltaa L^2_\xi}\leq e^{-0.8\<\f{\mathsf{T}_k}{A_*}\>^\sss}\sup_{\upsigma}\nnr{e^{\<\mrxi\>^\sss}\pa_\xi^\alpha\rr{\wh{f}(t'_{\upmu},\upsigma\bfe,\xi)-\int_{\mathsf{T}_k}^{t'_\upmu}\wh{Q}(s,\upsigma\bfe,\xi)ds}}_{\ltaa L^2_\xi}\\
		&+\rr{e^{-0.88\<\f{\mathsf{T}_k}{A_*}\>^\sss}+e^{-98.5A'^{0.34}}}Z(t),\\
		&\nnr{e^{\<\mrxi\>^\sss}\pa_\xi^\alpha\rr{\wh{f}(t,\upsigma\bfe,\xi)-\int_{t'_{\upmu}}^t\wh{Q}(\tau,\upsigma\bfe,\xi)d\tau-\wh{f}(t'_\upmu,\upsigma\bfe,\xi)}}_{\ltaa L^2_\xi}
		\leq\eps_0^{0.9}A'\rr{\nnr{\wh{f}(t'_{\upmu},\upmu\bfe,\upxi\bfe)\right.\right.\\
		&\left.\left.-\int_{\mathsf{T}_k}^{t'_\upmu}\wh{Q}(s,\upmu\bfe,\upxi\bfe)ds}_{L^2_\upxi}+Z(t)}+e^{-0.8\<\f{\mathsf{T}_k}{A_*}\>^\sss}\sup_{\upsigma\in\Z}\nnr{e^{\<\mrxi\>^\sss}\pa_\xi^\alpha\rr{\wh{f}(t'_{\upmu},\upsigma\bfe,\xi)-\int_{\mathsf{T}_k}^{t'_\upmu}\wh{Q}(s,\upsigma\bfe,\xi)ds}}_{\ltaa L^2_\xi}.
		\ea\]
		The first inequality is the refined form of \eqref{tya}, while the second is exactly \eqref{fjeif}.
		
		\noindent\textit{Step 2. Identification of the transfer operator.}
		Insert the Volterra representation \eqref{cwa} into the profile equation at $t'_{\upmu-1}$. This gives
		\[\ba
		&\wh{f}(t'_{\upmu-1},\upsigma\bfe,\xi)=\wh{f}(t'_\upmu,\upsigma\bfe,\xi)+\int_{t'_\upmu}^{t'_{\upmu-1}}\wh{Q}(\tau,\upsigma\bfe,\xi)d\tau-\int_\R\int_\R\rr{\wh{f}(t'_\upmu,\upmu\bfe,\upmu s\bfe)+\int_{t'_\upmu}^s\wh{Q}(s',\upmu\bfe,\upmu s\bfe)ds'}\\
		&\rr{\delta_0(\tau-s)-R(\tau-s,\upmu\bfe)\vv{1}_{\tau-s>0}}
		\f{\upxi-\upsigma \tau}{\upmu}\sF\rr{ f^0+\udl{f_\infty} }((\upsigma-\upmu)\bfe,\xi-\upmu\tau\bfe) dsd\tau+\wh{\sR}(t'_{\upmu-1},\upsigma\bfe,\xi)\\
		&=\wh{f}(t'_\upmu,\upsigma\bfe,\xi)+\int_{t'_\upmu}^{t'_{\upmu-1}}\wh{Q}(\tau,\upsigma\bfe,\xi)d\tau-\int_\R\int_\R\rr{\wh{f}(t'_\upmu,\upmu\bfe,\upeta\bfe)+\int_{t'_\upmu}^{\f\upeta\upmu}\wh{Q}(s',\upmu\bfe,\upeta\bfe)ds'}
		\left(\f{\upmu}{|\upmu|}\delta_0(\upzeta-\upeta)\right.\\&\left.-\f1{\upmu}R(\f{\upzeta-\upeta}{\upmu},\upmu\bfe)\vv{1}_{\upzeta-\upeta>0}\right)\rr{\f{\upxi-\upzeta}{\upmu^2}-\f{\upsigma-\upmu}{\upmu^3}\upzeta}
		\sF\rr{ f^0+\udl{f_\infty} }((\upsigma-\upmu)\bfe,\xi-\upzeta\bfe) d\upeta d\upzeta+\wh{\sR}(t'_{\upmu-1},\upsigma\bfe,\xi).
		\ea\]
		After the inverse Fourier transform, the double integral is precisely $\cD_{\upsigma,\upmu}$ applied to the renormalized resonant amplitude. This proves \eqref{rec41}; the remainder and interval bounds follow from Step 1.
	\end{proof}
	
	We next iterate the recurrence through the low modes. The growth condition on $Z$ is chosen to absorb the possible one-step amplification $\max\{1,\eps_0A'/\upmu^3\}$, while
	\[\mathtt{G}_\theta(\upmu)=\prod_{\upnu=\upmu}^{\min\{k'_0,k'_2\}-1}\Etnu\rr{1+\f{10^2}{\upnu^2}}.\]
	collects the harmless resolvent and same-mode losses. The statement distinguishes whether the iteration has crossed the source mode $k'$. At that one interval, the source enters the resonant amplitude directly, which enlarges only the numerical constants and not the form of the estimate.
	
	\begin{prop}\label{4.5}
		Assume that, the envelope satisfies that for fixed $A'>0$ and every $1\leq\upmu\leq\min\{k'_0,k'_2\}-1$, 
		\begin{equation}\label{growth1}
			\max\{1,\f{\eps_0A'}{\upmu^3}\}\rr{1+\eps_0^{10}}Z(t'_\upmu)\leq Z(t'_{\upmu-1})\rr{1+\f{10^2}{\upmu^2}}.
		\end{equation}
		For any $0\leq\upmu\leq\min\{k'_0,k'_2\}-1$, if $k'\notin[\upmu+1,\min\{k'_0,k'_2\}-1]$, then 
		\begin{equation}\label{rpe}
			\nnr{\gamma(t'_\upmu,\upsigma,\fv)}_{L^2_\fv}\leq A'Z(t'_\upmu)\mathtt{G}_\theta(\upmu+1)\left\{\begin{aligned}
				&10^4\rr{3-\<\upsigma-\upmu-1\>^{-7}},\quad&\upsigma\geq\upmu+1, \\
				&1,\quad&\upsigma=\upmu, \\
				&{A'}^{-1}+50\eps_{*}\<\upmu-\upsigma\>^{-7},\quad&\upsigma\leq\upmu-1.
			\end{aligned}\right.
		\end{equation}
		If $k'\in[\upmu+1,\min\{k'_0,k'_2\}-1]$, then 
		\begin{equation}\label{rpe1}
			\nnr{\gamma(t'_\upmu,\upsigma,\fv)}_{L^2_\fv}\leq A'Z(t'_\upmu)\mathtt{G}_\theta(\upmu+1)\left\{\begin{aligned}
				&10^5\rr{3-\<\upsigma-\upmu-1\>^{-7}},\quad&\upsigma\geq\upmu+1, \\
				&2,\quad&\upsigma=\upmu, \\
				&{A'}^{-1}+10^2\eps_{*}\<\upmu-\upsigma\>^{-7},\quad&\upsigma\leq\upmu-1.
			\end{aligned}\right.
		\end{equation}	
	\end{prop}
	\begin{proof}
		We split the proof into three steps.
		
		\noindent\textit{Step 1. One-step estimate.}
		Integrating \eqref{rec41} in the transverse velocity variables gives
		\[\gamma(t'_{\upmu-1},\upsigma,\fv)=\gamma(t'_{\upmu},\upsigma,\fv)+\fD_{\upsigma,\upmu}\rr{\gamma(t'_{\upmu},\upmu,\fv)+\sF_{\upzeta}^{-1}\rr{\int_{\mathsf{T}_k}^{\f{\upzeta}{\upmu}}\wh{Q}(\tau,\upmu\bfe,\upzeta\bfe)d\tau }}+\fR_{\upsigma,\upmu}(\fv),\]
		with the error term satisfying
		\[\nnr{\fR_{\upsigma,\upmu}(\fv)}_{L^2_\fv}\leq e^{-0.8\<\f{\mathsf{T}_k}{A_*}\>^\sss}\sup_{\upsigma\in\Z}\nnr{\gamma(t'_{\upmu},\upsigma,\fv)}_{ L^2_\fv}+\rr{e^{-0.8\<\f{\mathsf{T}_k}{A_*}\>^\sss}+e^{-98A'^{0.34}}}Z(t'_{\upmu-1}).\]
		Consequently, if \eqref{rpe} or \eqref{rpe1} holds at $\upmu$, then
		\begin{equation}\label{key}
			\begin{aligned}
				\nnr{\gamma(t'_{\upmu-1},\upsigma,\fv)}_{L^2_\fv}&\leq\nnr{\gamma(t'_{\upmu},\upsigma,\fv)}_{L^2_\fv}+\nnr{\fD_{\upsigma,\upmu}\rr{\gamma(t'_{\upmu},\upmu,\fv)+\sF_\upzeta^{-1}\rr{\int_{\mathsf{T}_k}^{\f\upzeta\upmu} \wh{Q}(\tau,\upmu\bfe,\upzeta\bfe)d\tau}}}_{L^2_\fv} \\
				&\quad +\rr{e^{-0.7\<\f{\mathsf{T}_k}{A_*}\>^\sss}+e^{-97A'^{0.34}}}Z(t'_{\upmu-1}).
			\end{aligned}
		\end{equation}
		Moreover we have 
		\[\ba&\nnr{\fD_{\upsigma,\upmu}\rr{\gamma(t'_{\upmu},\upmu,\fv)+\sF_\upzeta^{-1}\rr{\int_{\mathsf{T}_k}^{\f\upzeta\upmu} \wh{Q}(\tau,\upmu\bfe,\upzeta\bfe)d\tau}}}_{L^2_\fv}\leq e^{-99A'^{0.34}}Z(t'_\upmu)\\
		&+\nnr{\gamma(t'_{\upmu},\upmu,\fv)+\sF_\upzeta^{-1}\rr{\int_{\mathsf{T}_k}^{\f\upzeta\upmu} \wh{Q}(\tau,\upmu\bfe,\upzeta\bfe)d\tau}}_{L^2_\fv}\Etmu\left\{\ba&\f{\eps_0 A'}{\upmu^3}\rr{1+20\eps_*},\quad&\upsigma=\upmu\pm1,\\
		&\f{50\eps_0 A'}{\upmu^3}\eps_*\<\upsigma-\upmu\>^{-8},\quad&\upsigma\neq\upmu\pm1,\upmu,\\
		&\f{10^2\theta^{-2}}{\upmu^2},\quad&\upsigma=\upmu.\ea
		\right.
		\ea\]
		This operator bound is proved exactly as \eqref{Dmnup}: the main wave produces the two nearest-neighbor factors $\eps_0A'/\upmu^3$, the localized correction produces the algebraic tail, and the resolvent contributes the diagonal $\upmu^{-2}$ loss.
		
		\noindent\textit{Step 2. Iteration before the source mode is crossed.}
		Set $\upmu^\dagger=\min\{0\leq\upmu\leq\min\{k'_0,k'_2\}-1:k'\notin[\upmu+1,\min\{k'_0,k'_2\}-1]\}$. We prove \eqref{rpe} for $\upmu^\dagger\leq\upmu\leq\min\{k'_0,k'_2\}-1$. Propositions \ref{4.3} and \ref{4.4} provide the terminal case $\upmu=\min\{k'_0,k'_2\}-1$. Assume the estimate holds at all indices from $\upmu$ up to this terminal index. We prove it at $\upmu-1$ by separating the four relative positions of $\upsigma$. First, for $\upsigma=\upmu-1$,
		\[\ba
		&\nnr{\gamma(t'_{\upmu-1},\upmu-1,\fv)}_{L^2_\fv}\leq\nnr{\gamma(t'_{\upmu},\upmu-1,\fv)}_{L^2_\fv}+\nnr{\fD_{\upmu-1,\upmu}\gamma(t'_{\upmu},\upmu,\fv)}_{L^2_\fv}+\rr{e^{-0.7\<\f{\mathsf{T}_k}{A_*}\>^\sss}+e^{-97A'^{0.34}}}Z(t'_{\upmu-1})\\
		&\leq A'\mathtt{G}_\theta(\upmu+1)Z(t'_\upmu)\rr{{A'}^{-1}+50\eps_{*}+\f{\eps_0 A'}{\upmu^3}\rr{1+20\eps_*}\Etmu}+\rr{e^{-0.7\<\f{\mathsf{T}_k}{A_*}\>^\sss}+e^{-96A'^{0.34}}}Z(t'_{\upmu-1})\\
		&\leq A'\mathtt{G}_\theta(\upmu+1)Z(t'_\upmu)\max\{\f{\eps_0A'}{\upmu^3},1\}(1+\eps_0^{20})\Etmu+\rr{e^{-0.7\<\f{\mathsf{T}_k}{A_*}\>^\sss}+e^{-96A'^{0.34}}}Z(t'_{\upmu-1})\leq A'\mathtt{G}_\theta(\upmu)Z(t'_{\upmu-1}).
		\ea\]
		Next we consider the case $\upsigma\leq\upmu-2$,
		\[\ba
		&\nnr{\gamma(t'_{\upmu-1},\upsigma,\fv)}_{L^2_\fv}\leq\nnr{\gamma(t'_{\upmu},\upsigma,\fv)}_{L^2_\fv}+\nnr{\fD_{\upsigma,\upmu}\gamma(t'_{\upmu},\upmu,\fv)}_{L^2_\fv}+\rr{e^{-0.7\<\f{\mathsf{T}_k}{A_*}\>^\sss}+e^{-96A'^{0.34}}}Z(t'_{\upmu-1})\\
		&\leq A'\mathtt{G}_\theta(\upmu+1)Z(t'_\upmu)\rr{{A'}^{-1}+50\eps_{*}\<\upsigma-\upmu\>^{-7}+\f{50\eps_0 A'}{\upmu^3}\eps_*\<\upsigma-\upmu\>^{-8}\Etmu}+\rr{e^{-0.7\<\f{\mathsf{T}_k}{A_*}\>^\sss}+e^{-96A'^{0.34}}}Z(t'_{\upmu-1})\\
		&\leq A'\mathtt{G}_\theta(\upmu)Z(t'_{\upmu-1})\left({A'}^{-1}+50\eps_{*}\<\upsigma-\upmu\>^{-7}
		+50\eps_*\<\upsigma-\upmu\>^{-8}\right)\leq A'\mathtt{G}_\theta(\upmu)Z(t'_{\upmu-1})\rr{{A'}^{-1}+50\eps_{*}\<\upsigma-\upmu+1\>^{-7}}.	
		\ea\]
		Then we consider the case $\upsigma=\upmu$.
		\[\nnr{\gamma(t'_{\upmu-1},\upmu,\fv)}_{L^2_\fv}\leq A'\mathtt{G}_\theta(\upmu+1)Z(t'_\upmu)\rr{1+\f{10^2\theta^{-2}}{\upmu^2}\Etmu}
		\leq2\times10^4A'\mathtt{G}_\theta(\upmu)Z(t'_{\upmu-1}).\]	
		Then we consider the case $\upsigma\geq\upmu+1$.
		\[\ba
		&\nnr{\gamma(t'_{\upmu-1},\upsigma,\fv)}_{L^2_\fv}\leq 10^4A'\mathtt{G}_\theta(\upmu+1)Z(t'_\upmu)\rr{3-\<\upsigma-\upmu-1\>^{-7}+0.4\f{\eps_0A'}{\upmu^3}\<\upsigma-\upmu\>^{-8}\Etmu}+\rr{e^{-0.7\<\f{\mathsf{T}_k}{A_*}\>^\sss}\right.\\&\left.+e^{-96A'^{0.34}}}Z(t'_{\upmu-1})\leq 10^4A'\mathtt{G}_\theta(\upmu)Z(t'_{\upmu-1})\rr{3-\<\upsigma-\upmu-1\>^{-7}+\<\upsigma-\upmu\>^{-8}}\leq10^4 A'\mathtt{G}_\theta(\upmu)Z(t'_{\upmu-1})\rr{3-\<\upsigma-\upmu\>^{-7}}.	
		\ea\]
		
		\noindent\textit{Step 3. Crossing the source mode.}
		It remains to prove \eqref{rpe1} for $0\leq\upmu\leq\upmu^\dagger-1$. If this range is nonempty, then $\upmu^\dagger=k'$. The estimate \eqref{rpe} is already known at $\upmu=k'$, and on the crossing interval $[t'_{k'},t'_{k'-1}]$ the direct source contribution must be included in the transfer bound:
		\[\ba&\nnr{\fD_{\upsigma,k'}\rr{\gamma(t'_{k'},k',\fv)+\sF_\zeta^{-1}\rr{\int_{\mathsf{T}_k}^{\f{\zeta}{k'}} \wh{Q}(\tau,k'\bfe,\zeta\bfe)d\tau}}}_{L^2_\fv}\leq e^{-99A'^{0.34}}Z(t'_{k'})\\
		&+\rr{\nnr{\gamma(t'_{k'},{k'},\fv)}_{L^2_\fv}+Z(t'_{k'-1})}\mathsf{E}_\theta({k'})\left\{\ba&\f{\eps_0 A'}{{k'}^3}\rr{1+20\eps_*},\quad&\upsigma={k'}\pm1,\\
		&\f{50\eps_0 A'}{{k'}^3}\eps_*\<\upsigma-{k'}\>^{-8},\quad&\upsigma\neq{k'}\pm1,k',\\
		&\f{10^2\theta^{-2}}{{k'}^2},\quad&\upsigma={k'}.\ea
		\right.
		\ea\]
		We first consider the case $\upsigma=k'-1$.
		\[\ba
		&\nnr{\gamma(t'_{k'-1},k'-1,\fv)}_{L^2_\fv}\leq A'\mathtt{G}_\theta(k'+1)Z(t'_{k'})\rr{{A'}^{-1}+50\eps_{*}}+\rr{A'\mathtt{G}_\theta(k'+1)Z(t'_{k'})+Z(t'_{k'-1})}\f{\eps_0 A'}{{k'}^3}\rr{1+20\eps_*}\mathsf{E}_\theta(k')\\
		&+\rr{e^{-0.7\<\f{\mathsf{T}_k}{A_*}\>^\sss}+e^{-96A'^{0.34}}}Z(t'_{k'-1})
		\leq A'\mathtt{G}_\theta(k'+1) Z(t'_{k'})\rr{{A'}^{-1}+50\eps_{*}+\f{\eps_0 A'}{{k'}^3}\rr{1+20\eps_*}\mathsf{E}_\theta(k')}+A'Z(t'_{k'-1})\\
		&\leq 2A'\mathtt{G}_\theta(k')Z(t'_{k'-1}).
		\ea\]
		Next we consider the case $\upsigma\leq k'-2$.
		\[\ba
		&\nnr{\gamma(t'_{k'-1},\upsigma,\fv)}_{L^2_\fv}\leq A'\mathtt{G}_\theta(k'+1)Z(t'_{k'})\rr{{A'}^{-1}+50\eps_{*}\<k'-\upsigma\>^{-7}}+A'\mathtt{G}_\theta(k'+1)Z(t'_{k'})\f{50\eps_0 A'}{{k'}^3}\eps_*\<\upsigma-k'\>^{-8}\mathsf{E}_\theta(k')\\
		&+A'Z(t'_{k'-1})\eps_*\<\upsigma-k'\>^{-8}+\rr{e^{-0.7\<\f{\mathsf{T}_k}{A_*}\>^\sss}+e^{-96A'^{0.34}}}Z(t'_{k'-1})
		\leq A'\rr{{A'}^{-1}+50\eps_{*}\<k'-\upsigma\>^{-7}+70\eps_{*}\<k'-\upsigma\>^{-8}}\\
		&\times\mathtt{G}_\theta(k') Z(t'_{k'-1})\leq A'\mathtt{G}_\theta(k')Z(t'_{k'-1})\rr{{A'}^{-1}+10^2\eps_{*}\<k'-\upsigma-1\>^{-7}}.
		\ea\]
		Then we consider the case $\upsigma=k'$.
		\[\ba&\nnr{\gamma(t'_{k'-1},k',\fv)}_{L^2_\fv}\leq A'\mathtt{G}_\theta(k'+1)Z(t'_{k'})\rr{1+\f{10^2\theta^{-2}}{{k'}^2}\mathsf{E}_\theta(k')}+Z(t'_{k'-1})\f{200\theta^{-2}}{{k'}^2}\mathsf{E}_\theta(k')\leq10^5A'\mathtt{G}_\theta(k')Z(t'_{k'-1}).\ea\]
		Then we consider the case $\upsigma\geq k'+1$.
		\[\ba
		&\nnr{\gamma(t'_{k'-1},\upsigma,\fv)}_{L^2_\fv}\leq 10^5A'\mathtt{G}_\theta(k'+1)Z(t'_{k'})\rr{3-\<\upsigma-k'-1\>^{-7}+0.4\f{\eps_0 A'}{{k'}^3}\<\upsigma-k'\>^{-8}\mathsf{E}_\theta(k')}\\
		&+A'Z(t'_{k'-1})\<\upsigma-k'\>^{-8}+\rr{e^{-0.7\<\f{\mathsf{T}_k}{A_*}\>^\sss}+e^{-96A'^{0.34}}}Z(t'_{k'-1})
		\leq 10^5A'\rr{3-\<\upsigma-k'-1\>^{-7}+\<k'-\upsigma\>^{-8}}\\
		&\times\mathtt{G}_\theta(k') Z(t'_{k'-1})\leq 10^5A'\mathtt{G}_\theta(k')Z(t'_{k'-1})\rr{3-\<\upsigma-k'\>^{-7}}.
		\ea\]
		Thus \eqref{rpe1} holds at $\upmu=k'-1$. Below this index the source has already been crossed, so the induction is identical to Step 2 with the enlarged constants in \eqref{rpe1}. This proves the estimate for every $0\leq\upmu\leq k'-1$.
	\end{proof}

	The recurrence bounds control the reduced amplitude at the endpoints of the critical intervals. The next proposition combines them with the within-interval estimate \eqref{fjeif}, thereby returning to the full profile throughout the complete forward echo chain.
    
	\begin{prop}\label{79}
		Assume $A'>\mathsf{T}_k/2$. For any $t\in[\max\{t'_{k'_0-1},t'_{k'_2-1}\},t'_0]$,
		\begin{equation}\label{yvd}
			\nnr{e^{\<\na_{\mrv}\>^\sss}v^{\alpha}\rr{\sF_zf(t,\upsigma\bfe,v)-\int_{\mathsf{T}_k}^t\sF_zQ(\tau,\upsigma\bfe,v)d\tau}}_{\lta L^2_v}\leq \f12A'^2Z(t).
		\end{equation}
	\end{prop}
    \begin{proof}
	    	Set $\mathbf{T}=\max\{t'_{k'_0-1},t'_{k'_2-1}\}$. Propositions \ref{4.2} and \ref{4.3} provide \eqref{yvd} at $t=\mathbf{T}$. Assume it has been propagated through the critical intervals up to $t'_\upmu$. On $[t'_\upmu,t'_{\upmu-1}]$, equations \eqref{rec41} and \eqref{fjeif} decompose the profile into the entrance value, all completed transfers, the current partial transfer, and the remainder:
    	\[\ba
    	&\nnr{e^{\<\na_{\mrv}\>^\sss}v^{\alpha}\rr{\sF_zf(t,\upsigma\bfe,v)-\int_{\mathsf{T}_k}^t\sF_zQ(\tau,\upsigma\bfe,v)d\tau}}_{\ltaa L^2_v}\leq \nnr{e^{\<\na_{\mrv}\>^\sss}v^{\alpha}\rr{\sF_zf(\mathbf{T},\upsigma\bfe,v)-\int_{\mathsf{T}_k}^{\mathbf{T}}\sF_zQ(\tau,\upsigma\bfe,v)d\tau}}_{\ltaa L^2_v}\\
    	&+\sum_{\upnu=\upmu+1}^{\min\{k'_0-1,k'_2-1\}}\nnr{e^{\<\na_{\mrv}\>^\sss}v^{\alpha}\cD_{\upsigma,\upnu}\rr{\gamma(t'_{\upnu},\upnu,\fv)+\sF_\upzeta^{-1}\rr{\int_{\mathsf{T}_k}^{\f\upzeta\upnu} \wh{Q}(\tau,\upnu\bfe,\upzeta\bfe)d\tau}}}_{\ltaa L^2_v}
    	+\nnr{e^{\<\na_{\mrv}\>^\sss}v^{\alpha}\rr{\sF_zf(t,\upsigma\bfe,v)\right.\right.\\
    	&\left.\left.-\int_{t'_\upmu}^t\sF_zQ(\tau,\upsigma\bfe,v)d\tau-\sF_zf(t'_\upmu,\upsigma\bfe,v)}}_{\ltaa L^2_v}
    	+e^{-95A'^{0.34}}Z(t)\leq A'Z(t) 	+\sum_{\upnu=\upmu+1}^{\min\{k'_0-1,k'_2-1\}}\f{\eps_0^{0.5}A'}{\upnu^3}\<\upsigma-\upnu\>^{-8}\\
    	&\times\rr{\nnr{\gamma(t'_{\upnu},\upnu,\fv)}_{L^2_\fv}+Z(t)}
    	+\eps_0^{0.5}A'\rr{\nnr{\gamma(t'_{\upmu},\upmu,\fv)}_{L^2_\fv}+Z(t)}
    	+e^{-94A'^{0.34}}Z(t)\leq \f12A'^2Z(t).
    	\ea\]
	    	The last inequality uses \eqref{rpe}--\eqref{rpe1} and the summability of the transfer tails. This completes the induction and proves \eqref{yvd} on $[\mathbf{T},t'_0]$.
	    \end{proof} 
	
	Finally, after $\mathsf{T}=\max\{t'_0,\mathsf{T}_k\}$ every critical time lies in the past. The density created by a packet near $(k',A')$ is then exponentially small, so the subsequent evolution differs from the direct source response only by a perturbative tail.
	
	\begin{prop}\label{4.6}
		We denote $\mathsf{T}=\max\{t'_0,\mathsf{T}_k\}$. It holds that
		\[\nnr{e^{\<\na_{\mrv}\>^\sss}v^{\alpha}\rr{\sF_zf(t,\upsigma\bfe,v)-\int_{\mathsf{T}}^t\sF_zQ(\tau,\upsigma\bfe,v)d\tau-\sF_zf(\mathsf{T},\upsigma\bfe,v)}}_{\ltaa L^2_v}\leq\rr{e^{-\f12\<\f{\mathsf{T}_k}{A_*}\>^\sss}+e^{-98A'^{0.34}}}Z(t)\]
		on $[\max\{t'_0,\mathsf{T}_k\},\infty)$.
	\end{prop}
	\begin{proof}
		Because $t\geq\mathsf{T}\geq t'_0$, the phase relation $\upnu t\approx A'$ cannot hold for every nonzero $\upnu$. The density is therefore a remainder, and the profile satisfies
		\[\ba
		&\wh{\rho}(t,\upnu\bfe)=\wh{\sR}(t,\upnu\bfe)\quad\text{for any}\ \upnu\in\Z,\\
		&\wh{f}(t,\upsigma\bfe,\xi)-\int_{\mathsf{T}}^t\wh{Q}(\tau,\upsigma\bfe,\xi)d\tau=\wh{f}(\mathsf{T},\upsigma\bfe,\xi)-\sum_{\upnu\neq0}\int_{\mathsf{T}}^{t}\wh{\udl{\rho}}(\tau,\upnu\bfe)\f{\upxi-\upsigma \tau}{\upnu}\wh{f}(\tau,(\upsigma-\upnu)\bfe,\xi-\upnu\tau\bfe) d\tau+\wh{\sR}(t,\upsigma\bfe,\xi),
		\ea\]
		with the error terms satisfying
		\[\nr{\wh{\sR}(t,\upnu\bfe)}\leq e^{-100A'^{0.34}}e^{-0.124\<\f{\upnu-k',\upnu t-A'}{A_*}\>^\sss}Z(t),\quad\nnr{e^{0.124\<\f{\na_\fv-A'}{A_*}\>^\sss+\<\na_{\mrv}\>^{\sss}}\rr{v^{\alpha}\sR(t,\upsigma\bfe,v)}}_{\ltaa L^2_v}\leq e^{-100A'^{0.34}}Z(t).\]
		Substitution into the profile equation produces a small multiple of the profile norm. The continuity argument from Proposition \ref{prop:far-frequency-response} then proves the result.
	\end{proof}
	
	\subsection{Backward response on $[\mathfrak{T}_k,\mathsf{T}_k]$}
	We now evolve from the zero terminal value at $\mathsf{T}_k$ toward earlier times. The structure mirrors the forward argument, but the critical indices are traversed in the opposite order. We first control any nonresonant segment adjacent to $\mathsf{T}_k$, then derive and iterate a backward recurrence on complete critical intervals, and finally estimate the segment adjacent to $\mathfrak{T}_k$.

	Suppose first that $A'\leq\mathsf{T}_k/2$. Then $t'_0\leq\mathsf{T}_k$, and the interval between $\max\{\mathfrak{T}_k,t'_0\}$ and $\mathsf{T}_k$ lies after the last critical time.
	\begin{prop}\label{4.7}
		Assume $A'\leq\mathsf{T}_k/2$. On $[\max\{\mathfrak{T}_k,t'_0\},\mathsf{T}_k]$,
		\[\ba\nnr{e^{\<\na_{\mrv}\>^\sss}v^{\alpha}\rr{\sF_zf(t,\upsigma\bfe,v)+\int_{t}^{\mathsf{T}_k}\sF_zQ(\tau,\upsigma\bfe,v)d\tau}}_{\lta L^2_v}\leq\rr{e^{-\f12\<\f{\mathfrak{T}_k}{A_*}\>^\sss}+e^{-98A'^{0.34}}}Z(t)\ea\]
		for any $\upsigma\in\Z$.
	\end{prop}
	\begin{proof}
		On $[\max\{\mathfrak{T}_k,t'_0\},\mathsf{T}_k]$, every density mode is separated from its resonant line. Reversing the time orientation in the proof of Proposition \ref{4.6} gives
		\[\ba
		&\wh{\rho}(t,\upnu\bfe)=\wh{\sR}(t,\upnu\bfe)\quad\text{for any}\ \upnu\in\Z,\\
		&\wh{f}(t,\upsigma\bfe,\xi)+\int_{t}^{\mathsf{T}_k}\wh{Q}(\tau,\upsigma\bfe,\xi)d\tau=\sum_{\upnu\neq0}\int_{t}^{\mathsf{T}_k}\wh{\udl{\rho}}(\tau,\upnu\bfe)\f{\upxi-\upsigma \tau}{\upnu}\wh{f}(\tau,(\upsigma-\upnu)\bfe,\xi-\upnu\tau\bfe) d\tau+\wh{\sR}(t,\upsigma\bfe,\xi),
		\ea\]
		with the error terms satisfying
		\[\nr{\wh{\sR}(t,\upnu\bfe)}\leq e^{-100A'^{0.34}}e^{-0.124\<\f{\upnu-k',\upnu t-A'}{A_*}\>^\sss}Z(t),\quad\nnr{e^{0.124\<\f{\na_\fv-A'}{A_*}\>^\sss+\<\na_{\mrv}\>^{\sss}}\rr{v^{\alpha}\sR(t,\upsigma\bfe,v)}}_{\ltaa L^2_v}\leq e^{-100A'^{0.34}}Z(t).\]
		The remainder bounds follow from the same two frequency regions as in Proposition \ref{prop:far-frequency-response}. Since $Z$ is decreasing on the backward interval, the source accumulated between $t$ and $\mathsf{T}_k$ is controlled by $Z(t)$. The resulting small profile term is absorbed by continuity.
	\end{proof}
	
	If $A'>\mathsf{T}_k/2$, $\mathsf{T}_k$ lies inside the critical-time range. When $k'_0<k'_2$, it lies in the low-mode regime, and the partial interval $[t'_{k'_0},\mathsf{T}_k]$ contains only the mode $k'_0$.
	
	\begin{prop}\label{4.8}
		Assume $A'>\mathsf{T}_k/2$ and $k'_0<k'_2$. On $[t'_{k'_0},\mathsf{T}_k]$,
		\[\ba
		&\nnr{e^{\<\na_{\mrv}\>^\sss}v^{\alpha}\rr{\sF_zf(t,\upsigma\bfe,v)+\int_{t}^{\mathsf{T}_k}\sF_zQ(\tau,\upsigma\bfe,v)d\tau}}_{\lta L^2_v}\\
		&\leq\f12 A'Z(t)\rr{\vv{1}_{\upsigma=k'_0\pm1}+\eps_{*}\<\upsigma-k'_0\>^{-8}}+\rr{e^{-98A'^{0.34}}+e^{-0.8\<\f{\mathsf{T}_k}{A_*}\>^\sss} }Z(t)
		\ea\]
		for any $\upsigma\in\Z$.
	\end{prop}
	\begin{proof}
		Reverse the time orientation in Proposition \ref{4.3}. After subtracting $-\int_t^{\mathsf{T}_k}Q(\tau)d\tau$, all modes other than $k'_0$ are off resonant on $[t'_{k'_0},\mathsf{T}_k]$. The remaining scalar density equation is inverted with the negative-frequency resolvent from \eqref{eR}; its transfer through $\eps_0\vf\cos\fz$ produces the modes $k'_0\pm1$, while $f_{*,\infty}$ gives the algebraically decaying tail. The remainder estimates are unchanged, and the monotonicity of $Z$ in the backward direction gives the stated endpoint bound.
	\end{proof}
	
	If $A'\leq\mathfrak{T}_k/2$, Proposition \ref{4.7} already covers the entire backward interval. Hence assume $A'>\mathfrak{T}_k/2$ and recall the index $k'_1$ defined by $t'_{k'_1}\leq\mathfrak{T}_k<t'_{k'_1-1}$. Set
	\[
		m_{\mathrm b}:=k'_0+1,
		\qquad
		T_{\mathrm b}:=t'_{k'_0}.
	\]
	Thus $T_{\mathrm b}$ is the endpoint reached after the segment adjacent to $\mathsf{T}_k$, and $m_{\mathrm b}$ is the first complete critical interval encountered when moving backward from that endpoint. If $k'_0\geq k'_2$, the critical ranges below are empty and Proposition \ref{4.11} connects directly to $\mathsf{T}_k$. Otherwise, on the complete intervals between $T_{\mathrm b}$ and $\max\{t'_{k'_1-1},t'_{k'_2-1}\}$, define
	\[\gamma(t,\upmu,\fv)=\int_{\R^2}\sF_z f(t,\upmu\bfe,v)d\mrv+\int_{t}^{\mathsf{T}_k}\int_{\R^2}\sF_z Q(\tau,\upmu\bfe,v)d\mrv d\tau. \]
	This is the backward analogue of the renormalized amplitude used in the forward iteration.
	
	\begin{prop}\label{4.9}
		Assume $A'>\mathfrak{T}_k/2$. For any $m_{\mathrm b}\leq\upmu\leq\min\{k'_1,k'_2\}-1$, we have the following approximate recurrence relation on $[t'_\upmu,t'_{\upmu-1}]$:
		\[\ba&\sF_zf(t'_{\upmu},\upsigma\bfe,v) =\sF_zf(t'_{\upmu-1},\upsigma\bfe,v)-\int_{t'_\upmu}^{t'_{\upmu-1}}\sF_z Q(\tau,\upsigma\bfe,v)d\tau\\
		&+\cE_{\upsigma,\upmu}\rr{\gamma(t'_{\upmu-1},\upmu,\fv)-\sF_\zeta^{-1}\rr{\int_{\f{\zeta}{\upmu}}^{\mathsf{T}_k} \wh{Q}(\tau,\upmu\bfe,\zeta\bfe)d\tau}}+\sR_{\upsigma,\upmu}(v)\ea\]
		with the error term satisfying
		\[\ba
		&\nnr{e^{\<\na_{\mrv}\>^\sss}v^{\alpha}\sR_{\upsigma,\upmu}(v)}_{\ltaa L^2_v}\leq e^{-0.8\<\f{\mathsf{T}_k}{A_*}\>^\sss}\sup_{\upsigma\in\Z}\nnr{e^{\<\na_{\mrv}\>^\sss}v^\alpha\rr{\sF_z f(t'_{\upmu-1},\upsigma\bfe,v)+\int^{\mathsf{T}_k}_{t'_{\upmu-1}}\sF_z Q(\tau,\upsigma\bfe,v)d\tau}}_{\ltaa L^2_v}\\
		&+\rr{e^{-0.8\<\f{\mathsf{T}_k}{A_*}\>^\sss}+e^{-98A'^{0.34}}}Z(t'_{\upmu})\ea\]
		and
		\[\nnr{\int_{\R^2}\sR_{\upsigma,\upmu}(v)d\mrv}_{ L^2_\fv}\leq e^{-0.8\<\f{\mathsf{T}_k}{A_*}\>^\sss}\sup_{\upsigma\in\Z}\nnr{\gamma(t'_{\upmu-1},\upsigma,\fv)}_{ L^2_\fv}+\rr{e^{-0.8\<\f{\mathsf{T}_k}{A_*}\>^\sss}+e^{-98A'^{0.34}}}Z(t'_{\upmu}).\]
		Moreover the following estimate holds on $[t'_\upmu,t'_{\upmu-1}]$:
		\[\ba&\nnr{e^{\<\na_{\mrv}\>^\sss}v^{\alpha}\rr{\sF_zf(t,\upsigma\bfe,v)-\sF_zf(t'_{\upmu-1},\upsigma\bfe,v)+\int^{t'_{\upmu-1}}_{t} \sF_z Q(\tau,\upsigma\bfe,v)d\tau}}_{\ltaa L^2_v}\leq \eps_0^{0.8}A'\rr{\nnr{\gamma(t'_{\upmu-1},\upmu,\fv)}_{L^2_\fv}+Z(t)} \\
		&+ e^{-0.8\<\f{\mathsf{T}_k}{A_*}\>^\sss}\sup_{\upsigma\in\Z}\nnr{e^{\<\na_{\mrv}\>^\sss}v^\alpha\rr{\sF_z f(t'_{\upmu-1},\upsigma\bfe,v)+\int^{\mathsf{T}_k}_{t'_{\upmu-1}}\sF_z Q(\tau,\upsigma\bfe,v)d\tau}}_{\ltaa L^2_v}.\ea\]
	\end{prop}
	\begin{proof}
		The proof follows the two-stage argument of Proposition \ref{4.4} with time reversed. On $[t'_\upmu,t'_{\upmu-1}]$, isolate the density mode $\upmu$ and place every other density mode into the off-resonant remainder. Invert the scalar Volterra equation from the right endpoint $t'_{\upmu-1}$; after the change of variables $\upzeta=\upmu t$, the leading double integral is the backward transfer operator $\cE_{\upsigma,\upmu}$. The renormalized amplitude now contains the source accumulated from $\upzeta/\upmu$ to $\mathsf{T}_k$, which explains the sign and limits in the statement. The weighted remainder and within-interval estimates follow from the same frequency separation as in Proposition \ref{4.4}, using that $Z$ increases when time is moved backward.
	\end{proof}
	
	The recurrence must now be iterated toward increasing mode indices. Accordingly, the backward growth condition compares $Z(t'_{\upmu-1})$ with $Z(t'_\upmu)$, the reverse of \eqref{growth1}. As in the forward argument, the statement distinguishes whether the iteration has crossed the source mode $k'$.
	
	\begin{prop}\label{4.10}
		Assume $A'>\mathfrak{T}_k/2$ and that, the envelope satisfies that for fixed $A'>0$ and every $m_{\mathrm b}\leq\upmu\leq\min\{k'_1,k'_2\}-1$, 
		\begin{equation}\label{growth2}
			\max\{1,\f{\eps_0A'}{\upmu^3}\}\rr{1+\eps_0^{10}}Z(t'_{\upmu-1})\leq Z(t'_{\upmu})\rr{1+\f{100}{\upmu^2}}.
		\end{equation}
		For any $m_{\mathrm b}-1\leq\upmu\leq\min\{k'_1,k'_2\}-1$, if $k'\notin[m_{\mathrm b},\upmu]$, then
		\[\nnr{\gamma(t'_\upmu,\upsigma\bfe,\fv)}_{L^2_\fv}\leq A'Z(t'_\upmu)\rr{
			\prod_{\upnu=m_{\mathrm b}}^{\upmu}\Etnu\rr{1+\f{10^2}{\upnu^2}}}\left\{\ba&{A'}^{-1}+50\eps_{*}\<\upsigma-\upmu-1\>^{-7},\quad&\upsigma\geq\upmu+2,\\&1,\quad&\upsigma=\upmu+1,\\&10^4\rr{3-\<\upmu-\upsigma\>^{-7}},\quad&\upsigma\leq\upmu.\ea\right.
		\]
		If $k'\in[m_{\mathrm b},\upmu]$, then
		\[\nnr{\gamma(t'_\upmu,\upsigma\bfe,\fv)}_{L^2_\fv}\leq A'Z(t'_\upmu)\rr{
			\prod_{\upnu=m_{\mathrm b}}^{\upmu}\Etnu\rr{1+\f{10^2}{\upnu^2}}}\left\{\ba&{A'}^{-1}+10^2\eps_{*}\<\upsigma-\upmu-1\>^{-7},\quad&\upsigma\geq\upmu+2,\\&2,\quad&\upsigma=\upmu+1,\\&10^5\rr{3-\<\upmu-\upsigma\>^{-7}},\quad&\upsigma\leq\upmu.\ea\right.
		\]
		Here and below, a product over an empty index set is understood to be $1$.
	\end{prop}
	\begin{proof}
		Start from the bound supplied by Proposition \ref{4.7} or \ref{4.8}, according to the value of $T_{\mathrm b}$. Apply the recurrence in Proposition \ref{4.9} successively for $\upmu=m_{\mathrm b},m_{\mathrm b}+1,\ldots$. The nearest-neighbor term is bounded by $\eps_0A'/\upmu^3$, the diagonal term by $C_\theta/\upmu^2$, and the remaining modes by the summable tail $\eps_*\<\upsigma-\upmu\>^{-8}$. Condition \eqref{growth2} absorbs the first factor, while the product in the statement records the diagonal losses. Before the index $k'$ is reached, the first set of constants closes the induction. At $\upmu=k'$, the direct source term enters once and changes these constants to the second set. The subsequent induction is otherwise unchanged. This is the backward counterpart of the three steps in Proposition \ref{4.5}.
	\end{proof}
	
	We now transfer the endpoint estimates for $\gamma$ back to the full profile on every critical interval.
	\begin{prop}\label{89}
		Assume $A'>\mathfrak{T}_k/2$. For any $t\in[\max\{t'_{k'_1-1},t'_{k'_2-1}\},T_{\mathrm b}]$,
		\[\nnr{e^{\<\na_{\mrv}\>^\sss}v^{\alpha}\rr{\sF_z f(t,\upsigma\bfe,v)+\int_{t}^{\mathsf{T}_k}\sF_z Q(\tau,\upsigma\bfe,v)d\tau}}_{\lta L^2_v}\leq\f12 A'^2Z(t).\]
	\end{prop}
    \begin{proof}
	    	Begin at $T_{\mathrm b}$ with Proposition \ref{4.7} or \ref{4.8}. On each interval $[t'_\upmu,t'_{\upmu-1}]$, Proposition \ref{4.9} decomposes the profile into the endpoint value, the completed backward transfers, the current partial transfer, and an exponentially small remainder. Proposition \ref{4.10} bounds every resonant amplitude appearing in this decomposition. Summing the algebraically decaying transfer tails and using the monotonicity of $Z$ gives the factor $A'^2Z(t)$ uniformly until $\max\{t'_{k'_1-1},t'_{k'_2-1}\}$ is reached.
    \end{proof}
	
	It remains to connect the iterated critical region to the initial time $\mathfrak{T}_k$. As in the forward high-mode estimate, modes $\upnu\geq k'_2$ can first be controlled simultaneously because their transfer factor is perturbative. We write $\mathtt{T}=\min\{t'_{k'_2-1},\mathsf{T}_k\}$ for the right endpoint of this segment.
	
	\begin{prop}\label{4.11}
		Whenever $\mathfrak{T}_k\leq\min\{t'_{k'_2-1},\mathsf{T}_k\}$, on $[\mathfrak{T}_k,\min\{t'_{k'_2-1},\mathsf{T}_k\}]$,
		\[\nnr{e^{\<\na_{\mrv}\>^\sss}v^{\alpha}\rr{\sF_z f(t,\upsigma\bfe,v)+\int_{t}^{\min\{t'_{k'_2-1},\mathsf{T}_k\}}\sF_z Q(\tau,\upsigma\bfe,v)d\tau-\sF_z f(\min\{t'_{k'_2-1},\mathsf{T}_k\},\upsigma\bfe,v)}}_{\lta L^2_v}\leq {A'}^2Z(t)\]
		for any $\upsigma\in\Z$.
	\end{prop}
	\begin{proof}
		Reverse the argument of Proposition \ref{4.2}. On $[\mathfrak{T}_k,\mathtt{T}]$, retain the density modes $\upnu\geq k'_2$ and put the modes below $k'_2$ into the off-resonant remainder. This gives
		\[\ba
		&\wh{\rho}(t,\upnu\bfe)=\wh{\sR}(t,\upnu\bfe)\quad\text{for any}\ \upnu<k'_2,\\
		&\wh{\rho}(t,\upmu\bfe)=\wh{f}(\mathtt{T},\upmu\bfe,\upmu t\bfe)-\int^{\mathtt{T}}_t\wh{Q}(\tau,\upmu\bfe,\upmu t\bfe)d\tau+\int^{\mathtt{T}}_{t}\wh{\rho}(\tau,\upmu\bfe)(t-\tau)\sF\rr{f^0+\udl{f_\infty}}(0,\upmu(t-\tau)\bfe) d\tau\\
		&+\sum_{\substack{\upnu\geq k'_2\\\upnu\neq \upmu}}\int^{\mathtt{T}}_{t}\wh{\rho}(\tau,\upnu\bfe)\f{\upmu(t-\tau)}{\upnu}\wh{\udl{f_\infty}}((\upmu-\upnu)\bfe,(\upmu t-\upnu\tau)\bfe) d\tau
		+\wh{\sR}(t,\upmu\bfe )\quad\text{for any}\ \upmu\geq k'_2,\\
		&\wh{f}(t,\upsigma\bfe,\xi)=\wh{f}(\mathtt{T},\upsigma\bfe,\xi)-\int^{\mathtt{T}}_t\wh{Q}(\tau,\upsigma\bfe,\xi)d\tau\\
		&+\sum_{\upnu\geq k'_2}\int^{\mathtt{T}}_{t}\wh{\rho}(\tau,\upnu\bfe)\f{\upxi-\upsigma \tau}{\upnu}\sF\rr{f^0+\udl{f_\infty}}((\upsigma-\upnu)\bfe,\xi-\upnu\tau\bfe) d\tau
		+\wh{\sR}(t,\upsigma\bfe,\xi).
		\ea\]
		The error terms satisfy the estimates of Proposition \ref{4.2} with the time integrals reversed. Volterra inversion controls the active densities, and substitution into the profile equation gives
		\[\ba&\nnr{e^{\<\mrxi\>^\sss}\pa_\xi^\alpha\rr{\wh{f}(t,\upsigma\bfe,\xi)+\int_t^{\mathtt{T}}\wh{Q}(\tau,\upsigma\bfe,\xi)d\tau-\wh{f}(\mathtt{T},\upsigma\bfe,\xi)}}_{\ltaa L^2_\xi}\leq \eps_0^{0.97}\sup_{\upnu\geq k'_2}\rr{\nnr{\wh{\rho}(\tau,\upnu\bfe)}_{L^2_\tau[t,\mathtt{T}]}\upnu^\f12}\\
		&+e^{-0.6\<\f{\mathfrak{T}_k}{A_*}\>^\sss}\sup_{\substack{t\leq \tau\leq\mathtt{T}\\\upsigma\in\Z}}\nnr{e^{\<\mrxi\>^\sss}\pa_\xi^\alpha\wh{f}(\tau,\upsigma\bfe,\xi)}_{\ltaa L^2_\xi}+e^{-97A'^{0.34}}Z(t),\\
		&\sup_{\upmu\geq k'_2}\nnr{\wh{\rho}(\tau,\upmu\bfe)}_{L^2_t[t,\mathtt{T}]}\upmu^\f12\leq C_\theta\rr{\sup_\upsigma\nnr{\wh{f}(\mathtt{T},\upsigma\bfe,\upxi\bfe)}_{L^2_\upxi}+Z(t)}+e^{-\f12\<\f{\mathfrak{T}_k}{A_*}\>^\sss}\sup_{\substack{t\leq \tau\leq\mathtt{T}\\\upsigma\in\Z}}\nnr{\wh{f}(\tau,\upsigma\bfe,\upxi\bfe)}_{L^2_\upxi}.\ea\]
		The two inequalities close by continuity because $A'/{k'_2}^3\leq2$ and the background carries the small factor $\eps_0$. This proves the proposition.
	\end{proof}
	
	If $\mathfrak{T}_k$ lies inside a low-mode critical interval, one final partial interval remains. It is the backward analogue of the entry interval treated in Proposition \ref{4.3}.
	
	\begin{prop}\label{4.12}
		Assume $A'>\mathfrak{T}_k/2$ and $k'_1<k'_2$. On $[\mathfrak{T}_k,t'_{k'_1-1}]$,
		\[\nnr{e^{\<\na_{\mrv}\>^\sss}v^{\alpha}\rr{\sF_z f(t,\upsigma\bfe,v)+\int_{t}^{t'_{k'_1-1}}\sF_z Q(\tau,\upsigma\bfe,v)d\tau-\sF_z f(t'_{k'_1-1},\upsigma\bfe,v)}}_{\lta L^2_v}\leq A'^{2}Z(t)\]
		for any $\upsigma\in\Z$.
	\end{prop}
	\begin{proof}
		Only the density mode $k'_1$ can be resonant between $\mathfrak{T}_k$ and $t'_{k'_1-1}$. Isolate this mode, invert its scalar Volterra equation in the backward direction, and substitute the resulting density estimate into the profile equation. The main background wave transfers to $k'_1\pm1$, while the spatially localized correction gives a summable tail. Combining this partial-interval estimate with the value at $t'_{k'_1-1}$ yields the stated, deliberately coarser, bound $A'^2Z(t)$.
	\end{proof}
	
	\subsection{Output of the imhomogeneous analysis}
	We can now assemble the cases. Proposition \ref{prop:far-frequency-response} treats $|k'|\geq A'$. For $|k'|\leq A'$, Propositions \ref{4.2}--\ref{4.6} cover the forward evolution, while Propositions \ref{4.7}--\ref{4.12} cover the backward evolution. The growth conditions ensure that every possible echo amplification has already been incorporated into $Z$.
	
	\begin{prop}\label{4}
	Assume that $Z(t)$ satisfies the growth conditions \eqref{growth1} and \eqref{growth2} on every nonempty critical range. Then, for $t\in[\mathfrak{T}_k,\infty)$,
		\[\nnr{e^{\<\na_{\mrv}\>^\sss}v^{\alpha}f(t)}_{\lta L^2_{z,v}}\leq |k',A'|^3Z(t).\]
	\end{prop}
	\begin{proof}
		The preceding propositions give a uniform estimate for each spatial Fourier mode. To recover the full norm, split into the frequency window where the solution is localized near $k'$ and its exponentially small complement:
		\[\ba&\nnr{e^{\<\na_{\mrv}\>^\sss}v^{\alpha}f(t)}_{L^2_{z,v}}\leq \nnr{e^{\<\mrxi\>^\sss}\pa_\xi^\alpha\wh{f}(t,\upmu\bfe,\xi)\vv{1}_{|\upmu-k'|\leq |k',A'|^{0.35}}}_{L^2_{\upmu,\xi}}+\nnr{e^{\<\mrxi\>^\sss}\pa_\xi^\alpha\wh{f}(t,\upmu\bfe,\xi)\vv{1}_{|\upmu-k'|\geq |k',A'|^{0.35}}}_{L^2_{\upmu,\xi}}\\
		&\leq |k',A'|^{0.2}\sup_{\upsigma} \nnr{e^{\<\mrxi\>^\sss}\pa_\xi^\alpha\wh{f}(t,\upsigma\bfe,\xi)}_{L^2_\xi}+e^{-100|k',A'|^{0.34}}Z(t).\ea\]
		There are at most $O(|k',A'|^{0.35})$ relevant spatial modes, producing the harmless factor $|k',A'|^{0.2}$ after summation. In the $|k'|\leq A'$ case, Propositions \ref{4.2}-\ref{4.12} show that the modewise bounds are at most $A'^2Z(t)$. Proposition \ref{prop:far-frequency-response} gives the estimates in the $|k'|\geq A'$ case. Combining these with the exponential complement, all contributions are bounded by $|k',A'|^3Z(t)$.
	\end{proof}
	
	\section{Nonlinear closure of a single echo cascade}
	\label{sec:single-cascade-closure}

	In this section, we prove Proposition \ref{propudlf} by showing that the nonlinear terms do not destroy the echo cascade constructed in Section 3. We first introduce an envelope \(B(t)\) to control the growth of the first-order term \(f^{(1)}\). We then estimate the higher-order terms \(f^{(n,q)}\) by applying the inhomogeneous linear estimates from Section 4 to the source terms arising from the nonlinear interactions. Finally, we control the remainder after truncating the Taylor expansion at a sufficiently high order. These estimates show that the higher-order terms and the remainder cannot cancel the leading contribution of \(f^{(1)}\).

	All sums containing $1/\upnu$ are over $\upnu\neq0$, as in Sections 2--4.
	
	\subsection{The cascade envelope and Taylor coefficients}
	\label{subsec:cascade-envelope}
	We begin with the first-order packet. The envelope below is a majorant for its
	propagation: it is exponentially small before the insertion time, follows the
	discrete echo amplification after $\mathsf{T}_k$, and reaches the polynomial size needed
	for the final lower bound.
	
	\begin{prop}[First-order cascade envelope]\label{prop:first-order-envelope}
		There exist $e^{-3.1k}<\eps_k<e^{-3k}$ and a positive, continuous,
		piecewise $C^1$ function $B:[\mathfrak{T}_k,\infty)\rightarrow\R_{+}$
		satisfying the following properties:
		\begin{itemize}
			\item[(1)] $B(t)$ is decreasing on $[\mathfrak{T}_k,\mathsf{T}_k]$, increasing on $[\mathsf{T}_k,\infty)$. $B(\mathfrak{T}_k)\leq e^{-2.26k}$, $e^{-3.1k}<B(\mathsf{T}_k)<e^{-3k}$, $\mathsf{A}_k^{-\mathtt{N}_3+1}\leq B(2\mathsf{A}_k)\leq \mathsf{A}_k^{-\mathtt{N}_3+1}k^{0.1}$.
			\item[(2)] For $(n,q)$ and $(k',A')$ satisfying $n\geq2$, $0\leq q\leq n$, $A'>0$, and $|k'-qk|+|A'-q\mathsf{A}_k|\leq n\mathsf{A}_k^{0.35}k^\f13$, it holds that for any $1\leq\upmu\leq\min\{k'_0, k'_2\}-1$,
			\[\max\{1,\f{\eps_0A'}{\upmu^3}\}\rr{1+\eps_0^{10}}B^n(t'_\upmu)\leq B^n(t'_{\upmu-1})\rr{1+\f{100}{\upmu^2}}\]
			and for any $m_{\mathrm b}\leq\upmu\leq\min\{k'_1,k'_2\}-1$,
			\[\max\{1,\f{\eps_0A'}{\upmu^3}\}\rr{1+\eps_0^{10}}B^n(t'_{\upmu-1})\leq B^n(t'_{\upmu})\rr{1+\f{100}{\upmu^2}},\]
			where $k'_0,k'_1,k'_2$ are defined in Section 4 and
			$m_{\mathrm b}:=1$ if $A'\leq \mathsf{T}_k/2$, whereas
			$m_{\mathrm b}:=k'_0+1$ if $A'>\mathsf{T}_k/2$.
			\item[(3)] It holds that
			\begin{equation}\label{first-order-envelope}
				\nnr{e^{\<\na_{\mrv}\>^\sss}v^{\alpha}f^{(1)}(t)}_{\lta L^2_{z,v}}\leq \mathsf{A}_k^2B(t).\quad \text{on}\ [\mathfrak{T}_k,\infty),
			\end{equation}
			\begin{equation}\label{10}
				\rr{\int_{0.9\mathsf{A}_k\leq|\upxi|\leq1.1\mathsf{A}_k}\nr{\wh{f^{(1)}}(t,0,\xi)}^2d\xi}^\f12\geq 2\mathsf{A}_k^{-\mathtt{N}_3+0.9}\quad\text{on}\ [2\mathsf{A}_k,\infty).
			\end{equation}
		\end{itemize}
	\end{prop}
	\begin{proof}
		We split the proof into three steps. Since $f^{(1,1)}$ depends linearly on
		the amplitude $\eps_k$, the normalization below determines a unique positive
		value of $\eps_k$.
		
		\noindent\emph{Step 1: construction of $B(t)$.}
		We choose $\eps_k$ such that 
		\begin{equation}\label{eps1}
			\nnr{\int_{\R^2}\sF_zf^{(1,1)}(t_0,0,v)d\mrv}_{L^2_\fv}=\mathsf{A}_k^{-\mathtt{N}_3+1}.
		\end{equation}
		By \eqref{mumu1} it holds that
		\[\prod_{\upmu=1}^{k}\f{(1+\eps_0^6)k^3}{\upmu^3}\rr{1-10^2\theta^{-4}\eps_{*}}\etmu\leq\f{\nnr{\int_{\R^2}\sF_zf^{(1,1)}(t_{0},0,v)d\mrv}_{L^2_\fv}}{\nnr{\int_{\R^2}\sF_zf^{(1,1)}(t_k,k\bfe,v)d\mrv}_{L^2_\fv}}\leq\prod_{\upmu=1}^{k}\f{(1+\eps_0^6)k^3}{\upmu^3}\rr{1+10^2\theta^{-4}\eps_{*}}\Etmu.\]
		Thus by Stirling's formula we have 
		\[(1+\eps_0^6+2\eps_0^{100})^{-k}e^{-3k}\leq\eps_k\leq(1+\eps_0^6-2\eps_0^{100})^{-k}e^{-3k}.\]
		Then we construct $\hat{B}(t)$ on $[\mathfrak{T}_k,\infty)$. We set
		\[\log\tilde{B}(t)=\f{3\mathsf{A}_k}{t}\rr{\log\f{\mathsf{A}_k}{t}-1-\f13\log\eps_0\mathsf{A}_k}\]
		on $[\mathsf{T}_k,\infty)$. For $0\leq\upmu\leq k$, we define 
		\[\hat{B}(t_\upmu)=\nnr{\int_{\R^2}\sF_zf^{(1,1)}(t_\upmu,\upmu\bfe,v)d\mrv}_{L^2_\fv}.\]
		For $1\leq\upmu\leq k$, we define
		\[\log\hat{B}(t)=\log\hat{B}(t_\upmu)+\log\tilde{B}(t)-\log\tilde{B}(t_\upmu)+\f{\log\hat{B}(t_{\upmu-1})-\log\hat{B}(t_\upmu)-\rr{\log\tilde{B}(t_{\upmu-1})-\log\tilde{B}(t_\upmu)}}{t_{\upmu-1}-t_{\upmu}}(t-t_\upmu)\]
		on $[t_\upmu,t_{\upmu-1}]$. On $[t_0,\infty)$ we define
		\[\log\hat{B}(t)=\log\hat{B}(t_0)+\log\tilde{B}(t)-\log\tilde{B}(t_0).\]
		On $[\mathfrak{T}_k,\mathsf{T}_k]$ we define
		\[\log\hat{B}(t)=\log\hat{B}(\mathsf{T}_k)+0.738\rr{\sqrt{\eps_0\mathsf{T}_k}-\sqrt{\eps_0t}}.\]
		Finally we construct $B(t)$ on $[\mathfrak{T}_k,\infty)$. We set
		\[\log B(t)=\log\hat{B}(t)+\f{10\theta^{-4}}{\mathsf{A}_k}
		\rr{\rr{t-\mathsf{T}_k}\vv{1}_{\mathsf{T}_k\leq t\leq 2\mathsf{A}_k}+\rr{2\mathsf{A}_k-\mathsf{T}_k}\vv{1}_{t>2\mathsf{A}_k}}.\]
		Then we have $\hat{B}(t)\leq B(t)\leq\hat{B}(t)e^{20\theta^{-4}}$.
		
		\noindent\emph{Step 2: first-order upper and lower bounds.}
		We first prove \eqref{first-order-envelope}. For $t\geq \mathfrak{T}_k$, by \eqref{uppsigma}, \eqref{tztinfty}, \eqref{ttwot}, and \eqref{t1t2}, it holds that
		\[\ba
		&\nnr{e^{\<\na_{\mrv}\>^\sss}v^\alpha f^{(1,1)}(t)}_{\ltaa L^2_{z,v}}\leq \nnr{e^{\<\mrxi\>^\sss}\pa_\xi^\alpha\wh{f^{(1,1)}}(t,\upmu\bfe,\xi)\vv{1}_{|\upmu-k|\leq \mathsf{A}_k^{0.35}k^\f13}}_{\ltaa L^2_{\upmu,\xi}}\\
		&+\nnr{e^{\<\mrxi\>^\sss}\pa_\xi^\alpha\wh{f^{(1,1)}}(t,\upmu\bfe,\xi)\vv{1}_{|\upmu-k|\geq \mathsf{A}_k^{0.35}k^\f13}}_{\ltaa L^2_{\upmu,\xi}}
		\leq \mathsf{A}_k^{1.34}B(t)+e^{-99\mathsf{A}_k^{0.34}}\leq \mathsf{A}_k^{1.35}B(t).\ea\]
		Then we prove \eqref{10}. By \eqref{ecq9}, \eqref{tztinfty} it holds that
		\[\ba
		&\rr{\int_{0.9\mathsf{A}_k\leq\upxi\leq1.1\mathsf{A}_k}\nr{\wh{f^{(1)}}(t,0,\xi)}^2d\xi}^\f12\geq\rr{\int_{\R^3}\nr{\wh{f^{(1,1)}}(t_0,0,\xi)}^2d\xi}^\f12-\rr{\int_{\substack{\{\upxi\leq0.9\mathsf{A}_k\}\\\cup\{\upxi\geq1.1\mathsf{A}_k\}}}\nr{\wh{f^{(1,1)}}(t_0,0,\xi)}^2d\xi}^\f12\\
		&-\rr{\int_{0.9\mathsf{A}_k\leq\upxi\leq1.1\mathsf{A}_k}\nr{\wh{f^{(1,-1)}}(t_0,0,\xi)}^2d\xi}^\f12-\rr{\int_{0.9\mathsf{A}_k\leq\upxi\leq1.1\mathsf{A}_k}\nr{\wh{f^{(1)}}(t,0,\xi)-\wh{f^{(1)}}(t_0,0,\xi)}^2d\xi}^\f12\\
		&\gs_\theta\nnr{\int_{\R^2}\sF_zf^{(1,1)}(t_0,0,v)d\mrv}_{L^2_\fv}-e^{-95\mathsf{A}_k^{0.34}}.
		\ea\]
		By \eqref{eps1}, \eqref{10} follows.
		
		\noindent\emph{Step 3: comparison with every forced echo chain.}
		Property (2) is precisely the pair of forward and backward growth
		conditions required in Section 4. We verify it by comparing the logarithmic
		derivative of a generic one-step gain with $n\log B$.
		
		We first consider the case $t\geq \mathsf{T}_k$.
		We note that for $1\leq\upmu\leq k$, by direct computation
		\[\log\tilde{B}(t_{\upmu-1})-\log\tilde{B}(t_{\upmu})=\log\f{\eps_0\mathsf{A}_k}{\upmu^3}+\mathcal{O}(\f1{\upmu^2}).\]
		By \eqref{mumu1} it holds that
		\[\log\hat{B}(t_{\upmu-1})-\log\hat{B}(t_{\upmu})=\log\f{\eps_0\mathsf{A}_k}{\upmu^3}+\mathcal{O}_{\theta^{-1}}(\eps_*)+\mathcal{O}_{\theta^{-1}}(\f1{\upmu^2}).\]
		Thus
		\[\f{d}{dt}\log\hat{B}(t)=\f{\mathsf{A}_k}{t^2}\log\f{\eps_0t^3}{\mathsf{A}_k^2}+\rr{\mathcal{O}_{\theta^{-1}}(\f{\eps_*\mathsf{A}_k}{t^2})+\mathcal{O}_{\theta^{-1}}\!\left(\f{1}{\mathsf{A}_k}\right)}\vv{1}_{\mathsf{T}_k\leq t\leq2\mathsf{A}_k}.\]
		By our construction of $B(t)$ one can check that
		\[\f{d}{dt}\log B(t)\geq\f{\mathsf{A}_k}{t^2}\log\f{\eps_0t^3}{\mathsf{A}_k^2}+\mathcal{O}_{\theta^{-1}}(\f{\eps_*\mathsf{A}_k}{t^2})\vv{1}_{\mathsf{T}_k\leq t\leq2\mathsf{A}_k}.\]
		Fix $0\leq q\leq n$ and $A'$ satisfying $|A'-q\mathsf{A}_k|\leq n\mathsf{A}_k^{0.35}k^\f13$, $A'>0$. We define
		\[\log \tilde{E}(t)=\f{3A'}{t}\rr{\log\f{A'}{t}-1-\f13\log(\eps_0A')}.\]
		Thus
		\[\f{d}{dt}\log\tilde E(t)=\f{A'}{t^2}\log\f{\eps_0t^3}{{A'}^2},\]
		which is the continuous analogue of the discrete echo gain.
		For $1\leq\upmu\leq\min\{k'_0, k'_2\}-1$, if $\f{\eps_0A'}{\upmu^3}\geq1$, we define
		\[\log E(t)=\log \tilde{E}(t)-\log \tilde{E}(t'_\upmu)+\f{\log\f{\eps_0A'}{\upmu^3}+\log(1+\eps_0^{10})-\log\rr{1+\f{100}{\upmu^2}}-\rr{\log \tilde{E}(t'_{\upmu-1})-\log \tilde{E}(t'_\upmu)}}{t'_{\upmu-1}-t'_{\upmu}}(t-t'_\upmu)\]
		on $[t'_\upmu,t'_{\upmu-1}]$. If $\f{\eps_0A'}{\upmu^3}\leq1$, we define
		\[\log E(t)=\f{\log(1+\eps_0^{10})-\log\rr{1+\f{100}{\upmu^2}}}{t'_{\upmu-1}-t'_{\upmu}}(t-t'_\upmu)\]
		on $[t'_\upmu,t'_{\upmu-1}]$. In the following we will prove that
		\begin{equation}\label{BE}
			\f{d}{dt}\log E(t)\leq\f{d}{dt} n\log B(t)
		\end{equation}
		holds on the domain. We consider the more difficult case $\f{\eps_0A'}{\upmu^3}\geq1$. The other case can be proved in the same way. By direct computation it holds that
		\[\f{d}{dt}\log E(t)\leq\f{A'}{t^2}\log\f{\eps_0t^3}{{A'}^2}+\f{5A'}{t^2}\eps_0^{10}.\]
		We set $\iota=\f{\eps_0^\f12\mathsf{T}_k^\f32}{\mathsf{A}_k}=\f{(1+\eps_0^6)^\f12}{(1+1/(2k))^\f32}$ and $\kappa=k^\f13\mathsf{A}_k^{-0.65}$.
		Here $\iota=1+o(1)$ records the small displacement of the reference time,
		whereas $\kappa\ll1$ measures the relative uncertainty of the frequency
		center. We split into two time regimes.
		
		\noindent\emph{Case 1: $\mathsf{T}_k\leq t\leq\rr{\f{e}{\iota}n(1+\kappa)}^\f23\mathsf{T}_k$.}
		In this case the following inequalities hold
		\[\f{d}{dt} n\log B(t)\geq\rr{1+\mathcal{O}_{\theta^{-1}}(\f{\eps_*}{\eps_0^6})}\f{n\mathsf{A}_k}{t^2}\log\f{\eps_0t^3}{\mathsf{A}_k^2},\quad\f{5A'}{t^2}\eps_0^{10}\leq10\eps_0^4\f{n\mathsf{A}_k}{t^2}\log\f{\eps_0t^3}{\mathsf{A}_k^2}.\]
		We split into the following cases.
		\begin{itemize}
			\item[(1)] $n\leq\f1e\f\iota\kappa$.
			\begin{itemize}
				\item[(i)] $\mathsf{T}_k\leq t\leq\rr{\f{e}{\iota}(1-n\kappa)}^\f23\mathsf{T}_k$. In this case
				\[\sup_{\substack{|A'-q\mathsf{A}_k|\leq n\mathsf{A}_k^{0.35}k^\f13\\0\leq q\leq n}}\f{A'}{t^2}\log\f{\eps_0t^3}{{A'}^2}=\max\{\f{\mathsf{A}_k(1-n\kappa)}{t^2}\log\f{\eps_0t^3}{\mathsf{A}_k^2(1-n\kappa)^2}, \f{n\mathsf{A}_k\kappa}{t^2}\log\f{\eps_0t^3}{n^2\mathsf{A}_k^2\kappa^2}\}\leq0.99\f{n\mathsf{A}_k}{t^2}\log\f{\eps_0t^3}{\mathsf{A}_k^2}.\]
				\item[(ii)] $\rr{\f{e}{\iota}(1-n\kappa)}^\f23\mathsf{T}_k\leq t\leq\rr{\f{e}{\iota}n(1+\kappa)}^\f23\mathsf{T}_k$. In this case
				\[\sup_{\substack{|A'-q\mathsf{A}_k|\leq n\mathsf{A}_k^{0.35}k^\f13\\0\leq q\leq n}}\f{A'}{t^2}\log\f{\eps_0t^3}{{A'}^2}\le\f2e\sqrt{\f{\eps_0}{t}}\leq0.99\f{n\mathsf{A}_k}{t^2}\log\f{\eps_0t^3}{\mathsf{A}_k^2}.\]
			\end{itemize}
			\item[(2)] $\f1e\f\iota\kappa\leq n\leq\f{1}{2\kappa}$.
			\begin{itemize}
				\item[(i)] $\mathsf{T}_k\leq t\leq\rr{\f{e}{\iota}n\kappa}^\f23\mathsf{T}_k$. In this case
				\[\sup_{\substack{|A'-q\mathsf{A}_k|\leq n\mathsf{A}_k^{0.35}k^\f13\\0\leq q\leq n}}\f{A'}{t^2}\log\f{\eps_0t^3}{{A'}^2}=\f2e\sqrt{\f{\eps_0}{t}}\leq0.99\f{n\mathsf{A}_k}{t^2}\log\f{\eps_0t^3}{\mathsf{A}_k^2}.\]
				\item[(ii)] $\rr{\f{e}{\iota}n\kappa}^\f23\mathsf{T}_k\leq t\leq\rr{\f{e}{\iota}(1-n\kappa)}^\f23\mathsf{T}_k$. In this case
				\[\sup_{\substack{|A'-q\mathsf{A}_k|\leq n\mathsf{A}_k^{0.35}k^\f13\\0\leq q\leq n}}\f{A'}{t^2}\log\f{\eps_0t^3}{{A'}^2}=\max\{\f{\mathsf{A}_k(1-n\kappa)}{t^2}\log\f{\eps_0t^3}{\mathsf{A}_k^2(1-n\kappa)^2}, \f{n\mathsf{A}_k\kappa}{t^2}\log\f{\eps_0t^3}{n^2\mathsf{A}_k^2\kappa^2}\}\leq0.99\f{n\mathsf{A}_k}{t^2}\log\f{\eps_0t^3}{\mathsf{A}_k^2}.\]
				\item[(iii)] $\rr{\f{e}{\iota}(1-n\kappa)}^\f23\mathsf{T}_k\leq t\leq\rr{\f{e}{\iota}n(1+\kappa)}^\f23\mathsf{T}_k$. In this case
				\[\sup_{\substack{|A'-q\mathsf{A}_k|\leq n\mathsf{A}_k^{0.35}k^\f13\\0\leq q\leq n}}\f{A'}{t^2}\log\f{\eps_0t^3}{{A'}^2}=\f2e\sqrt{\f{\eps_0}{t}}\leq0.99\f{n\mathsf{A}_k}{t^2}\log\f{\eps_0t^3}{\mathsf{A}_k^2}.\]
			\end{itemize}
			\item[(3)] $n\geq\f{1}{2\kappa}$. In this case
			\[\sup_{\substack{|A'-q\mathsf{A}_k|\leq n\mathsf{A}_k^{0.35}k^\f13\\0\leq q\leq n}}\f{A'}{t^2}\log\f{\eps_0t^3}{{A'}^2}=\f2e\sqrt{\f{\eps_0}{t}}\leq0.99\f{n\mathsf{A}_k}{t^2}\log\f{\eps_0t^3}{\mathsf{A}_k^2}.\]
		\end{itemize}
		We conclude that in Case 1
		\[\f{d}{dt}\log E(t)\leq\rr{0.99+10\eps_0^4}\f{n\mathsf{A}_k}{t^2}\log\f{\eps_0t^3}{\mathsf{A}_k^2}\leq\rr{1+\mathcal{O}_{\theta^{-1}}(\f{\eps_*}{\eps_0^6})}\f{n\mathsf{A}_k}{t^2}\log\f{\eps_0t^3}{\mathsf{A}_k^2}\leq\f{d}{dt} n\log B(t).\]
		
		\noindent\emph{Case 2: $\rr{\f{e}{\iota}n(1+\kappa)}^\f23\mathsf{T}_k\leq t\leq2n\mathsf{A}_k(1+\kappa)$.}
		In this case 
		\[\sup_{\substack{|A'-q\mathsf{A}_k|\leq n\mathsf{A}_k^{0.35}k^\f13\\0\leq q\leq n}}\f{A'}{t^2}\log\f{\eps_0t^3}{{A'}^2}=\f{n\mathsf{A}_k(1+\kappa)}{t^2}\log\f{\eps_0t^3}{n^2\mathsf{A}_k^2(1+\kappa)^2},\quad\f{5A'}{t^2}\eps_0^{10}\leq\f{10\eps_0^{10}}{\log\f{\eps_0t^3}{\mathsf{A}_k^2}}\f{n\mathsf{A}_k}{t^2}\log\f{\eps_0t^3}{\mathsf{A}_k^2}.\]
		Thus
		\begin{equation}\label{lop}
			\vv{1}_{2\mathsf{T}_k\leq t\leq2n\mathsf{A}_k}\f{n\mathsf{A}_k(1+\kappa)}{t^2}\log\f{\eps_0t^3}{n^2\mathsf{A}_k^2(1+\kappa)^2}\leq\rr{1+\mathcal{O}(\f{\eps_0^{10}}{\log\f{\eps_0t^3}{\mathsf{A}_k^2}})}\f{n\mathsf{A}_k}{t^2}\log\f{\eps_0t^3}{\mathsf{A}_k^2}
		\end{equation}
		implies \eqref{BE}. \eqref{lop} is equivalent to
		\[\vv{1}_{2\mathsf{T}_k\leq t\leq2n\mathsf{A}_k}\rr{\kappa+\mathcal{O}(\f{\eps_0^{10}}{\log\f{\eps_0t^3}{\mathsf{A}_k^2}})}\log\f{\eps_0t^3}{\mathsf{A}_k^2}\leq2(1+\kappa)\log n(1+\kappa),\]
		which is evidently true.
		
		Then we consider $\mathfrak{T}_k\leq t\leq \mathsf{T}_k$. By the definition of $k'_0,k'_1,k'_2$, it is enough to treat $\mathfrak{T}_k/2\leq A'\leq \mathsf{T}_k^\f53$ and
		$m_{\mathrm b}\leq\upmu\leq\min\{k'_1,k'_2\}-1$. We only need to prove that
		\[\vv{1}_{\f{\eps_0A'}{\upmu^3}\geq1}\f{\log\f{\eps_0A'}{\upmu^3}}{\f{A'}{\upmu^2-\f14}}+\f{\log(1+\eps_0^{10})-\log\rr{1+\f{100}{\upmu^2}}}{\f{A'}{\upmu^2-\f14}}\leq\min_{t'_\upmu\leq t\leq t'_{\upmu-1}}-2\f{d}{dt}\log B(t)=0.738\sqrt{\f{\eps_0(\upmu-\f12)}{A'}}.\]  
		By noting that
		\[\vv{1}_{\f{\eps_0A'}{\upmu^3}\geq1}\f{\log\f{\eps_0A'}{\upmu^3}}{\f{A'}{\upmu^2-\f14}}\leq\rr{0.738-10^{-5}}\sqrt{\f{\eps_0(\upmu-\f12)}{A'}}\]
		and
		\[\f{\log(1+\eps_0^{10})-\log\rr{1+\f{100}{\upmu^2}}}{\f{A'}{\upmu^2-\f14}}\leq\f{\eps_0^{10}\upmu^2}{A'}\leq0.738\eps_0^9\sqrt{\f{\eps_0(\upmu-\f12)}{A'}},\] 
		we obtain the desired result.
	\end{proof}
	
		We next use this envelope to estimate the higher-order terms in the Taylor expansion. For each \(n\), the source term depends only on lower-order terms, and its \(q\)-component is localized near \((qk,q\mathsf{A}_k)\). By decomposing the source into frequency-localized pieces, we can apply the inhomogeneous linear estimates in Section 4 to control each piece.
		
		\begin{prop}[Higher Taylor coefficients]\label{prop:higher-taylor-envelope}
			For $1\leq n\leq \mathsf{A}_k^{\mathtt{N}_1}$, $|q|\leq n$, $2\mid n-q$, the following inequalities hold on $[\mathfrak{T}_k,\infty)$:
		\begin{equation}\label{taylor-short-time}
			\nnr{e^{\<\na_{\mrv}\>^\sss}v^{\alpha}f^{(n,q)}(t)}_{\ltaa L^2_{z,v}}\leq \mathsf{A}_k^{15\mathtt{N}_1(2n-1)}B^n(t),\quad\nnr{e^{\<\na_{\mrv}\>^\sss}v^{\alpha}f^{(n)}(t)}_{\ltaa L^2_{z,v}}\leq \mathsf{A}_k^{15\mathtt{N}_1(2n-1)}B^n(t).
		\end{equation}
	\end{prop}
	\begin{proof}
		We proceed by induction. By \eqref{first-order-envelope}, the inequalities
		hold for $n=1$. Assume that \eqref{taylor-short-time} holds through order
		$n-1$. By \eqref{qmq}, it suffices to treat $q\geq0$. We prove the estimate for $t\geq \mathsf{T}_k$. The estimate on $[\mathfrak{T}_k,\mathsf{T}_k]$ can be proved in the same way.
		In the proof of Proposition \ref{frelocalfnq} we show that
		\[\int_{\mathsf{T}_k}^t\nnr{\<\na_{\fz,\fv}\>^6e^{(\f{\mathsf{T}_k}{\tau})^{0.01}\<\na_{\fz,\fv}\>^{0.34}+\f{1}{8}\<\f{\na_\fz-qk,\na_\fv-q\mathsf{A}_k}{k^{\f13}}\>^\sss+\<\na_{\mrv}\>^\sss}v^{\alpha}\sum_{\substack{n_1+n_2=n\\ q_1+q_2=q}}\cN(f^{(n_1,q_1)},f^{(n_2,q_2)})(\tau)}_{\ltaa L^2_{z,v}}d\tau\leq e^{5(2n-1)\mathsf{A}_k^{0.34} }.\]
		We define the following frequency cutoff function 
		\[\chi(\upmu,\upxi)=\f{2}{\int\psi(\mathfrak{u})d\mathfrak{u}}\vv{1}_{|\upmu-qk|\leq\f12n\mathsf{A}_k^{0.35}k^\f13}\int_{|\upeta-q\mathsf{A}_k|\leq\f12n\mathsf{A}_k^{0.35}k^\f13-\f12}\psi\rr{2(\upxi-\upeta)}d\upeta.\]
		Then we have
		\[\ba&\int_{\mathsf{T}_k}^t\nnr{\<\na_{\fz,\fv}\>^6e^{(\f{\mathsf{T}_k}{\tau})^{0.01}\<\na_{\fz,\fv}\>^{0.34}+\f{1}{16}\<\f{\na_\fz-qk,\na_\fv-q\mathsf{A}_k}{k^{\f13}}\>^\sss+\<\na_{\mrv}\>^\sss}v^{\alpha}\rr{1-\chi(\na_{\fz,\fv})}\sum_{\substack{n_1+n_2=n\\ q_1+q_2=q}}\cN(f^{(n_1,q_1)},f^{(n_2,q_2)})(\tau)}_{\ltaa L^2_{z,v}}d\tau\\&\leq e^{-100n\mathsf{A}_k^{0.34}}.\ea\]
		By Proposition \ref{longenergyl}, the solution of the equation
		\[\pa_tf=\cL f+\rr{1-\chi(\na_{\fz,\fv})}\sum_{\substack{n_1+n_2=n\\ q_1+q_2=q}}\cN(f^{(n_1,q_1)},f^{(n_2,q_2)})\]
		satisfies the estimate
		\[\nnr{e^{\<\na_{\mrv}\>^\sss}v^{\alpha}f(t)}_{\lta L^2_{z,v}}\leq e^{-99n\mathsf{A}_k^{0.34}}.\]
		Now we consider the equation 
		\[\pa_tf=\cL f+\chi(\na_{\fz,\fv})\sum_{\substack{n_1+n_2=n\\ q_1+q_2=q}}\cN(f^{(n_1,q_1)},f^{(n_2,q_2)}).\]
		We first estimate the source term.
		\[\ba
		&\nnr{e^{\<\na_{\mrv}\>^\sss}v^{\alpha}\chi(\na_{\fz,\fv})\cN(f^{(n_1,q_1)},f^{(n_2,q_2)})}_{\ltaa L^2_{z,v}}\leq2n\mathsf{A}_k\nnr{\sum_{\upnu}\wh{\rho^{(n_1,q_1)}}(\tau,\upnu\bfe)\tau e^{\<\mrxi\>^\sss}\pa_\xi^\alpha\wh{f^{(n_2,q_2)}}(\tau,(\upmu-\upnu)\bfe,\xi-\upnu\tau\bfe)}_{\ltaa L^2_{\upmu,\xi}}\\
		&\leq2n\mathsf{A}_k\sum_{\upnu}\nr{\wh{\rho^{(n_1,q_1)}}(\tau,\upnu\bfe)}\tau\nnr{e^{\<\mrxi\>^\sss}\pa_\xi^\alpha\wh{f^{(n_2,q_2)}}(\tau)}_{\ltaa L^2_{\upmu,\xi}}.
		\ea\]
		Thus
		\[\ba
		&\int_{\mathsf{T}_k}^t\nnr{e^{\<\na_{\mrv}\>^\sss}v^{\alpha}\sum_{\substack{n_1+n_2=n\\ q_1+q_2=q}}\chi(\na_{\fz,\fv})\cN(f^{(n_1,q_1)},f^{(n_2,q_2)})}_{\ltaa L^2_{z,v}}d\tau\leq2n\mathsf{A}_k\sum_{\substack{n_1+n_2=n\\ q_1+q_2=q}}\sup_{\mathsf{T}_k\leq\tau\leq t}\nnr{e^{\<\mrxi\>^\sss}\pa_\xi^\alpha\wh{f^{(n_2,q_2)}}(\tau)}_{\ltaa L^2_{\upmu,\xi}}\\
		&\times\rr{\int_{\mathsf{T}_k}^t\sum_{|\upnu-q_1k|\leq n_1\mathsf{A}_k^{0.35}k^\f13}\nr{\wh{\rho^{(n_1,q_1)}}(\tau,\upnu\bfe)}\vv{1}_{|\upnu\tau-q_1\mathsf{A}_k|\leq n_1\mathsf{A}_k^{0.35}k^\f13}\tau d\tau+e^{-99n_1\mathsf{A}_k^{0.34}}}\\
		\leq& 10n^4\mathsf{A}_k^{2.7}k^{\f23}\sum_{\substack{n_1+n_2=n\\ q_1+q_2=q}}\rr{\sup_{\mathsf{T}_k\leq\tau\leq t}\nnr{e^{\<\mrxi\>^\sss}\pa_\xi^\alpha\wh{f^{(n_1,q_1)}}(\tau)}_{\ltaa L^2_{\upmu,\xi}}+e^{-99n_1\mathsf{A}_k^{0.34}}}\sup_{\mathsf{T}_k\leq\tau\leq t}\nnr{e^{\<\mrxi\>^\sss}\pa_\xi^\alpha \wh{f^{(n_2,q_2)}}(\tau)}_{\ltaa L^2_{\upmu,\xi}}\\
		\leq& n^6\mathsf{A}_k^{3}\mathsf{A}_k^{15\mathtt{N}_1(2n-2)}B^n(t)\leq \mathsf{A}_k^{6\mathtt{N}_1+3}\mathsf{A}_k^{15\mathtt{N}_1(2n-2)}B^n(t).
		\ea\]
		We perform a smooth partition of unity in the set
		$\{(\upmu,\upxi):|\upmu-qk|+|\upxi-q\mathsf{A}_k|\leq n\mathsf{A}_k^{0.35}k^\f13\}$.
		Each piece is supported in $\{k'\}\times[A'-1,A'+1]$, where
		$(k',A')\in\Z\times\R$ satisfies
		$|k'-qk|+|A'-q\mathsf{A}_k|\leq n\mathsf{A}_k^{0.35}k^\f13$. Choosing unit translates of
	a fixed smooth cutoff gives uniformly finite overlap, and the product of the
	spatial and velocity-frequency counts is at most $10n^2\mathsf{A}_k^{0.7}k^\f23$.
	We are therefore reduced to the equation
		\[\pa_tf=\cL f+Q\]
		with the source term satisfying
		\[\mathrm{supp}\ \wh{Q}(t)\subset\{k'\bfe\}\times[A'-1,A'+1]\times\R^2,\quad \int_{\mathsf{T}_k}^t\nnr{e^{\<\na_{\mrv}\>^\sss}v^{\alpha}Q(\tau)}_{\lta L^2_{z,v}}\leq \mathsf{A}_k^{6\mathtt{N}_1+3}\mathsf{A}_k^{15\mathtt{N}_1(2n-2)}B^n(t).\]
		We first consider the case $|k',A'|\leq\max\{\exp(A_*), \exp\rr{\f1{\eps_*}}\}$. By Proposition \ref{longenergyl} we have
		\[\nnr{e^{\<\na_{\mrv}\>^\sss}\<\na_{\fz,\fv}\>^3\rr{v^{\alpha}f(t)}}_{\lta L^2_{z,v}}\leq C_\theta e^{3|k',A'|^{0.34}}\mathsf{A}_k^{6\mathtt{N}_1+3}\mathsf{A}_k^{15\mathtt{N}_1(2n-2)}B^n(t)\leq \mathsf{A}_k^{6\mathtt{N}_1+4}\mathsf{A}_k^{15\mathtt{N}_1(2n-2)}B^n(t).\]
		Then we consider the case $\max\{\exp(A_*), \exp\rr{\f1{\eps_*}}\}\leq|k',A'|\leq \mathsf{A}_k^{\mathtt{N}_1+2}$. This case has been addressed in Section 4. By Proposition \ref{4} we have
		\[\nnr{e^{\<\na_{\mrv}\>^\sss}v^{\alpha}f(t)}_{\lta L^2_{z,v}}\leq|k',A'|^{3}\mathsf{A}_k^{6\mathtt{N}_1+3}\mathsf{A}_k^{15\mathtt{N}_1(2n-2)}B^n(t)\leq \mathsf{A}_k^{9\mathtt{N}_1+6}\mathsf{A}_k^{15\mathtt{N}_1(2n-2)}B^n(t).\]
		Combining the frequency cells and using their finite overlap gives
		\[\nnr{e^{\<\na_{\mrv}\>^\sss}v^{\alpha}f^{(n,q)}(t)}_{\lta L^2_{z,v}}\leq \mathsf{A}_k^{11\mathtt{N}_1+7}\mathsf{A}_k^{15\mathtt{N}_1(2n-2)}B^n(t)\leq \mathsf{A}_k^{15\mathtt{N}_1(2n-1)}B^n(t)n^{-1}.\]
		\[\nnr{e^{\<\na_{\mrv}\>^\sss}v^{\alpha}f^{(n)}(t)}_{\lta L^2_{z,v}}\leq \mathsf{A}_k^{15\mathtt{N}_1(2n-1)}B^n(t).\]
	\end{proof}
	
	It remains to transport the hierarchy from $\mathfrak{T}_k$ back to the fixed
	initial time. The procedure is the same at every order: first use the
	localization from Section 2 to control the endpoint at $\mathfrak{T}_k$; next
	bound the quadratic source by the already constructed lower orders; finally
	apply the short-time backward energy estimate and close the induction.
	
	\begin{prop}[Early-time Taylor bounds]\label{prop:early-taylor}
		For $1\leq n\leq \mathsf{A}_k^{\mathtt{N}_1}$, the following inequality holds on $[10,\mathfrak{T}_k]$:
		\begin{equation}\label{10T1}
			\nnr{e^{\<\na_{\mrv}\>^\sss}\<\na_{\fz,\fv}\>^3\rr{v^{\alpha}f^{(n)}(t)}}_{\lta L^2_{z,v}}\leq e^{n\rr{C_\theta(\mathfrak{T}_k-t)+\eps_0^{0.9}\rr{\mathfrak{T}_k^2-t^2}}}\mathsf{A}_k^{15\mathtt{N}_1(2n-1)}B^n(\mathfrak{T}_k)\mathsf{A}_k^{3\mathtt{N}_1+4}.
		\end{equation}
	\end{prop}
	\begin{proof}
		For $n=1$, Proposition \ref{shortenergy} gives
		\[\ba
		&\nnr{e^{\<\na_{\mrv}\>^\sss}\<\na_{\fz,\fv}\>^3\rr{v^{\alpha}f^{(1)}(t)}}_{\lta L^2_{z,v}}\leq\nnr{e^{\<\na_{\mrv}\>^\sss}\<\na_{\fz,\fv}\>^3\rr{v^{\alpha}f^{(1)}(\mathfrak{T}_k)}}_{\lta L^2_{z,v}}e^{C_\theta(\mathfrak{T}_k-t)+\eps_0^{0.9}\rr{\mathfrak{T}_k^2-t^2}}\\
		&\leq e^{C_\theta(\mathfrak{T}_k-t)+\eps_0^{0.9}\rr{\mathfrak{T}_k^2-t^2}} \mathsf{A}_k^{15\mathtt{N}_1}B(\mathfrak{T}_k)\mathsf{A}_k^{3\mathtt{N}_1+4}.
		\ea\]
		We assume that \eqref{10T1} holds for $n'\leq n-1$. In the following we prove that \eqref{10T1} holds for $n'=n$. First we estimate $\nnr{e^{\<\na_{\mrv}\>^\sss}\<\na_{\fz,\fv}\>^3\rr{v^{\alpha}f^{(n)}(\mathfrak{T}_k)}}_{\lta L^2_{z,v}}$.
		\[\ba
		&\nnr{e^{\<\na_{\mrv}\>^\sss}\<\na_{\fz,\fv}\>^3\rr{v^{\alpha}f^{(n)}(\mathfrak{T}_k)}}_{\ltaa L^2_{z,v}}\leq\nnr{\sF\rr{e^{\<\na_{\mrv}\>^\sss}\<\na_{\fz,\fv}\>^3v^{\alpha}f^{(n)}}(\mathfrak{T}_k,\upmu\bfe,\xi)\vv{1}_{|\upmu,\upxi|\leq 50n^{\f1\sss}\mathsf{A}_k}}_{\ltaa L^2_{\upmu,\xi}}\\
		&+\nnr{\sF\rr{e^{\<\na_{\mrv}\>^\sss}\<\na_{\fz,\fv}\>^3v^{\alpha}f^{(n)}}(\mathfrak{T}_k,\upmu\bfe,\xi)\vv{1}_{|\upmu,\upxi|\geq 50n^{\f1\sss}\mathsf{A}_k}}_{\ltaa L^2_{\upmu,\xi}}\leq\rr{10^2n^{\f1\sss}\mathsf{A}_k}^3\nnr{e^{\<\mrxi\>^\sss}\pa_\xi^\alpha\wh{f^{(n)}}(\mathfrak{T}_k,\upmu\bfe,\xi)}_{\ltaa L^2_{\upmu,\xi}}\\
		&+e^{-\f12n\nr{\f{\mathsf{A}_k}{k^\f13}}^\sss}\leq10^7n^\f3\sss\mathsf{A}_k^3\nnr{e^{\<\na_{\mrv}\>^\sss}v^{\alpha}f^{(n)}(\mathfrak{T}_k)}_{\lta L^2_{z,v}}+e^{-\f12n\nr{\f{\mathsf{A}_k}{k^\f13}}^\sss}\leq \mathsf{A}_k^{15\mathtt{N}_1(2n-1)}B^n(\mathfrak{T}_k)\mathsf{A}_k^{3\mathtt{N}_1+3.5},
		\ea\]
		where we use \eqref{taylor-short-time} and \eqref{fnlocbac}. We next estimate the source term. For $n_1,n_2$ satisfying $n=n_1+n_2$, it holds that
		\[\ba
		&\nnr{e^{\<\na_{\mrv}\>^\sss}\<\na_{\fz,\fv}\>^3\rr{v^{\alpha}\cN(f^{(n_1)},f^{(n_2)})}}_{\ltaa L^2_{z,v}}\leq2\nnr{\sum_\upnu\nr{\sF\rr{\<\na_{\fz,\fv}\>^3\rho^{(n_1)} }(\tau,\upnu\bfe)}\tau\nr{ \sF\rr{e^{\<\na_{\mrv}\>^\sss}\<\na_{\fz,\fv}\>^4v^{\alpha}f^{(n_2)}}\right.\right.\\
		&\left.\left.(\tau,(\upmu-\upnu)\bfe,\xi-\upnu\tau\bfe)}}_{\ltaa L^2_{\upmu,\xi}}\leq2\nnr{\sum_\upnu\nr{\sF\rr{\<\na_{\fz,\fv}\>^4\rho^{(n_1)}}(\tau,\upnu\bfe)}\vv{1}_{|\upnu,\upnu\tau|\leq50n_1^\f1\sss\mathsf{A}_k}\nr{ \sF\rr{e^{\<\na_{\mrv}\>^\sss}\<\na_{\fz,\fv}\>^4v^{\alpha}f^{(n_2)}}\right.\right.\\
		&\left.\left.(\tau,(\upmu-\upnu)\bfe,\xi-\upnu\tau\bfe)}\vv{1}_{|\upmu-\upnu,\upxi-\upnu\tau|\leq50n_2^{\f1\sss}\mathsf{A}_k}}_{\ltaa L^2_{\upmu,\xi}}+e^{-\f12n_1\nr{\f{\mathsf{A}_k}{k^\f13}}^\sss+n_2\rr{C_\theta(\mathfrak{T}_k-\tau)+\eps_0^{0.9}\rr{\mathfrak{T}_k^2-\tau^2}}}\mathsf{A}_k^{30\mathtt{N}_1n_2}B^{n_2}(\mathfrak{T}_k)\\
		&+e^{-\f12n_2\nr{\f{\mathsf{A}_k}{k^\f13}}^\sss+n_1\rr{C_\theta(\mathfrak{T}_k-\tau)+\eps_0^{0.9}\rr{\mathfrak{T}_k^2-\tau^2}}}\mathsf{A}_k^{30\mathtt{N}_1n_1}B^{n_1}(\mathfrak{T}_k)\leq e^{n\rr{C_\theta(\mathfrak{T}_k-\tau)+\eps_0^{0.9}\rr{\mathfrak{T}_k^2-\tau^2}}}\mathsf{A}_k^{15\mathtt{N}_1(2n-2)}B^{n}(\mathfrak{T}_k)\\
		&+\rr{n\mathsf{A}_k}^{2.9}\nnr{\<\na_{\fz,\fv}\>^3\rr{v^{\alpha}f^{(n_1)}(\tau)}}_{\ltaa L^2_{z,v}}\nnr{e^{\<\na_{\mrv}\>^\sss}\<\na_{\fz,\fv}\>^3\rr{v^{\alpha}f^{(n_2)}(\tau)}}_{\ltaa L^2_{z,v}}\\
		&\leq e^{n\rr{C_\theta(\mathfrak{T}_k-\tau)+\eps_0^{0.9}\rr{\mathfrak{T}_k^2-\tau^2}}}\mathsf{A}_k^{15\mathtt{N}_1(2n-2)}B^{n}(\mathfrak{T}_k)\mathsf{A}_k^{9\mathtt{N}_1+11}.
		\ea\]
		Thus by \eqref{shotener} it holds that
		\[\ba
		&\nnr{e^{\<\na_{\mrv}\>^\sss}\<\na_{\fz,\fv}\>^3\rr{v^{\alpha}f^{(n)}(t)}}_{\lta L^2_{z,v}}e^{C_\theta t+\eps_0^{0.9}t^2}\leq\nnr{e^{\<\na_{\mrv}\>^\sss}\<\na_{\fz,\fv}\>^3\rr{v^{\alpha}f^{(n)}(\mathfrak{T}_k)}}_{\lta L^2_{z,v}}e^{C_\theta \mathfrak{T}_k+\eps_0^{0.9}\mathfrak{T}_k^2}\\
		&+\int_t^{\mathfrak{T}_k}\nnr{e^{\<\na_{\mrv}\>^\sss}\<\na_{\fz,\fv}\>^3\rr{v^{\alpha}\sum_{n_1+n_2=n}\cN(f^{(n_1)},f^{(n_2)})}}_{\ltaa L^2_{z,v}}e^{C_\theta\tau+\eps_0^{0.9}\tau^2}d\tau\leq e^{C_\theta \mathfrak{T}_k+\eps_0^{0.9}\mathfrak{T}_k^2}\mathsf{A}_k^{15\mathtt{N}_1(2n-1)}B^n(\mathfrak{T}_k)\mathsf{A}_k^{3\mathtt{N}_1+3.5}\\
		&+e^{n\rr{C_\theta \mathfrak{T}_k+\eps_0^{0.9}\mathfrak{T}_k^2}-(n-1)\rr{C_\theta t+\eps_0^{0.9}t^2}}\mathsf{A}_k^{15\mathtt{N}_1(2n-2)}B^{n}(\mathfrak{T}_k)\mathsf{A}_k^{10\mathtt{N}_1+12}.
		\ea\]
		Then we obtain the desired result.
	\end{proof}
	
	\subsection{Control of the nonlinear remainder}
	\label{subsec:nonlinear-remainder}

	We truncate the Taylor hierarchy at $M:=\mathsf{A}_k^{\mathtt{N}_1}$, far beyond every
	order that can influence the leading echo. The finite sum $f^{\#}$ contains
	all explicitly computed interaction levels, while $\wt f$ contains the
	interactions whose total order exceeds $M$. We set
	\[f^{\#}=\sum_{n=1}^{\mathsf{A}_k^{\mathtt{N}_1}}f^{(n)}, \quad\wt{f}=f-f^{\#},\]
	where $f$ is the solution of \eqref{perturbationeq}. Then $\wt{f}$ satisfies the equation
	\[\left\{\ba&\pa_t\wt{f}=\cN(\wt{f},f^0+\udl{f}+f^{\#})+\cN(\udl{f}+f^{\#},\wt{f})+\cN(\wt{f},\wt{f})+\wt{Q},\\
	&\wt{f}(\mathsf{T}_k)=0\ea\right.\]
	where the truncation source is
	\[\wt{Q}=\sum_{\substack{1\leq n_1,n_2\leq M\\n_1+n_2>M}}
	\cN(f^{(n_1)},f^{(n_2)}).\]
	Thus every term in $\wt Q$ contains at least $M+1$ copies of the small
	packet. We first convert the Taylor estimates into uniform weighted bounds
	for $f^{\#}$ and $\wt Q$.
	
	\begin{prop}[Truncated expansion and truncation source]\label{prop:truncated-source}
		The following inequalities hold on $[10,\infty)$
		\begin{equation}\label{96}
			\nnr{\<\na_{\fz,\fv}\>^{20}e^{5\<\f{\na_{\fz,\fv}}{\mathsf{A}_k}\>^\sss+\<\na_{\mrv}\>^\sss}\rr{v^{\alpha}f^{\#}(t)}}_{\lta L^2_{z,v}}\leq \mathsf{A}_k^{-\mathtt{N}_3+24},
		\end{equation}
		\begin{equation}\label{97}
			\int_{10}^\infty\nnr{\<\na_{\fz,\fv}\>^{16}e^{5\<\f{\na_{\fz,\fv}}{\mathsf{A}_k}\>^\sss+\<\na_{\mrv}\>^\sss}\rr{v^{\alpha}\wt{Q}(\tau)}}_{\lta L^2_{z,v}}d\tau\leq \mathsf{A}_k^{-\f{\mathtt{N}_3}{2}\mathsf{A}_k^{\mathtt{N}_1}}.
		\end{equation}
	\end{prop}
	\begin{proof}
		We first show that the localization estimate for $f^{(n)}$ similar to \eqref{fnlocbac} holds on $[\mathfrak{T}_k,\infty)$. By \eqref{fnqlocfor} and \eqref{fnqlocbac} the following holds on $[\mathfrak{T}_k,\infty)$
		\[\ba
		&\nnr{\<\na_{\fz,\fv}\>^3e^{0.14\<\f{\na_{\fz,\fv}}{k^\f13}\>^\sss+\<\na_{\mrv}\>^\sss}\rr{v^{\alpha}f^{(n)}(t)}}_{\ltaa L^2_{z,v}}\leq\sum_{|q|\leq n}\nnr{\<\na_{\fz,\fv}\>^3e^{0.14\<\f{\na_\fz-qk,\na_\fv-q\mathsf{A}_k}{k^\f13}\>^\sss+\<\na_{\mrv}\>^\sss}\rr{v^{\alpha}f^{(n,q)}(t)}}_{\ltaa L^2_{z,v}}\\
		&\times e^{0.14\<\f{qk,q\mathsf{A}_k}{k^\f13}\>^\sss}\leq e^{0.52\nr{\f{n\mathsf{A}_k}{k^\f13}}^\sss}.
		\ea\]
		Thus, by \eqref{taylor-short-time} and \eqref{10T1}, the following holds on $[10,\infty)$:
		\begin{equation}\label{yup}
			\begin{aligned}
				&\nnr{\<\na_{\fz,\fv}\>^{20}e^{5\<\f{\na_{\fz,\fv}}{\mathsf{A}_k}\>^\sss+\<\na_{\mrv}\>^\sss}\rr{v^{\alpha}f^{(n)}(t)}}_{\ltaa L^2_{z,v}} \\
				&\quad\leq\nnr{\<\upmu,\upxi\>^{20}e^{5\<\f{\upmu,\upxi}{\mathsf{A}_k}\>^\sss}e^{\<\mrxi\>^\sss}\pa_\xi^\alpha\wh{f^{(n)}}(t,\upmu\bfe,\xi)\vv{1}_{|\upmu,\upxi|\leq 100n^\f1\sss\mathsf{A}_k}}_{\ltaa L^2_{\upmu,\xi}} \\
				&\qquad +\nnr{\<\upmu,\upxi\>^{20}e^{5\<\f{\upmu,\upxi}{\mathsf{A}_k}\>^\sss}e^{\<\mrxi\>^\sss}\pa_\xi^\alpha\wh{f^{(n)}}(t,\upmu\bfe,\xi)\vv{1}_{|\upmu,\upxi|\geq 100n^\f1\sss\mathsf{A}_k}}_{\ltaa L^2_{\upmu,\xi}} \\
				&\quad\leq 10^{50}e^{500n}\nnr{e^{\<\na_{\mrv}\>^\sss}v^{\alpha}f^{(n)}(t)}_{\ltaa L^2_{z,v}}\mathsf{A}_k^{20(\f{\mathtt{N}_1}{\sss}+1)} \\
				&\qquad +e^{-n\nr{\f{\mathsf{A}_k}{k^\f13}}^\sss} \\
				&\quad\leq e^{500n}\mathsf{A}_k^{\rr{30\mathtt{N}_1-\mathtt{N}_3+2}n}\mathsf{A}_k^{5\mathtt{N}_1+20}+e^{-\nr{\f{n\mathsf{A}_k}{k^\f13}}^\sss} \\
				&\quad\leq \mathsf{A}_k^{-0.6\mathtt{N}_3 n}.
			\end{aligned}
		\end{equation}
		For $f^{(1)}$, by \eqref{first-order-envelope} we have
		\[\nnr{\<\na_{\fz,\fv}\>^{20}e^{5\<\f{\na_{\fz,\fv}}{\mathsf{A}_k}\>^\sss+\<\na_{\mrv}\>^\sss}\rr{v^{\alpha}f^{(1)}(t)}}_{\ltaa L^2_{z,v}}\leq \mathsf{A}_k^{-\mathtt{N}_3+23.5}.\]
		Thus we get \eqref{96}. Then we prove \eqref{97}.
		\[\ba
		&\nnr{\<\na_{\fz,\fv}\>^{16}e^{5\<\f{\na_{\fz,\fv}}{\mathsf{A}_k}\>^\sss+\<\na_{\mrv}\>^\sss}\rr{v^{\alpha}\cN(f^{(n_1)},f^{(n_2)})}}_{\ltaa L^2_{z,v}}\leq\sum_\upnu\nr{\sF\rr{\<\na_{\fz,\fv}\>^{17}e^{5\<\f{\na_{\fz,\fv}}{\mathsf{A}_k}\>^\sss}\rho^{(n_1)}}(\tau,\upnu\bfe)}\\
		&\times\nnr{ \<\na_{\fz,\fv}\>^{17}e^{5\<\f{\na_{\fz,\fv}}{\mathsf{A}_k}\>^\sss+\<\na_{\mrv}\>^\sss}\rr{v^{\alpha}f^{(n_2)}(\tau)}}_{\ltaa L^2_{z,v}}\leq\rr{\sum_\upnu\nr{\sF\rr{\<\na_{\fz,\fv}\>^{17}e^{5\<\f{\na_{\fz,\fv}}{\mathsf{A}_k}\>^\sss}\rho^{(n_1)}}(\tau,\upnu\bfe)}\vv{1}_{|\upnu,\upnu\tau|\leq100n^\f1\sss\mathsf{A}_k}\right.\\
		&\left.+e^{-\f1{16}\<\f{\tau}{k^\f13}\>^\sss}e^{-\nr{\f{n_1\mathsf{A}_k}{k^\f13}}^\sss}}\nnr{ \<\na_{\fz,\fv}\>^{17}e^{5\<\f{\na_{\fz,\fv}}{\mathsf{A}_k}\>^\sss+\<\na_{\mrv}\>^\sss}\rr{v^{\alpha}f^{(n_2)}(\tau)}}_{\ltaa L^2_{z,v}}.
		\ea\]
		Thus
		\[
		\ba&\int_{10}^\infty\nnr{\<\na_{\fz,\fv}\>^{16}e^{5\<\f{\na_{\fz,\fv}}{\mathsf{A}_k}\>^\sss+\<\na_{\mrv}\>^\sss}\rr{v^{\alpha}\cN(f^{(n_1)},f^{(n_2)})}}_{\ltaa L^2_{z,v}}d\tau\leq\left(\sum_\upnu\int_{10}^\infty\nr{\sF\rr{\<\na_{\fz,\fv}\>^{17}e^{5\<\f{\na_{\fz,\fv}}{\mathsf{A}_k}\>^\sss}\rho^{(n_1)}}(\tau,\upnu\bfe)}\right.\\
		&\left.\vv{1}_{|\upnu,\upnu\tau|\leq100n^\f1\sss\mathsf{A}_k}d\tau+e^{-\f12\nr{\f{n_1\mathsf{A}_k}{k^\f13}}^\sss}\right)\mathsf{A}_k^{-0.6\mathtt{N}_3 n_2}
		\leq\rr{\rr{n^\f1\sss\mathsf{A}_k}^2\mathsf{A}_k^{-0.6\mathtt{N}_3 n_1}+e^{-\f12\nr{\f{n_1\mathsf{A}_k}{k^\f13}}^\sss}}\mathsf{A}_k^{-0.6\mathtt{N}_3 n_2}\leq \mathsf{A}_k^{-0.52\mathtt{N}_3(n_1+n_2)}.
		\ea\]
		Thus \eqref{97} holds.
	\end{proof}
	
	We estimate $\wt f$ in three regimes. On $[10,\mathsf{T}_k]$ and $[\mathsf{T}_k,\mathsf{A}_k^5]$, a
	short-time energy estimate is sufficient. After $\mathsf{A}_k^5$, the additional
	time-decaying weight used in Section 2 is needed to recover an integrable
	density estimate.
	
	\begin{prop}[Short and intermediate times]\label{ftildeshort}
		On $[10,\mathsf{T}_k]$,
		\[\nnr{\<\na_{\fz,\fv}\>^{16}e^{5\rr{1-\f{1}{t}}\<\f{\na_{\fz,\fv}}{\mathsf{A}_k}\>^\sss+\<\na_{\mrv}\>^\sss}\rr{v^{\alpha}\wt{f}(t)}}_{\lta L^2_{z,v}}\leq e^{\mathsf{T}_k^2}\mathsf{A}_k^{-\f{\mathtt{N}_3}{2}\mathsf{A}_k^{\mathtt{N}_1}}.\]
		On $[\mathsf{T}_k,\mathsf{A}_k^5]$,
		\begin{equation}\label{700}
			\nnr{\<\na_{\fz,\fv}\>^{16}e^{5\rr{\f45+\f{1}{t}}\<\f{\na_{\fz,\fv}}{\mathsf{A}_k}\>^\sss+\<\na_{\mrv}\>^\sss}\rr{v^{\alpha}\wt{f}(t)}}_{\lta L^2_{z,v}}\leq e^{\mathsf{A}_k^{10}}\mathsf{A}_k^{-\f{\mathtt{N}_3}{2}\mathsf{A}_k^{\mathtt{N}_1}}.
		\end{equation}
	\end{prop}
	\begin{proof}
		We first treat $\mathsf{T}_k\leq t\leq \mathsf{A}_k^5$. Set
		\[\fE(t,\mu,\xi)=\<\upmu,\upxi\>^{16}
		e^{5\rr{\f45+\f{1}{t}}\<\f{\upmu,\upxi}{\mathsf{A}_k}\>^\sss+\<\mrxi\>^\sss}.\]
		A standard energy estimate, analogous to the proof
		of Proposition \ref{shortenergy}, gives
		\[\ba
		&\f12\f{d}{dt}\nnr{\fE(t,\na)\rr{v^{\alpha}\wt{f}(t)}}_{\lta L^2_{z,v}}^2
		+\f{5}{t^2}\nnr{\<\f{\na_{\fz,\fv}}{\mathsf{A}_k}\>^{\f\sss2}\fE(t,\na)\rr{v^{\alpha}\wt{f}(t)}}_{\lta L^2_{z,v}}^2
		\\=&-\sum_{|\alpha|\leq2}\sum_{\upmu}\int\wh{\wt{\rho}}(t,\upmu\bfe)\f1\upmu\pa_{\xi}^\alpha\rr{ (\upxi-\upmu t)\wh{f^0}((\xi-\upmu t)\bfe)}\fE(t,\upmu\bfe,\xi)
		\ov{\sF\rr{\fE(t,\na)v^{\alpha}\wt{f}}(t,\upmu\bfe,\xi)}d\xi\\
		&-\sum_{|\alpha|\leq2}\sum_{\upmu,\upnu}\int\wh{\wt{\rho}}(t,\upnu\bfe)\f{\upxi-\upmu t}\upnu\pa_{\xi}^\alpha\sF\rr{f^{\#}+\udl{f}}(t,(\upmu-\upnu)\bfe,\xi-\upnu t\bfe)\fE(t,\upmu\bfe,\xi)\\&\times\ov{\sF\rr{\fE(t,\na)v^{\alpha}\wt{f}}(t,\upmu\bfe,\xi)}d\xi
		-\sum_{|\alpha|\leq2}\sum_{\upmu,\upnu}\int\sF\rr{\rho^{\#}+\udl{\rho}}(t,\upnu\bfe)\f{\upxi-\upmu t}\upnu\pa_{\xi}^\alpha\wh{\wt{f}}(t,(\upmu-\upnu)\bfe,\xi-\upnu t\bfe)\\&\times\rr{\fE(t,\upmu\bfe,\xi)-\fE(t,(\upmu-\upnu)\bfe,\xi-\upnu t\bfe)}\ov{\sF\rr{\fE(t,\na)v^{\alpha}\wt{f}}(t,\upmu\bfe,\xi)}d\xi-\sum_{|\alpha|\leq2}\sum_{\upmu,\upnu}\int\wh{\wt{\rho}}(t,\upnu\bfe)\\
		&\times\f{\upxi-\upmu t}\upnu\pa_{\xi}^\alpha\wh{\wt{f}}(t,(\upmu-\upnu)\bfe,\xi-\upnu t\bfe)\rr{\fE(t,\upmu\bfe,\xi)-\fE(t,(\upmu-\upnu)\bfe,\xi-\upnu t\bfe)}\ov{\sF\rr{\fE(t,\na)v^{\alpha}\wt{f}}(t,\upmu\bfe,\xi)}d\xi\\
		-&\sum_{\substack{|\alpha|\leq2\\\upalpha>0}}C_\alpha\sum_{\upmu,\upnu}\int\wh{\wt{\rho}}(t,\upnu\bfe)\f{1}\upnu\pa_{\xi}^{\alpha-\bfe}\sF\rr{f^{\#}+\udl{f}}(t,(\upmu-\upnu)\bfe,\xi-\upnu t\bfe)\fE(t,\upmu\bfe,\xi)\ov{\sF\rr{\fE(t,\na)v^\alpha\wt{f}}(t,\upmu\bfe,\xi)}d\xi\\
		-&\sum_{\substack{|\alpha|\leq2\\\upalpha>0}}C_\alpha\sum_{\upmu,\upnu}\int\sF\rr{\rho^{\#}+\udl{\rho}+\wt{\rho}}(t,\upnu\bfe)\f1\upnu\pa_{\xi}^{\alpha-\bfe}\wh{\wt{f}}(t,(\upmu-\upnu)\bfe,\xi-\upnu t\bfe)\fE(t,\upmu\bfe,\xi)\ov{\sF\rr{\fE(t,\na)v^\alpha\wt{f}}(t,\upmu\bfe,\xi)}d\xi\\
		+&\sum_{|\alpha|\leq2}\sum_\upmu\int\sF\rr{\fE(t,\na)v^{\alpha}\wt{Q}}(t,\upmu\bfe,\xi)\ov{\sF\rr{\fE(t,\na)v^{\alpha}\wt{f}}(t,\upmu\bfe,\xi)}d\xi:=\sum_{i=1}^{7}\mathscr{K}_i.
		\ea\]
		We estimate each term individually. The methods used in these estimates are almost the same as those in the proof of Proposition \ref{shortenergy}.
		\[\mathscr{K}_1\ls_\theta \nnr{\fE(t,\na)\rr{v^{\alpha}\wt{f}(t)}}_{\ltaa L^2_{z,v}}^2.\]
		\[\mathscr{K}_2\ls\nnr{\fE(t,\na)\rr{v^{\alpha}\wt{f}(t)}}^2_{\ltaa L^2_{z,v}}\nnr{\<\na_{\fz,\fv}\>^2\fE(t,\na)v^{\alpha}\rr{f^{\#}+\udl{f}}(t)}_{\ltaa L^2_{z,v}}t\ls\eps_0t\nnr{\fE(t,\na)\rr{v^{\alpha}\wt{f}(t)}}_{\ltaa L^2_{z,v}}^2.\]
		For $\mathscr{K}_3$ we use a standard commutator estimate and get
		\[\mathscr{K}_3\ls\nnr{\<\na_{\fz,\fv}\>^4\fE(t,\na_{\fz,\fv})\rr{\rho^{\#}+\udl{\rho}}(t)}_{L^2_x}t^{-2}\nnr{\<\f{\na_{\fz,\fv}}{\mathsf{A}_k}\>^{\f\sss2}\fE(t,\na)\rr{v^{\alpha}\wt{f}(t)}}_{\lta L^2_{z,v}}^2.\]
		For the nonlinear term $\mathscr{K}_4$ we use the standard paraproduct decomposition and split it into two parts
		\[\ba&\mathscr{K}_4=-\sum_{|\alpha|\leq2}\sum_{\upmu,\upnu}\int\wh{\wt{\rho}}(t,\upnu\bfe)
		\f{\upxi-\upmu t}\upnu\rr{\vv{1}_{\<\upnu,\upnu t\>\leq\<\upmu-\upnu,\upxi-\upnu t\>}+\vv{1}_{\<\upnu,\upnu t\>\geq\<\upmu-\upnu,\upxi-\upnu t\>}}\pa_{\xi}^\alpha\wh{\wt{f}}(t,(\upmu-\upnu)\bfe,\xi-\upnu t\bfe)\\&\times\rr{\fE(t,\upmu\bfe,\xi)-\fE(t,(\upmu-\upnu)\bfe,\xi-\upnu t\bfe)}\ov{\sF\rr{\fE(t,\na)v^{\alpha}\wt{f}}(t,\upmu\bfe,\xi)}d\xi=\mathscr{K}_{4,LH}+\mathscr{K}_{4,HL}.\ea\]
		We first estimate $\mathscr{K}_{4,LH}$. We note that on the region $\<\upnu,\upnu t\>\leq\<\upmu-\upnu,\upxi-\upnu t\>$ it holds that $\<\upmu,\upxi\>\ls\<\upmu-\upnu,\upxi-\upnu t\>$.
		\[\ba
		&\mathscr{K}_{4,LH}\ls\nnr{\<\na_{\fz,\fv}\>^4e^{5\rr{\f45+\f{1}{t}}\<\f{\na_{\fz,\fv}}{\mathsf{A}_k}\>^\sss}\wt{\rho}(t)}_{L^2_x}t^{-2}\nnr{\<\f{\na_{\fz,\fv}}{\mathsf{A}_k}\>^{\f\sss2}\fE(t,\na)\rr{v^{\alpha}\wt{f}(t)}}_{\ltaa L^2_{z,v}}^2\\
		&\ls t^{-2}\nnr{\fE(t,\na)\rr{v^{\alpha}\wt{f}(t)}}_{\ltaa L^2_{z,v}}\nnr{\<\f{\na_{\fz,\fv}}{\mathsf{A}_k}\>^{\f\sss2}\fE(t,\na)\rr{v^{\alpha}\wt{f}(t)}}_{\ltaa L^2_{z,v}}^2.
		\ea\]
		Then we estimate $\mathscr{K}_{4,HL}$. We note that on the region $\<\upnu,\upnu t\>\geq\<\upmu-\upnu,\upxi-\upnu t\>$ it holds that $\<\upmu,\upxi\>\ls\<\upnu,\upnu t\>$.
		\[\ba
		&\mathscr{K}_{4,HL}\ls\nnr{\fE(t,\na)\wt{\rho}(t)}_{L^2_x}\nnr{\<\na_{\fz,\fv}\>^2e^{5\rr{\f45+\f{1}{t}}\<\f{\na_{\fz,\fv}}{\mathsf{A}_k}\>^\sss+\<\na_{\mrv}\>^\sss}\rr{v^{\alpha}\wt{f}(t)}}_{\ltaa L^2_{z,v}}\nnr{\fE(t,\na)\rr{v^{\alpha}\wt{f}(t)}}_{\ltaa L^2_{z,v}}t\\
		&\ls t\nnr{\fE(t,\na)\rr{v^{\alpha}\wt{f}(t)}}_{\ltaa L^2_{z,v}}^3.
		\ea\]
		The remaining terms can be estimated in the same way as before. Combining all the estimates we get
		\begin{equation}\label{800}
			\begin{aligned}
				&\f12\f{d}{dt}\nnr{\fE(t,\na)\rr{v^{\alpha}\wt{f}(t)}}_{\ltaa L^2_{z,v}}^2
		+\f{1}{t^2}\rr{4-10^2\nnr{\fE(t,\na)\rr{v^{\alpha}\wt{f}(t)}}_{\ltaa L^2_{z,v}}}\nnr{\<\f{\na_{\fz,\fv}}{\mathsf{A}_k}\>^{\f\sss2}\fE(t,\na)\rr{v^{\alpha}\wt{f}(t)}}_{\ltaa L^2_{z,v}}^2 \\
				&\leq\rr{C_\theta+\eps_0^{0.9}t}\nnr{\fE(t,\na)\rr{v^{\alpha}\wt{f}(t)}}_{\ltaa L^2_{z,v}}^2+10^2t\nnr{\fE(t,\na)\rr{v^{\alpha}\wt{f}(t)}}_{\ltaa L^2_{z,v}}^3 \\
				&\quad +\nnr{\fE(t,\na)\rr{v^{\alpha}\wt{Q}(t)}}_{\ltaa L^2_{z,v}}\nnr{\fE(t,\na)\rr{v^{\alpha}\wt{f}(t)}}_{\ltaa L^2_{z,v}}
			\end{aligned}
		\end{equation}
		In the following we apply a continuity argument to prove \eqref{700}. We assume that \eqref{700} holds true on $[\mathsf{T}_k,t^\dagger]$. Then by \eqref{800} the following holds on $[\mathsf{T}_k,t^\dagger]$
		\[\f{d}{dt}\nnr{\fE(t,\na)\rr{v^{\alpha}\wt{f}(t)}}_{\ltaa L^2_{z,v}}\leq\rr{C_\theta+\eps_0^{0.8}t}\nnr{\fE(t,\na)\rr{v^{\alpha}\wt{f}(t)}}_{\ltaa L^2_{z,v}}+\nnr{\fE(t,\na)\rr{v^{\alpha}\wt{Q}(t)}}_{\ltaa L^2_{z,v}}\]
		Solving the ordinary differential equation we get
		\[e^{-\rr{C_\theta t+\eps_0^{0.8}t^2}}\nnr{\fE(t,\na)\rr{v^{\alpha}\wt{f}(t)}}_{\ltaa L^2_{z,v}}\leq\int_{\mathsf{T}_k}^te^{-\rr{C_\theta \tau+\eps_0^{0.8}\tau^2}}\nnr{\fE(\tau,\na)\rr{v^{\alpha}\wt{Q}(\tau)}}_{\ltaa L^2_{z,v}}d\tau.\]
		Thus on $[\mathsf{T}_k,t^\dagger]$ it holds that
		\[\nnr{\fE(t,\na)\rr{v^{\alpha}\wt{f}(t)}}_{\ltaa L^2_{z,v}}\leq e^{C_\theta \rr{t-\mathsf{T}_k}+\eps_0^{0.8}\rr{t^2-\mathsf{T}_k^2}}\int_{\mathsf{T}_k}^t\nnr{\fE(\tau,\na)\rr{v^{\alpha}\wt{Q}(\tau)}}_{\ltaa  L^2_{z,v}}d\tau\leq e^{\eps_0^{0.7}\mathsf{A}_k^{10}}\mathsf{A}_k^{-\f{\mathtt{N}_3}{2}\mathsf{A}_k^{\mathtt{N}_1}}.\]
		Thus \eqref{700} holds on a longer time interval. By a standard continuity argument \eqref{700} holds on $[\mathsf{T}_k,\mathsf{A}_k^5]$.

		For $10\leq t\leq \mathsf{T}_k$, use instead
		\[\fE_-(t,\mu,\xi)=\<\upmu,\upxi\>^{16}
		e^{5\rr{1-\f1t}\<\f{\upmu,\upxi}{\mathsf{A}_k}\>^\sss+\<\mrxi\>^\sss}\]
		and integrate the remainder equation backward from $\mathsf{T}_k$, where $\wt f(\mathsf{T}_k)=0$.
		The seven terms above obey the same commutator and paraproduct bounds, with
		the time orientation reversed. Proposition \ref{shortenergy} and \eqref{97}
		therefore give
		\[\nnr{\fE_-(t,\na)v^\alpha\wt f(t)}_{\lta L^2_{z,v}}
		\leq e^{C_\theta(\mathsf{T}_k-t)+\eps_0^{0.8}(\mathsf{T}_k^2-t^2)}
		\int_t^{\mathsf{T}_k}\nnr{\fE_-(\tau,\na)v^\alpha\wt Q(\tau)}_{\lta L^2_{z,v}}d\tau
		\leq e^{\mathsf{T}_k^2}\mathsf{A}_k^{-\f{\mathtt{N}_3}{2}\mathsf{A}_k^{\mathtt{N}_1}},\]
		which proves the first estimate.
	\end{proof}
	
We now turn to the long-time regime. The key is to combine a Volterra estimate for the density with a weighted energy estimate for the profile.
	
	\begin{prop}[Long-time remainder bound]\label{escs}
		On $[\mathsf{A}_k^5,\infty)$ it holds that 
		\begin{equation}\label{jsp}
			\nnr{\<\na_{\fz,\fv}\>^{14}e^{\rr{\f{\mathsf{T}_k}{t}}^{0.01}\<\na_{\fz,\fv}\>^{0.34}}e^{4\rr{\f45+\f{1}{t}}\<\f{\na_{\fz,\fv}}{\mathsf{A}_k}\>^\sss+\<\na_{\mrv}\>^\sss}\rr{v^{\alpha}\wt{f}(t)}}_{\lta L^2_{z,v}}\leq e^{2\mathsf{A}_k^{10}}\mathsf{A}_k^{-\f{\mathtt{N}_3}{2}\mathsf{A}_k^{\mathtt{N}_1}}.
		\end{equation}
	\end{prop}
	\begin{proof}
		We combine the Volterra and energy estimates from Proposition
		\ref{longenergyl}. The proof has four steps.
		
		\noindent\emph{Step 1: preparation of the long-time weights.}
		By \eqref{700} we have
		\[\nnr{\<\na_{\fz,\fv}\>^{16}e^{\rr{\f{\mathsf{T}_k}{\mathsf{A}_k^5}}^{0.01}\<\na_{\fz,\fv}\>^{0.34}}e^{4\rr{\f45+\f{1}{\mathsf{A}_k^5}}\<\f{\na_{\fz,\fv}}{\mathsf{A}_k}\>^\sss+\<\na_{\mrv}\>^\sss}\rr{v^{\alpha}\wt{f}(\mathsf{A}_k^5)}}_{\lta L^2_{z,v}}\leq e^{\mathsf{A}_k+\mathsf{A}_k^{10}}\mathsf{A}_k^{-\f{\mathtt{N}_3}{2}\mathsf{A}_k^{\mathtt{N}_1}},\]
		where we use $\rr{\f{\mathsf{T}_k}{\mathsf{A}_k^5}}^{0.01}\<\upmu,\upxi\>^{0.34}\leq\f45\<\f{\upmu,\upxi}{\mathsf{A}_k}\>^\sss+\mathsf{A}_k$. Similarly by \eqref{97}
		\[\int_{\mathsf{A}_k^5}^\infty\nnr{\<\na_{\fz,\fv}\>^{16}e^{\rr{\f{\mathsf{T}_k}{\tau}}^{0.01}\<\na_{\fz,\fv}\>^{0.34}}e^{4\rr{\f45+\f{1}{\tau}}\<\f{\na_{\fz,\fv}}{\mathsf{A}_k}\>^\sss+\<\na_{\mrv}\>^\sss}\rr{v^{\alpha}\wt{Q}(\tau)}}_{\lta L^2_{z,v}}d\tau\leq e^{\mathsf{A}_k}\mathsf{A}_k^{-\f{\mathtt{N}_3}{2}\mathsf{A}_k^{\mathtt{N}_1}}.\]
		In the following we set $\fM_a(t,\mu,\xi)=\<\upmu,\upxi\>^{13+a}e^{\rr{\f{\mathsf{T}_k}{t}}^{0.01}\<\upmu,\upxi\>^{0.34}}e^{4\rr{\f45+\f{1}{t}}\<\f{\upmu,\upxi}{\mathsf{A}_k}\>^\sss+\<\mrxi\>^\sss}$, $\fM(t,\mu,\xi)=\fM_0(t,\mu,\xi)$.
		
		\noindent\emph{Step 2: density estimate.}
		The density $\wh{\wt{\rho}}(t,\upmu\bfe)$ satisfies
		\[\ba
		&\wh{\wt{\rho}}(t,\upmu\bfe)=H(t,\upmu\bfe)-\int_{\mathsf{A}_k^5}^tH(\tau,\upmu\bfe)R^0(t-\tau,\upmu\bfe)d\tau,\quad	H(t,\upmu\bfe)=\wh{\wt{f}}(\mathsf{A}_k^5,\upmu\bfe,\upmu t\bfe)+\int_{\mathsf{A}_k^5}^t\wh{\wt{Q}}(\tau,\upmu\bfe,\upmu t\bfe)d\tau\\
		&-\sum_\upnu\int_{\mathsf{A}_k^5}^t\wh{\wt{\rho}}(\tau,\upnu\bfe)\f{\upmu(t-\tau)}{\upnu}\sF\rr{f^{\#}+\udl{f}}(\tau,(\upmu-\upnu)\bfe,(\upmu t-\upnu\tau)\bfe)d\tau-\sum_\upnu\int_{\mathsf{A}_k^5}^t\sF\rr{\rho^{\#}+\udl{\rho}}(\tau,\upnu\bfe)\\
		&\times\f{\upmu(t-\tau)}{\upnu}\wh{\wt{f}}(\tau,(\upmu-\upnu)\bfe,(\upmu t-\upnu\tau)\bfe)d\tau-\sum_\upnu\int_{\mathsf{A}_k^5}^t\wh{\wt{\rho}}(\tau,\upnu\bfe)\f{\upmu(t-\tau)}{\upnu}\wh{\wt{f}}(\tau,(\upmu-\upnu)\bfe,(\upmu t-\upnu\tau)\bfe)d\tau.
		\ea\]
		By \eqref{eRz},
		\[\nnr{\<t\>^3\fM(t,\upmu\bfe,\upmu t\bfe)\wh{\wt{\rho}}(t,\upmu\bfe)}_{L^2_t([\mathsf{A}_k^5,\tilde T];L^2_\upmu)}\ls_\theta\nnr{\<t\>^3\fM(t,\upmu\bfe,\upmu t\bfe)H(t,\upmu\bfe)}_{L^2_t([\mathsf{A}_k^5,\tilde T];L^2_\upmu)}.\]
		Thus it remains to estimate $H$. Denote its five terms by
		$\mathscr{P}_1,\ldots,\mathscr{P}_5$. We estimate them
		with the same kernel bounds as in Proposition \ref{longenergyl}.
		\[\nnr{\<t\>^3\fM(t,\upmu\bfe,\upmu t\bfe)\mathscr{P}_1(t,\upmu\bfe)}_{L^2_t([\mathsf{A}_k^5,\tilde T];L^2_\upmu)}\leq\rr{\sum_\upmu\int_{\mathsf{A}_k^5}^{\tilde T}\nr{\sF\rr{\fM_{3}(t,\na)\wt{f}}(\mathsf{A}_k^5,\upmu\bfe,\upmu t\bfe)}^2dt}^\f12\leq e^{\mathsf{A}_k+\mathsf{A}_k^{10}}\mathsf{A}_k^{-\f{\mathtt{N}_3}{2}\mathsf{A}_k^{\mathtt{N}_1}}.\]
		\[\nnr{\<t\>^3\fM(t,\upmu\bfe,\upmu t\bfe)\mathscr{P}_2(t,\upmu\bfe)}_{L^2_t([\mathsf{A}_k^5,\tilde T];L^2_\upmu)}
		\leq\int_{\mathsf{A}_k^5}^{\tilde T}\rr{\sum_\upmu\int_{\mathsf{A}_k^5}^{\tilde T} \nr{\sF\rr{\fM_3(\tau,\na)\wt{Q}}(\tau,\upmu\bfe,\upmu t\bfe)}^2dt }^\f12d\tau \leq e^{\mathsf{A}_k}\mathsf{A}_k^{-\f{\mathtt{N}_3}{2}\mathsf{A}_k^{\mathtt{N}_1}}.\]
		For $\mathscr{P}_3$ we note that 
		\[\rr{\f{\mathsf{T}_k}{\tau}}^{0.01}\<\upmu,\upmu t\>^{0.34}\leq\rr{\f{\mathsf{T}_k}{\tau}}^{0.01}\<\upnu,\upnu\tau\>^{0.34}+\f14\<\f{\upmu-\upnu,\upmu t-\upnu\tau}{\mathsf{A}_k}\>^\sss\] 
		since $|\upnu,\upnu\tau|\geq \mathsf{A}_k^5$. Thus it holds that
		\[\ba
		&\<t\>^3\fM(t,\upmu\bfe,\upmu t\bfe)\mathscr{P}_3(t,\upmu\bfe)
		\leq\sup_{\mathsf{A}_k^5\leq\tau\leq t}\nnr{\<\na_{\fz,\fv}\>^{20}e^{5\<\f{\na_{\fz,\fv}}{\mathsf{A}_k}\>^\sss+\<\na_{\mrv}\>^\sss}v^{\alpha}\rr{f^{\#}+\udl{f}}(\tau)}_{\ltaa  L^2_{z,v}}\\
		&\times\sum_\upnu\int_{\mathsf{A}_k^5}^t\<\tau\>^3\nr{\sF\rr{\fM(\tau,\na)\wt{\rho}}(\tau,\upnu\bfe)}\f{\tau\vv{1}_{\upnu\neq\upmu}+1}{|\upnu|}\<\upmu-\upnu,\upmu t-\upnu\tau\>^{-6}e^{\rr{\rr{\f{\mathsf{T}_k}{t}}^{0.01}-\rr{\f{\mathsf{T}_k}{\tau}}^{0.01}}\<\upmu,\upmu t\>^{0.34}}\f{\<t\>^3}{\<\tau\>^3}d\tau.
		\ea\]
		By the Schur test and Lemma \ref{auxi} it holds that
		\[\nnr{\<t\>^3\fM(t,\upmu\bfe,\upmu t\bfe)\mathscr{P}_3(t,\upmu\bfe)}_{L^2_t([\mathsf{A}_k^5,\tilde T];L^2_\upmu)}\leq\eps_{0}^{0.94}\nnr{\<\tau\>^3\fM(\tau,\upnu,\upnu\tau)\wh{\wt{\rho}}(\tau,\upnu\bfe)}_{L^2_\tau([\mathsf{A}_k^5,\tilde T];L^2_\upnu)}.\]	
		For $\mathscr{P}_4$ we note that $|\upnu,\upnu\tau|\geq \mathsf{A}_k^5$ implies that 
		\[\rr{\f{\mathsf{T}_k}{\tau}}^{0.01}\<\upmu,\upmu t\>^{0.34}\leq\f14\rr{\f{\mathsf{T}_k}{\tau}}^{0.01}\<\f{\upnu,\upnu\tau}{\mathsf{A}_k}\>^{\sss}+\rr{\f{\mathsf{T}_k}{\tau}}^{0.01}\<\upmu-\upnu,\upmu t-\upnu\tau\>^{0.34}.\]
		Thus
		\[\ba
		&\<t\>^3\fM(t,\upmu\bfe,\upmu t\bfe)\mathscr{P}_4(t,\upmu\bfe)\leq\sup_{\mathsf{A}_k^5\leq\tau\leq t}\nnr{\<\na_{\fz,\fv}\>^{20}e^{5\<\f{\na_{\fz,\fv}}{\mathsf{A}_k}\>^\sss}\rr{\rho^{\#}+\udl{\rho}}(\tau)}_{L^2_{x}}\\
		&\sum_\upnu\int_{\mathsf{A}_k^5}^t\<\upnu,\upnu\tau\>^{-3}\nr{\sF\rr{\fM_1(\tau,\na)\wt{f}}(\tau,(\upmu-\upnu)\bfe,(\upmu t-\upnu\tau)\bfe)}d\tau.
		\ea\]
		By Minkowski inequality we have
		\[\nnr{\<t\>^3\fM(t,\upmu\bfe,\upmu t\bfe)\mathscr{P}_4(t,\upmu\bfe)}_{L^2_t([\mathsf{A}_k^5,\tilde T];L^2_\upmu)}\leq\eps_0^{0.94}\sup_{\mathsf{A}_k^5\leq\tau\leq\tilde T}\nnr{\fM_1(\tau,\na)\wt{f}}_{L^2_{z,v}}.\]
		For the nonlinear term $\mathscr{P}_5$ we apply the paraproduct decomposition 
		\[\ba&\mathscr{P}_5=-\sum_\upnu\int_{\mathsf{A}_k^5}^t\wh{\wt{\rho}}(\tau,\upnu\bfe)\f{\upmu(t-\tau)}{\upnu}\rr{\vv{1}_{\<\upnu,\upnu\tau\>\leq\<\upmu-\upnu,\upmu t-\upnu\tau\>}+\vv{1}_{\<\upnu,\upnu\tau\>\geq\<\upmu-\upnu,\upmu t-\upnu\tau\>}}\wh{\wt{f}}(\tau,(\upmu-\upnu)\bfe,(\upmu t-\upnu\tau)\bfe)d\tau\\
		&=\mathscr{P}_{5,LH}+\mathscr{P}_{5,HL}.\ea\]
		We first estimate $\mathscr{P}_{5,LH}$. We note that on the region $\<\upnu,\upnu\tau\>\leq\<\upmu-\upnu,\upmu t-\upnu\tau\>$ it holds that $\<\upmu,\upmu t\>\ls\<\upmu-\upnu,\upmu t-\upnu\tau\>$.
		\[\ba
		&\<t\>^3\fM(t,\upmu\bfe,\upmu t\bfe)\mathscr{P}_{5,LH}(t,\upmu\bfe)\ls\int_{\mathsf{A}_k^5}^t\<\tau\>^4\nr{\sF\rr{\fM_{-13}(\tau,\na)\wt{\rho}}(\tau,\upnu\bfe)}\nr{\sF\rr{\fM_1(\tau,\na)\wt{f}}(\tau,(\upmu-\upnu)\bfe,(\upmu t-\upnu\tau)\bfe)}d\tau\\	&\ls\sup_{\mathsf{A}_k^5\leq\tau\leq t}\nnr{\fM(\tau,\na)\rr{v^{\alpha}\wt{f}(\tau)}}_{\ltaa L^2_{z,v}}\sum_\upnu\int_{\mathsf{A}_k^5}^t\<\upnu,\upnu\tau\>^{-3}\nr{\sF\rr{\fM_1(\tau,\na)\wt{f}}(\tau,(\upmu-\upnu)\bfe,(\upmu t-\upnu\tau)\bfe)}d\tau.
		\ea\]
		Thus
		\[\nnr{\<t\>^3\fM(t,\upmu\bfe,\upmu t\bfe)\mathscr{P}_{5,LH}(t,\upmu\bfe)}_{L^2_t([\mathsf{A}_k^5,\tilde T];L^2_\upmu)}\ls\sup_{\mathsf{A}_k^5\leq\tau\leq \tilde T}\nnr{\fM_1(\tau,\na)\rr{v^{\alpha}\wt{f}(\tau)}}^2_{\ltaa L^2_{z,v}}.\]
		Then we estimate $\mathscr{P}_{5,HL}$. On the complementary region
		$\<\upnu,\upnu\tau\>\geq\<\upmu-\upnu,\upmu t-\upnu\tau\>$, it holds that
		$\<\upmu,\upmu t\>\ls\<\upnu,\upnu\tau\>$.
		\[\ba
		&\<t\>^3\fM(t,\upmu\bfe,\upmu t\bfe)\mathscr{P}_{5,HL}(t,\upmu\bfe)\ls\int_{\mathsf{A}_k^5}^t\<\tau\>^3\nr{\sF\rr{\fM(\tau,\na)\wt{\rho}}(\tau,\upnu\bfe)}\nr{\sF\rr{\fM_{-12}(\tau,\na)\wt{f}}(\tau,(\upmu-\upnu)\bfe,(\upmu t-\upnu\tau)\bfe)}\\
		&\times\f{\tau\vv{1}_{\upnu\neq\upmu}+1}{|\upnu|} e^{\rr{\rr{\f{\mathsf{T}_k}{t}}^{0.01}-\rr{\f{\mathsf{T}_k}{\tau}}^{0.01}}\<\upmu,\upmu t\>^{0.34}}\f{\<t\>^3}{\<\tau\>^3}d\tau
		\ls\sup_{\mathsf{A}_k^5\leq\tau\leq t}\nnr{\fM_1(\tau,\na)\rr{v^{\alpha}\wt{f}}(\tau)}_{\ltaa  L^2_{z,v}}\\
		&\times\sum_\upnu\int_{\mathsf{A}_k^5}^t\<\tau\>^3\nr{\sF\rr{\fM(\tau,\na)\wt{\rho}}(\tau,\upnu\bfe)}\f{\tau\vv{1}_{\upnu\neq\upmu}+1}{|\upnu|}\<\upmu-\upnu,\upmu t-\upnu\tau\>^{-6}e^{\rr{\rr{\f{\mathsf{T}_k}{t}}^{0.01}-\rr{\f{\mathsf{T}_k}{\tau}}^{0.01}}\<\upmu,\upmu t\>^{0.34}}\f{\<t\>^3}{\<\tau\>^3}d\tau.
		\ea\]
		By the Schur test and Lemma \ref{auxi} it holds that
		\[\ba\nnr{\<t\>^3\fM(t,\upmu\bfe,\upmu t\bfe)\mathscr{P}_{5,HL}(t,\upmu\bfe)}_{L^2_t([\mathsf{A}_k^5,\tilde T];L^2_\upmu)}\ls\sup_{\mathsf{A}_k^5\leq\tau\leq \tilde T}\nnr{\fM_1(\tau,\na)v^{\alpha}\wt{f}(\tau)}_{\ltaa  L^2_{z,v}}\nnr{\<\tau\>^3\fM(\tau,\na)\wt{\rho}(\tau)}_{L^2_\tau([\mathsf{A}_k^5,\tilde T];L^2_x)}.\ea\]
		Summing up all the estimates we get
		\begin{equation}\label{lerh9o}
			\begin{aligned}
				&\nnr{\<t\>^3\fM(t,\na)\wt{\rho}(t)}_{L^2_t([\mathsf{A}_k^5,\tilde T];L^2_x)} \\
				&\quad\leq C_\theta e^{\mathsf{A}_k+\mathsf{A}_k^{10}}\mathsf{A}_k^{-\f{\mathtt{N}_3}{2}\mathsf{A}_k^{\mathtt{N}_1}} \\
				&\qquad +\eps_0^{0.9}\sup_{\mathsf{A}_k^5\leq\tau\leq\tilde T}\nnr{\fM_1(\tau,\na)\wt{f}}_{L^2_{z,v}} \\
				&\qquad +10^2C_\theta\sup_{\mathsf{A}_k^5\leq\tau\leq \tilde T}\nnr{\fM_1(\tau,\na)\rr{v^{\alpha}\wt{f}(\tau)}}_{\ltaa L^2_{z,v}} \\
				&\qquad\quad\times\Bigl(\sup_{\mathsf{A}_k^5\leq\tau\leq \tilde T}\nnr{\fM_1(\tau,\na)\rr{v^{\alpha}\wt{f}(\tau)}}_{\ltaa L^2_{z,v}} \\
				&\qquad\qquad +\nnr{\<t\>^3\fM(t,\na)\wt{\rho}(t)}_{L^2_t([\mathsf{A}_k^5,\tilde T];L^2_x)}\Bigr).
			\end{aligned}
		\end{equation}
		
		\noindent\emph{Step 3: profile energy estimate.}
		By direct computation,
		\[\ba
		&\f12\f{d}{dt}\nnr{\fM_1(t,\na)v^{\alpha}\wt{f}(t)}_{\ltaa  L^2_{z,v}}^2+\f{0.01}{t}\rr{\f{\mathsf{T}_k}{t}}^{0.01}\nnr{\fM_{1.17}(t,\na)v^{\alpha}\wt{f}(t)}_{\ltaa  L^2_{z,v}}^2
		+\f{4}{t^2}\nnr{\<\f{\na_{\fz,\fv}}{\mathsf{A}_k}\>^{\f\sss2}\fM_{1}(t,\na)v^{\alpha}\wt{f}(t)}_{\ltaa  L^2_{z,v}}^2\\
		&=-\sum_{|\alpha|\leq2}\sum_{\upmu}\int\wh{\wt{\rho}}(t,\upmu\bfe)\f1\upmu\pa_{\xi}^\alpha\rr{ (\upxi-\upmu t)\wh{f^0}((\xi-\upmu t)\bfe)}\fM_1(t,\upmu\bfe,\xi)
		\ov{\sF\rr{\fM_1(t,\na)v^{\alpha}\wt{f}}(t,\upmu\bfe,\xi)}d\xi\\
		&-\sum_{|\alpha|\leq2}\sum_{\upmu,\upnu}\int\wh{\wt{\rho}}(t,\upnu\bfe)\f{\upxi-\upmu t}\upnu\pa_{\xi}^\alpha\sF\rr{f^{\#}+\udl{f}}(t,(\upmu-\upnu)\bfe,\xi-\upnu t\bfe)\fM_1(t,\upmu\bfe,\xi)\ov{\sF\rr{\fM_1(t,\na)v^{\alpha}\wt{f}}(t,\upmu\bfe,\xi)}d\xi\\
		&-\sum_{|\alpha|\leq2}\sum_{\upmu,\upnu}\int\sF\rr{\rho^{\#}+\udl{\rho}}(t,\upnu\bfe)\f{\upxi-\upmu t}\upnu\pa_{\xi}^\alpha\wh{\wt{f}}(t,(\upmu-\upnu)\bfe,\xi-\upnu t\bfe)\rr{\fM_1(t,\upmu\bfe,\xi)-\fM_1(t,(\upmu-\upnu)\bfe,\xi-\upnu t\bfe)}\\&\times\ov{\sF\rr{\fM_1(t,\na)v^{\alpha}\wt{f}}(t,\upmu\bfe,\xi)}d\xi-\sum_{|\alpha|\leq2}\sum_{\upmu,\upnu}\int\wh{\wt{\rho}}(t,\upnu\bfe)\f{\upxi-\upmu t}\upnu\pa_{\xi}^\alpha\wh{\wt{f}}(t,(\upmu-\upnu)\bfe,\xi-\upnu t\bfe)\\&\times\rr{\fM_1(t,\upmu\bfe,\xi)-\fM_1(t,(\upmu-\upnu)\bfe,\xi-\upnu t\bfe)}\ov{\sF\rr{\fM_1(t,\na)v^{\alpha}\wt{f}}(t,\upmu\bfe,\xi)}d\xi\\
		-&\sum_{\substack{|\alpha|\leq2\\\upalpha>0}}C_\alpha\sum_{\upmu,\upnu}\int\wh{\wt{\rho}}(t,\upnu\bfe)\f{1}\upnu\pa_\xi^{\alpha-\bfe}\sF\rr{f^{\#}+\udl{f}}(t,(\upmu-\upnu)\bfe,\xi-\upnu t\bfe)\fM_1(t,\upmu\bfe,\xi)\ov{\sF\rr{\fM_1(t,\na)v^\alpha\wt{f}}(t,\upmu\bfe,\xi)}d\xi\\
		-&\sum_{\substack{|\alpha|\leq2\\\upalpha>0}}C_\alpha\sum_{\upmu,\upnu}\int\sF\rr{\rho^{\#}+\udl{\rho}+\wt{\rho}}(t,\upnu\bfe)\f1\upnu\pa_\xi^{\alpha-\bfe}\wh{\wt{f}}(t,(\upmu-\upnu)\bfe,\xi-\upnu t\bfe)\fM_1(t,\upmu\bfe,\xi)\ov{\sF\rr{\fM_1(t,\na)v^\alpha\wt{f}}(t,\upmu\bfe,\xi)}d\xi\\
		+&\sum_{|\alpha|\leq2}\sum_\upmu\int\sF\rr{\fM_1(t,\na)v^{\alpha}\wt{Q}}(t,\upmu\bfe,\xi)\ov{\sF\rr{\fM_1(t,\na)v^{\alpha}\wt{f}}(t,\upmu\bfe,\xi)}d\xi:=\sum_{i=1}^{7}\mathscr{Q}_i.
		\ea\]
		The estimates below use the same multiplier bounds as in
		Proposition \ref{longenergyl}.
		\[\mathscr{Q}_1\ls_\theta\nnr{\fM(t,\na)\wt{\rho}(t)}_{L^2_x}\nnr{\fM_1(t,\na)\rr{v^{\alpha}\wt{f}(t)}}_{\ltaa  L^2_{z,v}}t\]
		\[\mathscr{Q}_2\ls\nnr{\fM(t,\na)\wt{\rho}(t)}_{L^2_x}\nnr{\<\na_{\fz,\fv}\>^{16}e^{5\<\f{\na_{\fz,\fv}}{\mathsf{A}_k}\>^\sss+\<\na_{\mrv}\>^\sss}v^{\alpha}\rr{f^{\#}+\udl{f}}(t)}_{\ltaa  L^2_{z,v}}\nnr{\fM_1(t,\na)\rr{v^{\alpha}\wt{f}(t)}}_{\ltaa  L^2_{z,v}}t^2.\]
		\[\mathscr{Q}_3\ls t^{-2}\nnr{\<\na_{\fz,\fv}\>^{18}e^{5\<\f{\na_{\fz,\fv}}{\mathsf{A}_k}\>^\sss}\rr{\rho^{\#}+\udl{\rho}}(t)}_{L^2_x}\rr{\nnr{\fM_{1.17}(t,\na)\rr{v^{\alpha}\wt{f}(t)}}_{\ltaa  L^2_{z,v}}^2+\nnr{\<\f{\na_{\fz,\fv}}{\mathsf{A}_k}\>^{\f\sss 2}\fM_1(t,\na)\rr{v^{\alpha}\wt{f}(t)}}^2_{\ltaa  L^2_{z,v}}}.\]
		For the nonlinear term $\mathscr{Q}_4$ we use the standard paraproduct decomposition as before
		\[\ba&\mathscr{Q}_4=-\sum_{|\alpha|\leq2}\sum_{\upmu,\upnu}\int\wh{\wt{\rho}}(t,\upnu\bfe)\f{\upxi-\upmu t}\upnu\rr{\vv{1}_{\<\upnu,\upnu t\>\leq\<\upmu-\upnu,\upxi-\upnu t\>}+\vv{1}_{\<\upnu,\upnu t\>\geq\<\upmu-\upnu,\upxi-\upnu t\>}}\pa_{\xi}^\alpha\wh{\wt{f}}(t,(\upmu-\upnu)\bfe,\xi-\upnu t\bfe)\\
		\times&\rr{\fM_1(t,\upmu\bfe,\xi)-\fM_1(t,(\upmu-\upnu)\bfe,\xi-\upnu t\bfe)}\ov{\sF\rr{\fM_1(t,\na)v^{\alpha}\wt{f}}(t,\upmu\bfe,\xi)}d\xi=\mathscr{Q}_{4,LH}+\mathscr{Q}_{4,HL}.\ea\]
		For $\mathscr{Q}_{4,LH}$, we have
		\[\mathscr{Q}_{4,LH}\ls t^{-2}\nnr{\fM_{-9}(t,\na)\wt{\rho}(t)}_{L^2_x}\rr{\nnr{\fM_{1.17}(t,\na)\rr{v^{\alpha}\wt{f}(t)}}_{\ltaa  L^2_{z,v}}^2+\nnr{\<\f{\na_{\fz,\fv}}{\mathsf{A}_k}\>^{\f\sss 2}\fM_1(t,\na)\rr{v^{\alpha}\wt{f}(t)}}^2_{\ltaa  L^2_{z,v}}}.\]
		For $\mathscr{Q}_{4,HL}$, we have
		\[\mathscr{Q}_{4,HL}\ls\nnr{\fM(t,\na)\wt{\rho}(t)}_{L^2_x}\nnr{\fM_1(t,\na)\rr{v^{\alpha}\wt{f}(t)}}^2_{\ltaa  L^2_{z,v}}t^2.\]
		The remaining terms can be estimated in the same way as before. By summing up all the estimates and integrating over the interval $[\mathsf{A}_k^5,t]$ we get
		\begin{equation}\label{w62}
			\begin{aligned}
				&\nnr{\fM_1(t,\na)v^{\alpha}\wt{f}(t)}_{\ltaa  L^2_{z,v}}^2 \\
				&\quad+\int_{\mathsf{A}_k^5}^t \rr{\f{0.01\mathsf{T}_k^{0.01}}{\tau^{1.01}}-\f{\eps_0^{0.9}}{\tau^2}-\f{10^2}{\tau^2}\nnr{\fM_{1}(\tau,\na)v^{\alpha}\wt{f}(\tau)}_{\ltaa  L^2_{z,v}}}\nnr{\fM_{1.17}(\tau,\na)v^{\alpha}\wt{f}(\tau)}_{\ltaa  L^2_{z,v}}^2d\tau \\
				&\quad +\int_{\mathsf{A}_k^5}^t\f{1}{\tau^2}\rr{4-\eps_0^{0.9}-10^2\nnr{\fM_{1}(\tau,\na)v^{\alpha}\wt{f}(\tau)}_{\ltaa  L^2_{z,v}}}\nnr{\<\f{\na_{\fz,\fv}}{\mathsf{A}_k}\>^{\f\sss2}\fM_{1}(\tau,\na)v^{\alpha}\wt{f}(\tau)}_{\ltaa  L^2_{z,v}}^2d\tau \\
				&\leq\nnr{\fM_1(\mathsf{A}_k^5,\na)v^{\alpha}\wt{f}(\mathsf{A}_k^5)}_{\ltaa  L^2_{z,v}}^2 \\
				&\quad+C_\theta\nnr{\<\tau\>^3\fM(\tau,\na)\wt{\rho}(\tau)}_{L^2_\tau([\mathsf{A}_k^5,t];L^2_x)}\sup_{\mathsf{A}_k^5\leq\tau\leq t}\nnr{\fM_1(\tau,\na)\rr{v^{\alpha}\wt{f}(\tau)}}_{\ltaa  L^2_{z,v}} \\
				&\quad +10^2\nnr{\<\tau\>^3\fM(\tau,\na)\wt{\rho}(\tau)}_{L^2_\tau([\mathsf{A}_k^5,t];L^2_x)}\sup_{\mathsf{A}_k^5\leq\tau\leq t}\nnr{\fM_1(\tau,\na)\rr{v^{\alpha}\wt{f}(\tau)}}^2_{\ltaa  L^2_{z,v}} \\
				&\quad +\int_{\mathsf{A}_k^5}^t\nnr{\fM_1(\tau,\na)\rr{v^{\alpha}\wt{Q}(\tau)}}_{\ltaa  L^2_{z,v}}d\tau \\
				&\qquad\times\sup_{\mathsf{A}_k^5\leq\tau\leq t}\nnr{\fM_1(\tau,\na)\rr{v^{\alpha}\wt{f}(\tau)}}_{\ltaa  L^2_{z,v}}
			\end{aligned}
		\end{equation}
		
		\noindent\emph{Step 4: bootstrap closure.}
		Assume that the following holds on $[\mathsf{A}_k^5,t^\dagger]$:
		\begin{equation}\label{bct}
			\left\{\begin{aligned}
				&\nnr{\<\tau\>^3\fM(\tau,\na)\wt{\rho}(\tau)}_{L^2_\tau([\mathsf{A}_k^5,t];L^2_x)}\leq2C_\theta e^{\mathsf{A}_k+\mathsf{A}_k^{10}}\mathsf{A}_k^{-\f{\mathtt{N}_3}{2}\mathsf{A}_k^{\mathtt{N}_1}}, \\
				&\nnr{\fM_1(t,\na)\rr{v^{\alpha}\wt{f}(t)}}_{\ltaa L^2_{z,v}}\leq 3C^2_\theta e^{\mathsf{A}_k+\mathsf{A}_k^{10}}\mathsf{A}_k^{-\f{\mathtt{N}_3}{2}\mathsf{A}_k^{\mathtt{N}_1}}.
			\end{aligned}\right.
		\end{equation}
		By \eqref{lerh9o} and \eqref{w62} the following holds on $[\mathsf{A}_k^5,t^\dagger]$
		\[\left\{\ba&\nnr{\<\tau\>^3\fM(\tau,\na)\wt{\rho}(\tau)}_{L^2_\tau([\mathsf{A}_k^5,t];L^2_x)}\leq 1.9C_\theta e^{\mathsf{A}_k+\mathsf{A}_k^{10}}\mathsf{A}_k^{-\f{\mathtt{N}_3}{2}\mathsf{A}_k^{\mathtt{N}_1}},\\
		&\nnr{\fM_1(t,\na)\rr{v^{\alpha}\wt{f}(t)}}_{\ltaa  L^2_{z,v}}\leq 2.9C^2_\theta e^{\mathsf{A}_k+\mathsf{A}_k^{10}}\mathsf{A}_k^{-\f{\mathtt{N}_3}{2}\mathsf{A}_k^{\mathtt{N}_1}}.
		\ea\right.\]
		Thus \eqref{bct} holds on a longer time interval. By a standard continuity argument \eqref{bct} holds on $[\mathsf{A}_k^5,\infty)$. 
	\end{proof}
	
	\subsection{Completion of the single-cascade construction}
	\label{subsec:single-cascade-completion}

	We now verify the four conclusions of Proposition \ref{propudlf}: uniform
	weighted control and scattering, persistence of the first-order output, and
	Gevrey-3 smallness at every time before the packet is inserted.
	
	\begin{proof}
		\noindent\emph{Step 1: global weighted control.}
		By \eqref{96}, Proposition \ref{ftildeshort} and Proposition \ref{escs} we have
		\[\ba&\nnr{\<\na_{\fz,\fv}\>^{20}e^{3\<\f{\na_{\fz,\fv}}{\mathsf{A}_k}\>^\sss+\<\na_{\mrv}\>^\sss}\rr{v^{\alpha}f(t)}}_{\ltaa  L^2_{z,v}}\leq\nnr{\<\na_{\fz,\fv}\>^{20}e^{3\<\f{\na_{\fz,\fv}}{\mathsf{A}_k}\>^\sss+\<\na_{\mrv}\>^\sss}\rr{v^{\alpha}f^{\#}(t)}}_{\ltaa  L^2_{z,v}}\\&+\nnr{\<\na_{\fz,\fv}\>^{20}e^{3\<\f{\na_{\fz,\fv}}{\mathsf{A}_k}\>^\sss+\<\na_{\mrv}\>^\sss}\rr{v^{\alpha}\wt{f}(t)}}_{\ltaa  L^2_{z,v}}\leq \mathsf{A}_k^{-\mathtt{N}_3+24.1}.
		\ea\]
		For $\wt f$, the derivative difference in Propositions \ref{ftildeshort}
		and \ref{escs} is absorbed using
		$\<\upmu,\upxi\>^m\ls \mathsf{A}_k^m e^{0.01\<\f{\upmu,\upxi}{\mathsf{A}_k}\>^\sss}$,
		with $m=4$ or $m=6$.

		\noindent\emph{Step 2: scattering decomposition.}
		We define 
		\[f_\infty=f(\mathsf{T}_k)+\int_{\mathsf{T}_k}^\infty
		\rr{\cN(f,f^0+\udl{f})+\cN(\udl{f},f)+\cN(f,f)}d\tau.\]
		Now we estimate the integral term.
		\[\ba
		&\nnr{\<\na_{\fz,\fv}\>^{17}e^{1.5\<\f{\na_{\fz,\fv}}{\mathsf{A}_k}\>^\sss+\<\na_{\mrv}\>^\sss}v^{\alpha}\cN(g,h)}_{\ltaa  L^2_{z,v}}\ls\nnr{\<\na_{\fz,\fv}\>^{18}e^{1.5\<\f{\na_{\fz,\fv}}{\mathsf{A}_k}\>^\sss}\rho_g (\tau)}_{L^2_x}\nnr{\<\na_{\fz,\fv}\>^{20}e^{1.5\<\f{\na_{\fz,\fv}}{\mathsf{A}_k}\>^\sss+\<\na_{\mrv}\>^\sss}v^{\alpha}h(\tau)}_{\ltaa  L^2_{z,v}}\\
		&\ls e^{-1.5\<\f{\tau}{\mathsf{A}_k}\>^\sss}\<\tau\>^{-2}\nnr{\<\na_{\fz,\fv}\>^{20}e^{3\<\f{\na_{\fz,\fv}}{\mathsf{A}_k}\>^\sss}\rho_g}_{L^2_x}\nnr{\<\na_{\fz,\fv}\>^{20}e^{2\<\f{\na_{\fz,\fv}}{\mathsf{A}_k}\>^\sss+\<\na_{\mrv}\>^\sss}v^{\alpha}h(\tau)}_{\ltaa  L^2_{z,v}}
		\ea\]
		Thus
		\[\ba
		&\nnr{\<\na_{\fz,\fv}\>^{17}e^{1.5\<\f{\na_{\fz,\fv}}{\mathsf{A}_k}\>^\sss+\<\na_{\mrv}\>^\sss}\rr{v^{\alpha}f_\infty}}_{\ltaa  L^2_{z,v}}\leq\nnr{\<\na_{\fz,\fv}\>^{17}e^{1.5\<\f{\na_{\fz,\fv}}{\mathsf{A}_k}\>^\sss+\<\na_{\mrv}\>^\sss}\rr{v^{\alpha}f(\mathsf{T}_k)}}_{\ltaa  L^2_{z,v}}+\int_{\mathsf{T}_k}^\infty 10e^{-1.5\<\f{\tau}{\mathsf{A}_k}\>^\sss}\<\tau\>^{-2}d\tau\\
		&\times \sup_{\tau\geq \mathsf{T}_k}\nnr{\<\na_{\fz,\fv}\>^{20}e^{3\<\f{\na_{\fz,\fv}}{\mathsf{A}_k}\>^\sss+\<\na_{\mrv}\>^\sss}\rr{v^{\alpha}f(\tau)}}_{\ltaa  L^2_{z,v}}\left(\nnr{\<\na_{\fz,\fv}\>^{20}e^{3\<\f{\na_{\fz,\fv}}{\mathsf{A}_k}\>^\sss+\<\na_{\mrv}\>^\sss}\rr{v^{\alpha}f^0(\tau)}}_{\ltaa  L^2_{v}}+\sup_{\tau\geq \mathsf{T}_k}\nnr{\<\na_{\fz,\fv}\>^{20}\right.\right.\\
			&\left.\left.e^{3\<\f{\na_{\fz,\fv}}{\mathsf{A}_k}\>^\sss+\<\na_{\mrv}\>^\sss}\rr{v^{\alpha}\udl{f}(\tau)}}_{\ltaa  L^2_{z,v}}+\sup_{\tau\geq \mathsf{T}_k}\nnr{\<\na_{\fz,\fv}\>^{20}e^{3\<\f{\na_{\fz,\fv}}{\mathsf{A}_k}\>^\sss+\<\na_{\mrv}\>^\sss}\rr{v^{\alpha}f(\tau)}}_{\ltaa  L^2_{z,v}}\right)\leq \mathsf{A}_k^{-\mathtt{N}_3+24.2}.
		\ea\]
		Since $f(t)=f(\mathsf{T}_k)+\int_{\mathsf{T}_k}^t\partial_\tau f(\tau)d\tau$, the decomposition
		$f(t)=f_\infty+f_{err}(t)$ requires
		\[f_{err}(t)=-\int_t^\infty\rr{\cN(f,f^0+\udl{f})+\cN(\udl{f},f)+\cN(f,f)}d\tau.\]
		By a similar argument we have
		\[\nnr{\<\na_{\fz,\fv}\>^{17}e^{1.5\<\f{\na_{\fz,\fv}}{\mathsf{A}_k}\>^\sss+\<\na_{\mrv}\>^\sss}\rr{v^{\alpha}f_{err}(t)}}_{\ltaa  L^2_{z,v}}\leq\int_{t}^\infty 10e^{-1.5\<\f{\tau}{\mathsf{A}_k}\>^\sss}\<\tau\>^{-2}d\tau \mathsf{A}_k^{-\mathtt{N}_3+24.15}\leq e^{-1.5\<\f{t}{\mathsf{A}_k}\>^\sss}\mathsf{A}_k^{-\mathtt{N}_3+24.2}.\]
		Thus
		\[\ba
		&\nnr{\<\na_{\fz,\fv}\>^{20}e^{\<\f{\na_{\fz,\fv}}{\mathsf{A}_k}\>^\sss+\<\na_{\mrv}\>^\sss}\rr{v^{\alpha}f_\infty}}_{\ltaa  L^2_{z,v}}\ls \mathsf{A}_k^3\nnr{\<\na_{\fz,\fv}\>^{17}e^{1.5\<\f{\na_{\fz,\fv}}{\mathsf{A}_k}\>^\sss+\<\na_{\mrv}\>^\sss}\rr{v^{\alpha}f_\infty}}_{\ltaa  L^2_{z,v}}\leq \mathsf{A}_k^{-\mathtt{N}_3+27.3},\\
		&\nnr{\<\na_{\fz,\fv}\>^{20}e^{\<\f{\na_{\fz,\fv}}{\mathsf{A}_k}\>^\sss+\<\na_{\mrv}\>^\sss}\rr{v^{\alpha}f_{err}(t)}}_{\ltaa  L^2_{z,v}}\ls \mathsf{A}_k^3\nnr{\<\na_{\fz,\fv}\>^{17}e^{1.5\<\f{\na_{\fz,\fv}}{\mathsf{A}_k}\>^\sss+\<\na_{\mrv}\>^\sss}\rr{v^{\alpha}f_{err}(t)}}_{\ltaa  L^2_{z,v}}\leq \mathsf{A}_k^{-\mathtt{N}_3+27.3}e^{-1.5\<\f{t}{\mathsf{A}_k}\>^\sss}.
		\ea\]
		\noindent\emph{Step 3: survival of the echo output.}
		By \eqref{10} and \eqref{yup},
		\[\ba
		&\rr{\int_{0.9\mathsf{A}_k\leq|\upxi|\leq1.1\mathsf{A}_k}\nr{\wh{f}(t,0,\xi)}^2d\xi}^\f12\geq\rr{\int_{0.9\mathsf{A}_k\leq|\upxi|\leq1.1\mathsf{A}_k}\nr{\wh{f^{(1)}}(t,0,\xi)}^2d\xi}^\f12-\sum_{n=2}^{\mathsf{A}_k^{\mathtt{N}_1}}\nnr{f^{(n)}(t)}_{L^2_{z,v}}-\nnr{\wt{f}(t)}_{L^2_{z,v}}\\&\geq	2\mathsf{A}_k^{-\mathtt{N}_3+0.9}-\mathsf{A}_k^{-1.1\mathtt{N}_3}\geq \mathsf{A}_k^{-\mathtt{N}_3+0.9}.
		\ea\]
		\noindent\emph{Step 4: Gevrey-3 smallness before insertion.}
		By \eqref{taylor-short-time} and \eqref{10T1}, on $[10,\mathsf{T}_k]$ it holds that
		\[\nnr{e^{\<\na_{\mrv}\>^\sss}\<\na_{\fz,\fv}\>^3\rr{v^{\alpha}f^{(n)}(t)}}_{\ltaa L^2_{z,v}}\leq e^{-2.25nk}\]
		Thus
		\[\ba&\nnr{\<\na_{\fz,\fv}\>^3e^{1.1\<\eps_0\na_{\fz,\fv}\>^\f13+\<\na_{\mrv}\>^\sss}\rr{v^{\alpha}f^{(n)}(t)}}_{\ltaa L^2_{z,v}}\leq\nnr{\<\upmu,\upxi\>^3e^{1.1\<\eps_0|\upmu,\upxi|\>^\f13+\<\mrxi\>^\sss}\pa_{\xi}^\alpha\wh{f^{(n)}}(t,\upmu\bfe,\xi)\vv{1}_{|\upmu,\upxi|\leq5n^{\f1\sss}\mathsf{A}_k}}_{\ltaa L^2_{\upmu,\xi}}\\
		&+\nnr{\<\upmu,\upxi\>^3e^{1.1\<\eps_0|\upmu,\upxi|\>^\f13+\<\mrxi\>^\sss}\pa_{\xi}^\alpha\wh{f^{(n)}}(t)\vv{1}_{|\upmu,\upxi|\geq5n^{\f1\sss}\mathsf{A}_k}}_{\ltaa L^2_{\upmu,\xi}}
		\leq e^{1.89n^{\f{1}{3\sss}}k}\nnr{e^{\<\na_{\mrv}\>^\sss}\<\na_{\fz,\fv}\>^3\rr{v^{\alpha}f^{(n)}(t)}}_{\ltaa L^2_{z,v}}\\
		&+e^{-0.1\nr{\f{n\mathsf{A}_k}{k^\f13}}^\sss}\leq e^{-0.3nk}.
		\ea\]
		Moreover we have 
		\[\nnr{\<\na_{\fz,\fv}\>^3e^{1.1\<\eps_0\na_{\fz,\fv}\>^\f13+\<\na_{\mrv}\>^\sss}\rr{v^{\alpha}\wt{f}(t)}}_{\ltaa L^2_{z,v}}\leq\nnr{\<\na_{\fz,\fv}\>^3e^{1.1\<\f{\na_{\fz,\fv}}{\mathsf{A}_k}\>^\sss+\<\na_{\mrv}\>^\sss}\rr{v^{\alpha}\wt{f}(t)}}_{\ltaa L^2_{z,v}}e^{\mathsf{A}_k}\leq e^{3\mathsf{A}_k^{10}}\mathsf{A}_k^{-\f{\mathtt{N}_3}{2}\mathsf{A}_k^{\mathtt{N}_1}}. \]
		Thus 
		\[\nnr{\<\na\>^3e^{1.1\<\eps_0\na\>^\f13}v^{\alpha}f(t)}_{\ltaa L^2_{z,v}}\leq\sum_{n=1}^{\mathsf{A}_k^{\mathtt{N}_1}}\nnr{\<\na\>^3e^{1.1\<\eps_0\na\>^\f13}v^{\alpha}f^{(n)}(t)}_{\ltaa L^2_{z,v}}+\nnr{\<\na\>^3e^{1.1\<\eps_0\na\>^\f13}v^{\alpha}\wt{f}(t)}_{\ltaa  L^2_{z,v}}\leq e^{-0.25k}.\]
	\end{proof}

	\section{Infinite echo cascades and global instability}
	\label{sec:global-instability}

	This section proves Theorem \ref{T1} by combining infinitely many echo cascades constructed in Proposition \ref{propudlf}. We first construct a low-frequency background solution and then successively introduce cascades at rapidly increasing frequency scales. The resulting series converges to a global solution. Moreover, the frequency separation ensures that the contributions from different cascades cannot cancel the dominant contribution on each output band. Since the lower bound for each cascade persists for all sufficiently large times, we obtain the divergence of the high-order Sobolev norm as \(t\to\infty\).

	\subsection{The background solution and scale hierarchy}
	
    We construct the background wave by solving the final-state problem with asymptotic profile $\delta\vf(v)\cos\fz$. The proof adapts the Gevrey estimates of \cite[Sections 2--4]{waveop2} to the equilibrium assumptions in \eqref{Penrose}.

	\begin{lem}[Background wave]\label{lem:background-wave}
	There exists $\delta_0=\delta_0(\theta)>0$ such that, for every $0<\delta<\delta_0$, equation \eqref{1} has a real-valued solution $f_{[1]}$ on $[10,\infty)$ of the form
	\[f_{[1]}(t,z,v)=\delta\vf(v)\cos\fz+f_{(1)}(t,z,v),\qquad
	f_{(1)}(t,z,v)=f_{(1)}(t,\fz\bfe,v),\]
	which satisfies
	\[\nnr{\<\na\>^{20}e^{1.5\<\na\>^\sss}v^{\alpha}f_{(1)}(t)}_{\lta L^2_{z,v}}
	\leq C_{\theta}\delta e^{-\<t/4\>^\sss}.\]
	The associated distribution $f^0+f_{[1]}$ is nonnegative, and the corresponding density and electric field converge to zero as $t\to\infty$.
	\end{lem}

	The next observation records all scale separations used in the
	infinite iteration. Its proof is rather elementary.

	\begin{lem}[Scale separation]\label{lem:scale-separation}
	Let
	\[
	\eps_1:=\f14\rr{\exp\exp(1/\eps_0)}^{-1},\qquad
	A_1:=4\exp\exp(1/\eps_0),
	\]
	and define
	\[
	k_2:=2\lfloor\exp\exp(A_1)\rfloor,\qquad
	A_2:=\f{(1+\eps_0^6)k_2^3}{\eps_0}.
	\]
	For $j\geq2$, define recursively
	\[
	k_{j+1}:=\lfloor\exp\exp\exp(A_j)\rfloor,\qquad
	A_{j+1}:=\f{(1+\eps_0^6)k_{j+1}^3}{\eps_0},\qquad
	T_j:=\f{A_j}{k_j+\f12}.
	\]
	If $\eps_0$ is sufficiently small depending only on $\theta$, then:
	\begin{enumerate}[(i)]
	\item At each step of the iteration, the frequency parameter satisfies the condition
	
	$$
	k\geq \max\{\exp\exp(1/\eps_*),\exp\exp(A_*)\}
	$$
	
	required in Proposition \ref{propudlf}. More precisely, this holds with \((k,\eps_*,A_*)=(k_2,\eps_1,A_1)\) for the first insertion, and with \((k,\eps_*,A_*)=(k_{j+1},A_2^{-\mathtt N_3+29},A_j)\) for every \(j\geq2\);
	\item for every $j\geq2$,
	\[
	A_2^{-\mathtt N_3+28}\leq\f12\eps_0,\qquad
	A_{j+1}^{-\mathtt N_3+28}
	\leq2^{-j}A_2^{-\mathtt N_3+28},\qquad
	\sum_{\ell\geq j+1}A_\ell^{-\mathtt N_3+28}
	\leq A_{j+1}^{-\mathtt N_3+29};
	\]
	\item for every $J\geq3$,
	\[
	\sum_{j=2}^{J-1}
	e^{-\f12\<A_J/A_j\>^\sss}A_j^{-\mathtt N_3+28}
	\leq e^{-0.1A_J^{0.5}}.
	\]
	\end{enumerate}
	\end{lem}

	\subsection{Infinite cascade iteration and convergence}
	We now construct the infinite sequence of cascades and prove Theorem \ref{T1}.
	
	\begin{proof}
		\noindent\emph{Step 1: initialization and inductive construction.}
		Fix $0<\epsilon<\epsilon_{\rm th}$ and set $\eps_0:=\epsilon$. We choose
		$\epsilon_{\rm th}(\theta)$ below all the smallness thresholds in Lemma
		\ref{lem:background-wave} and Proposition \ref{propudlf}. Apply Lemma \ref{lem:background-wave} with
		$\delta=\eps_0$, and set
		\[f_{(1)}:=f_{[1]}-\eps_0\vf(v)\cos\fz.\]
		We use the decomposition $f_{(1),\infty}:=0$ and
		$f_{(1),err}:=f_{(1)}$.
		Then $f_{(1)}=0+f_{(1)}$, together with $(\eps_1,A_1)$ from Lemma
		\ref{lem:scale-separation}, satisfies the assumptions on $f_*$ in
		Proposition \ref{propudlf}.
		
		Let $k_2,A_2,T_2$ be as in Lemma \ref{lem:scale-separation}.
		Proposition \ref{propudlf} gives an increment $f_{(2)}(t,z,v)$ solving
		\[\ba\left\{\ba
		&\pa_tf=\cN(f,f^0+f_{[1]})+\cN(f_{[1]},f)+\cN(f,f),\\
		&f(T_2)=\eps_2\psi(k_2^\f13\fv)\phi(\mrv)\cos(k_2\fz+A_2\fv),
		\ea\right.\ea\]
		and satisfying all the conclusions of Proposition \ref{propudlf}. In
		particular, $f_{[2]}:=f_{[1]}+f_{(2)}$ solves \eqref{1}, because the
		perturbation equation is exactly the difference between the equations for
		$f_{[2]}$ and $f_{[1]}$.
		
		We construct $\{f_{[J]}\}_{J=2}^{\infty}$ by induction. Assume that
		$f_{[J]}$ is a real-valued, neutral solution of \eqref{1} on
		$[10,\infty)$ and
		\[f_{[J]}(t,z,v)=\eps_0\vf(v)\cos \fz+\sum_{j=1}^{J}f_{(j)}(t,z,v),\]
		where $f_{(j)}(t,z,v)=f_{(j)}(t,\fz\bfe,v)$ and
		\[\ba
		&f_{(j)}(t,z,v)=f_{(j),\infty}(z,v)+f_{(j),err}(t,z,v),\\
		&\nnr{\<\na_{\fz,\fv}\>^{20}e^{\<\f{\na_{\fz,\fv}}{A_{J}}\>^\sss+\<\na_{\mrv}\>^\sss}\rr{v^{\alpha}\sum_{j=1}^{J}f_{(j),\infty}}}_{\ltaa  L^2_{z,v}}\leq\rr{\sum_{j=1}^{J}2^{-j+1}}A_2^{-\mathtt N_3+28},\\
		&\nnr{\<\na_{\fz,\fv}\>^{20}e^{\<\f{\na_{\fz,\fv}}{A_{J}}\>^\sss+\<\na_{\mrv}\>^\sss}\rr{v^{\alpha}\sum_{j=1}^{J}f_{(j),err}(t)}}_{\ltaa  L^2_{z,v}}\leq C_\theta\eps_0e^{-\<\f{t}{A_{J}}\>^\sss}\rr{\sum_{j=1}^{J}2^{-j+1}},
		\ea\]
		where $\{A_j\}_{j=1}^J$ is strictly increasing. By the above construction, the inductive hypothesis holds for $J=2$. Moreover, the inductive bounds imply
		\[
		\nnr{\<\na_{\fz,\fv}\>^{20}
		e^{\<\na_{\fz,\fv}/A_J\>^\sss+\<\na_{\mrv}\>^\sss}
	\rr{v^\alpha\sum_{j=1}^Jf_{(j),\infty}}}_{\lta L^2_{z,v}}
		\leq \eps_0 A_2^{-\mathtt N_3+29}.
		\]
		Thus $f_{[J]}$ satisfies the
		background assumptions of Proposition \ref{propudlf} with
		$(\eps_*,A_*)=(A_2^{-\mathtt N_3+29},A_J)$.
		
		Let $k_{J+1},A_{J+1},T_{J+1}$ be given by Lemma
		\ref{lem:scale-separation}. Proposition \ref{propudlf} gives a solution
		$f_{(J+1)}$ to
		\[\ba\left\{\ba
		&\pa_tf=\cN(f,f^0+f_{[J]})+\cN(f_{[J]},f)+\cN(f,f),\\
		&f(T_{J+1})=\eps_{J+1}\psi(k_{J+1}^\f13\fv)\phi(\mrv)\cos(k_{J+1}\fz+A_{J+1}\fv),
		\ea\right.\ea\]
		that satisfies the conclusion of Proposition \ref{propudlf}. Set
		$f_{[J+1]}=f_{[J]}+f_{(J+1)}$. Expanding the equation shows that
		$f_{[J+1]}$ again solves \eqref{1}. It's easy to see that $f_{(J+1)}$ is neutral and satisfies $f_{(J+1)}(t,z,v)=f_{(J+1)}(t,\fz\bfe,v)$. By \eqref{finfty8} and
		\eqref{fperturb}, $f_{[J+1]}$ satisfies the inductive
		hypothesis for $J+1$.
		
		\noindent\emph{Step 2: convergence to a global solution.}
		Set
		\[f(t)=\lim_{J\rightarrow\infty}f_{[J]}(t).\]
		We justify this limit before using it, recalling the norm $\mathcal X^m$
		defined before Theorem \ref{T1}.
		By \eqref{finfty8} and \eqref{fperturb}, for every $j\geq2$,
		\[\sup_{t\geq10}\nnr{f_{(j)}(t)}_{\mathcal X^{20}}
		\leq2\mathsf{A}_{k_j}^{-\mathtt N_3+28}.\]
		The sequence $\{\mathsf{A}_{k_j}\}$ grows super-exponentially and $\mathtt N_3>30$; hence the series on the right is summable. Therefore $f_{[J]}$ is Cauchy in $C([10,R];\mathcal X^{20})$ for every $R>10$, and the convergence is uniform on $[10,\infty)$ in $\mathcal X^{20}$.

		For $t\in[10,R]$, Lemma \ref{lem:nonlinear-continuity} gives
		\[\nnr{\cN(g,h)(t)}_{\mathcal X^{17}}\leq C_R\nnr{g(t)}_{\mathcal X^{20}}\nnr{h(t)}_{\mathcal X^{20}}.\]
		Consequently the densities, electric fields, and nonlinear terms of $f_{[J]}$ converge in the corresponding Sobolev spaces, uniformly on compact time intervals. Passing to the limit in the integral formulation of \eqref{1} shows that
		\[f\in C([10,R];\mathcal X^{20})\cap C^1([10,R];\mathcal X^{17})\]
		and that $f$ solves \eqref{1}. The mass condition passes to the limit as
		well. Since $\mathcal X^{20}$ embeds into $C^0_{z,v}$, the convergence is
		also uniform in $(z,v)$; once nonnegativity of the approximants is proved
		below, it therefore passes to the limit.

		\noindent\emph{Step 3: scattering and damping.}
		The same summability identifies the scattering profile. Define
		\[f_\infty:=\eps_0\vf(v)\cos\fz+\sum_{j=2}^{\infty}f_{(j),\infty}.\]
		Lemma \ref{lem:background-wave}, \eqref{fperturb}, and dominated convergence
		yield
		\[\lim_{t\to\infty}\nnr{f(t)-f_\infty}_{\mathcal X^{18}}=0.\]
		Every increment has zero mass, and its limiting profile has zero mass by
		\eqref{fperturb}. Hence both $f_\infty$ and $f(t)-f_\infty$ are neutral.
		For a neutral $g\in\mathcal X^{18}$, Lemma \ref{lem:fourier-trace} and the
		identity $\wh{\rho}_g(t,\mu)=\wh{g}(\mu,\mu t)$ give
		\[\nnr{\rho_g(t)}_{H^{10}_x}+\nnr{E_g(t)}_{H^{11}_x}
		\leq C\<t\>^{-6}\nnr{g}_{\mathcal X^{18}}.\]
		Applying this estimate to $g=f_\infty$ and, at each time $t$, to
		$g=f(t)-f_\infty$, proves
		\[\lim_{t\to\infty}\rr{\nnr{\rho_f(t)}_{H^{10}_x}+\nnr{E_f(t)}_{H^{11}_x}}=0.\]
		\noindent\emph{Step 4: persistence and non-cancellation of the cascade outputs.}
		Fix $J\geq2$. On $[2A_J,\infty)$, \eqref{lowerbound} gives
		\[\rr{\int_{0.9A_{J}\leq|\upxi|\leq1.1A_{J}}\nr{\wh{f_{(J)}}(t,0,\xi)}^2d\xi}^\f12\geq A_{J}^{-\mathtt N_3+0.9}.\]
		For $2\leq j\leq J-1$, by \eqref{finfty8} and \eqref{fperturb} it holds that
		\[\rr{\int_{0.9A_{J}\leq|\upxi|\leq1.1A_{J}}\nr{\wh{f_{(j)}}(t,0,\xi)}^2d\xi}^\f12\leq e^{-\f12\<\f{A_{J}}{A_{j}}\>^\sss}\nnr{e^{\<\f{\na_{\fz,\fv}}{A_{j}}\>^\sss}f_{(j)}(t)}_{L^2_{z,v}}\leq 2e^{-\f12\<\f{A_{J}}{A_{j}}\>^\sss}A^{-\mathtt N_3+28}_{j}.\]
		For $j\geq J+1$, by \eqref{finfty8} and \eqref{fperturb} it holds that
		\[\rr{\int_{0.9A_{J}\leq|\upxi|\leq1.1A_{J}}\nr{\wh{f_{(j)}}(t,0,\xi)}^2d\xi}^\f12\leq 2A_{j}^{-\mathtt N_3+28}.\]
		The persistent part of $f_{[1]}$ has spatial frequency $\pm1$ and hence
		does not contribute at spatial mode zero. Lemma \ref{lem:background-wave}
		therefore gives
		\[\rr{\int_{0.9A_J\leq|\upxi|\leq1.1A_J}
		\nr{\wh{f_{[1]}}(t,0,\xi)}^2d\xi}^\f12
		\leq e^{-0.1A_J^\sss}.\]
		Using Lemma \ref{lem:scale-separation}(ii)--(iii), we conclude that
		\[\ba
		&\nnr{\<\na\>^{\mathtt N_3}f(t)}_{L^2_{z,v}}\gs A_{J}^{\mathtt N_3} \rr{\int_{0.9A_{J}\leq|\upxi|\leq1.1A_{J}}\nr{\wh{f}(t,0,\xi)}^2d\xi}^\f12\geq A_{J}^{\mathtt N_3}\left\{\rr{\int_{0.9A_{J}\leq|\upxi|\leq1.1A_{J}}\nr{\wh{f_{(J)}}(t,0,\xi)}^2d\xi}^\f12\right.\\
		&-\sum_{j\geq J+1}\rr{\int_{0.9A_{J}\leq|\upxi|\leq1.1A_{J}}\nr{\wh{f_{(j)}}(t,0,\xi)}^2d\xi}^\f12-\sum_{j=2}^{J-1}\rr{\int_{0.9A_{J}\leq|\upxi|\leq1.1A_{J}}\nr{\wh{f_{(j)}}(t,0,\xi)}^2d\xi}^\f12\\
		&\left.-\rr{\int_{0.9A_{J}\leq|\upxi|\leq1.1A_{J}}\nr{\wh{f_{[1]}}(t,0,\xi)}^2d\xi}^\f12\right\}\geq A_{J}^{\mathtt N_3}\rr{A_{J}^{-\mathtt N_3+0.9}-2A_{J+1}^{-\mathtt N_3+29}-2e^{-0.1A^{0.5}_{J} }-e^{-0.1A_{J}^\sss}}\geq0.5A^{0.9}_{J}
		\ea\]
		for every $t\geq2A_J$. Since $\mathtt N_3=1000$ and $A_J\to\infty$,
		\[\lim_{t\rightarrow\infty}\nnr{\<\na\>^{1000}f(t)}_{L^2_{z,v}}=\infty.\]

		\noindent\emph{Step 5: initial regularity and nonnegativity.}
		By \eqref{initialdata},
		\[\nnr{\<\na\>^3e^{1.1\<\eps_0\na\>^\f13}\rr{v^{\alpha}f_{(j)}(10)}}_{\lta L^2_{z,v}}\leq e^{-\f14k_j}\]
		for any $j\geq2$. Thus 
		\[\nnr{\<\na\>^3e^{1.1\<\eps_0\na\>^\f13}\rr{v^{\alpha}f(10)}}_{\lta L^2_{z,v}}\leq C_\theta\eps_0.\]
		Here the background contribution follows from Lemma
		\ref{lem:background-wave}, and the series of packet contributions is
		summable by Lemma \ref{lem:scale-separation}. By standard local well-posedness and Gevrey energy estimates, the solution extends uniquely backward to \([0,10]\). Moreover, the initial Gevrey radius can be chosen slightly smaller than that at \(t=10\), yielding
		\[\nnr{\<\na\>^3e^{\<\eps_0\na\>^\f13}\rr{v^{\alpha}f(0)}}_{\lta L^2_{z,v}}\leq C_\theta\eps_0=C_\theta\epsilon.\]

		It remains to verify nonnegativity. Lemma \ref{lem:background-wave} gives
		$f^0+f_{[1]}\geq0$. Suppose inductively that $f^0+f_{[J]}\geq0$. At the
		insertion time,
		\[f^0(v)+f_{[J+1]}(T_{J+1},z,v)=f^0(v)+f_{[J]}(T_{J+1},z,v)+\eps_{J+1}\psi(k_{J+1}^\f13\fv)\phi(\mrv)\cos(k_{J+1}\fz+A_{J+1}\fv).\]
		On $\{|\fv|\leq k_{J+1}^{-\f13}\}\cap\{|\mrv|\leq 1\}$ it holds that 
		\[f^0(v)+f_{[J]}(T_{J+1},z,v)+\eps_{J+1}\psi(k_{J+1}^\f13\fv)\phi(\mrv)\cos(k_{J+1}\fz+A_{J+1}\fv)\geq\theta e^{-3}-\nnr{f_{[J]}(T_{J+1})}_{L^\infty_{z,v}}-\eps_{J+1}\nnr{\psi}_{L^\infty}\nnr{\phi}_{L^\infty}\geq\f12\theta e^{-3}.\]
		On $\{|\fv|\geq k_{J+1}^{-\f13}\}\cup\{|\mrv|\geq 1\}$ it holds that
		\[f^0(v)+f_{[J+1]}(T_{J+1},z,v)=f^0(v)+f_{[J]}(T_{J+1},z,v)\geq0.\]
		The total distribution satisfies a transport equation, so nonnegativity at
		$T_{J+1}$ propagates both forward and backward on $[0,\infty)$. Thus every
		$f^0+f_{[J]}$ is nonnegative there. Uniform convergence gives
		$f^0+f\geq0$ on $[0,\infty)$.
	\end{proof}

	\appendix
	\section{Invariant subspaces and the Fourier formulation}
	\label{app:invariant-subspace}

	We first record the Fourier equation and justify the reduction to solutions
	that depend only on the first spatial coordinate. Taking the Fourier transform
	of \eqref{1} gives
	\begin{equation}\label{xv}
		\pa_t\wh{f}(t,\mu,\xi)+\sum_{\nu\in\Z^3\backslash\{0\}}
	\wh{f}(t,\nu,\nu t)\f{\nu\cdot(\xi-\mu t)}{|\nu|^2}
	\sF\rr{f^0+f}(t,\mu-\nu,\xi-\nu t)=0.
	\end{equation}

	\begin{prop}[Invariant one-dimensional spatial subspace]
	\label{propreductionto1d}
	Let $I$ be an interval containing $t_0$, and let $f$ be a unique solution of
	\eqref{1} on $I$ in a class embedded in
	$C(I;H^6_{z,v}(\<v\>^4dzdv))$. If
	\[f(t_0,z,v)=f(t_0,\fz\bfe,v),\]
	then
	\[f(t,z,v)=f(t,\fz\bfe,v)\qquad\text{for every }t\in I.\]
	Moreover, the mass $\int f(t,z,v)dzdv$ is constant on $I$.
	\end{prop}
	\begin{proof}
	For $a\in\T^2$, define
	$f_a(t,\fz,\mathring z,v):=f(t,\fz,\mathring z+a,v)$. The equilibrium and
	the Vlasov--Poisson equation are invariant under this translation, so $f_a$
	and $f$ solve the same equation and have the same value at $t_0$. Uniqueness
	gives $f_a=f$ for every $a$, proving independence of $\mathring z$. Equivalently,
	$\wh f(t,\mu,\xi)$ is supported on $\mathring\mu=0$, and \eqref{xv} then
	reduces to the one-dimensional spatial sum used in the main text. Finally,
	evaluating \eqref{xv} at $(\mu,\xi)=(0,0)$ gives
	$\partial_t\wh f(t,0,0)=0$, which proves conservation of mass.
	\end{proof}

	\section{Auxiliary Sobolev and trace estimates}
	\label{app:auxiliary-estimates}

	Recall the norm
	\[
	\nnr{g}_{\mathcal X^m}:=
	\nnr{\<\na\>^m\rr{v^\alpha g}}_{\ell^2(|\alpha|\leq2)L^2_{z,v}}.
	\]

	\begin{lem}[Continuity of the nonlinearity]
	\label{lem:nonlinear-continuity}
	For every $R\geq10$ and $t\in[0,R]$,
	\[
	\nnr{\cN(g,h)(t)}_{\mathcal X^{17}}
	\leq C_R\nnr{g(t)}_{\mathcal X^{20}}
	\nnr{h(t)}_{\mathcal X^{20}}.
	\]
	\end{lem}
	\begin{proof}
	For nonzero spatial modes, the Poisson multiplier gains one spatial
	derivative. Sobolev embedding and the velocity-moment bounds give
	\[
	\nnr{E_g(t)}_{W^{17,\infty}_z}
	\leq C\nnr{g(t)}_{\mathcal X^{20}}.
	\]
	The operator $\na_v-t\na_z$ costs one derivative and contributes at most a
	factor $1+R$. Applying the Sobolev product estimate to
	$E_g(t,z+vt)\cdot(\na_v-t\na_z)h$ and distributing the two velocity moments
	between the factors proves the result.
	\end{proof}

	\begin{lem}[Fourier trace estimate]\label{lem:fourier-trace}
	Suppose $g\in\mathcal X^{18}$ and
	$\int_{\T^3\times\R^3}g(z,v)dzdv=0$. Then, for every $t\geq0$,
	\[
	\nnr{\rho_g(t)}_{H^{10}_x}+\nnr{E_g(t)}_{H^{11}_x}
	\leq C\<t\>^{-6}\nnr{g}_{\mathcal X^{18}}.
	\]
	\end{lem}
	\begin{proof}
	The zero spatial mode vanishes by the neutrality assumption. For
	$\mu\neq0$, Sobolev embedding $H^2_\xi(\R^3)\hookrightarrow L^\infty_\xi$
	and the identity $\wh\rho_g(t,\mu)=\wh g(\mu,\mu t)$ give
	\[
	\nr{\wh g(\mu,\mu t)}
	\leq C\<\mu,\mu t\>^{-18}
	\sum_{|\alpha|\leq2}
	\nnr{\<\mu,\xi\>^{18}\pa_\xi^\alpha\wh g(\mu,\xi)}_{L^2_\xi}.
	\]
	Since $\mu\neq0$,
	$\<\mu,\mu t\>\gtrsim\<\mu\>\<t\>$. Squaring, summing in $\mu$, and
	using Plancherel yield
	\[
	\nnr{\rho_g(t)}_{H^{10}_x}
	\leq C\<t\>^{-6}\nnr{g}_{\mathcal X^{18}}.
	\]
	Finally, the multiplier of $E_g=\na_x\Delta_x^{-1}\rho_g$ satisfies
	$\<\mu\>^{11}/|\mu|\lesssim\<\mu\>^{10}$ for $\mu\neq0$, which proves
	the electric-field estimate.
	\end{proof}

	{{\bf Acknowledgments.}} The author would like to thank Prof. Deng Yu and Prof. He Ling-Bing for fruitful discussions.

	{\bf {Data availability statement.}} Data sharing is not applicable to this article as no datasets were generated or analysed during the current study.

	\textbf{AI declaration} This work was completed without the assistance of artificial intelligence (AI) tools. The author subsequently used ChatGPT solely for language editing and polishing of the manuscript. The author takes full responsibility for the correctness and validity of all mathematical results presented in this paper.

\end{document}